\documentclass{article}

\usepackage{amsmath,amssymb,graphics}
\usepackage[all]{xy}

\usepackage{fullpage}

\def\wedgebullet {{\mbox{$\bigwedge^\bullet$}}}

\newtheorem{definition}{Definition}
\newtheorem{proposition}{Proposition}
\newtheorem{corollary}{Corollary}
\newtheorem{theorem}{Theorem}

\newtheorem{desired statement}{Desired statement}

\def\gg {\mathfrak{g}}
\def\ff {\mathfrak{f}}
\def\hh {\mathfrak{h}}
\def\uu {\mathfrak{u}}

\def\( {\begin{equation}}
\def\) {\end{equation}}

\def\adress#1{\gdef\@adress{#1}}

\renewcommand{\(}{\begin{equation}}
\renewcommand{\)}{\end{equation}}
\newcommand{\bea}{\begin{eqnarray}}
\newcommand{\eea}{\end{eqnarray}}

\def\proof {{Proof.}\hspace{7pt}}

\def\endofproof {\hfill{$\Box$}\\}

\begin{document}

\title{$L_\infty$-algebra connections and applications\\
 to String- and Chern-Simons $n$-transport}

\author{
  Hisham Sati\thanks{hisham.sati@yale.edu}, 
  Urs Schreiber\thanks{schreiber@math.uni-hamburg.de} 
  and 
  Jim Stasheff\thanks{jds@math.upenn.edu}}

\maketitle

\begin{abstract}
  We give a generalization of the notion of a Cartan-Ehresmann connection
  from Lie algebras to $L_\infty$-algebras and use it to study
  the obstruction theory of lifts through higher String-like extensions of 
  Lie algebras. We find (generalized) Chern-Simons and BF-theory 
  functionals this way and describe aspects of their parallel transport
  and quantization.

  It is known 
  that over a D-brane the Kalb-Ramond background field of the string
  restricts to a 2-bundle with
  connection (a gerbe) which can be seen as the obstruction to lifting the
  $PU(H)$-bundle on the D-brane to a $U(H)$-bundle. We discuss how this 
  phenomenon generalizes
  from the ordinary central extension $U(1) \to U(H) \to PU(H)$ to higher
  categorical central extensions, like the String-extension
  $\mathbf{B} U(1) \to \mathrm{String}(G) \to G$. Here the obstruction to the lift
  is a 3-bundle with connection (a 2-gerbe): the Chern-Simons 3-bundle classified
  by the first Pontrjagin class.
  For $G = \mathrm{Spin}(n)$ this obstructs the existence of a String-structure.
  We discuss how to describe this obstruction problem in terms of Lie $n$-algebras 
  and their corresponding
  categorified Cartan-Ehresmann connections. Generalizations even beyond 
  String-extensions are then straightforward. For $G = \mathrm{Spin}(n)$ the
  next step is  ``Fivebrane structures'' whose existence is obstructed by
  certain generalized Chern-Simons 7-bundles classified by the second
  Pontrjagin class.
\end{abstract}

\thispagestyle{empty}

\newpage

\tableofcontents

\newpage

\section{Introduction}

The study of extended $n$-dimensional relativistic objects which arise in string 
theory has shown that these couple to background fields which can naturally be
thought of as $n$-fold categorified generalizations of fiber bundles with connection.
These structures, or various incarnations of certain special cases of them, 
are probably most commonly known as (bundle-)$(n-1)$-gerbes with connection.
These are known to be equivalently described by Deligne cohomology, by abelian
gerbes with connection (``and curving'') and by Cheeger-Simons
differential characters. Following \cite{Bartels,BS} we address them 
as $n$-bundles with connection.

\begin{table}[h]
  \begin{center}
  \begin{tabular}{c|c}
     \begin{tabular}{c}
       {\bf fundamental}
       \\
       {\bf object}
     \end{tabular}
     &
     \begin{tabular}{c}
       {\bf background}
       \\
       {\bf field}
     \end{tabular}     
     \\
     \hline
     \\
     $n$-particle & $n$-bundle
     \\
     \\
     $(n-1)$-brane & $(n-1)$-gerbe
  \end{tabular}
  \end{center}
  \caption{
    {\bf The two schools of counting} higher dimensional 
    structures. Here $n$ is in $\mathbb{N} = \{0,1,2,\cdots\}$.
  }
\end{table}
 
In string theory, the first departure from bundles with connections 
to higher bundles with connection
occured with the fundamental 
(super)string coupling to the Neveu-Schwarz
(NS) $B$-field. Locally, the $B$-field is just an ${\mathbb R}$-valued two-form. 
However, the study of the path integral, which amounts to `exponentiation', reveals that the 
$B$-field can be thought of as an abelian gerbe with connection whose curving corresponds to the $H$-field 
$H_3$ or as a Cheeger-Simons differential character, whose  holonomy \cite{FW} 
can be described \cite{CJM} in the language of bundle gerbes \cite{Murray}. 

\vspace{3mm}
The next step up occurs with the M-theory (super)membrane which couples
to the $C$-field \cite{BST}. In supergravity, this is viewed locally as an $\mathbb{R}$-valued
differential three-form. However, the study of the path integral has shown that
this field is quantized in a rather nontrivial way \cite{Wi2}. This makes the $C$-field
not precisely a 2-gerbe or degree 3 Cheeger-Simons differential character but 
rather a shifted version \cite{DFM} that can also be modeled using the Hopkins-Singer
description of differential characters \cite{HS}. Some aspects of the description in terms of 
Deligne cohomology is given in \cite{Clingher}. 

\vspace{3mm}
From a purely formal point of view, the need of higher connections for the
description of higher dimensional branes is not a surprise: 
$n$-fold categorified bundles with connection should be precisely
those objects that allow us to define a consistent assignment of  
``phases'' to $n$-dimensional paths in their 
base space. We address such an assignment as {\bf parallel $n$-transport}.
This is in fact essentially the definition of Cheeger-Simons
differential characters \cite{CS} as these are consistent 
assignments of phases to chains. 
However, abelian bundle gerbes, Deligne cohomology and Cheeger-Simons
differential characters all have one major restriction: they only
know about assignments of elements in $U(1)$.

\vspace{3mm}
While the group of phases that enter the
path integral is usually abelian, more general $n$-transport
is important nevertheless. For instance, the latter plays a role at intermediate stages. 
This is well understood for $n=2$: over a $D$-brane the 
abelian bundle gerbe corresponding to
the NS field has the special property that it measures the obstruction to 
lifting a $PU(H)$-bundle to a $U(H)$-bundle, i.e. lifting a bundle with structure
group the infinite projective unitary group on a Hilbert space $H$ to the 
corresponding unitary group  \cite{BM} \cite{BCMMS}.
Hence, while itself an abelian 2-structure, it is crucially related to a nonabelian 
1-structure.  

\vspace{3mm}
That this phenomenon deserves special attention becomes
clear when we move up the dimensional ladder: The Green-Schwarz 
anomaly cancelation \cite{GS} in the heterotic string leads to 
a 3-structure with the special property that, over the
target space, it measures the
obstruction to lifting an $E_8 \times \mathrm{Spin}(n)$-bundle
to a certain nonabelian principal 2-bundle, called a
\emph{String 2-bundle}.  Such a 3-structure is also known as 
a Chern-Simons 2-gerbe \cite{CJMSW}. By itself this is abelian, but its
structure is constrained by certain nonabelian data.
Namely this string 2-bundle with connection, from which the Chern-Simons 
3-bundle arises, is itself an instance of a structure that yields parallel 2-transport. 
It can be described neither by abelian bundle gerbes, nor by Cheeger-Simons 
differential characters, nor by Deligne cohomology.

\vspace{3mm}
In anticipation of such situations, previous works have considered
nonabelian gerbes and nonabelian bundle gerbes with connection. 
However, it turns out that care is needed in order to find the right setup. 
For instance, the kinds of nonabelian gerbes with connection studied in 
\cite{BrM} \cite{AJ}, although very interesting, are not sufficiently general to 
capture String 2-bundles. Moreover, it is not easy to see how to obtain the 
parallel $2$-transport assignment from these structures. For the application to 
string physics, it would be much more suitable to have a nonabelian
generalization of the notion of a Cheeger-Simons differential character,
and thus a structure which, by definition, knows how to assign
generalized phases to $n$-dimensional paths.

\vspace{3mm}
The obvious generalization that is needed is that of a parallel 
transport $n$-functor. Such a notion was described in \cite{BS} \cite{SW}: 
a structure defined by the very fact that it labels $n$-paths by algebraic
objects that allow composition in $n$ different directions, such
that this composition is compatible with the gluing of $n$-paths.
One can show that such transport $n$-functors encompass 
abelian and nonabelian gerbes with connection as special cases
\cite{SW}. However, these $n$-functors are more general. For instance,
String 2-bundles with connection are given by parallel transport
$2$-functors.
 Ironically, the strength of the latter --  namely their knowledge about general phase 
assignments to higher dimensional paths -- is to some degree also 
a drawback: for many computations, a description
\emph{entirely} in terms of differential form data would be more tractable.
However, the passage from parallel $n$-transport to the corresponding differential
structure is more or less straightforward: a parallel transport $n$-functor
is essentially a morphism of Lie $n$-groupoids. As such, it can be sent,
by a procedure generalizing the passage from Lie groups to Lie algebras,
to a morphism of Lie $n$-algebroids.

\vspace{3mm}
The aim of this paper is to describe two topics: First, to set up a formalism for higher 
bundles with connections entirely in terms of $L_{\infty}$-algebras,
which may be thought of as a categorification of the theory of 
Cartan-Ehresmann connections.  This is supposed to be
the differential version of the theory of parallel transport $n$-functors,
but an exhaustive discussion of the differentiation procedure 
is not given here. Instead we discuss a couple
of examples and then show how the lifting problem has a nice description in
this language.
To do so, we present a family of $L_{\infty}$-algebras that 
govern the gauge structure of $p$-branes, as above, and discuss the lifting problem for them.
By doing so, we characterize Chern-Simons 3-forms as local connection data
on 3-bundles with connection which arise as the obstruction to lifts of
ordinary bundles to the corresponding String 2-bundles, governed by the
String Lie 2-algebra.

\vspace{3mm}
The formalism immediately allows the generalization of this situation to 
higher degrees. Indeed we indicate how certain 7-dimensional generalizations of 
Chern-Simons 3-bundles obstruct the lift of ordinary bundles to certain 6-bundles
governed by the Fivebrane Lie 6-algebra. The latter correspond to what we define 
as the fivebrane structure, for which the degree seven NS field $H_7$ plays the role
that the degree three dual NS field $H_3$ plays for the $n=2$ case. 

\vspace{3mm}
The paper is organized in such a way that section \ref{plan} serves more or less as 
a self-contained description of the basic ideas and construction, 
with the rest of the document having all the details and all the proofs.

\vspace{3mm}
In this paper we make use of the homotopy algebras usually referred to as
 $L_\infty$-algebras. These algebras also go by other names
 such as sh-Lie algebras \cite{lada-jds}. In our context we may also call such algebras 
 Lie $\infty$-algebras which we think of as the abstract concept of an 
 $\infty$-vector space with an antisymmetric  bracket 
$\infty$-functor on it, which satisfies a Jacobi identity up to coherent equivalence,  
whereas ``$L_\infty$-algebra'' is concretely a 
codifferential coalgebra of sorts. In this paper we will nevertheless follow 
the standard notation of $L_{\infty}$-algebra.

\section{The Setting and Plan}
\label{plan}
 
 We set up a useful framework for describing higher order bundles with connection
 entirely in terms of Lie $n$-algebras, which can be thought of as arising
 from a categorification of the concept of an Ehresmann connection on a principal
 bundle.
 Then we apply this to the study of Chern-Simons $n$-bundles with connection as
 obstructions to lifts of principal $G$-bundles through higher String-like
 extensions of their structure Lie algebra.

 \subsection{$L_{\infty}$-algebras and their String-like central extensions}

   A Lie group has all the right properties to  locally describe the
   phase change of a charged particle as it traces out a worldline.
   A Lie $n$-group is a higher structure with precisely all the right
   properties to  describe locally the phase change of a charged 
   $(n-1)$-brane  as it traces out an $n$-dimensional
   worldvolume. 
   
  \subsubsection{$L_{\infty}$-algebras}   

   \label{Linfty in plan}
   
   Just as ordinary Lie groups have Lie algebras,
   Lie $n$-groups have Lie $n$-algebras. If the Lie $n$-algebra is what
   is called \emph{semistrict}, these are \cite{BaC} precisely 
   $L_\infty$-algebras
   \cite{lada-jds} which have come to play a significant role in cohomological physics.
   A (``semistrict'' and finite dimensional) 
   Lie $n$-algebra is any of the following three equivalent structures:
   \begin{itemize}
     \item  an $L_\infty$-algebra structure on a graded vector space $\gg$ concentrated
       in the first $n$ degrees $(0,...,n-1)$;
     \item
       a quasi-free differential graded-commutative algebra 
       (``qDGCA'': free as a graded-commutative)
       algebra on the dual of that vector space: this is the Chevalley-Eilenberg algebra
       $\mathrm{CE}(\gg)$ of $\gg$;
     \item
       an $n$-category internal to the category of {\it graded} vector spaces and equipped with
       a skew-symmetric linear bracket functor which satisfies a Jacobi identity
       up to higher coherent equivalence.
   \end{itemize}

  For every  $L_\infty$-algebra $\gg$, we have the following three qDGCAs:
  \begin{itemize}
    \item
      the {\bf Chevalley-Eilenberg algebra} $\mathrm{CE}(\gg)$
    \item
      the {\bf Weil algebra} $\mathrm{W}(\gg)$
    \item
      the algebra of {\bf invariant polynomials} or {\bf basic forms}
      $\mathrm{inv}(\gg)$.
  \end{itemize}

 These sit in a sequence
  \(
    \xymatrix{
      \mathrm{CE}(\gg)
      &&
      \mathrm{W}(\gg)
      \ar@{->>}[ll]
      &&
      \mathrm{inv}(\gg)
      \ar@{_{(}->}[ll]
    }
    \,,
  \)
  where all morphisms are morphisms of dg-algebras. This sequence plays the role 
  of the sequence of differential forms on the
  ``universal $\gg$-bundle''.

\begin{figure}[h]
 $$
  \hspace{-1.8cm}
  \xymatrix@C=0pt@R=7pt{
    &
    \mbox{\begin{tabular}{c}{\bf dg-algebras}\end{tabular}}
    &
    \ar@{-}[dddddd]    
    &
    &
    \mbox{\begin{tabular}{c}{\bf  $L_\infty$-algebras}\end{tabular}}
    &
    \ar@{-}[dddddd]
    &
    &
    \mbox{{\bf groupoids}}
    &
    \ar@{-}[dddddd]
    &&
    \mbox{\begin{tabular}{c}{\bf (pointed) }\\{\bf topological }\\ {\bf spaces}\end{tabular}}
    \\
    \ar@{-}[rrrrrrrrrrr]
    &&&&&&&&&&&&
    \\
    \mbox{\begin{tabular}{c} 
      Chevalley-
     \\Eilenberg\\ algebra 
    \end{tabular}}
    &
    \mathrm{CE}(\gg) 
    &&
    \mbox{Lie algebra}
    &
    \gg
    \ar@{^{(}->}[dd]
    &
    &
    \mbox{\begin{tabular}{c}structure \\group as \\ monoidal \\ set \end{tabular}}
    &
     G
     \ar@{^{(}->}[dd]
    && 
    \mbox{
      \begin{tabular}{c}
        structure \\ group
      \end{tabular}
    }
    &
    G \ar@{^{(}->}[dd]^i
    \\
    \\
    \mbox{Weil algebra}
    &\mathrm{W}(\gg) 
    \ar@{->>}[uu]_{i^*}
    &&
     \mbox{
       \begin{tabular}{c}
         inner \\ derivation \\
         Lie \\ 2-algebra
       \end{tabular}
     }
    &
    \mathrm{inn}(\gg)
    &&
    \mbox{\begin{tabular}{c}inner \\ automorphism \\ 2-group as \\ monoidal \\ groupoid \end{tabular}}
    &
    \mathrm{INN}(G)
    \ar@{->>}[dd]
    &
    & 
      \mbox{\begin{tabular}{c}universal \\ $G$-bundle\end{tabular}}
      &
    EG \ar@{->>}[dd]^p
    \\
    \\
    \mbox{\begin{tabular}{c}
       algebra of \\ invariant \\ polynomials
     \end{tabular}}
    &\mathrm{inv}(\gg)
    \ar@{^{(}->}[uu]_{p^*}
    &&
    &&
    &\mbox{\begin{tabular}{c} structure group\\ as one-object \\ groupoid \end{tabular}}
    &
     \mathbf{B} G
    &&
    \mbox{\begin{tabular}{c}
      classifying \\ space
      \\ for $G$
    \end{tabular}}
    &
    BG
    \\
    &&&&
    \ar@{|->}[lll]|{(\cdot)^*}^{\mbox{\tiny regard as codifferential coalgebra and dualize }}
    &&&
    \ar@{|->}[lll]|{\mathrm{Lie}}^{\mbox{\tiny differentiate}}  
    \ar@{|->}[rrr]|{|\cdot|}_{\mbox{\tiny take nerves and geometrically realize}}&&&&&&
  }
  $$
  \caption{
    {\bf The universal $G$-bundle} in its various incarnations.
    That the ordinary universal $G$ bundle is the realization of the nerve of the groupoid
    which we denote here by $\mathrm{INN}(G)$ is an old result by Segal (see 
    \cite{RS} for a review and a discussion of the situation for 2-bundles). 
    This groupoid $\mathrm{INN}(G)$ is in fact a 2-group. The corresponding Lie 2-algebra
    (2-term $L_\infty$-algebra) we denote by $\mathrm{inn}(\gg)$. Regarding this
    as a codifferential coalgebra and then dualizing that to a differential algebra
    yields the Weil algebra of the Lie algebra $\gg$. This plays the role of
    differential forms on the universal $G$-bundle, as already known to Cartan.
    The entire table is expected to admit an $\infty$-ization. Here we concentrate
    on discussing $\infty$-bundles with connection in terms just of $L_\infty$-algebras
    and their dual dg-algebras. An integration of this back to the world of 
    $\infty$-groupoids should proceed along the lines of 
    \cite{Getzler,Henriques}, but is
    not considered here.
    \label{universal G-bundle}
  }
\end{figure}
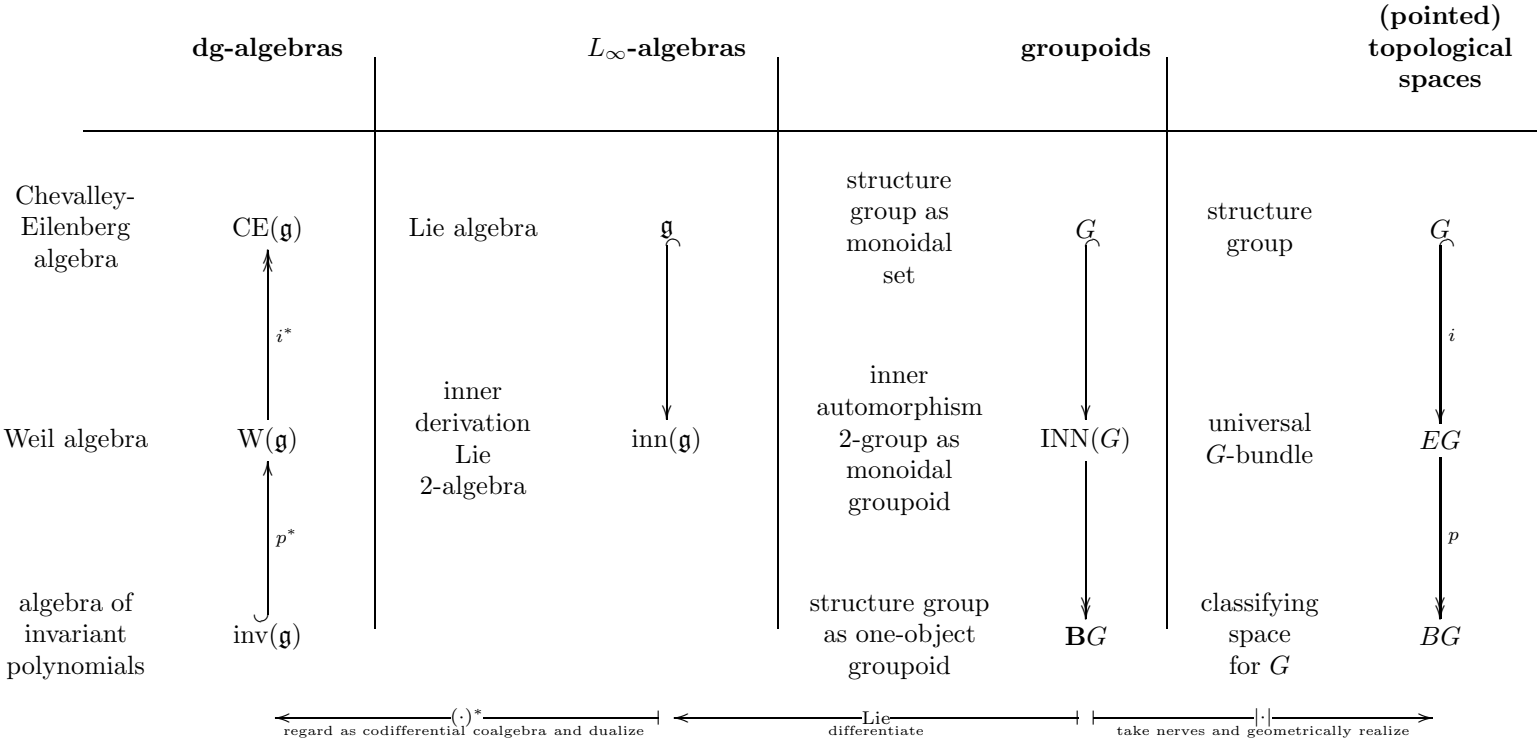

  \subsubsection{$L_\infty$-algebras from cocycles: String-like extensions}

   A simple but important source of examples for higher Lie $n$-algebras
comes from the abelian Lie algebra $\uu(1)$ which may be shifted into higher
categorical degrees. We write $b^{n-1}\uu(1)$ for the Lie $n$-algebra which 
is entirely trivial except in its $n$th degree, where it looks like $\uu(1)$.
Just as $\uu(1)$ corresponds to the Lie group $U(1)$ , so $b^{n-1}\uu(1)$ corresponds 
to the iterated classifying space $B^{n-1}U(1)$, realizable as the topological group 
given by the Eilenberg-MacLane space
$K({\mathbb{Z}},n)$.
Thus an important source for interesting Lie $n$-algebras comes from
extensions
\(
  0 \to b^{n-1}\uu(1) \to \hat \gg \to \gg \to 0
\)
of an ordinary Lie algebra $\gg$ by such a shifted abelian Lie $n$-algebra
$b^{n-1}\uu(1)$.
We find that, for each $(n+1)$-cocycle $\mu$ in the Lie algebra cohomology of $\gg$,
we do obtain such a central extension, which we describe by
\(
  0 \to b^{n-1}\uu(1) \to \gg_\mu \to \gg \to 0
  \,.
\)
Since, for the case when $\mu = \langle \cdot, [\cdot, \cdot]\rangle$ 
is the canonical 3-cocycle on a semisimple Lie algebra $\gg$, this $\gg_\mu$
is known (\cite{BCSS} and \cite{Henriques}) 
to be the Lie 2-algebra of the String 2-group, we call these central
extensions \emph{String-like} central extensions. (We also refer to these
as Lie $n$-algebras ``of Baez-Crans type'' \cite{BaC}.)
Moreover, whenever the cocycle $\mu$ is related by transgression to an
invariant polynomial $P$ on the Lie algebra, we find that $\gg_\mu$
fits into a short \emph{homotopy} exact sequence of Lie $(n+1)$-algebras
\(
  0 \to \gg_\mu \to \mathrm{cs}_P(\mu) \to \mathrm{ch}_P(\mu) \to 0
  \,.
\)
Here $\mathrm{cs}_P(\gg)$ is a Lie $(n+1)$-algebra governed by the 
Chern-Simons term corresponding to the  
transgression element interpolating between $\mu$ and $P$. In a similar
fashion $\mathrm{ch}_P(\gg)$ knows about the characteristic (Chern) class 
associated with $P$.

\vspace{3mm}
In summary, from elements of the cohomology of $\mathrm{CE}(\gg)$
together with related elements in $\mathrm{W}(\gg)$ we obtain the String-like extensions of Lie algebras 
to Lie $2n$-algebras and the associated Chern- and Chern-Simons Lie $(2n-1)$-algebras:

        
 \vspace{2mm}
   \begin{tabular}{cc|cc}
       \hline
       Lie algebra cocycle & $\mu$ & Baez-Crans Lie $n$-algebra & $\gg_{\mu}$
       \\
       invariant polynomial & $P$ & Chern Lie $n$-algebra 
          & $\mathrm{ch}_P(\gg)$
       \\
       transgression element & $\mathrm{cs}$ & Chern-Simons Lie $n$-algebra
         & $\mathrm{cs}_P(\gg)$
         \\
         \hline
     \end{tabular}
    \medskip     


\subsubsection{ $L_\infty$-algebra differential forms}

\label{Lie oo-algebra valued forms in plan}

For $\gg$ an ordinary Lie algebra and $Y$ some manifold, 
one finds that dg-algebra morphisms $\mathrm{CE}(\gg) \to \Omega^\bullet(Y)$
from the Chevally-Eilenberg algebra of $\gg$ to the DGCA of differential forms
on $Y$
   are in bijection with $\gg$-valued 1-forms $A \in \Omega^1(Y,\gg)$ 
   whose ordinary curvature 2-form 
   \(
     F_A
     = 
     d A + [A\wedge A]
   \) 
    vanishes. Without the flatness, the correspondence is with 
   algebra morphisms \emph{not} respecting the differentials. 
   But dg-algebra morphisms $A: \mathrm{W}(\gg) \to \Omega^\bullet(Y)$
   are in bijection with arbitrary $\gg$-valued 1-forms. 
   These are flat precisely if $A$ factors through $\mathrm{CE}(\gg)$. 
    This situation is depicted in the following diagram:
   \(
    \raisebox{30pt}{
     \xymatrix{
       \mathrm{CE}(\gg) 
       \ar@{..>}[d]_{(A,F_A = 0)}
       && 
       \mathrm{W}(\gg) \ar@{->>}[ll]
       \ar[d]_{(A,F_A)}
       \\
       \Omega^\bullet(Y)
       \ar@{-}[rr]^=
       &&
       \Omega^\bullet(Y)
     }
     }
     \,.
   \)
 This has an obvious generalization for $\gg$ an arbitrary $L_\infty$-algebra.
 For $\gg$ any $L_\infty$-algebra, we write
 \(
   \Omega^\bullet(Y,\gg) = \mathrm{Hom}_{\mathrm{dg-Alg}}(\mathrm{W}(\gg), \Omega^\bullet(X))
 \)
 for the collection of {\bf $\gg$-valued differential forms} and
 \(
   \Omega^\bullet_{\mathrm{flat}}(Y,\gg) 
    = \mathrm{Hom}_{\mathrm{dg-Alg}}(\mathrm{CE}(\gg), \Omega^\bullet(X))
 \)
 for the collection of {\bf flat $\gg$-valued differential forms}.

\subsection{ $L_\infty$-algebra Cartan-Ehresmann connections}

\subsubsection{$\gg$-Bundle Descent data}

\label{principal n-bundles in plan}

 A {\it descent object} for an ordinary  principal $G$-bundle on $X$ is a surjective submersion
 $\pi : Y \to X$ 
 together with a functor $g : Y \times_X Y \to \mathbf{B}G$
 from the groupoid whose morphisms are pairs of points in the
 same fiber of $Y$, to the groupoid 
 $\mathbf{B}G$ which is the one-object groupoid 
 corresponding to the group $G$.  Notice that the groupoid $\mathbf{B}G$ is not
 itself the classifying space $B G$ of $G$, but the geometric realization of its nerve,
 $|\mathbf{B}G|$, is: $|\mathbf{B}G | = BG$.

\vspace{3mm}
 We may take $Y$ to be the disjoint union of some open subsets $\{U_i\}$ of $X$ that 
 form a good open cover of $X$. Then $g$ is the familiar concept of a transition
 function decribing a bundle that has been locally trivialized over the 
 $U_i$.
 But one can also use more general surjective submersions. For instance,
 for $P \to X$ any principal $G$-bundle, it is sometimes useful to take
 $Y = P$. In this case one obtains a canonical choice for the cocycle
 \(
   g : Y \times_X Y = P \times_X P \to \mathbf{B}G
 \)
 since $P$ being principal means that
 \(
   P \times_X P \simeq_{\mathrm{diffeo}} P \times G
   \,.
 \)
 This reflects the fact that every principal bundle canonically trivializes
 when pulled back to its own total space.
 The choice $Y = P$ differs from that of a good cover crucially in the following
 aspect: if the group $G$ is connected, then also the fibers of $Y = P$ are connected.
  Cocycles over surjective submersions with connected fibers have special properties,
 which we will utilize:
{\it  When the fibers of $Y$ are connected, we may 
 think of the assignment of group elements to pairs of points in one fiber as 
 arising from the parallel transport with respect to a flat vertical 
 1-form $A_{\mathrm{vert}} \in \Omega^1_{\mathrm{vert}} (Y,\gg)$, flat along the fibers.} 
 As we shall see, this can be thought of as the vertical part of a 
 Cartan-Ehresmann connection 1-form. 
 This provides a morphism
 \(
   \xymatrix{
     \Omega^\bullet_{\mathrm{vert}}(Y)
     &
     \mathrm{CE}(\gg)
     \ar[l]_{A_{\mathrm{vert}}}
   }
 \)
 of differential graded algebras 
 from the Chevalley-Eilenberg algebra of $\gg$ to the 
 vertical differential forms on $Y$.
 
 \vspace{3mm}
Unless otherwise specified, \emph{morphism} will always mean 
\emph{homomorphism of differential graded algebra}.
$A_{\mathrm{vert}}$ has an obvious generalization: for $\gg$ any Lie $n$-algebra, 
 we say that a $\gg$-bundle descent object 
 for a $\gg$-$n$-bundle on $X$ is a surjective submersion
 $\pi : Y \to X$ together with a morphism
 $
   \xymatrix{
     \Omega^\bullet_{\mathrm{vert}}(Y)
     &
     \mathrm{CE}(\gg)
     \ar[l]_{A_{\mathrm{vert}}}
   }
   \,.
 $
 Now $A_{\mathrm{vert}} \in \Omega^\bullet_{\mathrm{vert}}(Y,\gg)$ 
 encodes a collection of vertical $p$-forms on $Y$,
 each taking values in the degree $p$-part of $\gg$ and all together
 satisfying a certain flatness condition, controlled by the nature
 of the differential on $\mathrm{CE}(\gg)$.

\subsubsection{Connections on $n$-bundles: the extension problem}

\label{connection in plan}

  Given a descent object 
  $\xymatrix{
       \Omega^\bullet_{\mathrm{vert}}(Y)
       &&
       \mathrm{CE}(\gg)
       \ar[ll]_{A_{\mathrm{vert}}}
    }$
  as above, a {\bf flat connection} on it is an extension of
  the morphism $A_{\mathrm{vert}}$ to a morphism $A_{\mathrm{flat}}$
  that factors through differential forms on $Y$
  \(
    \xymatrix{
       \Omega^\bullet_{\mathrm{vert}}(Y)
       &&
       \mathrm{CE}(\gg)
       \ar[ll]_{A_{\mathrm{vert}}}
       \ar@{..>}[ddll]^{A_{\mathrm{flat}}}
       \\
       \\
       \Omega^\bullet(Y)
       \ar@{->>}[uu]_{i^*}
    }\,.
  \)
  In general, such an extension does not exist. 
  A general {\bf connection} on a $\gg$-descent object $A_{\mathrm{vert}}$
  is a morphism
  \(
    \xymatrix{
      \Omega^\bullet(Y)
      &&
      \mathrm{W}(\gg)
      \ar[ll]_{(A,F_A)}
    }
  \)
  from the Weil algebra of $\gg$ to the differential forms on $Y$
  together with a morphism
  \(
    \xymatrix{
      \Omega^\bullet(Y)
      &&
      \mathrm{inv}(\gg)
      \ar[ll]_{\{K_i\}}
    }
  \)
  from the invariant polynomials on $\gg$, as in \ref{Linfty in plan}, 
  to the differential forms on $X$, such that
  the following two squares commute:
  \(
   \raisebox{40pt}{
    \xymatrix{
       \Omega^\bullet_{\mathrm{vert}}(Y)
       &&
       \mathrm{CE}(\gg)
       \ar[ll]_{A_{\mathrm{vert}}}
       \\
       \\
       \Omega^\bullet(Y)
       \ar@{->>}[uu]_{i^*}
       &&
       W(\gg)
       \ar@{->>}[uu]
       \ar[ll]_{(A,F_A)}^<{\ }="s"
       \\
       \\
       \Omega^\bullet(X)
       \ar@{^{(}->}[uu]_{\pi^*}
       &&
       \mathrm{inv}(\gg)
       \ar@{^{(}->}[uu]
       \ar[ll]^{\{K_i\}}_>{\ }="t"
       %
    }
    }
    \,.
  \)

  Whenever we have such two commuting squares, we say
  \begin{itemize}
    \item
      $A_{\mathrm{vert}} \in \Omega^\bullet_{\mathrm{vert}}(Y,\gg)$
      is a $\gg$-bundle {\bf descent object} (playing the role of a 
      {\bf transition function});
    \item
      $A \in \Omega^\bullet(Y,\gg)$
      is a (Cartan-Ehresmann) {\bf  connection} with values in the $L_\infty$-algebra
      $\gg$ on the total space of the 
      surjective submersion;
    \item
      $F_A \in \Omega^{\bullet + 1}(Y,\gg)$
      are the corresponding {\bf curvature} forms;
    \item
      and the set $\{K_i \in \Omega^\bullet(X)\}$ are the
      corresponding {\bf characteristic forms}, whose classes $\{[K_i]\}$
      in deRham cohomology 
      \(
        \xymatrix{
          \Omega^\bullet(X)
          &&
          \mathrm{inv}(\gg)
          \ar[ll]_{\{K_i\}}
          \\
          H_{\mathrm{deRham}}^\bullet(X)
          &&
          H^\bullet(\mathrm{inv}(\gg))
          \ar[ll]_{\{[K_i]\}}          
        }
      \)
      are the corresponding {\bf characteristic classes}
      of the given descent object $A_{\mathrm{vert}}$.
  \end{itemize}
  
\begin{figure}[h]
  \[
    \xymatrix@R=6pt{
       \Omega^\bullet_{\mathrm{vert}}(Y)
       &&
       \mathrm{CE}(\gg)
       \ar[ll]_{A_{\mathrm{vert}}}
       &&
       \mbox{ \bf
         \begin{tabular}{c}
           descent \\ data
         \end{tabular}
       }
       \\
       &&&&
       \mbox{
         \begin{tabular}{c}
           first \\ Cartan-Ehresmann \\ condition
         \end{tabular}
       }
       \\
       \Omega^\bullet(Y)
       \ar@{->>}[uu]_{i^*}
       &&
       W(\gg)
       \ar@{->>}[uu]
       \ar[ll]_{(A,F_A)}^<{\ }="s"
       &&
       \mbox{
         \bf
         \begin{tabular}{c}
           connection \\ data
         \end{tabular}
       }
       \\
       &&&&
       \mbox{
         \begin{tabular}{c}
           second \\ Cartan-Ehresmann \\ condition
         \end{tabular}
       }
       \\
       \Omega^\bullet(X)
       \ar@{^{(}->}[uu]_{\pi^*}
       &&
       \mathrm{inv}(\gg)
       \ar@{^{(}->}[uu]
       \ar[ll]|{\{K_i\}}_>{\ }="t"
       &&
       \mbox{
         \bf 
         \begin{tabular}{c}
           characteristic \\ forms
         \end{tabular}
       }       
       \\
       H^\bullet_{\mathrm{dR}}(X)
       &&
       H^\bullet(\mathrm{inv}(\gg))
       \ar[ll]^{\{[K_i]\}}
       &&
       \mbox{
         \begin{tabular}{c}
           Chern-Weil \\ homomorphism
         \end{tabular}
       }
       %
    }
  \]
  \caption{
    {\bf A $\gg$-connection descent object and its interpretation.}
    For $\gg$-any $L_\infty$-algebra and $X$ a smooth space, a $\gg$-connection
    on $X$ is an equivalence class of pairs $(Y,(A,F_A))$ consisting of 
    a surjective submersion $\pi: Y \to X$ and dg-algebra morphisms forming
    the above commuting diagram. The equivalence relation is concordance of such
    diagrams. The situation for ordinary Cartan-Ehresmann (1-)connections is
    described in \ref{examples for connection descent objects}.
  }
\end{figure}
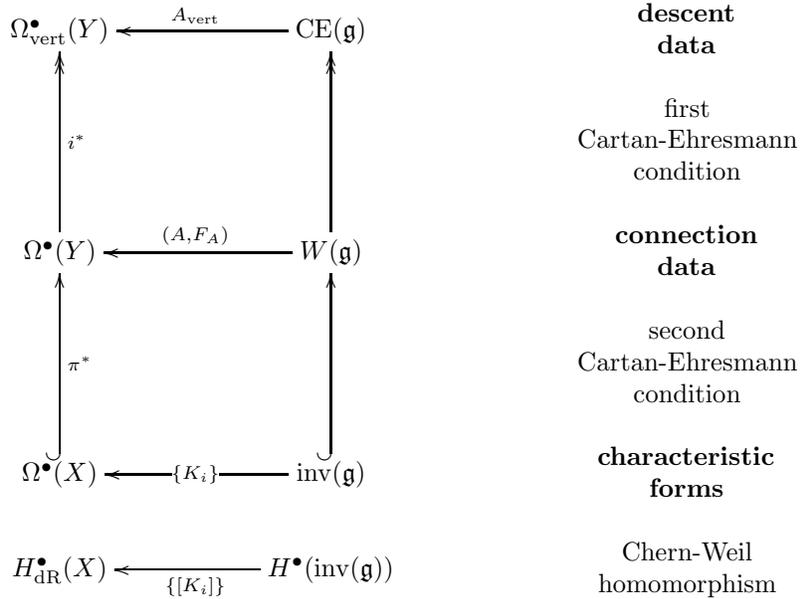

Hence we realize the curvature of a $\gg$-connection as the \emph{obstruction} to
extending a $\gg$-descent object to a \emph{flat} $\gg$-connection.

\subsection{Higher String and Chern-Simons $n$-transport: the lifting problem}
\label{plan: the lifting problem}

 Given a $\gg$-descent object
 \(
   \xymatrix{
     &
     \mathrm{CE}(\gg)
     \ar[dl]^{A_{\mathrm{vert}}}
     \\
     \Omega^\bullet_{\mathrm{vert}}(Y)
   },
 \)
 and given an extension of $\gg$ by a String-like $L_\infty$-algebra
 \(
   \xymatrix{
     \mathrm{CE}(b^{n-1}\uu(1))
     &&
     \mathrm{CE}(\gg_\mu)
     \ar@{->>}[ll]_i
     &&
     \mathrm{CE}(\gg)     
     \ar@{_{(}->}[ll]
   },
 \)
 we ask if it is possible to \emph{lift the descent object} through this 
 extension, i.e. to find a dotted arrow in
 \(
   \xymatrix@C=8pt{
     \mathrm{CE}(b^{n-1}\uu(1))
     &&
     \mathrm{CE}(\gg_\mu)
     \ar@{->>}[ll]
     \ar@{..>}[dr]
     &&
     \mathrm{CE}(\gg)     
     \ar@{_{(}->}[ll]
      \ar[dl]^{A_{\mathrm{vert}}}
     \\
     &&&\Omega^\bullet_{\mathrm{vert}}(Y)
  }
  \,.
 \)
 In general this is not possible. We seek a straightforward way to 
 compute the obstruction to the existence of the lift.
The strategy is to form the \emph{weak} (homotopy) kernel of 
 \(
   \xymatrix@C=8pt{
     \mathrm{CE}(b^{n-1}\uu(1))
     &&
     \mathrm{CE}(\gg_\mu)
     \ar@{->>}[ll]_i
   }
 \)
 which we denote by 
  $\mathrm{CE}(b^{n-1}\uu(1) \hookrightarrow \gg_\mu)$
  and realize as a mapping cone of qDGCAs. 
This comes canonically
  with a morphism $f$ from $\mathrm{CE}(\gg)$ which happens to 
  have a \emph{weak} inverse 
 \(
   \xymatrix@C=8pt{
   &&&&&&
   \mathrm{CE}(b^{n-1}\uu(1) \hookrightarrow \gg_\mu)
   \ar[dllll]
   \ar@/^1pc/[dll]|{f^{-1}}
   \\
     \mathrm{CE}(b^{n-1}\uu(1))
     &&
     \mathrm{CE}(\gg_\mu)
     \ar@{->>}[ll]_i
     &&
     \mathrm{CE}(\gg)     
     \ar@{_{(}->}[ll]
     \ar[urr]|f
   }
   \,.
 \)
 Then we see that, while the lift to a $\gg_\mu$-cocycle may not always exist, the
 lift to a $(b^{n-1}\uu(1)\hookrightarrow \gg_\mu)$-cocycle 
 does always exist. We form $A_{\mathrm{vert}}\circ f^{-1}$:
 \(
   \xymatrix@C=8pt{
   &&&&&&
   \mathrm{CE}(b^{n-1}\uu(1) \hookrightarrow \gg_\mu)
   \ar[dllll]
   \ar[dll]|{f^{-1}}
   \\
     \mathrm{CE}(b^{n-1}\uu(1))
     &&
     \mathrm{CE}(\gg_\mu)
     \ar@{->>}[ll]_i
     \ar@{..>}[dr]
     &&
     \mathrm{CE}(\gg)     
     \ar@{_{(}->}[ll]
      \ar[dl]^{A_{\mathrm{vert}}}
     \\
     &&&\Omega^\bullet_{\mathrm{vert}}(Y)
   }
   \,.
 \)
  The failure of this lift to be a true lift to $\gg_\mu$ is measured by the component
  of $A_{\mathrm{vert}}\circ f^{-1}$ on $b^{n-1}\uu(1)[1] \simeq b^n\uu(1)$.
  Formally this is the composite $A'_{\mathrm{vert}} := A_{\mathrm{vert}} \circ f^{-1} \circ j$
  in
 \(
   \xymatrix@C=8pt{
   &&&&&&
   \mathrm{CE}(b^{n-1}\uu(1) \hookrightarrow \gg_\mu)
   \ar[dllll]
   \ar[dll]|{f^{-1}}
   &&
   \mathrm{CE}(b^n \uu(1))
   \ar@{_{(}->}[ll]_<<<j
   \ar[ddlllll]|{A'_{\mathrm{vert}}}
   \\
     \mathrm{CE}(b^{n-1}\uu(1))
     &&
     \mathrm{CE}(\gg_\mu)
     \ar@{->>}[ll]_i
     \ar@{..>}[dr]
     &&
     \mathrm{CE}(\gg)     
     \ar@{_{(}->}[ll]
      \ar[dl]|{A_{\mathrm{vert}}}
     \\
     &&&\Omega^\bullet_{\mathrm{vert}}(Y)
   }
   \,.
 \)
 The nontriviality of the $b^n\uu(1)$-descent object $A'_{\mathrm{vert}}$
 is the obstruction to constructing the desired lift.

\vspace{3mm}

   We thus find the following results, for any $\gg$-cocycle $\mu$
   which is in transgression with the the invariant polynomial $P$
   on $\gg$,
      
   \begin{itemize}
     \item
       The characteristic classes (in deRham cohomology) of $\gg_\mu$-bundles
       are those of the corresponding $\gg$-bundles modulo
       those coming from the invariant polynomial $P$.
     \item
       The lift of a $\gg$-valued connection to a $\gg_\mu$-valued connection
       is obstructed by a $b^n\uu(1)$-valued $(n+1)$-connection whose 
       $(n+1)$-form curvature is $P(F_A)$, i.e. the image under the Chern-Weil
       homomorphism of the invariant polynomial corresponding to $\mu$.
     \item
        Accordingly, the $(n+1)$-form connection of the obstructing 
        $b^n\uu(1)$ $(n+1)$-bundle is a Chern-Simons form for this characteristic class.
   \end{itemize}
   
   We call the obstructing $b^n \uu(1)$ $(n+1)$-descent object 
   the corresponding Chern-Simons $(n+1)$-bundle. For the case when
   $\mu = \langle \cdot, [\cdot,\cdot]\rangle$ is the canonical 3-cocycle
   on a semisimple Lie algebra $\gg$, this structure (corresponding to a 2-gerbe)
   has a  3-connection given by the ordinary Chern-Simons 3-form and has 
   a curvature 4-form given by the (image in deRham cohomology of) the
   first Pontrjagin class of the underlying $\gg$-bundle.

\section{Physical applications: String-, Fivebrane- and $p$-Brane structures}

We can now discuss physical applications of the formalism that
we have developed. What we describe is a useful way to handle obstructing $n$-bundles
of various kinds that appear in string theory. In particular, we can
describe generalizations of string structure in string theory. In the context 
of $p$-branes, such generalizations have been suggested based on 
$p$-loop spaces  \cite{DDS} \cite{BPSST} \cite{PS} and, more generally,
on the space of maps $\mathrm{Map(M,X)}$ from the brane worldvolume 
$M$ to spacetime $X$ \cite{MP}. The statements
in this section will be established in detail in \cite{SSS2}.

\vspace{3mm}
From the point of view of supergravity, all branes, called $p$-branes in that
setting, are a priori treated in a unified way. In tracing back to string theory, 
however, there is a distinction in the form-fields between the 
Ramond-Ramond (RR) and the Neveu-Schwarz (NS) forms. The former 
live in generalized cohomology and the latter play two roles: they
act as twist fields for the RR fields and they are also connected to the 
geometry and 
topology of spacetime. The $H$-field $H_3$ plays the role of a twist in K-theory 
for the RR fields \cite{Ka} \cite{BM} \cite{MSa}. The twist
for the degree seven dual field $H_7$ is observed in \cite{H7}
at the rational level.  

\vspace{3mm}
The ability to define fields and their corresponding partition functions 
puts constraints on the topology of the underlying spacetime. The 
most commonly understood example is that of fermions where the 
ability to define
them requires spacetime to be spin, and the ability to describe 
theories with chiral fermions requires certain restrictions coming 
from the index theorem. In the context of heterotic string theory,
the Green-Schwarz anomaly cancelation leads to the condition 
that the difference between the Pontrjagin classes of the tangent 
bundle and that of the gauge bundle be zero. This is called the
string structure, which can be thought of as a spin structure on 
the loop space of spacetime \cite{Kill} \cite{CP}. In M-theory, the ability to 
define the partition function leads to an anomaly 
given by the integral seventh-integral Steifel-Whitney class of spacetime
\cite{DMW} whose cancelation requires spacetime to be
orientable with respect to generalized cohomology theories beyond
K-theory  \cite{KS1} . 

\vspace{3mm}
In all cases, the corresponding structure is 
related to the homotopy groups of the orthogonal group: the spin 
structure amounts to killing the first homotopy group, the string structure
and -- to some extent-- the $W_7$ condition to killing the third homotopy 
group. Note that 
when we say that the $n$-th homotopy group is killed, we really mean
that all homotopy groups up to and including the $n$-th one are killed.
For instance,  a String structure requires killing everything up to and
including the third, hence everything through the sixth, since there are no
homotopy groups in degrees four, five or six. 

\vspace{3mm}
The Green-Schwarz anomaly cancelation condition for the heterotic 
string can be translated to the language of $n$-bundles as follows. 
We have two bundles, the spin bundle with structure group 
$G={\rm Spin}(10)$, and the gauge bundle with structure group 
$G'$ being either ${\rm SO}(32)/{\mathbb{Z}_2}$ or $E_8 \times E_8$.
Considering the latter, we have one copy of $E_8$ on each 
ten-dimensional boundary component, which can be viewed as
an end-of-the-world nine-brane, or M9-brane \cite{HW}. The structure of the 
four-form on the boundary 
which we write as 
\( 
G_4|_{\partial}  = d H_3
\)
implies that the 3-bundle (2-gerbe) becomes the trivializable 
lifting 2-gerbe of a ${\rm String}({\rm Spin}(10) \times E_8)$ bundle 
over the M9-brane. As the four-form contains the difference of the
Pontrjagin classes of the bundles with structure groups $G$ and
$G'$, the corresponding three-form will be a difference of 
Chern-Simons forms. The bundle aspect of this has been
studied in \cite{BCRS} and will be revisited in the current context
in \cite{SSS2}.

\vspace{3mm}
The NS fields play a special role in relation to the homtopy 
groups of the orthogonal group. The degree three class $[H_3]$ 
plays the role of a twist for a spin structure. Likewise, the 
degree seven class plays a role of a twist for a higher structure
related to $BO\langle 9\rangle$, the 8-connected cover of $BO$,
which we might call a 
{\it Fivebrane}-structure on spacetime. We can talk about such a structure 
once the spacetime already has a string structure. The obstructions 
are given in table \ref{ext objs and top structures}, where $A$ is the connection on 
the $G'$ bundle and $\omega$ is a connection on the $G$ bundle.

\begin{table}[h]
\begin{center}
  \begin{tabular}{c|cccc}
    $n$
    &
    \begin{tabular}{c}
      2
      \\
      $=4\cdot 0 + 2$
    \end{tabular}
    &
    \begin{tabular}{c}
      6
      \\
      $=4\cdot 1 + 2$
    \end{tabular}
    &
    \\
    \hline
    \begin{tabular}{c}
       fundamental object
       \\
       ($n-1$)-brane
       \\
       $n$-particle
    \end{tabular}
    &
    string
    &
    5-brane 
    \\
    \hline
    \begin{tabular}{l}
      target space
      \\
      structure
    \end{tabular}
    &
    \begin{tabular}{c}
      string structure
      \\
      $\mathrm{ch}_2(A) - p_1(\omega) = 0$
    \end{tabular}
    &
    \begin{tabular}{c}
      fivebrane structure
      \\
      $\mathrm{ch}_4(A) - \frac{1}{48}p_2(\omega) = 0$
    \end{tabular}
  \end{tabular}
\end{center}
  \caption{
    \label{ext objs and top structures}
    {\bf Higher dimensional extended objects} and the
    corresponding topological structures.
  }
\end{table}

\vspace{3mm}
In the above we alluded to how the brane structures are related to obstructions 
to having spacetimes with connected covers of the orthogonal groups as 
structures. The obstructing classes here may be regarded as 
classifying the corresponding obstructing $n$-bundles,
after we apply the general formalism that we outlined earlier. 
The main example of this general mechanism that will be of interest to us here 
is the case where $\gg$ is an ordinary  semisimple Lie algebra.
In particular, we consider $\gg = \mathfrak{spin}(n)$.
  For $\gg = \mathfrak{spin}(n)$ and $\mu$ a $(2n+1)$-cocycle on
  $\mathfrak{spin}(n)$, we call $\mathfrak{spin}(n)_\mu$ the
  (skeletal version of the) {\bf $(2n-1)$-brane Lie $(2n)-algebra$}.
Thus, the case of String structure and Fivebrane structure occurring in the 
fundamental string and NS fivebrane correspond to the cases $n=1$
and $n=3$ respectively. Now applying our formalism
for $\gg = \mathfrak{spin}(n)$, and $\mu_3, \mu_7$ the canonical
  3- and 7-cocycle, respectively, we have:
  \begin{itemize}
    \item
      the obstruction to lifting a $\gg$-bundle descent object to a 
      String 2-bundle (a $\gg_{\mu_3}$-bundle descent object)
      is the first Pontryagin class of the original $\gg$-bundle cocycle;
    \item
      the obstruction to lifting a String 2-bundle 
     descent object to a 
      Fivebrane 6-bundle cocycle (a $\gg_{\mu_7}$-bundle descent object)
      is the second Pontryagin class of the original $\gg$-bundle cocycle.
  \end{itemize}

The cocyles and invariant polynomials corresponding to the two structures
are given in the following table

  \begin{table}[h]
\begin{center}
\begin{tabular}{llll}
  {\bf $p$-brane} & {\bf cocycle} & {\bf invariant polynomial}
  \\
  \hline
  $p = 1 = 4\cdot 0 +1 $ & $\mu_3 = \langle\cdot, [\cdot, \cdot]\rangle$ & $P_1 = \langle \cdot, \cdot \rangle$
  & first Pontrjagin
  \\
  $p = 5 = 4\cdot 1 + 1$ &
   $\mu_7 = \langle \cdot, [\cdot,\cdot],[\cdot,\cdot],[\cdot,\cdot] \rangle$
   & $P_2 = \langle \cdot, \cdot , \cdot, \cdot \rangle$
    & second Pontrjagin
     
\end{tabular}
\end{center}
 \caption{
   {\bf Lie algebra cohomology governing NS $p$-branes.}
  }
\end{table}
  

In case of the fundamental string, the obstruction to lifting the $PU(H)$ bundles
to $U(H)$ bundles is measured by a gerbe or a line 2-bundle. In the language of 
$E_8$ bundles this corresponds to lifting the loop group $LE_8$ bundles to the central 
extension ${\hat L}E_8$ bundles \cite{MSa}. The obstruction for the case of the String structure
is a $2$-gerbe and that of a Fivebrane structure is a $6$-gerbe. The structures
are summarized in the following table

\begin{table}[h]
\begin{tabular}{ccccc}
  {\bf obstruction}
  &&
  {\bf $G$-bundle}
  &&
  {\bf $\hat G$-bundle}  
  \\
 \hline
 \begin{tabular}{c}
   1-gerbes / line 2-bundles
   \\
   2-gerbes / line 3-bundles
   \\
   6-gerbes / line 7-bundles
 \end{tabular}
 &obstruct the lift of&
 \begin{tabular}{c}
   $PU(H)$-bundles
   \\
   $\mathrm{Spin}(n)$-bundles
   \\
   $\mathrm{String}(n)$-2-bundles   
 \end{tabular}
 &to&
 \begin{tabular}{c}
   $U(H)$-bundles
   \\
   $\mathrm{String}(n)$-2-bundles
   \\
   $\mathrm{FiveBrane}(n)$-6-bundles   
 \end{tabular}
 \end{tabular}
 \caption{
   {\bf Obstructing line $n$-bundles} appearing in string theory.
 }
\end{table}

A description can also be given in terms of (higher) loop spaces, generalizing
the known case where a String structure on a space $X$ can be viewed as 
a Spin structure on the loop space $LX$.
A fuller discussion of the ideas of this section will be given in \cite{SSS2}.

\section{Statement of the main results}

We define, for any $L_\infty$-algebra $\gg$ and any smooth space $X$, a notion of 
\begin{itemize}
  \item $\gg$-descent objects over $X$;
\end{itemize}
and an extension of these to
\begin{itemize}
  \item $\gg$-connection descent objects over $X$\,.
\end{itemize}
These descent objects are to be thought of as the data obtained from 
locally trivializing an $n$-bundle (with connection) 
whose structure $n$-group has the Lie $n$-algebra $\gg$. 
Being differential versions of $n$-functorial descent data
of such $n$-bundles, they consist of morphisms of 
quasi free differential graded-commutative algebras
(qDGCAs).

We define for each $L_\infty$-algebra $\gg$ a dg-algebra $\mathrm{inv}(\gg)$
of \emph{invariant polynomials} on $\gg$.
We show that every $\gg$-connection descent object gives rise to a collection of deRham
classes on $X$: its \emph{characteristic classes}. These are images of the
elements of $\mathrm{inv}(\gg)$.
Two descent objects are taken to be equivalent if they are \emph{concordant}
in a natural sense. 

Our first main result is

\begin{theorem}[characteristic classes]
  Characteristic classes are indeed characteristic  
  of  $\gg$-descent objects (but do not necessarily fully characterize them) in the 
  following sense: 
  \begin{itemize}
    \item Concordant $\gg$-connection descent objects 
       have the same characteristic classes.
    \item If the $\gg$-connection descent objects differ just by a gauge transformation,
      they even have the same characteristic forms.
  \end{itemize}  
  \label{characteristic classes are indeed characteristic}
\end{theorem}

This is our proposition 
\ref{characteristic classes of connection descent are invariant} and corollary
\ref{same characteristic forms for equivalent gg-connections}.

\paragraph{Remark.}
  We expect that this result can be strengthened. Currently our characteristic
  classes are just in deRham cohomology. One would expect that these 
  are images of classes in integral cohomology. 
  While we do not
  attempt here to discuss integral characteristic classes in general,
  we discuss some aspects of this for the case of abelian 
  Lie $n$-algebras $\gg = b^{n-1}\uu(1)$ in \ref{descent objects examples} 
  by relating $\gg$-descent objects to Deligne cohomology. 

  The reader should also note that in our main examples to be discussed
  in section \ref{lifting problem} we start with an $L_\infty$-connection
  which happens to be an ordinary Cartan-Ehresmann connection on an ordinary
  bundle and is hence known to have integral classes. It follows from our
  results then that also the corresponding Chern-Simons 3-connections
  in particular have an integral class. 

\medskip

We define String-like extensions $\gg_\mu$ of $L_\infty$-algebras coming
from any $L_\infty$-algebra cocycle $\mu$: a closed element in the
Chevalley-Eilenberg dg-algebra  $\mathrm{CE}(\gg)$
corresponding to $\gg$: $\mu \in \mathrm{CE}(\gg)$.
These generalize the String Lie 2-algebra which governs the dynamics
of (heterotic) superstrings.

Our second main results is

\begin{theorem}[string-like extensions and their properties]
  Every degree $(n+1)$-cocycle $\mu$ on an $L_\infty$-algebra $\gg$
  we obtain the string-like extension $\gg_\mu$ which sits in an exact seqeuence
  $$
    0 \to b^{n-1}\uu(1) \to \gg_\mu \to \gg \to 0
    \,.
  $$
  When $\mu$ is in transgression with an invariant polynomial $P$ we furthermore
  obtain a weakly exact sequence 
  $$
    0 \to \gg_\mu \to \mathrm{cs}_P(\gg) \to \mathrm{ch}_P(\mu) \to 0
  $$
  of $L_\infty$-algebras, where $\mathrm{cs}_P(\gg) \simeq \mathrm{inn}(\gg_\mu)$
  is trivializable (equivalent to the trivial $L_\infty$-algebra).
  There is an algebra of invariant polynomials on $\gg$ associated with $\mathrm{cs}_P(\gg)$
  and we show that it is the algebra of invariant polynomials of $\gg$ modulo the ideal 
  generaled by $P$.  
\end{theorem}

This is proposition \ref{the string-like extension sequence}, 
proposition \ref{Chern and Chern-Simons} and proposition \ref{basic forms on Chern-Simons algebra}.

Our third main result is

\begin{theorem}[obstructions to lifts through String-like extensions]
  For $\mu \in \mathrm{CE}(\gg)$ any degree $n+1$ $\gg$-cocycle that 
  transgresses to an invariant
  polynomial $P \in \mathrm{inv}(\gg)$,
  the obstruction to lifting a $\gg$-descent object to a $\gg_\mu$-descent object
  is a $(b^n \uu(1))$-descent object whose single characteristic class is 
  the class corresponding to $P$ of the original $\gg$-descent object.

  This is reflected by the fact that the cohomology of the basic forms on the Chevalley-Eilenberg
  algebra of the corresponding
  Chern-Simons $L_\infty$-algebra $\mathrm{cs}_P(\gg)$ is that of  the algebra
  of basic forms on $\mathrm{inv}(\gg)$ modulo
  the ideal generated by $P$.
\end{theorem}
This is our proposition \ref{basic forms on Chern-Simons algebra} and proposition \ref{coc}.

We discuss the following {\bf applications}.

\begin{itemize}

\item
  For $\gg$ an ordinary semisimple Lie algebra and $\mu$ its canonical 3-cocycle, 
 the obstruction to lifting a $\gg$-bundle to a String 2-bundle
is a Chern-Simons 3-bundle with characteristic class the Pontrjagin class
 of the original bundle. This is a special case of
 our proposition \ref{coc} which is spelled out in detail in
 in \ref{obstruction examples}.

The vanishing of this obstruction is known as a String structure
\cite{Kill, Kuribayashi, MurrayStevenson}. In categorical language, this issue was 
first discussed in \cite{StolzTeichner}.

By passing from our Lie $\infty$-algebraic description to smooth spaces
along the lines of \ref{differential forms on spaces of maps} and then
forming fundamental $n$-groupoids of these spaces, one can see that 
our construction of obstructing $n$-bundles to lifts through String-like
extensions reproduces the construction \cite{BrylinskiMcLaughlin,BrylinskiMcLaughlinII}
of {\v C}ech cocycles representing characteristic classes.
This, however, will not be discussed here.

\item

  This result generalizes to all String-like extensions. Using the
  7-cocycle on $\mathfrak{so}(n)$ we obtain lifts through 
  extensions by a Lie 6-algebra, which we call the Fivebrane Lie 
  6-algebra. Accordingly, fivebrane structures 
  on string structures are obstructed by
  the second Pontrjagin class. 

  This pattern continues and one would expect our obstruction theory
  for lifts through string-like extensions with respect to the 
  11-cocycle on $\mathfrak{so}(n)$ to correspond to    
  \emph{Ninebrane} structure.

  The issue of $p$-brane structures for higher $p$ was discussed
  before in \cite{MP}. In contrast to the discussion there, we here
  see $p$-brane structures only for $p = 4n + 1$, corresponding to the
  list of invariant polynomials and cocycles for $\mathfrak{so}(n)$.
  While our entire obstruction theory applies to all cocycles
  on all Lie $\infty$-algebras, it is only for those on 
  $\mathfrak{so}(n)$ and maybe $\mathfrak{e}_8$ for which the
  physical interpretation in the sense of $p$-brane structures
  is understood.

\item
  We discuss how the action functional of the
  topological field theory known as BF-theory  
  arises from an invariant polynomial on a strict Lie
  2-algebra, in a generalization of the integrated
  Pontrjagin 4-form of the topological term in Yang-Mills
  theory. See proposition \ref{invariant polynomials lift to crossed module} 
  and 
  the example in \ref{examples for characteristic forms}.

  This is similar to but different from the Lie 2-algebraic
  interpretation of BF theory indicated in 
  \cite{GirelliPfeiffer, GirelliPfeifferPopescu}, where the
  ``cosmological'' bilinear in the connection 2-form is not
  considered and a constraint on the admissable strict
   Lie 2-algebras is imposed.

 \item
  We indicate in \ref{parallel transport} 
  the notion of parallel transport induced by a $\gg$-connection,
  relate it to the $n$-functorial parallel transport of
  \cite{BS, SW, SWII, SWIII} and point out how this leads to
  $\sigma$-model actions in terms of dg-algebra morphisms.
  See section \ref{partra and sigma model}.

 \item
  We indicate in \ref{example for oo-configuration spaces}
  how by forming configuration spaces by sending DGCAs
  to smooth spaces and then using the internal hom of smooth space, 
  we obtain for every
  $\gg$-connection descent object configuration
  spaces of maps equipped with an action functional induced by
  the transgressed $\gg$-connection. We show that the algebra of
  differential forms on these configuration spaces naturally
  supports the structure of the corresponding BRST-BV complex,
  with the iterated ghost-of-ghost structure inherited from the
  higher degree symmetries induced by $\gg$.
  
  This construction is similar in spirit to the one given in
  \cite{AKSZ}, reviewed in \cite{Roytenberg}, but also, at least 
  superficially, a bit different. 
  
 \item
    We indicate also in \ref{example for oo-configuration spaces} 
    how this construction of configuration spaces induces the
    notion of transgression of $n$-bundles on $X$ to $(n-k)$-bundles on
    spaces of maps from $k$-dimensional spaces into $X$. An analogous
    integrated description of transgression in terms of inner homs is in
    \cite{SWII}. We show in \ref{example for oo-configuration spaces} 
    in particular that this
    transgression process relates the concept of String-structures
    in terms of 4-classes down on $X$ with the corresponding 3-classes
    on $L X$, as discussed for instance in \cite{Kuribayashi}.
  Our construction immediately generalizes
  to fivebrane and higher classes.  
  
\end{itemize}

All of our discussion here pertains to \emph{principal} $L_\infty$-connections.
One can also discuss \emph{associated} $\gg$-connections induced by 
($\infty$-)representations of
$\gg$ (for instance as in \cite{LadaMarkl}) and then study the collections 
of ``sections'' or ``modules'' of 
such associated $\gg$-connections. 

The extended quantum field theory of
a $(n-1)$-brane charged under an $n$-connection 
(``a charged $n$-particle'', definition \ref{charged n-particle}) should 
(see for instance \cite{FreedII, FreedIII, StolzTeichner, Hopkins}) assign
to each $d$-dimensional part $\Sigma$ of the brane's parameter space (``worldvolume'') the 
collection (an $(n-d-1)$-category, really) of sections/modules of the transgression 
of the $n$-bundle to the configuration space of maps from $\Sigma$.

For instance, the space of sections of a Chern-Simons 3-connection trangressed
to maps from the circle should yield the representation category of the 
Kac-Moody extension of the corresponding loop group. 

Our last proposition \ref{transgression to get Kac-Moody cocycle} points in
this direction. But a more detailed discussion will not be given here.

\section{Differential graded-commutative algebra}

Differential $\mathbb{N}$-graded commutative algebras (DGCAs) 
play a prominent
role in our discussion. One way to understand what is 
special about DGCAs is to realize that every DGCA can be
regarded, essentially, as the algebra of 
\emph{differential forms on some generalized smooth space}.

We explain what this means precisely in 
\ref{differential forms on spaces of maps}. The underlying
phenomenon is essentially the familiar governing principle
of Sullivan models in rational homotopy theory \cite{Hess,Sullivan},
but instead of working with simplicial spaces, we 
here consider presheaf categories. This will not
become relevant, though, until the discussion of 
configuration spaces, parallel transport and action functionals
in \ref{partra and sigma model}.

\subsection{Differential forms on smooth spaces}
\label{differential forms on spaces of maps}

We can think of every differential graded commutative
algebra essentially as being the algebra of differential
forms on \emph{some} space, possibly a generalized space.

\begin{definition}
  Let $S$ be the category whose objects are the open subsets
  of $\mathbb{R} \cup \mathbb{R}^2 \cup \mathbb{R}^3 \cup \cdots$
  and whose morphisms are smooth maps between these.
  We write 
  \(
    S^\infty := \mathrm{Set}^{S^{\mathrm{op}}}
  \)
  for the category of set-valued presheaves on $S$. 
\end{definition}

So an object $X$ in $S^\infty$ is an assignment of sets
$U \mapsto X(U)$ to each open subset $U$, 
together with an assignment
\(
  (\xymatrix{U \ar[r]^\phi & V})
  \mapsto 
  (\xymatrix{X(U) \ar@{<-}[r]^{\phi_X^*} & X(V)})
\)
of maps of sets to maps of smooth subsets which respects composition.
A morphism
\(
  f : X \to Y
\)
of smooth spaces is an assignment 
$U \mapsto 
(
  \xymatrix{X(U) \ar[r]^{f_U} & Y(U)}
)
$
of maps of sets to open subsets, such that for all smooth maps
of subsets $\xymatrix{U \ar[r]^\phi & V}$ we have that the square
\(
  \xymatrix{
     X(V) 
     \ar[r]^{f_V} 
     \ar[d]^{\phi_X^*}
     & Y(V)
     \ar[d]^{\phi_Y^*}
     \\
     X(U) \ar[r]^{f_U} & Y(U)
  }
\)
commutes.
We think of the objects of $S^\infty$ smooth spaces.
The set $X(U)$ that such a smooth space $X$ assigns to an open
subset $U$ is to be thought of as the set of smooth maps
from $U$ into $X$. As opposed to manifolds which are 
\emph{locally isomorphic} to an object in $S$, smooth spaces
can hence be thought of as being objects which are just
required to have the property that they may be \emph{probed}
by objects of $S$.
Every open subset $V$ becomes a smooth space by setting
\(
  V : U \mapsto \mathrm{Hom}_{S^\infty}(U,V)
  \,.
\)
This are the \emph{representable} presheaves.
Similarly, every ordinary manifold $X$ becomes a smooth space by
setting
\(
  X : U \mapsto \mathrm{Hom}_{\mathrm{manifolds}}(U,X)
  \,.
\)
The special property of smooth spaces which we need here is
that they form a (cartesian) \emph{closed} category: 
\begin{itemize}
\item
for any two smooth spaces $X$ and $Y$ there is a cartesian
product $X \times Y$, which is again a smooth space,
given by the assignment
\(
  X \times Y : U \mapsto X(U) \times Y(U)
  \,;
\)
where the cartesian product on the right is that of sets;
\item
the collection $\mathrm{hom}(X,Y)$
of morphisms from one smooth space $X$ to another smooth space
$Y$ is again a smooth space, given by the assignment
\(
  \mathrm{hom}_{S^\infty}(X,Y) : U \mapsto \mathrm{Hom}_{S^\infty}(X \times U, Y)
  \,.
\)
\end{itemize}

A very special smooth space is the smooth space of 
differential forms.
\begin{definition}
  We write $\Omega^\bullet$ for the smooth space which
  assigns to each open subset the set of differential forms
  on it
  \(
    \Omega^\bullet : U \mapsto \Omega^\bullet(U)
    \,.
  \)
  Using this object we define
  the DGCA of differential forms on any smooth 
  space $X$ to be the set
  \(
    \Omega^\bullet(X) := \mathrm{Hom}_{S^\infty}(X,\Omega^\bullet)
  \)
  equipped with the obvious DGCA structure induced by the 
  local DGCA structure of each $\Omega^\bullet(U)$.
\end{definition}

Therefore the object $\Omega^\bullet$ is in a way both a smooth space
as well as a differential graded commutative algebra: it is a DGCA-valued
presheaf. Such objects are known as \emph{schizophrenic} \cite{Johnstone} 
or better \emph{ambimorphic} \cite{Trimble} objects: they relate two different
worlds by duality.
In fact, the process of mapping into these objects provides an adjunction between 
the dual categories:

\begin{definition}
  \label{functors between forms and spaces}
  There are contravariant functors from smooth spaces to
  DGCAs given by
  \begin{eqnarray}
    \Omega^\bullet : S^\infty &\to& \mathrm{DGCA}s
  \nonumber\\
    X &\mapsto& \Omega^\bullet(X)
  \end{eqnarray}
  and
  \begin{eqnarray}
    \mathrm{Hom}(-,\Omega^\bullet(-)) : \mathrm{DGCA}s &\to& S^\infty
  \nonumber\\
    A &\mapsto& X_A
  \end{eqnarray}
\end{definition}

These form an adjunction of categories.
  The unit
  \(
    \xymatrix{
      \mathrm{DGCAs}
      \ar@/^2pc/[rrrr]^{\mathrm{Id}}_{\ }="s"
      \ar[rr]_{\mathrm{Hom}(-,\Omega^\bullet(-))}
      &&
      S^\infty 
      \ar[rr]_{\Omega^\bullet}
      &&
      \mathrm{DGCAs}
      \ar@{=>}^{u} "s"; "s"+(0,-5)
    }
  \)
  of this adjunction is a natural transformation 
  whose component map embeds each DGCA $A$ into the 
  algebra of differential forms on the smooth space it defines
  \(
    \xymatrix{
       A \ar@{^{(}->}[rr] && \Omega^\bullet(X_A)
    }
  \)
  by sending every $a \in A$ to the map of presheaves
  \(
    (f \in \mathrm{Hom}_{\mathrm{DGCAs}}(A,\Omega^\bullet(U)))
    \mapsto
    (f(a) \in \Omega^\bullet(U))
    \,.
  \)
This way of obtaining forms on $X_A$ from elements of $A$ 
will be crucial for our construction of differential forms
on spaces of maps, $\mathrm{hom}(X,Y)$, 
used in \ref{configuration spaces}.

Using this adjunction, we can ``pull back'' the internal hom of $S^\infty$ to
$\mathrm{DGCAs}$. Since the result is not literally the internal hom in 
$\mathrm{DGCA}s$ (which does not exist since DGCAs are not \emph{profinite}
as opposed to codifferential coalgebras \cite{Getzler}) we call it 
``$\mathrm{maps}$'' instead of ``$\mathrm{hom}$''.

\begin{definition}[forms on spaces of maps]
  \label{forms on maps functor}
  Given any two DGCAs $A$ and $B$, we define the DGCA
  of ``maps'' from $B$ to $A$
  \(
    \mathrm{maps}(B,A)
    :=
    \Omega^\bullet(\mathrm{hom}_{S^\infty}(X_A, X_B))
    \,.
  \)
  \label{space of maps between two DGCAs}
\end{definition}
  This is a functor
  \(
    \mathrm{maps}
    :
    \mathrm{DGCAs}^{\mathrm{op}}
    \times
    \mathrm{DGCAs}  
    \to
    \mathrm{DGCAs}\,.
  \)

Notice the fact (for instance corollary 35.10 in \cite{KolarSlovakMichor}
and theorem 2.8 in \cite{MoerdijkReyes}) 
that for any two smooth spaces $X$ and $Y$, algebra homomorphisms 
$C^\infty(X) \stackrel{\phi^*}{\leftarrow} C^\infty(Y)$ and hence DGCA morphisms
$\Omega^\bullet(X) \stackrel{\phi^*}{\leftarrow} \Omega^\bullet(Y)$
are in bijection with smooth maps $\phi : X \to Y$.

It follows that an element of $\mathrm{hom}(X_A,X_B)$ is, over test domains 
$U$ and $V$ a natural map of sets
\(
  \mathrm{Hom}_{\mathrm{DGCAs}}(A,\Omega^\bullet(V))
  \times
  \mathrm{Hom}_{\mathrm{DGCAs}}(\Omega^\bullet(U),\Omega^\bullet(V))
  \to
  \mathrm{Hom}_{\mathrm{DGCAs}}(B, \Omega^\bullet(V))
  \,.
\)
One way to obtain such maps is from pullback along algebra homomorphisms
$$
  B \to A \otimes \Omega^\bullet(U)
  \,.
$$
This will be an important source of DGCAs of maps for the case that $A$ is 
the Chevalley-Eilenberg algebra of an an $L_\infty$-algebra, as described in 
\ref{examples for forms on smooth spaces}.

\subsubsection{Examples}
\label{examples for forms on smooth spaces}

\paragraph{Diffeological spaces}

Particularly useful are smooth spaces $X$ which, while not
quite manifolds, have the property that there is a set
$X_s$ such that
\(
  X : U \mapsto X(U) \subset \mathrm{Hom}_{\mathrm{Set}}(U,X_s)
\)
for all $U \in S$. These are the Chen-smooth or diffeological
spaces used in \cite{BS,SWII,SWIII}. In particular, all
spaces of maps $\mathrm{hom}_{S^\infty}(X,Y)$ for $X$ and $Y$
manifolds are of this form. This includes in particular 
loop spaces.

\paragraph{Forms on spaces of maps.}

When we discuss parallel transport and its transgression
to configuration spaces in \ref{configuration spaces}, 
we need the following construction of differential
forms on spaces of maps.

\begin{definition}[currents]
  \label{currents}
  For $A$ any DGCA, we say that a current on $A$ is a smooth linear map
  \(
    c : A \to \mathbb{R}
    \,.
  \)
\end{definition}
For $A = \Omega^\bullet(X)$ this reduces to the ordinary notion of currents.

\begin{proposition}
  \label{forms from currents}
  Let $A$ be a quasi free DGCAs in positive degree (meaning that the
underlying graded commutative algebras are freely generated from some graded
 vector space in positive degree).
  For each element $b \in B$ and current $c$ on $A$, we get an element  in 
  $
    \Omega^\bullet(\mathrm{Hom}_{\mathrm{DGCA}s}(B,A \otimes \Omega^\bullet(-)))
  $
  by
  mapping, for each $U \in S$ 
  \begin{eqnarray}
    \mathrm{Hom}_{\mathrm{DGCA}s}(B,A \otimes \Omega^\bullet(U))
    &\to&
    \Omega^\bullet(U)
  \nonumber\\
    f^* &\mapsto& c(f^*(b))
    \,.
  \end{eqnarray}
\end{proposition}

If $b$ is in degree $n$ and $c$ in degree $m \leq n$, then this differential form is in degree
$n - m$.

\paragraph{The superpoint.}

  Most of the DGCAs we shall consider here are non-negatively graded or even
  positively graded. These can be thought of as Chevalley-Eilenberg algebras
  of Lie $n$-algebroids and Lie $n$-algebras, respectively, as discussed in more
  detail in \ref{Inf}.
 However, DGCAs of arbitrary degree do play an important role, too. 
  Notice that a DGCA of non-positive degree is in particular a cochain complex
  of non-positive degree. But that is the same as a chain complex of non-negative
  degree. 
  
  The following is a very simple but important example of a DGCA in non-positive degree.

   \begin{definition}[superpoint]
     \label{superpoint}
     The ``algebra of functions on the superpoint'' is the DGCA
     \(
       C(\mathbf{pt}) :=
        (\mathbb{R} \oplus \mathbb{R}[-1] ,d_{\mathbf{pt}})
     \)
     where the product on $\mathbb{R} \oplus \mathbb{R}[-1]$ 
     is the tensor product over $\mathbb{R}$, and 
     where the differential $d_{\mathbf{pt}} : \mathbb{R}[-1] \to \mathbb{R}$ is the canonical 
     isomorphism.
    \end{definition}   
    
    The smooth space associated to this algebra according to definition 
    \ref{functors between forms and spaces} is just the ordinary point, because
    for any test domain $U$ the set
    \(
      \mathrm{Hom}_{\mathrm{DGCAs}}(C(\mathbf{pt}),\Omega^\bullet(U))
    \)
    contains only the morphism which sends $1 \in \mathbb{R}$ to the constant unit function
    on $U$, and which sends $\mathbb{R}[-1]$ to 0.
   However, as is well known from the theory of supermanifolds, the algebra $C(\mathbf{pt})$
    is important in that morphisms from any other DGCA $A$ into it compute 
    the (shifted) \emph{tangent space}
    corresponding to $A$.
    From our point of view here this manifests itself in particular 
    by the fact that for $X$ any manifold, we have a canonical injection
    \(
      \Omega^\bullet(TX) \hookrightarrow \Omega^\bullet( \mathrm{maps}(C^\infty(X),C(\mathbf{pt})) )
    \)
    of the differential forms on the tangent bundle of $X$ into the differential forms
    on the smooth space of algebra homomorphisms of $C^\infty(X)$ to $C(\mathbf{pt})$:
    
    for every test domain $U$ an element in 
    $\mathrm{Hom}_{\mathrm{DGCAs}}(C^\infty(X),C(\mathbf{pt}\otimes \Omega^\bullet(U)))$
    comes from a pair consisting of a smooth map $f : U \to X$ and a vector field $v \in \Gamma(TX)$.
    Together this constitutes a smooth map $\hat f : U \to T X$ and hence for every form
    $\omega \in \Omega^\bullet(TX)$ we obtain a form on 
   $\mathrm{maps}(C^\infty(X), C(\mathbf{pt}))$ by the assignment
   \(
     ((f,v) \in \mathrm{Hom}_{\mathrm{DGCAs}}(C^\infty(X),C(\mathbf{pt}\otimes \Omega^\bullet(U))))
     \mapsto 
     (\hat f^* \omega \in \Omega^\bullet(U))
   \)
   over each test domain $U$.

   In \ref{examples for Loo algebras} we discuss how in the analogous fashion we obtain the Weil algebra
   $\mathrm{W}(\gg)$ of any $L_\infty$-algebra $\gg$ from its Chevalley-Eilenberg algebra
   $\mathrm{CE}(\gg)$ by mapping that to $C(\mathbf{pt})$. 
   This says that the Weil alghebra is like the space of functions on the shifted tangent
   bundle of the ``space'' that the Chevalley-Eilenberg algebra is the space of functions
   on. See also figure \ref{a remarkable coincidence of concepts}.

\subsection{Homotopies and inner derivations}

\label{homotopies and inner derivations}

When we forget the algebra structure of DGCAs, they are simply
cochain complexes. As such they naturally live in a 2-category
$\mathrm{Ch}^\bullet$ whose objects are cochain complexes 
$(V^\bullet,d_V)$, whose morphisms 
\(
 \xymatrix{
   (V^\bullet,d_V)
   &&
   (W^\bullet,d_W)
   \ar[ll]_{f^*}
  }
\)
are degree preserving linear maps 
$
  \xymatrix{
    V^\bullet
    &
    W^\bullet
    \ar[l]_{f^*}
  }
$
that do respect the differentials,
\(
  [d,f^*] := d_V \circ f^* - f^* \circ d_W = 0
  \,,
\)
and whose 2-morphisms 
\(
 \xymatrix{
   (V^\bullet,d_V)
   &&
   (W^\bullet,d_W)
   \ar@/_2pc/[ll]_{f^*}^{\ }="s"
   \ar@/^2pc/[ll]^{g^*}_{\ }="t"
   \ar@{=>}^{\rho} "s"; "t"
  }
\)
are cochain homotopies, namely linear degree -1 maps
$\rho : W^\bullet \to V^\bullet$ with the property that
\(
  g^* = f^* + [d,\rho]
  = f^* + d_V \circ \rho + \rho \circ d_W
  \,.
\)
Later in \ref{homotopies and concordances} we will also look at morphisms that do 
preserve the algebra structure, and homotopies of these. 
Notice that we can compose a 2-morphism from left and
right with 1-morphisms, to obtain another 2-morphism
\(
 \xymatrix{
   (U^\bullet,d_U)
   &&
   (V^\bullet,d_V)
   \ar[ll]_{h^*}
   &&
   (W^\bullet,d_W)
   \ar@/_2pc/[ll]_{f^*}^{\ }="s"
   \ar@/^2pc/[ll]^{g^*}_{\ }="t"
   &&
   (X^\bullet, d_X)
   \ar[ll]_{j^*}
   \ar@{=>}^{\rho} "s"; "t"
  }
\)
whose component map now is 
\(
  h^* \circ \rho \circ j^*
  :
  \xymatrix{
    X^\bullet 
    \ar[r]^{j^*}
    &
    W^\bullet
    \ar[r]^\rho
    &
    V^\bullet
    \ar[r]^{h^*}
    &
    U^\bullet
  }
  \,.
\)
This will be important for the interpretation of the
diagrams we discuss, 
of the type \ref{condition on vertical derivations}
and
\ref{diagram for invariance of basic forms under vertical derivations}
below.

Of special importance are linear endomorphisms  $\xymatrix{V^\bullet & V^\bullet \ar[l]_{\rho}}$
of DGCAs which are algebra derivations. Among them, the \emph{inner} derivations
in turn play a special role:

\begin{definition}[inner derivations]
  On any DGCA $(V^\bullet,d_V)$, a degree 0 endomorphism
  \(
    \xymatrix{
      (V^\bullet,d_V)
      &&
      (V^\bullet,d_V)
      \ar[ll]_L
    }
  \)
  is called an \emph{inner derivation} if
  \begin{itemize}
    \item
      it is an algebra derivation of degree 0;
    \item
      it is connected to the 0-derivation, i.e. there is a 2-morphims
  \(
    \xymatrix{
      (V^\bullet,d_V)
      &&
      (V^\bullet,d_V)
      \ar@/_2pc/[ll]_0^{\ }="s"
      \ar@/^2pc/[ll]^{L = [d_V,\rho]}_{\ }="t"
      \ar@{=>}^\rho "s"; "t"
    }
    \,,
  \)
  where $\rho$ comes from an algebra derivation of degree -1.      
  \end{itemize}
\end{definition}

\paragraph{Remark.}
Inner derivations generalize the notion of a Lie derivative on differential forms, and
hence they encode the notion of vector fields in the context of DGCAs.

\subsubsection{Examples}

\paragraph{Lie derivatives on ordinary differential forms.}

The formula somtetimes known as ``\emph{Cartan's magic formula}'', which says
that on a smooth space $Y$ the Lie derivative $L_v \omega $ of a differential form
$\omega \in \Omega^\bullet(Y)$ along a vector field $v \in \Gamma(T Y)$ is 
given by
\(
  L_v \omega = [d , \iota_v]
  \,,
\)
where $\iota_v : \Omega^\bullet(Y) \to \Omega^\bullet(Y)$, says that Lie derivatives on
differential forms are inner derivations, in our sense.
When $Y$ is equipped with a smooth projection $\pi : Y \to X$, it is of importance to 
distinguish the vector fields vertical with respect to $\pi$. The abstract formulation
of this, applicable to arbitrary DGCAs, is given in \ref{vertical flows and basic forms} below.

\subsection{Vertical flows and basic forms}

 \label{vertical flows and basic forms}

We will prominently be dealing with surjections
\(
  \xymatrix{
     A & B \ar@{->>}[l]_{i^*} 
  }
\)
of differential graded commutative algebras that play the role of the dual of an injection
\(
  \xymatrix{
     F \ar@{^{(}->}[r]^i & P
  }
\)
of a fiber into a bundle. We need a way to identitfy in the context of DGCAs which 
inner derivations of $P$ are \emph{vertical} with respect to $i$. Then we can find the
algebra corresponding to the \emph{basis} of $P$ as those elements of $B$ which are annihilated
by all vertical derivations.

\begin{definition}[vertical derivations]
  \label{vertical derivations}
  Given any surjection of differential graded algebras
  \(
    \xymatrix{
       F
       &
       P
       \ar@{->>}[l]_{i^*}
    }
  \)
  we say that the vertical inner derivations 
  \(
    \xymatrix{
       P
       &&
       P
       \ar@/_2pc/[ll]_{0}^{\ }="s"
       \ar@/^2pc/[ll]^{[d_P, \rho]}_{\ }="t"
       \ar@{=>}^\rho "s"; "t"
    }
  \)
  (this diagram is in the category of cochain complexes, 
  compare the beginning of \ref{homotopies and concordances})
  on $P$ with
  respect to $i^*$ are those inner derivations 
     \begin{itemize}
     \item
       for which there
     exists an inner derivation of $F$
  \(
    \xymatrix{
       F
       &&
       F
       \ar@/_2pc/[ll]_{0}^{\ }="s"
       \ar@/^2pc/[ll]^{[d_P, \rho']}_{\ }="t"
       \ar@{=>}^{\rho'} "s"; "t"
    }
  \)
  such that
\(
  \label{condition on vertical derivations}
  \raisebox{40pt}{
  \xymatrix{
      F
      &&
      F
      \ar@/^1.5pc/[ll]^{[d,\rho']}_{\ }="t1"
      \ar@/_1.5pc/[ll]_{0}^{\ }="s1"
      \\
      \\
      P
      \ar@{->>}[uu]
      &&
      P
      \ar@{->>}[uu]
      \ar@/^1.5pc/[ll]^{[d,\rho]}_{\ }="t2"
      \ar@/_1.5pc/[ll]_{0}^{\ }="s2"
      \ar@{=>}^{\rho'} "s1"; "t1"
      \ar@{=>}^{\rho} "s2"; "t2"
  }}
  \,;
\)
\item
  and where $\rho'$ is a \emph{contraction}, $\rho' = \iota_x$, i.e. a derivation
  which sends indecomposables to degree 0.
 \end{itemize}
\end{definition}

\begin{definition}[basic elements]
  \label{basic elements}
  Given any surjection of differential graded algebras
  \(
    \xymatrix{
       F
       &
       P
       \ar@{->>}[l]_{i^*}
    }
  \)
  we say that the algebra 
  \(
    P_{\mathrm{basic}} = \bigcap\limits_{\rho \,\mathrm{vertical}} 
      \mathrm{ker}(\rho )\cap \mathrm{ker}(\rho \circ d_p)
  \)
  of basic elements of $P$ (with respect to the surjection $i^*$) 
  is the subalgebra of $P$ of all those
  elements $a \in P$ which are annihilated by all 
  $i^*$-vertical derivations $\rho$, in that
  \begin{eqnarray}
    \rho (a) &=& 0
  \\
    \rho (d_P a) &=& 0
    \,.
  \end{eqnarray}
\end{definition}
We have a canonical inclusion
\(
  \xymatrix{
    P & P_{\mathrm{basic}} \ar@{_{(}->}[l]_{p^*}
  }
  \,.
\)
Diagrammatically the above condition says that 
\(
  \label{diagram for invariance of basic forms under vertical derivations}
  \xymatrix{
      F
      &&
      F
      \ar@/^1.5pc/[ll]^{[d,\rho']}_{\ }="t1"
      \ar@/_1.5pc/[ll]_{0}^{\ }="s1"
      \\
      \\
      P
      \ar@{->>}[uu]_{i^*}
      &&
      P
      \ar@{->>}[uu]_{i^*}
      \ar@/^1.5pc/[ll]^{[d,\rho]}_{\ }="t2"
      \ar@/_1.5pc/[ll]_{0}^{\ }="s2"
      \\
      \\
      P_{\mathrm{basic}}
      \ar@{^{(}->}[uu]_{p^*}
      &&
      P_{\mathrm{basic}}
      \ar@{^{(}->}[uu]_{p^*}
      \ar[ll]_0
      \ar@{=>}^{\rho'} "s1"; "t1"
      \ar@{=>}^{\rho} "s2"; "t2"
  }
  \,.
\)

\subsubsection{Examples}

\label{surj subm examples}

\paragraph{Basic forms on a bundle}
 \label{surjective submersions and differential forms}

As a special case of the above general defintion, we reobtain the standard notion of
basic differential forms on a smooth surjective submersion $\pi : Y \to X$ with connected
fibers.

\begin{definition}
  \label{vertical deRham complex}
   Let $\pi : Y \to X$ be a smooth map.
    The {\bf vertical deRham complex} ,
    $\Omega^\bullet_{\mathrm{vert}}(Y)$, with respect to $Y$
    is the deRham complex of $Y$ modulo those forms that vanish 
    when restricted in all arguments to 
    vector fields in the kernel of $\pi_* : \Gamma(TY )\to \Gamma(TX)$, 
    namely to vertical vector fields.
\end{definition}
The induced differential on $\Omega^\bullet_{\mathrm{vert}}(Y)$
sends $\omega_{\mathrm{vert}} = i^* \omega$ to
\(
  d_{\mathrm{vert}} : i^* \omega \mapsto i^* d\omega
  \,.
\)
\begin{proposition}
  This is well defined.
  The quotient $\Omega^\bullet_{\mathrm{vert}}(Y)$ with the differential
  induced from $\Omega^\bullet(Y)$ is indeed a dg-algebra, and the projection
    \(
      \xymatrix{
        \Omega^\bullet_{\mathrm{vert}}(Y)
        &&
        \Omega^\bullet(Y)
        \ar[ll]_{i^*}        
      }
    \)
    is a homomorphism of dg-algebras (in that it does respect the differential).
\end{proposition}
\proof
Notice that if $\omega \in \Omega^\bullet(Y)$ vanishes when evalutated on
vertical vector fields then obviously so does $\alpha\wedge \omega$, 
 for any
$\alpha \in \Omega^\bullet(Y)$. Moreover, due to the formula
\(
  d \omega (v_1,\cdots, v_{n+1})
  = 
  \sum\limits_{\sigma \in \mathrm{Sh}(1,n+1)} 
    \pm
  v_{\sigma_1} 
  \omega(v_{\sigma_2}, \cdots, v_{\sigma_{n+1}})
  +
  \sum\limits_{\sigma \in \mathrm{Sh}(2,n+1)}
  \pm
  \omega([v_{\sigma_1},v_{\sigma_2}],v_{\sigma_3}, \cdots, v_{\sigma_{n+1}})
\)
and the fact that for $v,w$ vertical so is $[v,w]$ and hence 
$d \omega$ is also vertical. This gives that vertical differential forms on $Y$
form a dg-subalgebra of the algebra of all forms on $Y$. 
Therefore if $i^* \omega = i^* \omega'$ then
\(
  d i^* \omega' = i^* d \omega' = i^* d (\omega + (\omega' - \omega)) 
  = i^* d \omega + 0 = d i^* \omega
  \,.
\)
Hence the differential is well defined and $i^*$ is then, by construction,
a morphism of dg-algebras.

\endofproof

Recall the following standard definition of basic differential forms.
\begin{definition}[basic forms]
  \label{basic forms, standard definition}
  Given a surjective submersion $\pi : Y \to X$, the \emph{basic forms} on $Y$
  are those with the property that they and their differentials are annihilated by
  all vertical vector fields
  \(
    \omega \in \Omega^\bullet(Y)_{\mathrm{basic}}
    \hspace{5pt}    
      \Leftrightarrow
    \hspace{5pt}    
    \forall v \in \mathrm{ker}(\pi) : 
    \iota_v \omega = \iota_v d \omega = 0
    \,.
  \)
\end{definition}
It is a standard result that
\begin{proposition}
  If $\pi : Y \to X$ is locally trivial and has connected fibers, then 
  the basic forms are precisely those coming from pullback along $\pi$
  \(
    \Omega^\bullet(Y)_{\mathrm{basic}} \simeq \Omega^\bullet(X)
    \,.
  \)
\end{proposition}

\paragraph{Remark.} The reader should compare this situation with the
definition of invariant polynomials in \ref{Lie infty-algebra cohomology}.

The next proposition asserts that these statements about ordinary basic differential forms
are indeed a special case of the general definition of basic elements with respect
to a surjection of DGCAs, definition
\ref{basic elements}.

\begin{proposition}
  Given a surjective submersion $\pi: Y \to X$ with connected fibers, then
  \begin{itemize}
   \item the inner derivations of $\Omega^\bullet(Y)$ which are vertical with respect
     to 
    $
     \xymatrix{
      \Omega^\bullet_{\mathrm{vert}}(Y)
      &
      \Omega^\bullet(Y)
      \ar@{->>}[l]_{i^*}
     }
    $
    according to the general definition \ref{vertical derivations},
    come precisely from contractions $\iota_v$ with vertical vector fields
    $v \in \mathrm{ker}(\pi_*) \subset \Gamma(TY)$;
    \item
    the basic differential forms on $Y$ according to definition \ref{basic forms, standard definition}
    conincide with the basic elements of $\Omega^\bullet(Y)$ relative to the above surjection
    \(
     \xymatrix{
      \Omega^\bullet_{\mathrm{vert}}(Y)
      &
      \Omega^\bullet(Y)
      \ar@{->>}[l]_{i^*}
     }
    \)
    according to the general definition \ref{basic elements}.
  \end{itemize}
\end{proposition}
\proof
      We first show that if $\xymatrix{\Omega^\bullet(Y) &\Omega^\bullet(Y) \ar[l]_{\rho}}$ is
      a vertical algebra derivation, then $\rho$ has to annihilate all forms in the
      image of $\pi_*$.
       Let $\alpha \in \Omega^\bullet(Y)$ be any 1-form and $\omega = \pi^* \beta$
      for $\beta \in \Omega^1(X)$.
      Then the wedge product $\alpha \wedge \omega$ is annihilated by
      the projection to $\Omega^\bullet_{\mathrm{vert}}(Y)$ and we find
      \(
        \raisebox{50pt}{
        \xymatrix{
             \rho(\omega)\wedge \alpha
             &&
             0
             \ar@{|->}[ll]
             \\
             \\
           \mbox{
             \begin{tabular}{c}
                 $\rho(\omega) \wedge \alpha $
                 \\
                 $+ \rho(\alpha)\wedge \omega$
             \end{tabular}
           }         
             \ar@{|->}[uu]_{i^*}
             &&
             \alpha\wedge\omega
             \ar@{|->}[uu]_{i^*}
             \ar@{|->}[ll]_{\rho}
        }
        }
        \,.
      \)
  We see that $\rho(\omega)\wedge \alpha$ has to vanish for all $\alpha$. Therefore $\rho(\omega)$
  has to vanish for all $\omega$ pulled back along $\pi^*$. Hence $\rho$ must be
  contraction with a vertical vector field.
   It then follows from the condition \ref{condition on vertical derivations} that
   a basic form is one annihilated by all such $\rho$ and all such $\rho \circ d$.  
\endofproof

Possibly the most familiar kinds of surjective submersions are
\begin{itemize}
  \item Fiber bundles. 

    Indeed, the standard Cartan-Ehresmann theory
    of connections of principal bundles is obtained in our context by
    fixing a Lie group $G$ and a principal $G$-bundle $p : P \to X$
    and then using $Y = P$ itself as the surjective submersion.
     The definition of a connection on $P$ in terms of a $\gg$-valued
    1-form on $P$ can be understood as the descent data for a connection
    on $P$ obtained with respect to canonical trivialization of the pullback
    of $P$ to $Y = P$.
    Using for the surjective submersion $Y$ a principal $G$-bundle
    $P \to X$ is also most convenient for studying all kinds
    of higher $n$-bundles obstructing lifts of the given
    $G$-bundle. This is why we will often make use of this
    choice in the following.

  \item
    Covers by open subsets.
    
    The disjoint union of all sets in a cover of $X$ by open subsets of $X$
    forms a surjective submersion $\pi : Y \to X$.  In large parts of the
    literature on descent (locally trivialized bundles), these are the
    only kinds of surjective submersions that are considered. 
    We will find here,
    that in order to characterize principal $n$-bundles
    entirely in terms of $L_\infty$-algebraic data, open covers are too restrictive
    and the full generality of surjective submersions is needed.
    The reason is that, for $\pi : Y \to X$ a cover by open subsets, there are no
    nontrivial vertical vector fields
    \(
      \mathrm{ker}(\pi) = 0,
    \)
    hence
    \(
      \Omega^\bullet_{\mathrm{vert}}(Y) = 0
      \,.
    \)
    With the definition of $\gg$-descent objects in \ref{gg-descent object}
    this implies that all $\gg$-descent objects over a cover by open subsets
    are trivial.
\end{itemize}

There are two important subclasses of surjective submersions $\pi : Y \to X$:

\begin{itemize}
  \item
    those for which $Y$ is (smoothly) contractible;
  \item
    those for which the fibers of $Y$ are connected.
\end{itemize}

We say $Y$ is (smoothly) contractible if the identity map $\mathrm{Id} : Y \to Y$
is (smoothly) homotopic to a map $Y \to Y$ which is constant on each connected component.
Hence $Y$ is a disjoint union of spaces that are each (smoothly) contractible to a point.
In this case the Poincar{\'e} lemma says that the dg-algebra 
$
  \Omega^\bullet(Y)
$
of differential forms
on $Y$ is contractible; each closed form is exact:
\(
  \xymatrix{
    \Omega^\bullet(Y)
    &&
    \Omega^\bullet(Y)
    \ar@/_2pc/[ll]_{0}^{\ }="s"
    \ar@/^2pc/[ll]^{[d,\tau]}_{\ }="t"
    \ar@{=>}^\tau "s"; "t"
  }
  \,.
\)
Here $\tau$ is the familiar homotopy operator that appears in the proof of the 
Poincar{\'e} lemma.
In practice, we often make use of the best of both worlds: 
surjective submersions that
are (smoothly) contractible to a discrete set but still have a sufficiently rich collection of vertical 
vector fields. 
The way to obtain these is by refinement: starting with any surjective submersion 
$\pi : Y \to X$ which has good vertical vector fields but might not be contractible,
we can cover $Y$ itself with open balls, whose disjoint union, $Y'$, then forms 
a surjective submersion $Y' \to Y$ over $Y$. The composite $\pi'$
\(
  \raisebox{20pt}{
  \xymatrix{
    Y' \ar[dr]_{\pi}\ar[rr] && Y \ar[dl]    
    \\
    &
    X
  }
  }
\)
is then a contractible surjective submersion of $X$. We will see that all our
descent objects can be pulled back along refinements of surjective submersions
this way, so that it is possible, without restriction of generality, to always
work on contractible surjective submersions.
Notice that for these the structure of
\(
  \xymatrix{
    \Omega^\bullet_{\mathrm{vert}}(Y)
    &&
    \Omega^\bullet(Y)
    \ar@{->>}[ll]
    &&
    \Omega^\bullet(X)
    \ar@{_{(}->}[ll]
  }
\)
is rather similar to that of
\(
  \xymatrix{
    \mathrm{CE}(\gg)
    &&
    \mathrm{W}(\gg)
    \ar@{->>}[ll]
    &&
    \mathrm{inv}(\gg)
    \ar@{_{(}->}[ll]
  }
  \,,
\)
since $\mathrm{W}(\gg)$ is also contractible,  according to 
proposition \ref{properties of Weil algebra}.

\paragraph{Vertical derivations on universal $\gg$-bundles.}
The other important exmaple of vertical flows, those on DGCAs modelling universal
$\gg$-bundle for $\gg$ an $L_\infty$-algebra, is discussed at the beginning of 
\ref{Lie infty-algebra cohomology}.

\section{$L_\infty$-algebras and their String-like extensions}

\label{Inf}

$L_\infty$-algebras are a generalization of Lie algebras,
where the Jacobi identity is demanded to hold only up
to higher coherent equivalence, as the category theorist would say,
or ``strongly homotopic'', as the homotopy theorist would say.

\subsection{$L_\infty$-algebras}

\label{Lie infty algebras}

\begin{definition} Given a graded vector space $V$, the 
{\em tensor space} $T^\bullet (V): = \bigoplus_{n=0} 
V^{\otimes n}$
with $V^0$ being the ground field. We will denote by 
$T^a (V)$ the {\em tensor algebra} with  the concatenation product on  $T^\bullet (V)$:
\(
x_1 \otimes  x_2 \otimes \cdots \otimes  x_p \bigotimes x_{p+1}\otimes \cdots \otimes  
x_n \mapsto x_1 \otimes  x_2 \otimes \cdots \otimes  x_n  
\)
and by
$T^c (V)$ the {\em tensor coalgebra} with  the deconcatenation product 
on  $T^\bullet (V)$:
\(
x_1 \otimes  x_2 \otimes \cdots \otimes  x_n \mapsto \sum_{p+q=n}x_1 \otimes  
x_2 \otimes \cdots \otimes  x_p \bigotimes x_{p+1}\otimes \cdots \otimes  x_n.
\)
The {\em graded symmetric algebra} 
$\wedge^\bullet(V)$ is the quotient of the tensor algebra 
$T^a (V)$ by the graded action of the symmetric groups
$\mathbf{S}_n$ on the components $V^{\otimes n}.$
The {\em graded symmetric coalgebra} 
$\vee^\bullet(V)$ is the sub-coalgebra of the tensor coalgebra 
$T^c (V)$ fixed by the graded action of the symmetric groups
$\mathbf{S}_n$ on the components $V^{\otimes n}.$
\end{definition}

\paragraph{Remark.} $\vee^\bullet(V)$ is spanned by graded symmetric tensors
\(
 x_1 \vee  x_2 \vee \cdots \vee  x_p
 \)
for  $x_i \in V$ and $p \geq 0,$
where we use $\vee$ rather than $\wedge$ to emphasize the coalgebra aspect,
e.g. 
\(
x \vee y = x \otimes y \pm y \otimes x.
\)

In characteristic zero, the graded symmetric algebra can be identified with a sub-algebra 
of $T^a (V)$ but that is unnatural and we will try to avoid doing so.
The coproduct on $\vee^\bullet(V)$ is given by
\(
  \Delta(x_1 \vee x_2 \cdots \vee x_n)
  =
  \sum_{p+q = n} \sum_{\sigma \in \mathrm{Sh}(p,q)}
   \epsilon(\sigma) (x_{\sigma(1)}\vee x_{\sigma(2)} \cdots x_{\sigma(p)})
   \otimes
   (x_{\sigma(p+1)} \vee \cdots x_{\sigma(n)})
   \,.
\)
The notation here means the following:
\begin{itemize}
\item $\mathrm{Sh}(p,q)$ is the subset of all those bijections 
(the ``unshuffles'') of 
$\{1,2, \cdots, p+q\}$ that have the property that $\sigma(i) < \sigma(i+1)$ 
whenever $i \neq p$;
\item
$\epsilon(\sigma)$, which is shorthand for
$\epsilon(\sigma,x_1 \vee x_2, \cdots x_{p+q})$, 
the Koszul sign, defined by
\(
  x_1 \vee \cdots \vee x_n 
  =
  \epsilon(\sigma) x_{\sigma(1)} \vee \cdots x_{\sigma(n)}
  \,.
\)
\end{itemize}

\begin{definition}[$L_\infty$-algebra]
  \label{Loo algebra}
  An $L_\infty$-algebra $\gg = (\gg,D)$ is a $\mathbb{N}_+$-graded vector space $\gg$
  equipped with a degree -1 coderivation
  \(
    D : \vee^\bullet \gg \to \vee^\bullet \gg
  \)
  on the graded co-commutative coalgebra generated by $\gg$, such that
  $D^2 = 0$.
  This induces a differential
  \(
    d_{\mathrm{CE}(\gg)} : \mathrm{Sym}^\bullet(\gg) \to \mathrm{Sym}^{\bullet+1}(\gg)
  \)
  on graded-symmetric multilinear functions on $\gg$. When $\gg$ is finite dimensional
  this yields a degree +1 differential
  \(
    d_{\mathrm{CE}(\gg)} : \wedge^\bullet \gg^* \to \wedge^\bullet \gg^*
  \)
  on the graded-commutative algebra generated from $\gg^*$. This is the Chevalley-Eilenberg
  dg-algebra corresponding to the $L_\infty$-algebra $\gg$.
\end{definition}

\paragraph{Remark.}
That the original definition of $L_\infty$-algebras in terms of 
multibrackets yields a codifferential coalgebra
as above was shown in 
\cite{lada-jds}. That every such codifferential comes from a collection
of multibrackets this way is due to \cite{LadaMarkl}.

\paragraph{Example}
For $(\gg[-1], [\cdot,\cdot])$ an ordinary Lie algebra
(meaning that we regard the vector space $\gg$ to be in degree 1),
the corresponding Chevalley-Eilenberg qDGCA is 
\(
  \mathrm{CE}(\gg) = (\wedge^\bullet \gg^*, d_{\mathrm{CE}(\gg)})
\)
with
\(
  d_{\mathrm{CE}(\gg)} : 
     \xymatrix{
       \gg^*
       \ar[r]^<<<<<{[\cdot,\cdot]^*}
       &
       \gg^*\wedge \gg^*
     } 
     \,.
\)
If we let $\{t_a\}$ be a basis of $\gg$ and $\{C^a{}_{bc}\}$ 
the corresponding structure constants of the Lie bracket $[\cdot, \cdot]$, 
and if we denote by $\{t^a\}$ the corresponding basis of
$\gg^*$, then we get
\(
  d_{\mathrm{CE}(\gg)} t^a = - \frac{1}{2}C^a{}_{bc}t^b \wedge t^c
  \,.
\)
If $\gg$ is concentrated in degree $1,\dots,n$, we also say that $\gg$ is a
{\bf Lie $n$-algebra}.
Notice that built in we have a shift of degree for convenience, which makes
ordinary Lie 1-algebras be in degree 1 already. In much of the literature
a Lie $n$-algebra would be based on a vector space concentratred in degrees 0 to $n-1$.
 An ordinary Lie algebra is a Lie 1-algebra. Here the
coderivation
differential $D = [\cdot,\cdot]$ is just the Lie bracket, extended as a
coderivation to $\vee^\bullet \gg$, with $\gg$ regarded as being in
degree 1.

\vspace{3mm}
In the rest of the paper we assume, just for simplicity and since it is 
sufficient for our applications,
all $\gg$  to be finite-dimensional. Then, by the above,
these $L_\infty$-algebras are equivalently conceived of in terms of their
dual Chevalley-Eilenberg algebras, $\mathrm{CE}(\gg)$,  as indeed every quasi-free differential graded commutative algebra 
(``qDGCA'', meaning that it is free as a graded commutative algebra) corresponds
to an $L_\infty$-algebra. We will find it convenient to work entirely in terms
of qDGCAs, which we will usually denote as $\mathrm{CE}(\gg).$

While not very interesting in themselves, truly free differential algebras
are a useful tool for handling quasi-free differential algebras.
\begin{definition}
  \label{free DGCAs}
  We say a qDGCA is \emph{free} (even as a differential algebra)
  if it is of the form
  \(
    \mathrm{F}(V) := (\wedge^\bullet (V^* \oplus V^*[1]), d_{\mathrm{F}(V)})
  \)
  with 
  \(
    d_{\mathrm{F}(V)}|_{V^*}  = \sigma : V^* \to V^*[1]
  \)
  the canonical isomorphism and 
  \(
    d_{\mathrm{F}(V)}|_{V^*[1]} = 0
    \,.
  \)
\end{definition}

\paragraph{Remark.}
Such algebras are indeed free in that they satisfy the universal property:
given any linear map $V\to W$, it uniquely extends to a morphism of 
qDGCAs $F(V)\to (\wedgebullet(W^*), d)$ for any choice of differential $d$.

\paragraph{Example.} The free qDGCA on a 1-dimensional vector space in degree 0
is the graded commutative algebra freely generated by two generators,
$t$ of degree 0 and $dt$ of degree 1, with the differential acting as $d : t \mapsto dt$
and $d : dt \mapsto 0$. In rational homotopy theory, this models the
interval $I = [0,1]$. The fact that the qDGCA is free corresponds
to the fact that the interval is homotopy equivalent to the point.

We will be interested in  qDGCAs that arise as mapping cones of morphisms
of $L_\infty$-algebras.

\begin{definition}[``mapping cone'' of qDGCAs]
  Let 
  \( \xymatrix{
       \mathrm{CE}(\hh)
       &
       \mathrm{CE}(\gg)
       \ar[l]_{t^*}
    }
    \)
  be a morphism of qDGCAs. 
    The mapping cone of $t^*$, which we write
  $\mathrm{CE}(\hh \stackrel{t}{\to} g)$, 
  is the qDGCA whose underlying graded algebra is
  \(
    \wedge^\bullet(  \gg^* \oplus \hh^*[1])
  \)
  and whose differential $d_{t^*}$ is such that it acts 
  as
  \(
     d_{t^*} =
     \left(
       \begin{array}{cc}
         d_{\gg} & 0
         \\
         t^* & d_{\hh}
       \end{array}
     \right)
     \,.
  \)
  \label{mapping cone of qDGCAs}
 \end{definition}
We postpone a more detailed definition and discussion to 
\ref{weak cokernels of Lie infty-algebras}; see definition
\ref{mapping cone, detailed def} and proposition 
\ref{nilpotency for mapping cone}.
Strictly speaking, the more usual
 notion of mapping cones of chain complexes applies to $t:\hh\to\gg$, but then is  extended as a 
derivation differential to the entire qDGCA.

\begin{definition}[Weil algebra of an $L_\infty$-algebra]
  The mapping cone of the identity on $\mathrm{CE}(\gg)$ is the Weil algebra
  \(
    \mathrm{W}(\gg) := \mathrm{CE}(\gg \stackrel{\mathrm{Id}}{\to} \gg)
  \)
  of $\gg$.
  \label{Weil algebra definition}
\end{definition}

\begin{proposition}
  For $\gg$ an ordinary Lie algebra this does coincide with the
  ordinary Weil algebra of $\gg$.
\end{proposition}
\proof
See the example in \ref{examples for Loo algebras}.
\endofproof

The Weil algebra has two important properties.

\begin{proposition}
  \label{properties of Weil algebra}
  The Weil algebra $\mathrm{W}(\gg)$ of any $L_\infty$-algebra $\gg$
  \begin{itemize}
    \item is isomorphic to a free differential algebra
      \(
        \mathrm{W}(\gg) \simeq \mathrm{F}(\gg)\,,
      \)
      and hence is contractible;
    \item
      has a canonical surjection
      \(
        \xymatrix{
          \mathrm{CE}(\gg)
          &
          \mathrm{W}(\gg)
          \ar@{->>}[l]_{i^*}
        }
        \,.
      \)
  \end{itemize}
\end{proposition}
\proof
  Define a morphism
  \(
    f : \mathrm{F}(\gg) \to \mathrm{W}(\gg)
  \)
  by setting
  \begin{eqnarray}
    f &:& a \mapsto a
  \\
    f &:& (d_{\mathrm{F}(V)}a = \sigma a) \mapsto (d_{\mathrm{W}(\gg)} a = 
      d_{\mathrm{CE}(\gg)}a + \sigma a)
  \end{eqnarray}
  for all $a \in \gg^*$ and extend as an algebra homomorphism. This
  clearly respects the differentials: for all $a \in V^*$ we find
  \(
    \raisebox{20pt}{
    \xymatrix{
      a 
      \ar@{|->}[r]^{d_{F(\gg)}}
      \ar@{|->}[d]^{f}
      &
      \sigma a
      \ar@{|->}[d]^f
      \\
      a
      \ar@{|->}[r]_<<<<<{d_{\mathrm{W}(\gg)}}
      &
      d_{\mathrm{CE}(\gg)} a + \sigma a
    }
    }
    \hspace{10pt}
    \mbox{and}
    \hspace{10pt}
    \raisebox{20pt}{
    \xymatrix{
      \sigma a 
      \ar@{|->}[r]^{d_{F(\gg)}}
      \ar@{|->}[d]^{f}
      &
      0
      \ar@{|->}[d]^f
      \\
      d_{\mathrm{W}(\gg)}a
      \ar@{|->}[r]_<<<<<{d_{\mathrm{W}(\gg)}}
      &
      0
    }}\,.
  \)
  One checks that the strict inverse exists and is given by
  \begin{eqnarray}
    f^{-1}|_{\gg^*} &:& a \mapsto a
  \\
    f^{-1}|_{\gg^*[1]} &:& \sigma a \mapsto d_{F(\gg)} a - d_{\mathrm{CE}(\gg)} a
    \,.
  \end{eqnarray}
  Here $\sigma : \gg^* \to \gg^*[1]$ is the canonical  isomorphism that shifts the
  degree.
  The surjection      
  $
        \xymatrix{
          \mathrm{CE}(\gg)
          &
          \mathrm{W}(\gg)
          \ar@{->>}[l]_{i^*}
        }
      $ simply projects out all elements in the shifted copy
  of $\gg$:
  \begin{eqnarray}
    i^*|_{\wedge^\bullet \gg^*} &=& \mathrm{id}
  \\
    i^*|_{\gg^*[1]} &=& 0
    \,.
  \end{eqnarray}
  This is an algebra homomorphism that respects the differential.
\endofproof

As a corollary we obtain
\begin{corollary}
  \label{homotopy operator}
  For $\gg$ any $L_\infty$-algebra, the cohomology of
  $\mathrm{W}(\gg)$ is trivial.
\end{corollary}

\begin{proposition}
  \label{functoriality of W}
  The step from a Chevalley-Eilenberg algebra to the corresponding Weil algebra is
  functorial: for any morphism
  \(
    \xymatrix{
      \mathrm{CE}(\hh) && \mathrm{CE}(\gg) \ar[ll]_{f^*}
    }
  \)
  we obtain a morphism
  \(
    \xymatrix{
      \mathrm{W}(\hh) && \mathrm{W}(\gg) \ar[ll]_{\hat f^*}
    }
  \)
  and this respects composition.
\end{proposition}
\proof
  The morphism $\hat f^*$ acts as for all generators $a \in \gg^*$ as
  \(
    \hat f^* : a \mapsto f^*(a)
  \)
  and
  \(
    \hat f^* : \sigma a \mapsto \sigma f^*(a)
    \,.
  \)
  We check that this does repect the differentials
  \(
    \raisebox{15pt}{
    \xymatrix{
      a \ar@{|->}[d]^{\hat f^*} 
       \ar@{|->}[r]^<<<<<<<<<<<{d_{\mathrm{W}(\gg)}} & d_{\mathrm{CE}(\gg)} a + \sigma a
        \ar@{|->}[d]^{\hat f^*}
      \\
      f^*(a)
      \ar[r]^<<<<<{d_{\mathrm{W}(\hh)}}
      &
      d_{\mathrm{CE}(\hh)} f^*(a) + \sigma f^*(a)
    }}
    \hspace{9pt}
    \raisebox{15pt}{
    \xymatrix{
      \sigma a \ar@{|->}[d]^{\hat f^*} 
       \ar@{|->}[r]^<<<<<<<<<<<{d_{\mathrm{W}(\gg)}} 
       & -\sigma(d_{\mathrm{CE}(\gg)} a) 
        \ar@{|->}[d]^{\hat f^*}
      \\
      \sigma f^*(a)
      \ar[r]^<<<<<{d_{\mathrm{W}(\hh)}}
      &
      -\sigma (d_{\mathrm{CE}(\hh)} a)
    }}\,,
  \)
\endofproof

\paragraph{Remark.} As we will shortly see, 
$\mathrm{W}(\gg)$
plays the role of the algebra of differential forms on the universal
$\gg$-bundle. The surjection $\xymatrix{ CE(\gg) & W(\gg) \ar@{->>}[l]_{i^*} }$
plays the role of the restriction to the differential forms on the
fiber of the universal $\gg$-bundle.

\subsubsection{Examples}

\label{examples for Loo algebras}

In section \ref{String-like extensions} 
we construct large families of examples of  $L_\infty$-algebras,
based on the first two of the following examples:

\medskip
\noindent {\bf 1. Ordinary Weil algebras as Lie 2-algebras.} 
What is ordinarily called the
Weil algebra $\mathrm{W}(\gg)$ of a Lie algebra $(\gg[-1], [\cdot,\cdot])$ 
can, since it is again a DGCA, also be interpreted as the
Chevalley-Eilenberg algebra of a Lie 2-algebra. This Lie 2-algebra we 
call $\mathrm{inn}(\gg)$. It corresponds to the Lie 2-group
$\mathrm{INN}(G)$ discussed in \cite{RS}:
\(
  \mathrm{W}(\gg) = \mathrm{CE}(\mathrm{inn}(\gg))
  \,.
\)
We have
\(
  \mathrm{W}(\gg) = (\wedge^\bullet( \gg^* \oplus \gg^*[1]), d_{\mathrm{W}(\gg)})
  \,.
\)
Denoting by $\sigma : \gg^* \to \gg^*[1]$ the canonical isomorphism, extended
as a derivation to all of $\mathrm{W}(\gg)$, we have
\(
  d_{\mathrm{W}(\gg)} : \xymatrix{
    \gg^*
    \ar[rr]^<<<<<<<<<<<{[\cdot,\cdot]^* + \sigma}
    &&
    \gg^* \wedge \gg^* \oplus \gg^*[1]
  } 
\)
and
\(
  d_{\mathrm{W}(\gg)} : 
  \xymatrix{
     \gg^*[1]
     \ar[rr]^{- \sigma \circ d_{\mathrm{CE}(\gg)} \circ \sigma^{-1}}
     &&
     \gg^* \otimes \gg^*[1]
  }
  \,.
\)
With $\{t^a\}$ a basis for $\gg^*$ as above, and $\{\sigma t^a\}$ the corresponding
basis of $\gg^*[1]$ we find
\(
  d_{\mathrm{W}(\gg)} : t^a \mapsto -\frac{1}{2}C^a{}_{bc} t^b \wedge t^c + \sigma t^a
\)
and
\(
  \label{differential of Weil algebra}
  d_{\mathrm{W}(\gg)} : \sigma t^a \mapsto - C^a{}_{bc}t^b \sigma t^c
  \,.
\)
The Lie 2-algebra $\mathrm{inn}(\gg)$ is, in turn, nothing but
the strict Lie 2-algebra as in the third example below, which
comes from the infinitesimal crossed module
$(\gg \stackrel{\mathrm{Id}}{\to} \gg 
  \stackrel{\mathrm{ad}}{\to} \mathrm{der}(\gg))$.

\medskip
\noindent {\bf 2. Shifted $\uu(1)$.} By the above, the qDGCA corresponding to the Lie
algebra $\uu(1)$ is simply
\(
  \mathrm{CE}(\uu(1)) = (\wedge^\bullet \mathbb{R}[1], d_{\mathrm{CE}(\uu(1))} = 0)
  \,.
\)
We write
\(
  \mathrm{CE}(b^{n-1}\uu(1)) = (\wedge^\bullet \mathbb{R}[n], d_{\mathrm{CE}(b^n\uu(1))} = 0)
\)
for the Chevalley-Eilenberg algebras corresponding to the Lie $n$-algebras
$b^{n-1} \uu(1)$.

\medskip
\noindent {\bf 3. Infinitesimal crossed modules and strict Lie 2-algebras.}
 An \emph{infinitesimal crossed module} is a diagram
    \(
      (
       \xymatrix{
         \hh \ar[r]^t & \gg \ar[r]^\alpha & \mathrm{der}(\hh)
       }
      )
   \)
  of Lie algebras  
  where $t$ and $\alpha$ satisfy two compatibility
  conditions. These conditions are equivalent to the nilpotency of the differential
  on
  \(
    \mathrm{CE}(\hh \stackrel{t}{\to} \gg)
    :=
    (\wedge^\bullet (\gg^* \oplus \hh^*[1]), d_{t})
  \)
  defined by
  \begin{eqnarray}
    \label{differential in CE of strict Lie 2-algebra}
    d_t|_{\gg^*} &=& [\cdot,\cdot]^*_{\gg} + t^*
  \\
     d_t|_{\hh^*[1]} &=& \alpha^*
     \,,
  \end{eqnarray}
  where we consider the vector spaces underlying both $\gg$ and $\hh$ to be in degree 1.
  Here in the last line we regard $\alpha$ as a linear map 
  $\alpha : \gg \otimes \hh \to \hh$.
  The Lie 2-algebras $(\hh \stackrel{t}{\to} \gg)$ thus defined are
  called strict Lie 2-algebras: these are precisely those Lie 2-algebras
  whose Chevalley-Eilenberg differential contains at most co-binary components.

\medskip
\noindent {\bf 4. Inner derivation $L_\infty$-algebras.}
In straightforward generalization of the first exmaple we find:
for $\gg$ any $L_\infty$-algebra, 
its Weil algebra $\mathrm{W}(\gg)$ is again a DGCA, hence the
Chevalley-Eilenberg algebra of some other $L_\infty$-algbera.
This we address as the $L_\infty$-algebra of inner derivations
and write
\(
  \label{Weil and inner}
  \mathrm{CE}(\mathrm{inn}(\gg)) :=
  \mathrm{W}(\gg)
  \,.
\)
This identification is actually useful for identifying the 
Lie $\infty$-groups that correspond to an integrated picture
underlying our differential discussion. In \cite{RS}
the Lie 3-group corresponding to $\mathrm{inn}(\gg)$ for 
$\gg$ the strict Lie 2-algebra of any strict Lie 2-group is discussed. This 3-group is
in particular the right codomain for incorporating the
the non-fake flat nonabelian gerbes with connection considered
in \cite{BrM} into the integrated version of the picture 
discussed here. This is indicated in \cite{SWIII}
and should be discussed elsewhere.

\begin{figure}[h]
  $$
    \xymatrix@C=5pt{
      \mbox{tangent category}
      &
      \mbox{\begin{tabular}{c} inner automorphism \\ $(n+1)$-group \end{tabular}}
      &
      \mbox{\begin{tabular}{c} inner derivation  \\ Lie $(n+1)$-algebra \end{tabular}}
      &
      \mbox{Weil algebra}
      &
      \mbox{\begin{tabular}{c} shifted\\ tangent bundle \end{tabular}}
      \\
      \mathrm{CE}(\mathrm{Lie}(T \mathbf{B}G))
      \ar@{=}[r]
      &
      \mathrm{CE}(\mathrm{Lie}(\mathrm{INN}(G)))
      \ar@{=}[r]
      &
      \mathrm{CE}(\mathrm{inn}(\gg))
      \ar@{=}[r]
      &
      \mathrm{W}(\gg)
      \ar@{=}[r]
      &
      C^\infty(T[1] \gg)
    }
  $$
  \caption{
    \label{a remarkable coincidence of concepts}
    {\bf A remarkable coincidence of concepts} relates the notion of tangency to the
    notion of universal bundles. The leftmost equality
    is discussed in \cite{RS}. The second one from the right is
    the identification \ref{Weil and inner}. The rightmost
    equality is equation  \ref{Weil algebra as functions on tangent space}.
  }
\end{figure}

\begin{proposition}
  For $\gg$ any finite dimensional $L_\infty$-algebra,
  the differential forms 
  on the smooth space of morphisms from the Chevalley-Eilenberg algebra
  $\mathrm{CE}(\gg)$ to the algebra of ``functions on the superpoint'',
  definition \ref{superpoint}, i.e. the elements in
  $\mathrm{maps}(\mathrm{CE}(\gg),C(\mathbf{pt}))$, which
  come from currents as in definition \ref{currents},
  form the Weil algebra $\mathrm{W}(\gg)$ of $\gg$:
  \(
    \mathrm{W}(\gg) \subset \mathrm{maps}(\mathrm{CE}(\gg),C(\mathbf{pt}))
    \,.
  \)
\end{proposition}
\proof
  For any test domain $U$, an element in 
  $\mathrm{Hom}_{\mathrm{DGCAs}}(\mathrm{CE}(\gg),C(\mathbf{pt})\otimes \Omega^\bullet(U))$
  is specified by a degree 0 algebra homomorphism
  \(
    \lambda : \mathrm{CE}(\gg) \to \Omega^\bullet(U)
  \)
  and a degree +1 algebra morphism
  \(
    \lambda : \mathrm{CE}(\gg) \to \Omega^\bullet(U)
  \)
  by 
  \(
    a \mapsto \lambda(a) + c \wedge \omega(a)
  \) 
  for all $a \in \gg^*$ and for $c$ denoting the canonical degree -1 generator of
  $C(\mathbf{pt})$; such that the equality in the bottom right corner of the diagram
\(
  \xymatrix{
    a \ar@{|->}[rr]^{d_{\mathrm{CE}(\gg)}}
    \ar@{|->}[dd]
    && d_{\mathrm{CE}(\gg)} a
    \ar@{|->}[dd]
    \\
    \\
    \lambda(a) 
    +
    c \wedge \omega(a)
    \ar@{|->}[rr]^{d_{\mathbf{pt}} + d_U}
    &&
    {
        {\lambda(d_{\mathrm{CE}(\gg)}a) + c \wedge \omega(d_{\mathrm{CE}(\gg)}a)}
        \atop
     {=d (\lambda(a))
    +
    \omega(a)
    -
    c \wedge d(\omega(a))
    }        
    }
   }
\)
holds. Under the two canonical currents on $C(\mathbf{pt})$ of degree 0 and degree 1, respectively, 
this gives rise for each $a \in \gg^*$ of degree $|a|$ to an $|a|$-form and an $(|a|+1)$ form
on $\mathrm{maps}(\mathrm{CE}(\gg),C(\mathbf{pt}))$ whose values on a given plot are
$\lambda(a)$ and $\omega(a)$, respectively.

By the above diagram, the differential of these forms satisfies
\(
  d \lambda(a) = \lambda(d_{\mathrm{CE}(\gg)} a) + \omega(a)
\)
and
\(
  d \omega(a) = -\omega(d_{\mathrm{CE}(\gg)} a)
  \,.
\)
But this is precisely the structure of $\mathrm{W}(\gg)$.
\endofproof
To see the last step, it may be helpful to consider this for a simple case
in terms of a basis:

let $\gg$ be an ordinary Lie algebra, $\{t^a\}$ a basis of $\gg^*$ and $\{C^a{}_{bc}\}$
the corresponding structure constants. Then, using the fact that, since we are dealing with
algebra homomorphisms, we have
\(
  \lambda (t^a \wedge t^b) = \lambda(t^a) \wedge \lambda(t^b)
\)
and
\(
  \omega (t^a \wedge t^b) = c \wedge (\omega(t^a) \wedge \lambda(t^b) - \lambda(t^a)\wedge \omega(t^b))
\)
we find
\(
  d \lambda(t^a) = -\frac{1}{2}C^a{}_{bc} \lambda(t^b)\wedge \lambda(t^c) + \omega(t^a)
\)
and
\(
  d \omega(t^a) = - C^a{}_{bc} \lambda(t^b) \wedge \omega(t^c)
  \,.
\)
This is clearly just the structure of $\mathrm{W}(\gg)$.

\paragraph{Remark.} As usual, we may think of the superpoint as an ``infinitesimal interval''.
The above says that the algebra of inner derivations of an $L_\infty$-algebra consists
of the maps from the infinitesimal interval to the supermanifold on which 
$\mathrm{CE}(\gg)$ is the ``algebra of functions''. On the one hand this tells us that
\(
  \label{Weil algebra as functions on tangent space}
  \mathrm{W}(\gg) = C^\infty(T[1]\gg)
\)
in supermanifold language. On the other hand, this construction is clearly analogous to the
corresponding discussion for Lie $n$-groups given in \cite{RS}: there the 3-group 
$\mathrm{INN}(G)$ of inner automorphisms of a strict 2-group $G$ was obtained by
mapping the ``fat point'' groupoid $\mathbf{pt} = \{\xymatrix{\bullet \ar[r] & \circ}\}$
into $G$. As indicated there, this is a special  case of a construction of ``tangent categories''
which mimics the relation between $\mathrm{inn}(\gg)$ and the shifted tangent bundle
$T[1]\gg$ in the integrated world of Lie $\infty$-groups. 
This relation between these concepts is summarized in figure \ref{a remarkable coincidence of concepts}.

\subsection{$L_\infty$-algebra homotopy and concordance}

\label{homotopies and concordances}

Like cochain complexes, differental graded algebras can
be thought of as being objects in a higher categorical
structure, which manifests itself in the fact that there
are not only morphisms between DGCAs, but also higher
morphisms between these morphisms. It turns out that
we need to consider a couple of slightly differing notions
of morphisms and higher morphisms for these. While differing, 
these concepts are closely related among each other, as we
shall discuss.

In \ref{homotopies and inner derivations} we had already
considered 2-morphisms of DGCAs obtained after forgetting
their algebra structure and just remembering their 
differential structure. The 2-morphisms we present
now crucially do know about the algebra structure. 

\begin{table}[h]
\begin{center}
\begin{tabular}{c|cc}
  & {\bf name} & {\bf nature}
  \\
  \hline
  {\bf infinitesimal} & transformation & \begin{tabular}{c} chain homotopy
                                      \\
     $$
      \xymatrix{
        \mathrm{CE}(\gg)
        &
        \mathrm{CE}(\hh)
        \ar@/^1.4pc/[l]^{g^*}_{\ }="t"
        \ar@/_1.4pc/[l]_{f^*}^{\ }="s"
        \ar@{=>}^\eta "s"; "t"
      }
    $$
  \end{tabular}
  \\
  \hline
  {\bf finite} & homotopy/concordance & 
  \begin{tabular}{c}
   extension over interval
   \\
  $\xymatrix{
     \mathrm{CE}(\gg)
     &
     \mathrm{CE}(\gg) \otimes \Omega^\bullet(I)
     \ar@<+2pt>[l]^{\mathrm{Id}\otimes t^*}
     \ar@<-2pt>[l]_{\mathrm{Id}\otimes s^*}
     &
     \mathrm{CE}(\hh)
     \ar[l]|<<<<{\eta^*}
     \ar@/_2pc/[ll]_{g^*}
     \ar@/^2pc/[ll]^{f^*}
   }
  $
  \end{tabular}
\end{tabular}
\end{center}
\caption{
  The two different notions of {\bf higher morphisms} of qDGCAs.
  \label{table with higher morphisms}
}
\end{table}

\paragraph{Infinitesimal homotopies between dg-algebra homomorphisms.}

When we restrict attention to cochain maps between qDGCAs which
respect not only the differentials but also the free graded commutative
algebra structure, i.e. to qDGCA homomorphisms,
it becomes of interest to express the cochain homotopies in terms
of their action on generators of the algebra.
We now define transformations (2-morphisms) between morphisms of 
qDGCAs by
first defining them for the case when the domain is a Weil algebra,
and then extending the definition to arbitrary qDGCAs.

\begin{definition}[transformation of morphisms of $L_\infty$-algebras]
  \label{transformation of qDGCA morphisms}
  We define transformations between qDGCA morphisms in two steps
  \begin{itemize}
    \item
    A 2-morphism
    \(
      \xymatrix{
        \mathrm{CE}(\gg)
        &&
        \mathrm{F}(\hh)
        \ar@/^2pc/[ll]^{g^*}_{\ }="t"
        \ar@/_2pc/[ll]_{f^*}^{\ }="s"
        \ar@{=>}^\eta "s"; "t"
      }
    \)
    is defined by a degree -1 map $\eta : \hh^*\oplus \hh^*[1] 
       \to \mathrm{CE}(\gg)$ which is extended to 
     a linear degree -1 map 
     $\eta : \wedge^\bullet(\hh^* \oplus \hh^*[1]) 
        \to \mathrm{CE}(\gg)$
    by defining it on all monomials of generators by the formula
  $$
    \eta : 
       x_1 \wedge \cdots \wedge x_n
       \mapsto
  $$
       \(
       \frac{1}{n!} 
       \sum\limits_{\sigma}
       \epsilon(\sigma)
       \sum\limits_{k=1}^n
       (-1)^{\sum\limits_{i=1}^{k-1} |x_{\sigma(i)}|}
       g^*(x_{\sigma(1)} \wedge \cdots \wedge x_{\sigma(k-1)})
       \wedge \eta(x_{\sigma(k)})
       \wedge
       f^*(x_{\sigma(k+1)} \wedge \cdots \wedge x_{\sigma(n)})
       \label{formula for chain homotopy}
       \)
  for all $x_1, \cdots, x_n \in \hh^* \oplus \hh^*[1]$,
  such that this is a chain homotopy from $f^*$ to $g^*$:
  \(
    g^* = f^* + [d,\eta]
    \,.
  \)

  \item
    A general 2-morphism
    \(
      \xymatrix{
        \mathrm{CE}(\gg)
        &&
        \mathrm{CE}(\hh)
        \ar@/^2pc/[ll]^{g^*}_{\ }="t"
        \ar@/_2pc/[ll]_{f^*}^{\ }="s"
        \ar@{=>}^\eta "s"; "t"
      }
    \)
    is a 2-morphism 
       \(
         \label{2-morphism by pull-back to free qDGCA}
         \raisebox{30pt}{
         \xymatrix{
           & \mathrm{CE}(\hh)
           \ar@/_1pc/[dl]_{g^*}
           \\
			\mathrm{CE}(\gg)
            && 
            \mathrm{W}(\hh)
            \ar@{->>}[ul]_{i^*}^{\ }="s"
            \ar@{->>}[dl]^{i^*}
            &
            \mathrm{F}(\hh)  
            \ar[l]_{\simeq}
		   \\
		   & \mathrm{CE}(\hh)
           \ar@/^1pc/[ul]^{f^*}_{\ }="t"
		   \ar@{=>} "s"; "t"
         }
         }
       \)
      of the above kind
      that vanishes on the shifted generators, i.e. such that
       \(
         \label{vanishing condition on homotopies}
         \raisebox{30pt}{
         \xymatrix{
           & \mathrm{CE}(\hh)
           \ar@/_1pc/[dl]_{g^*}
           \\
			\mathrm{CE}(\gg)
            && 
            \mathrm{W}(\hh)
            \ar@{->>}[ul]_{i^*}^{\ }="s"
            \ar@{->>}[dl]^{i^*}
            &
            \hh^*[1]
            \ar@{_{(}->}[l]
		   \\
		   & \mathrm{CE}(\hh)
           \ar@/^1pc/[ul]^{f^*}_{\ }="t"
		   \ar@{=>} "s"; "t"
         }
         }
       \)
      vanishes.
  \end{itemize}
\end{definition}

\begin{proposition}
  Formula \ref{formula for chain homotopy} is consistent
  in that $g^*|_{\hh^* \oplus \hh^*[1]} = 
    (f^* + [d,\eta])|_{\hh^* \oplus \hh^*[1]}$ 
    implies that $g^* = f^* + [d,\eta]$ on all elements
    of $F(\hh)$. 
\end{proposition}

\paragraph{Remark.}
Definition \ref{transformation of qDGCA morphisms}, 
  which may look ad hoc at this point, 
  has a practical and a deep conceptual 
  motivation. 
\begin{itemize}
  \item
     {\bf Practical motivation.}
      While it is clear that 2-morphisms of qDGCAs should
      be chain homotopies, it is not straightforward, in general,
      to characterize these by their action on generators.
      Except when the domain qDGCA is free, in which case
      our formula \ref{transformation of qDGCA morphisms}
      makes sense. 
      The prescription \ref{2-morphism by pull-back to free qDGCA}
      then provides a systematic algorithm for extending this 
      to arbitrary qDGCAs. 

      In particular, using the isomorphism
      $\mathrm{W}(\gg) \simeq \mathrm{F}(\gg)$ 
       from proposition \ref{properties of Weil algebra}, the above
      yields the usual explicit description of the homotopy
      operator $\tau : \mathrm{W}(\gg) \to \mathrm{W}(\gg)$
      with $\mathrm{Id}_{\mathrm{W}(\gg)} = 
      [d_{\mathrm{W}(\gg)},\tau]$. Among other things, this
      computes for us the transgression elements
      (``Chern-Simons elements'') for $L_\infty$-algbras
      in \ref{Lie infty-algebra cohomology}.

   \item 
     {\bf Conceptual motivation.} As we will see in       
     \ref{Lie infty-algebra cohomology}
     and \ref{Lie infty-algebra valued forms}, the qDGCA
      $\mathrm{W}(\gg)$ plays an important twofold role:
      it is both the algebra of differential forms
       on the total space of the universal $\gg$-bundle
       -- while $\mathrm{CE}(\gg)$ is that of forms on the
       fiber --,
     as well as the domain for $\gg$-valued differential
     forms, where the shifted component, that in $\gg^*[1]$,
     is the home of the corresponding curvature.

     In the light of this, the above restriction
     \ref{vanishing condition on homotopies} can be 
     understood as saying either that
     \begin{itemize}
        \item \emph{vertical} transformations induce transformations on the fibers;
     \end{itemize}
     or
     \begin{itemize}
        \item gauge transformations of $\gg$-valued
         forms are transformations under which the
         curvatures transform covariantly.
     \end{itemize}
\end{itemize}

\paragraph{Finite transformations between qDGCA morphisms: concordances.}

We now consider the finite transformations of morphisms
of DGCAs.
What we called 2-morphisms or transformations for qDGCAs
above would in other contexts possibly be called a homotopy.
Also the following concept is a kind of homotopy, and
appears as such in \cite{StahseffSchlessinger} which goes back to
\cite{BousfieldGuggenheim}. Here we wish to
clearly distinguish these different kinds of homotopies
and address the following concept as
\emph{concordance} -- 
a finite notion of 2-morphism between dg-algebra morphisms.

\paragraph{Remark.}
  In the following the algebra of forms $\Omega^\bullet(I)$ on the interval
  $$
    I := [0,1]
  $$
  plays an important role. Essentially everything would also go through
  by instead using $F(\mathbb{R})$, the DGCA on a single degree 0 generator,
  which is the algebra of \emph{polynomial} forms on the interval. This is the model 
  used in \cite{StahseffSchlessinger}.  

\begin{definition}[concordance]

  We say that two qDGCA morphisms
  \(
    \xymatrix{
       \mathrm{CE}(\gg)
       &&
       \mathrm{CE}(\hh)
       \ar[ll]_{g^*}
    }
  \)  
  and
  \(
    \xymatrix{
       \mathrm{CE}(\gg)
       &&
       \mathrm{CE}(\hh)
       \ar[ll]_{h^*}
    }
  \)  
  are concordant, if there exists a dg-algebra homomorphism
  \(
    \xymatrix{
       \mathrm{CE}(\gg)
       \otimes \Omega^\bullet(I)
       &&
       \mathrm{CE}(\hh)
       \ar[ll]_{\eta^*} 
    }
  \)
  from the source $\mathrm{CE}(\hh)$ to the the target $\mathrm{CE}(\gg)$
  tensored with forms on the interval,
  which restricts to the two given homomorphisms 
  when pulled back along the two boundary inclusions
  \(
    \xymatrix{
       \{\bullet\} 
       \ar@<+3pt>[r]^s
       \ar@<-3pt>[r]_t
       &
       I
    }
    \,,
  \)
  so that the diagram of dg-algebra morphisms
  \(
  \xymatrix{
     \mathrm{CE}(\gg)
     &
     \mathrm{CE}(\gg) \otimes \Omega^\bullet(I)
     \ar@<+2pt>[l]^{\mathrm{Id}\otimes t^*}
     \ar@<-2pt>[l]_{\mathrm{Id}\otimes s^*}
     &
     \mathrm{CE}(\hh)
     \ar[l]_<<<<{\eta^*}
     \ar@/_2pc/[ll]_{g^*}
     \ar@/^2pc/[ll]^{f^*}
   }
  \)
  commutes. 
  \label{concordance of qDGCA morphisms}
\end{definition}

See also table \ref{table with higher morphisms}.
Notice that the above diagram is shorthand for two separate
commuting diagrams, one involving $g^*$ and $s^*$, the other
involving $f^*$ and $t^*$.

We can make precise the statement that definition \ref{transformation of qDGCA morphisms}
is the infinitesimal version of definition \ref{concordance of qDGCA morphisms}, as follows.

\begin{proposition}
  \label{concordance and chain homotopy}
  Concordances 
  \(
    \xymatrix{
      \mathrm{CE}(\gg)
      \otimes \Omega^\bullet(I)
      &&
      \mathrm{CE}(\hh)
      \ar[ll]_{\eta^*}
    }
  \)
  are in bijection with 1-parameter families 
  \(
    \alpha : [0,1] \to \mathrm{Hom}_{\mathrm{dg-Alg}}(\mathrm{CE}(\hh),\mathrm{CE}(\gg))
  \)
  of morphisms whose derivatives with respect to the parameter is a 
  chain homotopy, i.e. a 2-morphism
  \(
    \forall t \in [0,1]
     \hspace{5pt}
     :
     \hspace{5pt}
     \xymatrix{
       \mathrm{CE}(\gg)
       &&
       \mathrm{CE}(\hh)
       \ar@/^2pc/[ll]^{\frac{d}{dt}\alpha(t) = [d,\rho]}_{\ }="t"
       \ar@/_2pc/[ll]_{0}^{\ }="s"       
       \ar@{=>}^\rho "s"; "t"
     }
  \)
  in the 2-category of cochain complexes.
For any such $\alpha$, the morphisms $f^*$ and $g^*$ between which
it defines a concordance are defined by the value of $\alpha$ on the
boundary of the interval.
\end{proposition}

\proof
Writing $t : [0,1] \to \mathbb{R}$ for the canonical coordinate function
on the interval $I = [0,1]$ we can decompose 
the dg-algebra homomorphism 
$\eta^*$ as
\(
  \eta^* : \omega \mapsto (t \mapsto \alpha(\omega)(t) + dt \wedge \rho(\omega)(t) ) 
  \,.
\)
$\alpha$ is itself a degree 0 dg-algebra 
homomorphism, while $\rho$ is degree -1 map.
Then the fact that $\eta^*$ respects the differentials implies that for all
$\omega \in \mathrm{CE}(\hh)$ we have
\(
  \xymatrix{
    \omega
    \ar@{|->}[rr]^{d_\hh}
    \ar@{|->}[dd]^{\eta^*}
    &&
    d_\hh \omega
    \ar@{|->}[dd]^{\eta^*}
    \\
    \\
    (t \mapsto (\alpha(\omega)(t)
    +
    dt \wedge \rho(\omega)(t)))
    \ar@{|->}[rr]^{d_\gg + d_t}
    &&
    **[r]{\begin{array}{l}
      (t \mapsto (\alpha(d_\hh \omega)(t) + dt \wedge \rho(d_\hh\omega)(t))
      \\
      =
      ( t \mapsto (d_\gg (\alpha(\omega))(t) \\
      + dt \wedge (\frac{d}{dt}\alpha(\omega) - d_\gg \rho(\omega))(t))
    \end{array}}
  }
  \,.
\)
The equality in the bottom right corner says that
\(
  \alpha \circ d_\hh - d_\gg \circ\alpha = 0
\)
and
\(
  \forall \omega \in \mathrm{CE}(\gg)
  :
  \frac{d}{dt}\alpha(\omega)
  =
  \rho (d_\hh\omega) + d_\gg (\rho (\omega))
  \,.
\)
But this means that $\alpha$ is a chain homomorphism whose derivative
is given by a chain homotopy.
\endofproof

\subsubsection{Examples}

\paragraph{Transformations between DGCA morphisms.}

We demonstrate two examples for the application of the 
notion of transformations of DGCA morphisms from definition
\ref{transformation of qDGCA morphisms} which are
relevant for us.

\subparagraph{Computation of transgression forms.}
As an example for the transformation in 
definition \ref{transformation of qDGCA morphisms}, 
we show how the usual Chern-Simons transgression
form is computed using formula \ref{formula for chain homotopy}.
The reader may wish to first skip to our discussion of 
Lie $\infty$-algebra cohomology in \ref{Lie infty-algebra cohomology}
for more background.
So let $\gg$ be an ordinary Lie algebra with invariant bilinear
form $P$, which we regard as a $d_{\mathrm{W}(\gg)}$-closed
element $P \in \wedge^2 \gg^*[1] \subset \mathrm{W}(\gg)$. 
We would like to compute $\tau P$, where $\tau$ is the contracting
homotopy of $\mathrm{W}(\gg)$, such that
\(
  [d,\tau] = \mathrm{Id}_{\mathrm{W}(\gg)}
  \,,
\)
which according to 
proposition \ref{properties of Weil algebra} is given on generators
by
\begin{eqnarray}
  \tau &:& a \mapsto 0
  \\
  \tau &:& d_{\mathrm{W}(\gg)}a \mapsto a
\end{eqnarray}
for all $a \in \gg^*$. Let $\{t^a\}$ be a chosen basis of
$\gg^*$ and let $\{P_{ab}\}$ be the components of $P$ in that
basis, then
\(
  P = P_{ab} (\sigma t^a) \wedge (\sigma t^b)
  \,.
\)
In order to apply formula \ref{formula for chain homotopy}
we need to first rewrite this in terms of monomials in 
$\{t^a\}$ and $\{d_{\mathrm{W}(\gg)}t^a\}$. Hence, using
$\sigma t^a = d_{\mathrm{W}(\gg)}t^a + 
\frac{1}{2}C^a{}_{bc}t^b \wedge t^c$, we get
\(
  \tau P =
  \tau
  \left(
    P_{ab} (d_{\mathrm{W}(\gg)} t^a)
       \wedge (d_{\mathrm{W}(\gg)} t^a)
    - P_{ab} (d_{\mathrm{W}(\gg)} t^a)
       \wedge C^b{}_{cd} t^c \wedge t^d
    +
    \frac{1}{4}P_{ab}C^a{}_{cd} C^b_{ef}t^c \wedge t^d \wedge 
     t^c \wedge t^d
  \right)
  \,.
\)
Now equation \ref{formula for chain homotopy} can be applied to
each term. Noticing the combinatorial prefactor $\frac{1}{n!}$,
which depends on the number of factors in the above terms,
and noticing the sum over all permutations, we find
\begin{eqnarray}
  \tau \left(P_{ab} (d_{\mathrm{W}(\gg)} t^a)
       \wedge (d_{\mathrm{W}(\gg)} t^a) 
    \right)
  &=& P_{ab} (d_{\mathrm{W}(\gg)} t^a) \wedge t^b
\nonumber\\
  \tau \left(
    - P_{ab} (d_{\mathrm{W}(\gg)} t^a)
       \wedge C^b{}_{cd} t^c \wedge t^d
  \right)
  &=&
  \frac{1}{3!}\cdot 2\; P_{ab}C^b_{cd} t^b \wedge t^c \wedge t^d
  =
  \frac{1}{3} C_{abc}t^a \wedge t^b \wedge t^c
  \,,
\end{eqnarray}
where we write $C_{abc} := P_{ad}C^d{}_{bc}$ as usual.
Finally
$
   \tau\left( \frac{1}{4}P_{ab}C^a{}_{cd} C^b_{ef}t^c \wedge t^d \wedge 
     t^c \wedge t^d \right)
   = 0
  \,.
$
In total this yields
\(
  \tau P =
  P_{ab} (d_{\mathrm{W}(\gg)} t^a) \wedge t^b
  +
  \frac{1}{3} C_{abc} t^a \wedge t^b \wedge t^c
  \,.
\)
By again using $d_{\mathrm{W}(\gg)} t^a = 
-\frac{1}{2}C^a{}_{bc}t^b \wedge t^c + \sigma t^a$
together with the invariance of $P$ 
(hence the $d_{\mathrm{W}(\gg)}$-closedness of $P$ which
implies that the constants $C_{abc}$ are skew symmetric
in all three indices),
one checks that this does indeed satisfy
\(
  d_{\mathrm{W}(\gg)} \tau P = P
  \,.
\)
In \ref{Lie infty-algebra valued forms} we will see that
after choosing a $\gg$-valued connection on the space $Y$
the generators $t^a$ here will get sent to components of
a $\gg$-valued 1-form $A$, while the $d_{\mathrm{W}(\gg)} t^a$
will get sent to the components of $dA$. Under this map the element $\tau P \in \mathrm{W}(\gg)$
maps to the familiar Chern-Simons 3-form
\(
  \mathrm{CS}_P(A) :=
  P(A \wedge d A) + \frac{1}{3}
  P(A \wedge [A \wedge A])
\)
whose differential is the characteristic form of $A$
with respect to $P$:
\(
  d \mathrm{CS}_P(A) = P(F_A \wedge F_A)
  \,.
\)
Characteristic forms, for arbitrary Lie $\infty$-algebra
valued forms, are discussed 
further in \ref{characteristic forms}.

\subparagraph{2-Morphisms of Lie 2-algebras}

\begin{proposition}
  For the special case that $\gg$ is any Lie 2-algebra 
  (any $L_\infty$-algebra
  concentrated in the first two degrees) the 2-morphisms
  defined by definition \ref{transformation of qDGCA morphisms}
  reproduce the 2-morphisms of Lie 2-algebras 
  as stated in \cite{BaC} and used in \cite{BCSS}.
\end{proposition}
\proof
  The proof is given in the appendix.
\endofproof

This implies in particular that with the 1- and 
2-morphisms as defined above, Lie 2-algebras do form
a 2-category. There is an rather straightforward
generalization of definition \ref{transformation of qDGCA morphisms}
to higher morphisms, which one would expect yields 
correspondingly $n$-categories of Lie $n$-algebras. 
But this we shall not try to discuss here.

\subsection{$L_\infty$-algebra cohomology}

\label{Lie infty-algebra cohomology}

The study of ordinary Lie algebra cohomology and of invariant 
polynomials on the Lie algebra has a simple formulation in terms of 
the qDGCAs $\mathrm{CE}(\gg)$ and $\mathrm{W}(\gg)$. Furthermore, 
this has a straightforward
generalization to arbitrary $L_\infty$-algebras which we now state.

  For 
 $ \xymatrix{
      \mathrm{CE}(\gg)
      &&
      \mathrm{W}(\gg)
      \ar@{->>}[ll]_{i^*}
    }$
  the canonical morphism from proposition \ref{properties of Weil algebra},
  notice that 
  \(
    \mathrm{CE}(\gg) \simeq \mathrm{W}(\gg)/\mathrm{ker}(i^*)
  \)
  and that
  \(
     \mathrm{ker}(i^*) = \langle \gg^*[1]\rangle_{\mathrm{W}(\gg)}
     \,,
  \)
 \noindent
 the ideal in $\mathrm{W}(\gg)$ generated by $ \gg^*[1].$
  Algebra derivations
  \(
    \iota_X : \mathrm{W}(\gg) \to \mathrm{W}(\gg)
  \)
  for $X\in\gg$ 
  are like (contractions with) vector fields on the space   
  on which $\mathrm{W}(\gg)$
  is like differential forms. In the case of an ordinary Lie algebra $\gg$, 
  the corresponding inner
  derivations $[d_{\mathrm{W}(\gg)},\iota_X]$ for $X\in\gg$ are of degree
  -1 and are known as the Lie derivative $L_X.$ They generate 
  flows $\exp([d_{\mathrm{W}(\gg)},\iota_X]) : 
  \mathrm{W}(\gg) \to \mathrm{W}(\gg)$ along these vector
  fields.

  \begin{definition}[vertical derivations]
     \label{vertical derivations special}
    We say an algebra derivation $\tau : \mathrm{W}(\gg) \to \mathrm{W}(\gg)$
    is \emph{vertical} if it vanishes on the shifted copy
    $\gg^*[1]$ of $\gg^*$ inside $\mathrm{W}(\gg)$,
    \(
      \tau|_{\gg^*[1]} = 0
      \,.
    \)
  \end{definition}

 \begin{proposition}
   The vertical derivations are precisely those that come from 
   contractions
   \(
     \iota_X : \gg^* \mapsto \mathbb{R}
   \)
   for all $X \in \gg$, extended to 0 on $\gg^*[1]$ and extended
   as algebra derivations to all of 
   $\wedge^\bullet(\gg^* \oplus \gg^*[1])$.
 \end{proposition}

The reader should compare this and the following definitions
to the theory of vertical Lie derivatives and basic 
differential forms with respect to any surjective submersion
$\pi : Y \to X$. This is discussed in 
\ref{surjective submersions and differential forms}.

\begin{definition}[basic forms and invariant polynomials]
The algebra $\mathrm{W}(\gg)_{\mathrm{basic}}$ of  {\bf basic forms} 
in $\mathrm{W}(\gg)$ is the intersection of the kernels of all 
vertical derivations and Lie derivatives. i.e. all the 
contractions $\iota_X$ and Lie derivatives $L_X$ for $X\in \gg.$
Since $L_X = [d_{\mathrm{W}(\gg)}, \iota_X]$, it follows that in the 
kernel of $\iota_X$, the Lie derivative vanishes 
only if $\iota_Xd_{\mathrm{W}(\gg)}$ vanishes.
  \label{invariant polynomials}
\end{definition}

\begin{figure}[h]
  $$
    \xymatrix{
      \mathrm{CE}(\gg)
      &&
      \mathrm{CE}(\gg)
      \ar@/^{1.3pc}/[ll]^{[d,\iota_X]}_{\ }="t2"
      \ar@/_{1.3pc}/[ll]_{0}^{\ }="s2"
      &
      \mbox{
        \begin{tabular}{l}
          (co)adjoint action
          \\
          of $\gg$ on itself
        \end{tabular}
      }
      \\
      \mathrm{W}(\gg) 
      \ar@{->>}[u]
      &&
      \mathrm{W}(\gg)
      \ar@{->>}[u]
      \ar@/^{1.3pc}/[ll]^{[d,\iota_X]}_{\ }="t"
      \ar@/_{1.3pc}/[ll]_{0}^{\ }="s"
      &
      \mbox{
        \begin{tabular}{l}
          vertical derivation
          \\
          on $\mathrm{W}(\gg)$
        \end{tabular}
      }
      \\
      \mathrm{inv}(\gg)
      \ar@{^{(}->}[u]
      &&
      \mathrm{inv}(\gg)
      \ar@{^{(}->}[u]
      \ar@/^{1.3pc}/[ll]^{0}_{\ }="t1"
      \ar@/_{1.3pc}/[ll]_{0}^{\ }="s1"
      &
      \mbox{
        \begin{tabular}{l}
          leaves basic forms
          \\
          invariant
        \end{tabular}
      }
      \ar@{=>}^{\iota_X} "s"; "t"
      \ar@{=>}^{0} "s1"; "t1"
      \ar@{=>}^{\iota_X} "s2"; "t2"
    }
    \,.
  $$
  \caption{ 
    \label{interpretation of vertical derivations on Wg}
   {\bf Interpretation of vertical derivations on 
   $\mathrm{W}(\gg)$.}
   The algebra $\mathrm{CE}(\gg)$ plays the role of the
   algebra of differential forms on the Lie $\infty$-group
   that integrates the Lie $\infty$-algebra $\gg$. 
   The coadjoint action of $\gg$ on these forms corresponds to
   Lie derivatives along the fibers of the universal bundle.
   These vertical derivatives leave the forms on the base
   of this universal bundle invariant.
   The diagram displayed is in the 2-category $\mathbf{Ch}^\bullet$
   of cochain complexes, as described in the beginning of
   \ref{homotopies and concordances}.
  }
  \end{figure}
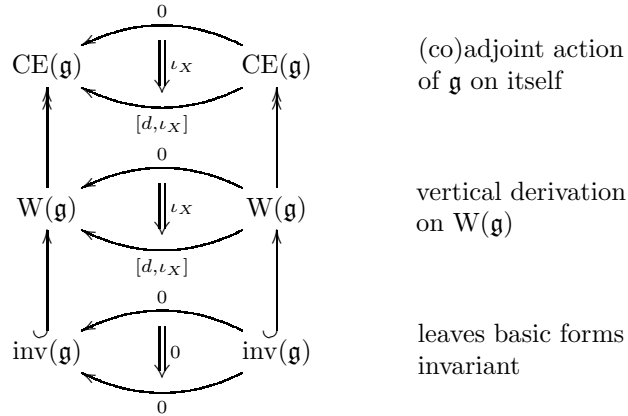

As will be discussed in a moment, basic forms in $\mathrm{W}(\gg)$ play the 
role of {\bf invariant polynomials} on the $L_\infty$-algebra
$\gg$. Therefore we often write $\mathrm{inv}(\gg)$ for 
$\mathrm{W}(\gg)_{\mathrm{basic}}$:
\(
  \mathrm{inv}(\gg) := \mathrm{W}(\gg)_{\mathrm{basic}}
  \,.
\)
Using the obvious inclusion 
$\xymatrix{\mathrm{W}(\gg) & \mathrm{inv}(\gg)\ar@{_{(}->}[l]_{p^*}}$
we obtain the sequence
\(
  \label{sequence which is dgc-version of universal bundle}
  \xymatrix{
     \mathrm{CE}(\gg)
     &&
     \mathrm{W}(\gg)
     \ar@{->>}[ll]_{i^*}
     &&
     \mathrm{inv}(\gg)
     \ar@{_{(}->}[ll]_{p^*}
  }
\)
of dg-algebras that plays a major role in our analysis: it
can be interpreted as coming from the universal bundle for the Lie 
$\infty$-algebra $\gg$. 
As shown in figure \ref{interpretation of vertical derivations on Wg},
we can regard vertical derivations on $\mathrm{W}(\gg)$
as derivations along the \emph{fibers} of the corresponding
dual sequence.

\begin{definition}[cocycles, invariant polynomials and transgression elements]
  Let $\gg$ be an $L_\infty$-algebra. Then
  \begin{itemize}
    \item
      An $L_\infty$-algebra {\bf cocycle} on $\gg$ 
       is a $d_{\mathrm{CE}(\gg)}$-closed element of $\mathrm{CE}(\gg)$.
      \(
        \mu \in \mathrm{CE}(\gg)\,, \hspace{20pt} d_{\mathrm{CE}(\gg)} \mu = 0
        \,.
      \)
    \item
      An $L_\infty$-algebra {\bf invariant polynomial} on $\gg$
      is an element 
      $P \in \mathrm{inv}(\gg) := {\mathrm{W}}(\gg)_{basic}$.
    \item
      An $L_\infty$-algebra 
      {\bf $\gg$-transgression element} for a given cocycle $\mu$ and
     an invariant polynomial $P$
      is an element $\mathrm{cs} \in \mathrm{W}(\gg)$ such that
      \(
        d_{\mathrm{W}(\gg)} \mathrm{cs} = p^*P
      \)
      \(
        i^* \mathrm{cs} = \mu
        \,.
      \)
  \end{itemize}
\end{definition}

\begin{figure}[h]
$$
  \xymatrix{
    \mbox{
      cocycle
    }
    &
    \mbox{transgression element}
    &
    \mbox{inv. polynomial}
    \\
    \mathrm{CE}(\gg)
    & 
    \mathrm{W}(\gg)
    \ar@{->>}[l]_<<<<{i^*} 
    & 
    \mathrm{inv}(\gg)
    \ar@{_{(}->}[l]_<<<<<{p^*}
    \\
    & 0
    \\
    0
    &
    p^* P
    \ar@{|->}[u]_{d_{\mathrm{W}(\gg)}}
    \ar@{|->}@/_1pc/[d]_\tau
    &
    P
    \ar@{|->}[l]^{p^*}
    \\
    \mu
    \ar@{|->}[u]_{d_{\mathrm{CE}(\gg)}}
    &
    \mathrm{cs}
    \ar@{|->}[l]_{i^*}
    \ar@{|->}[u]_{d_{\mathrm{W}(\gg)}}
    &
  }
$$
\caption{{\bf The homotopy operator $\tau$} is a contraction 
  homotopy for $\mathrm{W}(\gg)$. Acting with it
  on a closed invariant polynomial 
  $P \in \mathrm{inv}(\gg) \subset \wedge^\bullet \gg[1] \subset \mathrm{W}(\gg)$
  produces an element $\mathrm{cs} \in W(\gg)$ whose ``restriction to the fiber''
  $ \mu :=  i^* \mathrm{cs}$ is necessarily closed and hence a cocycle. 
  We say that $\mathrm{cs}$ induces the \emph{transgression} from $\mu$
  to $P$, or that $P$ \emph{suspends} to $\mu$.
\label{homotopy operator figure}}
\end{figure}
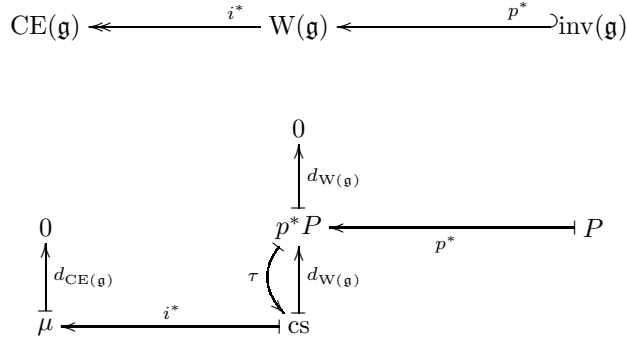

If a transgression element for $\mu$ and $P$ exists, we say
that \emph{$\mu$ {\bf transgresses} to $P$} and that
\emph{$P$ {\bf suspends} to $\mu$}. If $\mu = 0$ we say that
\emph{$P$ {\bf suspends to $0$}}.
The situation is illustrated diagrammatically in figure 
\ref{transgression} and figure \ref{homotopy operator figure}.

\begin{definition}[suspension to 0]
  An element $P \in \mathrm{inv}(\gg)$ is said to suspend to 0
  if under the inclusion
  \(
    \xymatrix{
      \mathrm{ker}(i^*)
      &&
      \mathrm{W}(\gg)
      \ar@{_{(}->}[ll]_{p^*}
    }
  \)
  it becomes a coboundary:
  \(
    p^*P = d_{\mathrm{ker}(i^*)} \alpha
  \)
  for some $\alpha \in \mathrm{ker}(i^*)$.
\end{definition}

\paragraph{Remark.} We will see that it is the intersection of $\mathrm{inv}(\gg)$ with
the cohomology of $\mathrm{ker}(i^*)$ that is a candidate, in general, for an algebraic
model of the classifying space of the object that integrates the $L_\infty$-algebra 
$\gg$. But at the moment we do not entirely clarify this relation to the integrated theory,
unfortunately.

\begin{figure}[h]
$$
  \xymatrix{
    \mbox{
      cocycle
    }
    &
    \mbox{transgression element}
    &
    \mbox{inv. polynomial}
    \\
    G
    \ar@{^{(}->}[r]^i 
    & 
    EG 
    \ar@{->>}[r]^p
    & 
    BG
    \\
    & 0
    \\
    0
    &
    p^* P
    \ar@{|->}[u]_d    
    &
    P
    \ar@{|->}[l]^{p^*}
    \\
    \mu
    \ar@{|->}[u]_d
    &
    \mathrm{cs}
    \ar@{|->}[l]^{i^*}
    \ar@{|->}[u]_d
    &    
  }
$$
\caption{{\bf Lie algebra cocycles, invariant polynomials
and transgression forms} in terms of cohomology of the 
universal $G$-bundle.
Let $G$ be a simply connected compact Lie group with Lie algebra
$\gg$. Then invariant polynomials $P$ on $\gg$ correspond to elements
in the cohomology $H^\bullet(BG)$ of the classifying space of $G$.
When pulled back to the total space of the universal $G$-bundle
$EG \to BG$,
these classes become trivial, due to the contractability of $EG$:
$p^* P = d(\mathrm{cs})$. Lie algebra cocycles, on the other hand,
correspond to elements in the cohomology $H^\bullet(G)$ of $G$ itself.
A cocycle $\mu \in H^\bullet(G)$ is in transgression with an invariant
polynomial $P \in H^\bullet(BG)$ if $\mu = i^* \mathrm{cs}$.
\label{transgression}}
\end{figure}

\begin{proposition}
For the case that $\gg$ is an ordinary Lie
algebra, the above definition reproduces the ordinary definitions 
of Lie algebra cocycles, invariant polynomials, and transgression
elements. Moreover, all elements in $\mathrm{inv}(\gg)$ are closed.
 \label{ordinary Lie cohomology is reproduced}
\end{proposition}
\proof
  That the definitions of Lie algebra cocycles and transgression elements
  coincides is clear. It remains to be checked that $\mathrm{inv}(\gg)$
  really contains the invariant polynomials. In the ordinary definition
  a $\gg$-invariant polynomial is a $d_{\mathrm{W}(\gg)}$-closed element
  in $\wedge^\bullet (\gg^*[1])$. Hence one only needs to check that
  all elements in $\wedge^\bullet (\gg^*[1])$ with the property that their
  image under $d_{\mathrm{W}(\gg)}$ is again in $\wedge^\bullet (\gg^*[1])$
  are in fact already closed. This can be seen for instance in components, 
  using the description of $\mathrm{W}(\gg)$ given in \ref{examples for Loo algebras}.
\endofproof

\paragraph{Remark.} For ordinary Lie algebras $\gg$ corresponding to a
simply connected compact Lie group $G$, the situation is often
discussed in terms of the cohomology of the universal $G$-bundle.
This is recalled in figure \ref{transgression} and in
\ref{examples for Lie 00-algebra cohomology}. The general definition
above is a precise analog of that familiar situation: 
$\mathrm{W}(\gg)$ plays the role of the algebra of 
(left invariant) differential forms
on the universal $\gg$-bundle and 
$\mathrm{CE}(\gg)$ plays the role of the
algebra of (left invariant) differential forms on its 
fiber. Then $\mathrm{inv}(\gg)$ plays the role of differential
forms on the base, $B G = EG/G$.
In fact, for $G$ a compact and
simply connected Lie group and $\gg$ its Lie algebra, we have
\(
  H^\bullet(\mathrm{inv}(\gg)) \simeq H^\bullet(BG, \mathbb{R})
  \,.
\)

In summary, the situation we thus obtain is that depicted in figure 
\ref{universal G-bundle}.
Compare this to the following fact.
\begin{proposition}
  For $p : P \to X$ a principal $G$-bundle, let 
  $\mathrm{vert}(P) \subset \Gamma(TP)$
  be the vertical vector fields on $P$. The horizontal differential
  forms on $P$ which are invariant under $\mathrm{vert}(P)$
  are precisely those that are pulled back along $p$ from $X$.
\end{proposition}
  These are called the {\bf basic differential forms} in 
  \cite{GHV}.

\paragraph{Remark.} We will see that, contrary to the situation for
ordinary Lie algebras, in general invariant polynomials of $L_\infty$ algebras
are not $d_{\mathrm{W}(\gg)}$-closed (the $d_{\mathrm{W}(\gg)}$-differential of
them is just horizontal). We will also see that those indecomposable
invariant polynomials in $\mathrm{inv}(\gg)$, i.e. those that become
exact in $\mathrm{ker}(i^*)$, are not characteristic for the corresponding
$\gg$-bundles. This probably means that the real cohomology of the
classifying space of the Lie $\infty$-group integrating $\gg$ is
spanned by invariant polynomials modulo those suspending to 0. But here
we do not attempt to discuss this further.

\begin{proposition}
  For every invariant polynomial $P \in \wedge^\bullet \gg[1] \subset \mathrm{W}(\gg)$ 
  on an $L_\infty$-algebra $\gg$ such that $d_{\mathrm{W}(\gg)} p^* P = 0$,
  there exists an $L_\infty$-algebra cocycle $\mu \in \mathrm{CS}(\gg)$ that
  transgresses to $P$.
  \label{every closed inv polynomial comes from transgression}
\end{proposition}
  
\proof
  This is a consequence of
  proposition \ref{properties of Weil algebra}
  and 
  proposition \ref{homotopy operator}.
  Let $P \in W(\gg)$ be an invariant polynomial.  By
  proposition \ref{properties of Weil algebra}, $p^* P$ 
   is in the kernel of the restriction homomorphism 
   $\xymatrix{\mathrm{CE}(\gg) & \mathrm{W}(\gg)\ar@{->>}[l]_{i^*}}$: $i^* P = 0$.
   By proposition \ref{homotopy operator}, $p^*P$ is the image under $d_{\mathrm{W}(\gg)}$
   of an element
   $\mathrm{cs} := \tau(p^*P)$
   and by the algebra homomorphism property of $i^*$ we know that 
   its restriction, $\mu := i^* \mathrm{cs}$,
   to the fiber is closed, because
   \(
     d_{\mathrm{CE}(\gg)} i^* \mathrm{cs} = i^* d_{\mathrm{W}(\gg)} \mathrm{cs}
     = i^* p^*P = 0
     \,.
   \)
   Therefore $\mu$ is an $L_\infty$-algebra cocycle for $\gg$ that
   transgresses to the invariant polynomial $P$.
\endofproof

\paragraph{Remark.} Notice that this statement is useful only for 
\emph{indecomposable} invariant polynomials. 
All others trivially suspend to the 0 cocycle.

\begin{proposition}
  An invariant polynomial which suspends to a 
  Lie $\infty$-algebra cocycle that is a coboundary also
  suspends to 0.
  \label{invariant polynomials and trivial cocycles}
\end{proposition}
\proof
  Let $P$ be an invariant polynomial, 
  $\mathrm{cs}$ the corresponding transgression
  element and $\mu = i^* \mathrm{cs}$ the corresponding
  cocycle, which is assumed to be a coboundary in that
  $\mu = d_{\mathrm{CE}(\gg)}b$ for some 
  $b \in \mathrm{CE}(\gg)$.   
  Then by the definition of $d_{\mathrm{W}(\gg)}$
  it follows that $\mu = i^* (d_{\mathrm{W}(\gg)} b)$.
  
  Now notice that
  \(
    \mathrm{cs}' := \mathrm{cs} - d_{\mathrm{W}(\gg)} b
  \)
  is another transgression element for $P$, since 
  \(
     d_{\mathrm{W}(\gg)} \mathrm{cs}' = p^* P
     \,.
  \)
  But now 
  \(
    i^*(\mathrm{cs}') = 
    i^*(\mathrm{cs} - d_{\mathrm{W}(\gg)} b )  = 0
    \,.
  \)
  Hence $P$ suspends to 0.
\endofproof

\subsubsection{Examples}
 \label{examples for Lie 00-algebra cohomology}

\paragraph{The cohomologies of $G$ and of $BG$ in terms of qDGCAs.}

To put our general considerations for $L_\infty$-algebras into
perspective, it is useful to keep the following classical results
for ordinary Lie algebras in mind.

A classical result of E. Cartan \cite{Cartan} \cite{Cartan2}
(see also \cite{Kalkman})
says that 
for a connected finite dimensional Lie group $G$, 
the cohomology $H^\bullet(G)$ of the group is isomorphic 
to that of the Chevalley-Eilenberg algebra $\mathrm{CE}(\gg)$
of its Lie algebra $\gg$:
\(
  H^\bullet(G) \simeq H^\bullet(\mathrm{CE}(\gg)) 
  \,,
\) 
namely to the algebra of Lie algebra cocycles on $\gg$.
If we denote by $Q_G$ the space of \emph{indecomposable} such cocycles, 
and form  the qDGCA $\wedge^\bullet Q_G = H^\bullet(\wedge^\bullet Q_G)$
with trivial differential, 
the above says that we have an isomorphism in cohomology
\(
  H^\bullet(G) \simeq H^\bullet(\wedge^\bullet Q_G) = \wedge^\bullet Q_G
\)
which is realized by the canonical inclusion
\(
 i : 
 \xymatrix{
    \wedge^\bullet Q_G
    \ar@{^{(}->}[r]
    &
    \mathrm{CE}(\gg)
 }
\)
of all cocycles into the Chevalley-Eilenberg algebra.

Subsequently, we have the classical result of Borel \cite{Borel}: 
For a connected finite dimensional Lie group $G$, the cohomology 
of its classifying space $BG$ is
a finitely generated polynomial algebra on
even dimensional generators:
\(
  H^\bullet(BG) \simeq \wedge^\bullet P_G\,.
\)
Here $P_G$ is the space of \emph{indecomposable} invariant polynomials on $\gg$, hence
\(
  H^\bullet(BG) \simeq H^\bullet(\mathrm{inv}(\gg))
  \,.
\)
In fact, $P_G$ and $Q_G$ are isomorphic after a shift:
\(
  P_G \simeq Q_G[1]
\)
and this isomorphism is induced by \emph{transgression} between indecomposable 
cocycles 
$\mu \in \mathrm{CE}(\gg)$ and indecomposable invariant polynomials 
$P \in \mathrm{inv}(\gg)$ via a transgression element 
$\mathrm{cs} = \tau P \in \mathrm{W}(\gg)$.

\paragraph{Cohomology and invariant polynomials of $b^{n-1}\uu(1)$}

\begin{proposition}
  \label{cohomology of shifted u(1)}
  For every integer $n \geq 1$, the Lie $n$-algebra
  $b^{n-1}\uu(1)$ (the $(n-1)$-folded shifted version of
  ordinary $\uu(1)$) from \ref{examples for Loo algebras})
  we have the following:
  \begin{itemize}
     \item
       there is, up to a scalar multiple, a single
       indecomposable Lie $\infty$-algebra cocycle
       which is of degree $n$ and \emph{linear},
       \(
         \mu_{b^{n-1}\uu(1)} \in \mathbb{R}[n]
               \subset \mathrm{CE}(b^{n-1}\uu(1))
         \,,
       \)
     \item
       there is, up to a scalar multiple, a single
       indecomposable Lie $\infty$-algebra invariant
       polynomial, which is of degree $(n+1)$
       \(
         P_{b^{n-1}\uu(1)} \in \mathbb{R}[n+1]
           \subset \mathrm{inv}(b^{n-1}\uu(1)) = \mathrm{CE}(b^n \uu(1))
         \,.
       \)
      \item
        The cocycle $\mu_{b^{n-1}\uu(1)}$ is in transgression with
        $P_{b^{n-1}\uu(1)}$.
  \end{itemize}
\end{proposition}
These statements are an obvious consequence of the
definitions involved, but they are important. The fact
that $b^{n-1}\uu(1)$ has a single invariant polynomial of degree $(n+1)$
will immediately imply, in \ref{Loo Cartan-Ehresmann connections},
that $b^{n-1}\uu(1)$-bundles have a single characteristic
class of degree $(n+1)$: known (at least for $n=2$, as the
Dixmier-Douady class).
Such a $b^{n-1}\uu(1)$-bundle classes appear in
\ref{lifting problem} as the obstruction classes for
lifts of $n$-bundles through string-like extensions of their
structure Lie $n$-algebra.

\paragraph{Cohomology and invariant polynomials of strict Lie 2-algebras.}

  Let $\gg_{(2)} = (\hh \stackrel{t}{\to} \gg \stackrel{\alpha}{\to} \mathrm{der}(\hh))$ 
  be a strict Lie 2-algebra
  as described in section \ref{Lie infty algebras}. 
    Notice that there is a canonical projection homomorphism
    \(
      \xymatrix{
        \mathrm{CE}(\gg)
        &&
        \mathrm{CE}(\hh \stackrel{t}{\to} \gg)
        \ar@{->>}[ll]_{j^*}
      }
    \)
    which, of course, extends to the Weil algebras
    \(
      \xymatrix{
        \mathrm{W}(\gg)
        &&
        \mathrm{W}(\hh \stackrel{t}{\to} \gg)
        \ar@{->>}[ll]_{j^*}
      }
      \,.
    \)
Here $j^*$ is simply the identity on $\gg^*$ and on $\gg^*[1]$ and 
    vanishes on $\hh^*[1]$ and $\hh^*[2]$.

  \begin{proposition}
    Every invariant polynomial $P \in \mathrm{inv}(\gg)$ 
    of the ordinary Lie algebra $\gg$ lifts to
    an invariant polynomial 
    on the Lie 2-algebra $(\hh \stackrel{t}{\to} \gg)$:
    \(
      \raisebox{40pt}{
      \xymatrix{
        \mathrm{W}(\hh \stackrel{t}{\to} \gg)
        \ar@{->>}[dd]^{i^*}
        \\
        \\
        \mathrm{W}(\gg)
        &&
        \mathrm{inv}(\gg)
        \ar@{_{(}->}[ll]
        \ar@{_{(}-->}[uull]
      }
      }
      \,.
    \)
     However, a closed invariant polynomial will not
     necessarily lift to a closed one.
    \label{invariant polynomials lift to crossed module}
  \end{proposition}
 \proof
   Recall that 
   $d_t := d_{\mathrm{CE}(\hh \stackrel{t}{\to} \gg)}$ acts on
   $\gg^*$ as
   \(
     d_t|_{\gg^*} = [\cdot,\cdot]_\gg^* + t^*
     \,.
   \)
   By definition \ref{mapping cone of qDGCAs} and definition \ref{Weil algebra definition}
   it follows that 
   $
     d_{\mathrm{W}(\hh \stackrel{t}{\to} \gg)}
   $
   acts on $\gg^*[1]$ as
   \(
     d_{\mathrm{W}(\hh \stackrel{t}{\to} \gg)}|_{\gg^*[1]}
     =
     -\sigma \circ [\cdot,\cdot]_\gg^*   - \sigma \circ t^* 
   \)
   and on $\hh^*[1]$ as
   \(
     d_{\mathrm{W}(\hh \stackrel{t}{\to} \gg)}|_{\hh^*[1]}
     =
     -\sigma \circ \alpha^*
     \,.
   \)
   Then notice that 
   \(
     (\sigma \circ t^*) : \gg^*[1] \to \hh^*[2]
     \,.
   \)
   But this means that $d_{\mathrm{W}(\hh \stackrel{t}{\to} \gg)}$
   differs from $d_{\mathrm{W}(\gg)}$ on $\wedge^\bullet (\gg^*[1])$
   only by elements that are annihilated by vertical $\iota_X$. 
   This proves the claim.
 \endofproof

  It may be easier to appreciate this proof by looking at what it does
  in terms of a chosen basis.

 \subparagraph{Same discussion in terms of a basis.}
  Let $\{t^a\}$ be a basis of $\gg^*$ and $\{b^i\}$ be a basis of $\hh^*[1]$.
  Let $\{C^a{}_{bc}\}$, $\{\alpha^i{}_{aj}\}$, and $\{t^a{}_i\}$,
  respectively,  
  be the components of $[\cdot,\cdot]_\gg$, $\alpha$ and $t$ 
  in that basis.
  Then corresponding to ${\mathrm{CE}(\gg)}$, ${\mathrm{W}(\gg)}$, 
 ${\mathrm{CE}(\hh \stackrel{t}{\to} \gg)}$, and 
 ${\mathrm{W}(\hh \stackrel{t}{\to} \gg)}$,
  respectively, we have the differentials
  \(
    d_{\mathrm{CE}(\gg)}
    : 
    t^a \mapsto -\frac{1}{2}C^a{}_{bc} t^b \wedge t^c, 
  \)
 
  \(
    d_{\mathrm{W}(\gg)}
    : 
    t^a \mapsto -\frac{1}{2}C^a{}_{bc} t^b \wedge t^c + \sigma t^a,
  \)
 
  \(
    d_{\mathrm{CE}(\hh \stackrel{t}{\to} \gg)}
    : 
    t^a \mapsto -\frac{1}{2}C^a{}_{bc} t^b \wedge t^c + t^a{}_i b^i,
     \)
  and
  \(
    d_{\mathrm{W}(\hh \stackrel{t}{\to} \gg)}
    : 
    t^a \mapsto -\frac{1}{2}C^a{}_{bc} t^b \wedge t^c + t^a{}_i b^i + \sigma t^a.
  \)
  Hence we get
  \(
    d_{\mathrm{W}(\gg)}    
    :
    \sigma t^a
    \mapsto
    -\sigma ( -\frac{1}{2}C^a{}_{bc} t^b \wedge t^c )
    =
    C^a{}_{bc} (\sigma t^b) \wedge t^c
  \)
  as well as
  \(
    d_{\mathrm{W}(\hh \stackrel{t}{\to} \gg)}    
    :
    \sigma t^a
    \mapsto
    -\sigma ( -\frac{1}{2}C^a{}_{bc} t^b \wedge t^c + t^a{}_i b^i )
    =
    C^a{}_{bc} (\sigma t^b) \wedge t^c + t^a{}_i \sigma b^i
    \,.
  \)
  Then if
  \( 
    P = P_{a_1 \cdots a_n} (\sigma t^{a_1}) \wedge \cdots \wedge (\sigma t^{a_n})
  \)
  is $d_{\mathrm{W}(\gg)}$-closed, i.e. an invariant polynomial on $\gg$,
  it follows that
  \(
    \label{differential of BF-like invariant polynomial}
    d_{\mathrm{W}(\hh \stackrel{t}{\to} \gg)}    P
    =
    n P_{a_1, a_2, \cdots a_n} (t^{a_1}{}_i \sigma b^i) \wedge (\sigma t^{a_2}) \wedge
     \cdots \wedge (\sigma t^{a_n})
     \,.
  \)  
  (all terms appearing are in the image of the
  shifting isomorphism $\sigma$), hence $P$ is also an invariant polynomial
  on $(\hh \stackrel{t}{\to} \gg)$.
\endofproof

  We will see a physical application of this fact in \ref{characteristic forms}.

 \paragraph{Remark.}
   Notice that the invariant polynomials $P$ lifted from $\gg$ to 
  $(\hh \stackrel{t}{\to} \gg)$ this way are no longer \emph{closed}, in general. 
  This is a new phenomenon we encounter for higher $L_\infty$-algebras.
  While, according to proposition \ref{ordinary Lie cohomology is reproduced},
  for $\gg$ an ordinary Lie algebra all elements in $\mathrm{inv}(\gg)$
  are  closed, this is no longer the case here: the lifted elements 
  $P$ above vanish only after we hit with them with both $d_{\mathrm{W}(\hh \stackrel{t}{\to} \gg)}$
  \emph{and} a vertical $\tau$.

 \begin{proposition}
   \label{transgressive elements from gg to hhtogg}
   Let $P$ be any invariant polynomial on the ordinary Lie algebra
   $\gg$ in transgression with the cocycle $\mu$ on $\gg$.
   Regarded both as elements of 
   $\mathrm{W}(\hh \stackrel{t}{\to}\gg)$ and 
   $\mathrm{CE}(\hh \stackrel{t}{\to}\gg)$ respectively. 
   Notice that $d_{\mathrm{CE}(\hh \stackrel{t}{\to}\gg)}\mu$ 
   in in general non-vanishing but is of course now an 
   exact cocycle on $(\hh \stackrel{t}{\to}\gg)$.

   We have : 
   the $(\hh \stackrel{t}{\to}\gg)$-cocycle 
   $d_{\mathrm{CE}(\hh \stackrel{t}{\to}\gg)}\mu$ 
   transgresses to 
   $d_{\mathrm{inv}(\hh \stackrel{t}{\to}\gg)} P$.
 \end{proposition}
 The situation is illustrated by the diagram in figure
 \ref{figure illustrating inv polynomials on strict Lie 2-algebra}.

 \begin{figure}
   $$
     \xymatrix{
        &&
        0
        \\
        \\
        0
        &&
        p^*
        d_{\mathrm{inv}(\hh \stackrel{t}{\to} \gg)} P
        \ar@/_1.3pc/@{|->}[dd]_{\tau_{\mathrm{W}(\hh \stackrel{t}{\to} \gg)}}
        \ar@{|->}[uu]^{d_{\mathrm{W}(\hh \stackrel{t}{\to} \gg)}}
        &&
        d_{\mathrm{inv}(\hh \stackrel{t}{\to}\gg)} P
        \ar@{|->}[ll]_{p^*_{\mathrm{W}(\hh \stackrel{t}{\to} \gg)}}
        \\
        \\
        d_{\mathrm{CE}(\hh \stackrel{t}{\to} \gg)} \mu
        \ar@{|->}[uu]^{d_{\mathrm{CE}(\hh \stackrel{t}{\to} \gg)}}
        &&
        \tau p^*
        d_{\mathrm{inv}(\hh \stackrel{t}{\to} \gg)} P
        \ar@{|->}[uu]_{d_{\mathrm{W}(\hh \stackrel{t}{\to} \gg)}}
        \ar@{|->}[ll]_{i^*} 
        &&
        P
        \ar@{|->}[uu]_{d_{\mathrm{inv}(\hh \stackrel{t}{\to}\gg)}}
     }
     \,.
   $$
   \caption{
     \label{figure illustrating inv polynomials on strict Lie 2-algebra}
     {\bf Cocycles and invariant polynomials on strict Lie 2-algebras}
     $(\hh \stackrel{t}{\to}\gg)$,   
     induced from cocycles and invariant polynomials on $\gg$.
     An invariant polynomial $P$ on $\gg$ in transgression with 
     a cocycle $\mu$ on $\gg$ lifts to a generally
     non-closed invariant polynomial on $(\hh \stackrel{t}{\to}\gg)$.
     The diagram says that its closure,
     $d_{\mathrm{inv}(\hh \stackrel{t}{\to}\gg)} P$, suspends to the
     $d_{\mathrm{CE}(\hh \stackrel{t}{\to}\gg)}$-closure of the 
     cocycle $\mu$.
     Since this $(\hh \stackrel{t}{\to} \gg)$-cocycle 
     $d_{(\hh \stackrel{t}{\to} \gg)} \mu$ is hence
     a coboundary, it follows from proposition   
     \ref{invariant polynomials and trivial cocycles}
     that $d_{\mathrm{inv}(\hh \stackrel{t}{\to}\gg)} P$ 
     suspends also to 0. Nevertheless the situation is of interest,
     in that it governs the topological field theory known as
     BF theory. 
     This is discussed in section \ref{examples for characteristic forms}.
   }
  \end{figure}

\paragraph{Concrete Example: $\mathfrak{su}(5) \to \mathfrak{sp}(5)$.}

It is known that the cohomology of the Chevalley-Eilenberg algebras for
$\mathfrak{su}(5)$ and $\mathfrak{sp}(5)$ are generated, respectively, by four
and five indecomposable cocycles,
\(
  H^\bullet(\mathrm{CE}(\mathfrak{su}(5)))
  =
  \wedge^\bullet[a ,b, c,d]
\)
and
\(
  H^\bullet(\mathrm{CE}(\mathfrak{sp}(5)))
  =
  \wedge^\bullet[v,w,x,y,z]
  \,,
\)
which have degree as indicated in the following table:
\begin{tabular}{c|cc}
    &$\mathbf{generator}$ & $\mathbf{degree}$
     \\
      &$a$ & 3
      \\
      &$b$ & 5
      \\
      $H^\bullet{\mathrm{CE}(\mathfrak{su}(5))}$
      &$c$ & 7
      \\
      &$d$ & 9
      \\
      \hline
      \\      
      &$v$	& 3
      \\
      &$w$ & 7
      \\
      $H^\bullet(\mathrm{CE}(\mathfrak{sp}(5)))$
      &$x$ & 11
      \\
      &$y$ & 15
      \\
      &$z$ & 19
\end{tabular}.

\noindent As discussed for instance in \cite{GHV}, the inclusion of groups
\(
  \mathrm{SU}(5) \hookrightarrow \mathrm{Sp}(5)
\)
is reflected in the morphism of DGCAs
\(
  \xymatrix{ 
    \mathrm{CE}(\mathfrak{su}(5))
    &&
    \mathrm{CE}(\mathfrak{sp}(5))
    \ar@{->>}[ll]_{t^*}
  }
\)
which acts, in cohomology, on $v$ and $w$ as 
\(
  \xymatrix@R=3pt{
    a & v \ar@{|->}[l] \\
    c & w \ar@{|->}[l] \\
  }
\)
and which sends $x$, $y$ and $z$ to wedge products of generators.

We would like to apply the above reasoning to this situation. Now, 
$\mathfrak{su}(5)$ is not normal in $\mathfrak{sp}(5)$ hence
$(\mathfrak{su}(5) \hookrightarrow \mathfrak{sp}(5))$ does not give
a Lie 2-algebra. But we can regard the cohomology complexes
$H^\bullet(\mathrm{CE}(\mathfrak{su}(5)))$ and
$H^\bullet(\mathrm{CE}(\mathfrak{sp}(5)))$ as Chevalley-Eilenberg
algebras of abelian $L_\infty$-algebras in their own right. Their inclusion
is normal, in the sense to be made precise below in definition 
\ref{normal subalgebras}. By useful abuse of notation, we write 
now $\mathrm{CE}(\mathfrak{su}(5) \hookrightarrow \mathfrak{sp}(5))$
for this inclusion at the level of cohomology.

Recalling from \ref{differential in CE of strict Lie 2-algebra}
that this means that in $\mathrm{CE}(\mathfrak{su}(5) \hookrightarrow \mathfrak{sp}(5))$
we have
\(
  \label{d v for su5 to sp5}
  d_{\mathrm{CE}(\mathfrak{su}(5) \hookrightarrow \mathfrak{sp}(5))} v := \sigma a
\)
and
\(
  \label{d w for su5 to sp5}
  d_{\mathrm{CE}(\mathfrak{su}(5) \hookrightarrow \mathfrak{sp}(5))}) w := \sigma c
\)
we see that the generators $\sigma a$ and $\sigma b$ drop out of the cohomology of
the Chevalley-Eilenberg algebra
\(
  \label{CE-algebra for su5 to sp5}
  \mathrm{CE}(\mathfrak{su}(5) \hookrightarrow \mathfrak{sp}(5))
  =
  (\wedgebullet (\mathfrak{sp}(5)^* \oplus \mathfrak{su}(5)^*[1]) , d_t)
\)
of the strict Lie 2-algebra coming from the infinitesimal crossed module
$(t : \mathfrak{su}(5) \hookrightarrow \mathfrak{sp}(5))$.

A simple spectral sequence argument shows that products 
are not killed in $H^\bullet(\mathrm{CE}(\mathfrak{su}(5) \hookrightarrow \mathfrak{sp}(5)))$,
but they may no longer be decomposable.
Hence 
\(
H^\bullet(\mathrm{CE}(\mathfrak{su}(5) \hookrightarrow \mathfrak{sp}(5)))
\)
is generated by classes in degrees 6 and 10 by $\sigma b$ and $\sigma d$, and in degrees
21 and 25, which are represented by products in 
\ref{CE-algebra for su5 to sp5} involving $\sigma a$ and $\sigma c$, with the only non zero product being 
\(
  6 \wedge 25 = 10 \wedge 21
  \,,
\)
where 31 is the dimension of the manifold $\mathrm{Sp}(5)/\mathrm{SU}(5)$.
Thus the strict Lie 2-algebra $(t : \mathfrak{su}(5) \hookrightarrow \mathfrak{sp}(5))$ 
plays the role of the quotient Lie 1-algebra $\mathfrak{sp}(5)/\mathfrak{su}(5)$.
We will discuss the general mechanism behind this phenomenon in \ref{weak cokernels of Lie infty-algebras}:
the Lie 2-algebra $(\mathfrak{su}(5)\hookrightarrow \mathfrak{sp}(5))$ is the \emph{weak cokernel},
i.e. the \emph{homotopy cokernel} of the inclusion $\mathfrak{su}(5)\hookrightarrow \mathfrak{sp}(5)$.

The Weil algebra of $(\mathfrak{su}(5) \hookrightarrow \mathfrak{sp}(5))$ is
\(
  \mathrm{W}(\mathfrak{su}(5) \hookrightarrow \mathfrak{sp}(5))
  =
  (\wedge^\bullet(
    \mathfrak{sp}(5)^* \oplus \mathfrak{su}(5)^*[1]
    \oplus
    \mathfrak{sp}(5)^*[1] \oplus \mathfrak{su}(5)^*[2]    
  ),
  d_{\mathrm{W}(\mathfrak{su}(5) \hookrightarrow \mathfrak{sp}(5))})
  \,.
\)
Recall the formula \ref{differential of Weil algebra} for 
the action of $d_{\mathrm{W}(\mathfrak{su}(5) \hookrightarrow \mathfrak{sp}(5))}$
on generators in $\mathfrak{sp}(5)^*[1] \oplus \mathfrak{su}(5)^*[2]$. 
By that formula, $\sigma v$ and $\sigma w$ are invariant polynomials 
on $\mathfrak{sp}(5)$ which lift to non-closed
invariant polynomials on $\mathfrak{su}(5) \hookrightarrow \mathfrak{sp}(5))$:
\(
  d_{\mathrm{W}(\mathfrak{su}(5) \hookrightarrow \mathfrak{sp}(5))})
  :
  \sigma v \mapsto 
  -\sigma ( d_{\mathrm{CE}(\mathfrak{su}(5) \hookrightarrow \mathfrak{sp}(5))} v )
  =
  - \sigma \sigma a
\)
by equation \ref{d v for su5 to sp5}; and
\(
  d_{\mathrm{W}(\mathfrak{su}(5) \hookrightarrow \mathfrak{sp}(5))})
  :
  \sigma w \mapsto 
  -\sigma ( d_{\mathrm{CE}(\mathfrak{su}(5) \hookrightarrow \mathfrak{sp}(5))} w )
  =
  - \sigma \sigma c
\)
by equation \ref{d w for su5 to sp5}.
Hence $\sigma v$ and $\sigma w$ are not closed in 
$\mathrm{CE}(\mathfrak{su}(5) \hookrightarrow \mathfrak{sp}(5)) $, but they are still
invariant polynomials according to definition \ref{invariant polynomials}, since their
differential sits entirely in the shifted copy $(\mathfrak{sp}(5)^* \oplus \mathfrak{su}(5)^*[1])[1]$.

On the other hand, notice that we do also have closed invariant polynomials on 
$(\mathfrak{su}(5) \hookrightarrow \mathfrak{sp}(5))$, for instance $\sigma \sigma b$
and $\sigma \sigma d$.

\subsection{$L_\infty$-algebras from cocycles: String-like extensions}

\label{String-like extensions}

We now consider the main object of interest here: families of 
$L_\infty$-algebras that are induced from $L_\infty$-cocycles
and invariant polynomials. First we need the following

\begin{definition}[String-like extensions of $L_\infty$-algebras]
  Let $\gg$ be an $L_\infty$-algebra. 
  \begin{itemize}
    \item
      For each degree $(n+1)$-cocycle $\mu$ on $\gg,$ let
      $\gg_\mu$
      be the $L_\infty$-algebra defined by
      \(
         \mathrm{CE}(\gg_\mu) = (\wedge^\bullet(\gg^* \oplus \mathbb{R}[n]), d_{\mathrm{CE}(\gg_\mu)})
      \)
      with differential given by
      \(
        d_{\mathrm{CE}(\gg_\mu)}|_{\gg^*} := d_{\mathrm{CE}(\gg)},
      \)
      and
      \(
        d_{\mathrm{CE}(\gg_\mu)})|_{\mathbb{R}[n]} : b \mapsto -\mu
        \,,
      \)
      where $\{b\}$ denotes the canonical basis of $\mathbb{R}[n]$.
      \emph{
        This we call the {\bf String-like extension} of $\gg$ with respecto to $\mu$,
        because, as described below in \ref{examples for String-like extensions}, 
        it generalizes the construction of the String Lie 2-algebra.
      }
    \item
      For each degree $n$ invariant polynomial $P$ on $\gg,$ let
       $\mathrm{ch}_P(\gg)$
      be the $L_\infty$-algebra defined by
      \(
         \mathrm{CE}(\mathrm{ch}_P(\gg)) 
          = (\wedge^\bullet(\gg^* \oplus \gg^*[1] \oplus \mathbb{R}[2n-1]), d_{\mathrm{CE}(\mathrm{ch}_P(\gg))})
      \)
      with the differential given by
      \(
        d_{\mathrm{CE}(\mathrm{ch}_P(\gg))}|_{\gg^*\oplus \gg^*[1]} := d_{\mathrm{W}(\gg)}
      \)
      and
      \(
        d_{\mathrm{CE}(\mathrm{ch}_P(\gg))})|_{\mathbb{R}[2n-1]} : c \mapsto P
        \,,
      \)
      where $\{c\}$ denotes the canonical basis of $\mathbb{R}[2n-1]$.
      \emph{
        This we call the {\bf Chern $L_\infty$-algebra} corresponding to the
        invariant polynomial $P$, because, as described below in 
        \ref{examples for Loo-valued forms}, 
         connections with values in it pick out the Chern-form
        corresponding to $P$.
      }

    \item
      For each degree $2n-1$ transgression element $\mathrm{cs},$ let
      $\mathrm{cs}_P(\gg)$
      be the $L_\infty$-algebra defined by
      \(
         \mathrm{CE}(\mathrm{cs}_P(\gg)) 
          = (\wedge^\bullet(\gg^* \oplus \gg^*[1] \oplus \mathbb{R}[2n-2]\oplus \mathbb{R}[2n-1]), d_{\mathrm{CE}(\mathrm{ch}_P(\gg))})
      \)
      with
      \(
        d_{\mathrm{CE}(\mathrm{cs}_P(\gg))}|_{\wedge^\bullet (\gg^* \oplus \gg^*[1])}
          = d_{\mathrm{W}(\gg)}
      \)
      \(
        d_{\mathrm{CE}(\mathrm{cs}_P(\gg))}|_{\mathbb{R}[2n-2]}
        :
        b \mapsto - \mathrm{cs} + c
      \)
      \(
        d_{\mathrm{CE}(\mathrm{ch}_p(\gg))}|_{\mathbb{R}[2n-1]}
        :
        c \mapsto P
        \,,
      \)
      where $\{b\}$ and $\{c\}$ denote the canonical bases of $\mathbb{R}[2n-2]$
      and $\mathbb{R}[2n-1]$, respectively.
      \emph{      
      This we call the {\bf Chern-Simons $L_\infty$-algebra} with respect to 
      the transgression element $\mathrm{cs}$, because, as described below in
      \ref{examples for Loo-valued forms}, 
      connections with values in these come from (generalized)
      Chern-Simons forms.}
      \end{itemize}

  \label{definition of string-like extension}
\end{definition}

The nilpotency of these differentials follows directly from the very
definition of $L_\infty$-algebra cocoycles and invariant polynomials.

\begin{proposition}[the string-like extensions]
  \label{the string-like extension sequence}
  For each $L_\infty$-cocycle $\mu \in \wedge^n (\gg^*)$ of degree $n$,
  the corresponding String-like extension
  sits in an exact sequence
  \(
    \xymatrix{
      0 
      &
      \mathrm{CE}(b^{n-1}\uu(1))
      \ar[l]
      &&
      \mathrm{CE}(\gg_\mu)
      \ar@{->>}[ll]
      &&
      \mathrm{CE}(\gg)
      \ar@{_{(}->}[ll]
      &
      0
      \ar[l]
    }
\)
\end{proposition}
\proof
  The morphisms are
  the canonical inclusion and
  projection.
\endofproof

\begin{proposition}
  For $\mathrm{cs} \in \mathrm{W}(\gg)$ any transgression element interpolating
  between the cocycle $\mu \in \mathrm{CE}(\gg)$ and the invariant
  polynomial $P \in \wedge^\bullet (\gg[1]) \subset \mathrm{W}(\gg)$, we obtain 
  a homotopy-exact sequence
  \(
    \xymatrix{
      \mathrm{CE}(\gg_\mu)
      &&
      \mathrm{CE}(\mathrm{cs}_P(\gg))
      \ar@{->>}[ll]
      \ar@{-}[d]^\simeq
      &&
      \mathrm{CE}(\mathrm{ch}_P(\gg))
      \ar@{_{(}->}[ll]
      \\
      &&
      \mathrm{W}(\gg_\mu)
    }
    \,.
  \)
  \label{Chern and Chern-Simons}
\end{proposition}

Here the isomorphism
\(
  f:
  \xymatrix{
    \mathrm{W}(\gg_\mu)
    \ar[r]^{\simeq}
    &
    \mathrm{CE}(\mathrm{cs}_P(\gg))
  }
\)
is the identity on $\gg^* \oplus \gg^*[1] \oplus \mathbb{R}[n]$
\(
  f|_{\gg^* \oplus \gg^*[1] \oplus \mathbb{R}[n]} = \mathrm{Id}
\)
and acts as
\(
  f|_{\mathbb{R}[n+1]} : b \mapsto c + \mu - \mathrm{cs}
\)
for $b$ the canonical basis of $\mathbb{R}[n]$ and $c$ that of $\mathbb{R}[n+1].$
We check that this does respect the differentials
\(
  \raisebox{40pt}{
  \xymatrix{
    b 
    \ar@{|->}[rr]^{d_{\mathrm{W}(\gg_\mu)}}
    \ar@{|->}[dd]^f
    &&
    -\mu + c
    \ar@{|->}[dd]^f
    \\
    \\
    b
    \ar@{|->}[rr]^{d_{\mathrm{CE}(\mathrm{cs}_P(\gg))}}
    &&
    -\mathrm{cs} + c
  }
  }
  \hspace{27pt}
  \raisebox{40pt}
  {
  \xymatrix{
    c 
    \ar@{|->}[rr]^{d_{\mathrm{W}(\gg_\mu)}}
    \ar@{|->}[dd]^f
    &&
    \sigma \mu
    \ar@{|->}[dd]^f
    \\
    \\
    c + \mu - \mathrm{cs}
    \ar@{|->}[rr]^{d_{\mathrm{CE}(\mathrm{cs}_P(\gg))}}
    &&
    \sigma \mu
  }
  }
  \,.
\)
Recall from definition \ref{mapping cone, detailed def}
that $\sigma$ is the canonical isomorphism 
$\sigma : \gg^* \to \gg^*[1]$
extended by 0 to $\gg^*[1]$ and then  
as a derivation to all of $\wedge^\bullet (\gg^*\oplus \gg^*[1])$.
In the above the morphism between the Weil algebra of $\gg_\mu$ 
and the Chevalley-Eilenberg algebra of $\mathrm{cs}_P(\gg)$
is indeed an \emph{iso}morphism (not just an equivalence).
This isomorphism exhibits one of the main points to be made here:
it makes manifest that the invariant polynomial $P$ that 
is related by transgression to the cocycle $\mu$ which 
induces $\gg_\mu$ becomes exact with respect to $\gg_\mu$.
This is the statement of proposition \ref{exactness of P in inv(g_mu)} below.

\paragraph{$L_\infty$-algebra cohomology and 
invariant polynomials of String-like extensions.}

The $L_\infty$-algebra $\gg_\mu$ obtained from an $L_\infty$-algebra
$\gg$ with an $L_\infty$-algebra cocycle $\mu \in H^\bullet(\mathrm{CE}(\gg))$
can be thought of as being obtained from $\gg$ by ``killing'' 
a cocycle $\mu$. 
This is familiar from Sullivan models in rational homotopy theory.

\begin{proposition}
  Let $\gg$ be an ordinary semisimple Lie algebra and $\mu$ a cocycle 
  on it. Then
  \(
    H^\bullet(\mathrm{CE}(\gg_\mu))
    =
    H^\bullet(\mathrm{CE}(\gg))/{\langle\mu\rangle}
    \,.
  \)
  \label{cohomology of g-mu}
\end{proposition}

Accordingly, one finds that in cohomology the invariant polynomials 
on $\gg_\mu$ are those of $\gg$ except that the polynomial
in transgression with $\mu$ now suspends to 0.

\begin{proposition}
  \label{exactness of P in inv(g_mu)}
  Let $\gg$ be an $L_\infty$-algebra and $\mu \in \mathrm{CE}(\gg)$
  in transgression with the invariant polynomial 
  $P \in \mathrm{inv}(\gg)$. Then with respect to 
  the String-like extension $\gg_\mu$ the polynomial $P$ suspends to 0.
\end{proposition}
\proof
  Since $\mu$ is a coboundary in $\mathrm{CE}(\gg_\mu)$,
  this is a corollary of 
  proposition \ref{invariant polynomials and trivial cocycles}.
\endofproof

\paragraph{Remark.} We will see in \ref{Loo Cartan-Ehresmann connections}
that those invariant polynomials which suspend to 0 do actually not 
contribute to the characteristic classes. As we will also see there,
this can be understood in terms of the invariant polynomials not
with respect to the projection 
$\xymatrix{
  \mathrm{CE}(\gg)
  &
  \mathrm{W}(\gg)
  \ar@{->>}[l]
}$
but with respect to the projection
$\xymatrix{
  \mathrm{CE}(\gg)
  &
  \mathrm{CE}(\mathrm{cs}_P(\gg_\mu))
  \ar@{->>}[l]
  \,,
}$
recalling from \ref{Chern and Chern-Simons} that $\mathrm{W}(\gg)$ is isomorphic to $\mathrm{cs}_P(\gg)$.

\begin{proposition}
  For $\gg$ any $L_\infty$-algebra with cocycle $\mu$ of degree 
  $2n+1$ in transgression with the 
  invariant polynomial $P$, denote by $\mathrm{cs}_P(\gg)_{\mathrm{basic}}$
  the DGCA of basic forms with respect to the canonical projection
  \(
  \xymatrix{
   \mathrm{CE}(\gg)
   &
   \mathrm{CE}(\mathrm{cs}_P(\gg_\mu))
   \ar@{->>}[l]
   \,,
  }
  \)
  according to the general definition \ref{basic elements}.
  
  Then the cohomology of $\mathrm{cs}_P(\gg)_{\mathrm{basic}}$ is that of $\mathrm{inv}(\gg)$
  modulo $P$:
  \(
    H^\bullet(\mathrm{cs}_P(\gg)_{\mathrm{basic}}) \simeq H^\bullet(\mathrm{inv}(\gg))/\langle P \rangle
    \,.
\)
  \label{basic forms on Chern-Simons algebra}
\end{proposition}
\proof
  One finds that the vertical derivations on $\mathrm{CE}(\mathrm{cs}_P(\gg))
  = \wedge^\bullet(\gg^* \oplus \gg^*[1] \oplus \mathbb{R}[n] \oplus \mathbb{R}[n+1])$ 
  are those that vanish on everything except the unshifted copy of $\gg^*$. 
  Therefore the basic forms are those in $\wedge^\bullet(\gg^*[1] \oplus \mathbb{R}[n] \oplus \mathbb{R}[n+1])$
  such that also their $d_{\mathrm{cs}_P(\gg)}$-differential is in that space.
  Hence all invariant $\gg$-polynomials are among them. But one of them now becomes exact, namely $P$. 
\endofproof

\begin{figure}[h]
$$
  \xymatrix@C=4pt{
    \mathrm{CE}(\gg_\mu)
    \ar@{-}[r]^{=}
    &
    \mathrm{CE}(\gg_\mu)
    \\
    \mathrm{W}(\gg_\mu) \ar@{-}[r]^\simeq 
    \ar@{->>}[u]
    & \mathrm{CE}(\mathrm{cs}_P(\gg))
    \ar@{->>}[u]
    \\
    **[l]\mathrm{inv}(\gg) = \mathrm{inv}(\gg_\mu) 
    \ar@{^{(}->}[r]
    \ar@{_{(}->}[u]
    & 
    \mathrm{cs}_P(\gg)_{\mathrm{basic}}
    \ar@{_{(}->}[u]
    \\
    & **[r]H^\bullet(\mathrm{cs}_P(\gg)_{\mathrm{basic}}) \simeq H^\bullet(\mathrm{inv}(\gg))/\langle P \rangle
  }
$$
\caption{
  The DGCA sequence playing the role of differential forms on 
  the {\bf universal (higher) String $n$-bundle}
  for a String-like extension $\gg_\mu$, definition \ref{definition of string-like extension},
  of an $L_\infty$-algebra $\gg$ by a cocycle $\mu$ of odd degree in transgression with an 
  invariant polynomial $P$. Compare with figure \ref{universal G-bundle}.
  In $H^\bullet(\mathrm{inv}(\gg_\mu)) = H^\bullet(\mathrm{W}(\gg_\mu)_{\mathrm{basic}})$ 
  the class of $P$ is still 
  contained, but suspends to 0,
  according to proposition \ref{exactness of P in inv(g_mu)}. In 
  $H^\bullet(\mathrm{cs}_P(\gg)_{\mathrm{basic}})$ the class of $P$ vanishes, according to
  proposition \ref{basic forms on Chern-Simons algebra}. 
  The isomorphism $\mathrm{W}(\gg_\mu) \simeq \mathrm{CE}(\mathrm{cs}_P(\gg))$ is
  from proposition \ref{Chern and Chern-Simons}.
  For $\gg$ an ordinary semisimple Lie algebra and $\gg_\mu$ the ordinary String extension 
  coming from the canonical 3-cocycle, this corresponds to the fact that the classifying 
  space of the String 2-group \cite{BCSS,Henriques} has the cohomology of the classifying
  space of the underlying group, modulo the first Pontrajagin class \cite{BaezStevenson}.
}
\end{figure}

\paragraph{Remark.} 
The first example below, definition \ref{string Lie 2-algebra}, introduces the 
String Lie 2-algebra of an ordinary semisimple Lie algebra $\gg$,
which gave all our String-like extensions its name. It is
known, corollary 2 in \cite{BaezStevenson} 
that the real cohomology of the classifying space of the 2-group integrating it
is that of $G = \exp(\gg)$, modulo the ideal generated by the class corresponding to
$P$. Hence $\mathrm{CE}(\mathrm{cs}_P(\gg))$ is an algebraic model for this space.

\subsubsection{Examples}

\label{examples for String-like extensions}

 \paragraph{Ordinary central extensions.}
 
  Ordinary central extensions coming from a 2-cocycle $\mu \in H^2(\mathrm{CE}(\gg))$
  of an ordinary Lie algebra $\gg$ are a special case of the ``string-like'' extensions
  we are considering:
  
  By definition \ref{definition of string-like extension} the $L_\infty$-algebra
  $\gg_\mu$ is the Lie 1-algebra whose Chevalley-Eilenberg algebra is
  \(
    \mathrm{CE}(\gg_\mu) = (\wedge^\bullet( \gg^* \oplus \mathbb{R}[1] ), d_{\mathrm{CE}(\gg_\mu)})
  \)
  where
  \(
    d_{\mathrm{CE}(\gg_\mu)}|_{\gg^*} = d_{\mathrm{CE}(\gg)}
  \)
  and
  \(
    d_{\mathrm{CE}(\gg_\mu)}|_{\mathbb{R}[1]} : b \mapsto \mu
  \)
  for $b$ the canonical basis of $\mathbb{R}[1]$. (Recall that in our conventions $\gg$ is
  in degree 1).  
 
  This is indeed the Chevalley-Eilenberg algebra corresponding to the Lie bracket
  \(
    [(x,c), (x',c')] = ([x,x'], \mu(x,x'))
  \)
  (for all $x,x' \in \gg$, $c,c' \in \mathbb{R}$)
  on the centrally extended Lie algebra.

 \paragraph{The String Lie 2-algebra.}

  \begin{definition}
    Let $\gg$ be a semisiple Lie algebra and $\mu = \langle \cdot , [\cdot, \cdot] \rangle$
    the canonical 3-cocycle on it. Then
    \(
      \mathrm{string}(\gg)
    \)
    is defined to be the strict Lie 2-algebra coming from the crossed module
    \(
      (\hat \Omega \gg \to P \gg)
      \,,
    \)
    where $P\gg$ is the Lie algebra of based paths in $\gg$ and $\hat \Omega \gg$
    the Lie algebra of based loops in $\gg$, with central extension induced by
    $\mu$. Details are in \cite{BCSS}. 
    \label{string Lie 2-algebra}
  \end{definition}

  \begin{proposition}[\cite{BCSS}]
     \label{string Lie 2-algebra as strict Lie 2-algebra}
    The Lie 2-algebra $\gg_\mu$ obtained from $\gg$ and
    $\mu$ as in definition \ref{definition of string-like extension} is equivalent
    to the strict string Lie 2-algebra
    \(
      \gg_\mu \simeq \mathrm{string}(\gg)
      \,.
    \)
    \label{equivalence of strict and skeletal string lie 2-algebra}
  \end{proposition}
  This means there are morphisms $\gg_\mu \to \mathrm{string}(\gg)$ and
  $\mathrm{string}(\gg) \to \gg_\mu$ whose composite is the identity
  only up to homotopy
  \(
    \xymatrix{
       \gg_\mu \ar[r] 
       \ar@/_2pc/[rr]_{\mathrm{Id}}^{\ }="t"
       & \mathrm{string}(\gg) 
       \ar[r]
       &
       \gg_\mu
       \ar@{=>}^= "t"+(0,+4); "t"
    }
    \hspace{17pt}
    \xymatrix{
       \mathrm{string}(\gg) \ar[r] 
       \ar@/_2pc/[rr]_{\mathrm{Id}}^{\ }="t"
       & 
       \gg_\mu 
       \ar[r]
       &
       \mathrm{string}(\gg)
       \ar@{=>}^\eta "t"+(0,+4); "t"
    }
  \)

  We call $\gg_\mu$ the \emph{skeletal} and $\mathrm{string}(\gg)$ the \emph{strict}
  version of the String Lie 2-algebra.

 \paragraph{The Fivebrane Lie 6-algebra}

  \begin{definition}
    Let $\gg = \mathrm{so}(n)$ and $\mu$
    the canonical 7-cocycle on it. Then
    \(
      \mathrm{fivebrane}(\gg)
    \)
    is defined to be the strict Lie 7-algebra which is equivalent to $\gg_\mu$
    \(
      \gg_\mu \simeq \mathrm{fivebrane}(\gg)
      \,.
    \)
  \end{definition}
      
  A Lie $n$-algebra is \emph{strict} if it corresponds to a 
  differential graded Lie algebra on a vector space in degree 1 to $n$.
  (Recall our grading conventions from \ref{Lie infty algebras}.)

  \paragraph{Remark.} It is a major open problem to identify 
     the strict $\mathrm{fivebrane}(\gg)$. 
     Proposition \ref{equivalence of strict and skeletal string lie 2-algebra}
     suggests that it might involve hyperbolic Kac-Moody algebras and/or
     the torus algebra of $\gg$, since these would seem to be
     what comes beyond the affine Kac-Moody algebras 
     relevant for $\mathrm{string}(n)$.

  \paragraph{The BF-theory Lie 3-algebra.}

   \begin{definition}
     For $\gg$ any ordinary Lie algebra with 
     bilinear invariant symmetric form 
     $\langle \cdot, \cdot\rangle \in \mathrm{inv}(\gg)$
     in transgression with the 3-cocycle $\mu$, and for
     $\hh \stackrel{t}{\to} \gg$ a strict Lie 2-algebra
     based on $\gg$, denote by 
     \(
       \hat \mu := d_{\mathrm{CE}(\hh \stackrel{t}{\to} \gg)} \mu
     \) 
     the corresponding exact 4-cocycle on $(\hh \stackrel{t}{\to} \gg)$
     discussed in \ref{examples for Lie 00-algebra cohomology}.
     Then we call the string-like extended Lie 3-algebra
     \(
       \mathfrak{bf}(\hh \stackrel{t}{\to} \gg)
       :=
       (\hh \stackrel{t}{\to} \gg)_{\hat \mu}
     \)
     the corresponding BF-theory Lie 3-algebra.
   \end{definition}

   The terminology here will become clear 
   once we describe in \ref{obstruction examples} 
   and \ref{examples for parallel transport} how the
   BF-theory action functional discussed in 
   \ref{examples for characteristic forms} arises as the
   parallel 4-transport given by the $b^3\uu(1)$-4-bundle
   which arises as the obstruction to lifting 
   $(\hh \stackrel{t}{\to}\gg)$-2-descent objects to 
   $\mathfrak{bf}(\hh \stackrel{t}{\to}\gg)$-3-descent objects.

\subsection{$L_\infty$-algebra valued forms}

  \label{Lie infty-algebra valued forms}

Consider an ordinary Lie algebra $\gg$ valued connection form $A$ regarded as a linear map
$\gg^*\to  \Omega^1(Y)$. Since $CE(\gg)$ is free as a graded commutative algebra,
this linear map extends uniquely to a morphism of graded commutative algebras, though 
not in general of differential graded commutative algebra.  In fact, the deviation is measured by
the \emph{curvature} $F_A$  of the connection. However, the differential in $\mathrm{W}(\gg)$ is
precisely such that the connection does extend to a morphism of 
differential graded-commutative algebras 
\(
  \mathrm{W}(\gg)\to  \Omega^\bullet(Y)\,.
\)
This implies that a good notion of a $\gg$-valued differential form on a smooth space $Y$, for 
$\gg$ any $L_\infty$-algebra,
  is a
  morphism of differential graded-commutative algebras from the
  Weil algebra of $\gg$ to the algebra of differential forms on $Y$.

 \begin{definition}[$\gg$-valued forms]
   For $Y$ a smooth space and $\gg$ an $L_\infty$-algebra,
   we call 
   \(
     \Omega^\bullet(Y,\gg)
     :=
     \mathrm{Hom}_{\mathrm{dgc-Alg}}(\mathrm{W}(\gg), \Omega^\bullet(Y))
   \)
   the space of {\bf $\gg$-valued differential forms on $X$}.
 \end{definition}

   \begin{definition}[curvature]
   We write $\gg$-valued differential forms as
   \(
     (\xymatrix{
       \Omega^\bullet(Y)
       &&
       \mathrm{W}(\gg)
       \ar[ll]_{(A,F_A)}
     }
     )
     \hspace{4pt}
       \in
     \hspace{4pt}
     \Omega^\bullet(Y,\gg)
     \,,
   \)
   where $F_A$ denotes the restriction to the 
   shifted copy $\gg^*[1]$ given by
   \(
     \mathrm{curv} : 
      (\xymatrix{
       \Omega^\bullet(Y)
       &
       \mathrm{W}(\gg)
       \ar[l]_{(A,F_A)}
     })
     \hspace{4pt}
     \mapsto
     \hspace{4pt}
      (\xymatrix{
       \Omega^\bullet(Y)
       &
       \mathrm{W}(\gg)
       \ar[l]_{(A,F_A)}
       &
       \gg^*[1]
       \ar@{_{(}->}[l]
       \ar@/_2pc/[ll]_{F_A}
     })
     \,.
   \)
   $F_A$ we call the curvature of $A$.
  \end{definition}

  \begin{proposition}
    The $\gg$-valued differential form
    $\xymatrix{
     \Omega^\bullet(Y) && \mathrm{W}(\gg)\ar[ll]_{(A,F_A)}
    }$
    factors through $\mathrm{CE}(\gg)$ precisely when its curvature 
    $F_A$ vanishes.
 \end{proposition}
 \(
   \xymatrix{
      \mathrm{CE}(\gg)
      \ar@{..>}[dd]^{(A,F_{A} = 0)}
       &&
      \mathrm{W}(\gg)
      \ar@{->>}[ll]
      \ar[dd]^{(A,F_A)}
      \\
      \\
      \Omega^\bullet(Y)
      \ar@{-}[rr]^=
       &&
      \Omega^\bullet(Y)
   }
   \,.
 \)
  In this case we say that $A$ is {\bf flat}.
  Hence the space of flat $\gg$-valued forms is 
   \(
     \Omega^\bullet_{\mathrm{flat}}(Y,\gg)
     \simeq
     \mathrm{Hom}_{\mathrm{dgc-Alg}}(\mathrm{CE}(\gg), \Omega^\bullet(Y))
     \,.
   \)

  \paragraph{Bianchi identity.}
    Recall from \ref{Lie infty algebras} that the Weil algebra $\mathrm{W}(\gg)$
    of an $L_\infty$-algebra $\gg$ is the same as the
    Chevalley-Eilenberg algebra $\mathrm{CE}(\mathrm{inn}(\gg))$ of the 
    $L_\infty$-algebra of inner derivation of $\gg$.
    It follows that $\gg$-valued differential forms on $Y$ are the same
    as \emph{flat} $\mathrm{inn}(\gg)$-valued differential forms on $Y$:
    \(
      \Omega^\bullet(Y,\gg) = \Omega^\bullet_{\mathrm{flat}}(\mathrm{inn}(\gg))
      \,.
    \)
    By the above definition of curvature, this says that the curvature $F_A$
    of a $\gg$-valued connection $(A,F_A)$ is itself a flat $\mathrm{inn}(\gg)$-valued
    connection. This is the generalization of the ordinary \emph{Bianchi identity}
    to $L_\infty$-algebra valued forms.

\begin{definition}
  \label{gauge transformation of g-valued forms}
  Two $\gg$-valued forms $A, A' \in \Omega^\bullet(Y,\gg)$ 
  are called {\bf (gauge) equivalent}
  precisely if they are related by a \emph{vertical}  
  concordance, i.e. by a concordance, such that the
  corresponding derivation $\rho$ from proposition 
  \ref{concordance and chain homotopy} is vertical, in the sense
  of definition \ref{vertical derivations}.  
\end{definition}

\subsubsection{Examples}

\label{examples for Loo-valued forms}

\medskip
\noindent  {\bf 1. Ordinary Lie-algebra valued 1-forms.} 
We have already mentioned ordinary Lie algebra valued 1-forms in this
general context in \ref{Lie oo-algebra valued forms in plan}.

\medskip 
\noindent {\bf 2. Forms with values in shifted $b^{n-1}\uu(1)$} 

A $b^{n-1}\uu(1)$-valued form is nothing but an ordinary $n$-form $A \in \Omega^n(Y)$:
\(
  \Omega^\bullet(b^{n-1}\uu(1),Y) \simeq \Omega^n(Y)
  \,.
\)
A flat $b^{n-1}\uu(1)$-valued form is precisely a closed $n$-form.

$$
  \xymatrix{
    \mathrm{CE}(b^{n-1}\uu(1))
    \ar@{..>}[dd]_<<<<<<<<<{(A)}|>>>>>>>{dA = 0}
    \ar@{<<-}[rr]
    &&
    \mathrm{W}(b^{n-1}\uu(1))
    \ar[dd]_<<<<<<<<{(A,F_A)}|>>>>>>>>>{A = dA}
    \\
    \\
    \Omega^\bullet(Y)
    \ar[rr]^{=}
    &&
    \Omega^\bullet(Y)
  }
$$

\medskip
\noindent {\bf 3. Crossed module valued forms.}
Let $\gg_{(2)} = (\hh \stackrel{t}{\to} \gg)$ be a strict Lie 
2-algebra coming from a crossed module. Then a 
$\gg_{(2)}$-valued form is an ordinary $\gg$-valued 1-form
$A$ and an ordinary $\hh$-valued 2-form $B$. The corresponding
curvature is an ordinary $\gg$-valued 2-form 
$\beta = F_A + t(B)$ and an ordinary $\hh$-valued 
3-form $H = d_A B$. This is denoted by the right 
vertical arrow in the following diagram.
\(
  \raisebox{30pt}{
  \xymatrix{
    \mathrm{CE}(\hh \stackrel{t}{\to} \gg)
    \ar@{..>}[dd]_<<<<<<<<{(A,B)}|>>>>>>>{F_A + t(B) = 0}
    &&
    \mathrm{W}(\hh \stackrel{t}{\to} \gg)
    \ar@{->>}[ll]
    \ar[dd]_<<<<<<<<{(A,B,\beta,H)}|>>>>>>>{\beta = F_A + t(B) \atop H = d_A B}
    \\
    \\
    \Omega^\bullet(Y)    
    \ar@{-}[rr]^=
    &&
    \Omega^\bullet(Y)    
  } 
  }
  \,.
\)
Precisely if the curvature components $\beta$ and $H$ vanish,
does this morphism on the right factor through 
$\mathrm{CE}(\hh \stackrel{t}{\to} \gg)$, which is
indicated by the left vertical arrow of the above diagram.

\medskip 
\noindent {\bf 4. String Lie $n$-algebra valued forms.} For $\gg$ an ordinary Lie algebra and $\mu$ a degree $(2n+1)$-cocycle on $\gg$
  the situation is captured by the following diagram 
       \(
  \xymatrix{
    &
    \mbox{String-like}
    &
    \mbox{Chern-Simons}
    &
    \mbox{Chern}
    \\
    \mbox{$1$}
    &
    \mbox{$2n$}
    &
    \mbox{$2n+1$}
    &
    \mbox{$2n+1$}
    \\
    \mathrm{CE}(\gg) 
    & 
    \mathrm{CE}(\gg_{\mu})
    \ar@{<-_{)}}[l] 
    \ar@{<<-}[r] 
    &
    \mathrm{CE}(\mathrm{cs}_P(\gg))    
    \ar@{<-^{)}}[r]
    &
    \mathrm{CE}(\mathrm{ch}_P(\gg))    
    \\
    \\
    \Omega^\bullet(Y)
    \ar@{-}[r]^=
    \ar@{<-}[uu]^>>>>>>>>{(A) }|<<<<<<<<{F_A = 0}
    &
    \Omega^\bullet(Y)
    \ar@{-}[r]^=
    \ar@{<..}[uu]^>>>>>>>>{(A,B) }
      |<<<<<<<{F_A = 0 \atop dB + \mathrm{CS}_k(A) = 0}  
    &
    \Omega^\bullet(Y)
    \ar@{-}[r]^=
    \ar@{<-}[uu]^>>>>>>>>{(A,B,C)}
      |<<<<<<<<{ C = dB + \mathrm{CS}_P(A) }
    &
    \Omega^\bullet(Y)
    \ar@{<-}[uu]^>>>>>>>>{(A,C)}
      |<<<<<<<{ dC = k( (F_A)^{n+1})}
  }
  \,.
\)

Here $\mathrm{CS}_P(A)$ denotes the Chern-Simons form such that 
$d \mathrm{CS}_P(A) = P(F_A)$, given by the specific contracting homotopy.

The standard example is that corresponding to the ordinary String-extension.
\(
 \raisebox{120pt}{
  \xymatrix{
    \mathrm{CE}(\gg) 
    \ar@{=}[d]
    & 
    \mathrm{CE}(\mathrm{string}(\gg)) 
    \ar@{<-_{)}}[l] 
    \ar@{<<-}[r] 
    \ar@{=}[d]^\simeq
    &
    \mathrm{W}(\mathrm{string}_k(\gg))
    \\
    \mathrm{CE}(\gg) 
    & 
    \mathrm{CE}(\gg_\mu) 
    \ar@{<-_{)}}[l] 
    \ar@{<<-}[r] 
    &
    \mathrm{CE}(\mathrm{cs}_k(\gg))    
    \ar@{=}[u]_\simeq
    \ar@{<-_{)}}[r]
    &
    \mathrm{CE}(\mathrm{ch}_P(\gg))    
    \\
    \Omega^\bullet(Y)
    \ar@{-}[r]^=
    \ar@{<-}[u]^>>>{(A) }|<<<{F_A = 0}
    &
    \Omega^\bullet(Y)
    \ar@{-}[r]^=
    \ar@{<..}[u]^>>>{(A,B) }
      |<<<{F_A = 0 \atop dB + \mathrm{CS}_P(A) = 0}  
    &
    \Omega^\bullet(Y)
    \ar@{-}[r]^=
    \ar@{<-}[u]^>>>{(A,B,C)}
      |<<<{ C = dB + \mathrm{CS}_P(A) }
    &
    \Omega^\bullet(Y)
    \ar@{<-}[u]^>>>{(A,C)}
      |<<<{ dC = \langle F_A \wedge F_A \rangle }
  }
  }
\)
In the above, $\gg$ is semisimple with invariant bilinear form $P = \langle \cdot, \cdot \rangle$
related by transgression to the 3-cocycle $\mu = \langle\cdot, [\cdot,\cdot]\rangle$.
Then the Chern-Simons 3-form for any $\gg$-valued 1-form $A$ is
\(
  \mathrm{CS}_{\langle\cdot,\cdot\rangle}(A) = 
  \langle A \wedge dA\rangle + \frac{1}{3} \langle A \wedge [A \wedge A]\rangle
  \,.
\) 

\noindent {\bf 5. Fields of 11-dimensional supergravity.}

While we shall not discuss it in detail here, it is clear that the entire
discussion we give has a straightforward generalization to \emph{super $L_\infty$-algebras},
obtained simply by working entirely within the category of super vector spaces
(the category of $\mathbb{Z}_2$-graded vector spaces equipped with the unique non-trivial
symmetric braiding on it, which introduces a sign whenever two odd-graded vector spaces
are interchanged).

A glance at the definitions shows that, up to mere differences in terminology,
the theory of ``FDA''s (``free differential algebras'') considered in 
\cite{dAuriaFre,CastellanidAuriaFre} is nothing but
that of what we call qDGCAs here: quasi-free differential graded commutative algebras.

Using that and the interpretation of qDGCAs in terms of $L_\infty$-algebras,
one can translate everything said in \cite{dAuriaFre,CastellanidAuriaFre}
into our language here to obtain the following statement:

The field content of 11-dimensional supergravity is 
nothing but a $\gg$-valued form, for 
\(
  \gg = \mathfrak{sugra}(10,1)
\)
the Lie 3-algebra which is the string-like extension
\(
  0 \to b^2\uu(1) \to \mathfrak{sugra}(10,1) \to s\mathrm{iso}(10,1) \to 0
\)
of the \emph{super-Poincar{\'e}} Lie algebra in 10+1 dimensions, coming from
a certain 4-cocycle on that.

\(
  \raisebox{20pt}{
  \xymatrix{
    \mathrm{CE}(\mathfrak{sugra}(10,1))
    \ar@{..>}[dd]_{\mbox{\small \begin{tabular}{c} 
      (flat connection, \\ parallel gravitino,\\ vanishing supertorsion) \end{tabular}}}    
    &&
    \mathrm{W}(\mathfrak{sugra}(10,1))
    \ar@{->>}[ll]
    \ar[dd]^{\mbox{\small \begin{tabular}{c} (graviton $g$, \\ gravitino $\Psi$, \\ 3-form $C_3$) \end{tabular}}}
    \\
    \\
    \Omega^\bullet(Y)
    \ar@{-}[rr]^=
    &&
    \Omega^\bullet(Y)
  }
  }
\)


While we shall not further pursue this here, this implies the following two 
interesting issues.

\begin{itemize}
  \item
    It is known in string theory \cite{DFM} that the supergravity 3-form in fact
    consists of three parts: two Chern-Simons parts for an $\mathfrak{e}_8$ and for a
    $\mathfrak{so}(10,1)$-connection, as well as a further fermionic part, coming
    precisely from the 4-cocycle that governs $\mathfrak{sugra}(10,1)$. 
    As we discuss in \ref{lifting problem} and \ref{partra and sigma model}, the
    two Chern-Simons components can be understood in terms of certain Lie 3-algebra
    connections coming from the Chern-Simons Lie 3-algebra $\mathrm{cs}_P(\gg)$ from
    definition \ref{definition of string-like extension}. It hence seems that there should
    be a Lie $n$-algebra which nicely unifies $\mathrm{cs}_P(\mathfrak{e}_8)$,
    $\mathrm{cs}_P(\mathfrak{\mathfrak{so}(10,1)}_8)$ and $\mathfrak{sugra}(10,1)$.
    This remains to be discussed.
    
  \item
    The discussion in \ref{Loo Cartan-Ehresmann connections} shows how to obtain from
    $\gg$-valued forms globally defined connections on possibly nontrivial $\gg$-$n$-bundles.
    Applied to $\mathfrak{sugra}(10,1)$ this should yield a \emph{global} description
    of the supergravity field content, which extends the local field content considered
    in \cite{dAuriaFre,CastellanidAuriaFre} in the way a connection in a possibly nontrivial
    Yang-Mills bundle generalizes a Lie algebra valued 1-form. This should for instance allow
    to discuss supergravity \emph{instanton} solutions.
\end{itemize}

\subsection{ $L_\infty$-algebra characteristic forms}

\label{characteristic forms}

\begin{definition}
  For $\gg$ any $L_\infty$ algebra and 
  \(
    \xymatrix{
      \Omega^\bullet(Y)
      &
      \mathrm{W}(\gg)
      \ar[l]_{(A,F_A)}
    }
  \)
  any $\gg$-valued differential form, we call the composite
  \(
    \xymatrix{
      \Omega^\bullet(Y)
      &
      \mathrm{W}(\gg)
      \ar[l]_{(A,F_A)}
      &
      \mathrm{inv}(\gg)
      \ar@{_{(}->}[l]      
      \ar@/_2pc/[ll]_{\{P(F_A)\}}
    }    
  \)
  the collection of {\bf invariant forms} of the $\gg$-valued form $A$.
  We call the deRham classes $[P(F_A)]$ of the characteristic forms arising as the image of
  closed invariant polynomials
  \(
    \xymatrix{
      \Omega^\bullet(Y)
      &
      \mathrm{W}(\gg)
      \ar[l]_{(A,F_A)}
      &
      \mathrm{inv}(\gg)
      \ar@{_{(}->}[l]      
      \ar@/_2pc/[ll]_{\{P_i(F_A)\}}
      \\
      H_{\mathrm{dR}}^\bullet(Y)
      &&
      H^\bullet(\mathrm{inv}(\gg))
      \ar[ll]_{\{[P(F_A)]\}}
    }    
  \)
 the collection of {\bf characteristic classes} of the $\gg$-valued form $A$.  
\end{definition}

Recall from \ref{Lie infty-algebra cohomology} that for ordinary Lie algebras all 
invariant polynomials are
closed, while for general $L_\infty$-algebras it is only true that their 
$d_{\mathrm{W}(\gg)}$-differential is horizontal.
Notice that $Y$ will play the role of a cover of some space
$X$ soon, and that characteristic forms really live down on
$X$. We will see shortly a constraint imposed which makes the
characteristic forms descend down from the $Y$ here to such 
an $X$.

\begin{proposition}
  Under gauge transformations as in definition 
  \ref{gauge transformation of g-valued forms}, 
  characteristic classes are invariant. 
  \label{characteristic classes are invariant under gauge transformations}
\end{proposition}
\proof
  This follows from proposition \ref{concordance and chain homotopy}:
 By that proposition,
  the derivative of the concordance form $\hat A$ along the interval 
  $I = [0,1]$ is a chain homotopy
  \(
    \frac{d}{dt}\hat A(P) = [d,\iota_X] P = d \tau(P) + \iota_X (d_{\mathrm{W}(\gg)} P)
    \,.
  \)
  By definition of gauge-transformations, $\iota_X$ is vertical. 
  By definition of basic forms, $P$ is both in the kernel of $\iota_X$
  as well as in the kernel of $\iota_X \circ d$.
   Hence the right hand
  vanishes.
\endofproof

\subsubsection{Examples}

 \label{examples for characteristic forms}

\paragraph{Characteristic forms of $b^{n-1}\uu(1)$-valued forms.}

\begin{proposition}
  A $b^{n-1}\uu(1)$-valued form 
  $
    \xymatrix{
      \Omega^\bullet(Y)
      &&
      \mathrm{W}(b^{n-1}\uu(1))
      \ar[ll]_{A}
    }
  $
  is precisely an $n$-form on $Y$:
  \(
    \Omega^\bullet(Y,b^{n-1}\uu(1))
    \simeq
    \Omega^n(Y)
    \,.
  \)
  If two such $b^{n-1}\uu(1)$-valued forms are gauge equivalent according to 
  definition \ref{gauge transformation of g-valued forms}, then their curvatures
  coincide
  \(
    (\xymatrix{
      \Omega^\bullet(Y)
      &
      \mathrm{W}(b^{n-1}\uu(1))
      \ar[l]_{A}
    })
    \sim
    (\xymatrix{
      \Omega^\bullet(Y)
      &
      \mathrm{W}(b^{n-1}\uu(1))
      \ar[l]_{A'}
    })
    \hspace{10pt}
    \Rightarrow
    \hspace{10pt}
    dA = d A'
    \,.
  \)
\end{proposition}

\paragraph{BF-theory.}
  We demonstrate that the expression known in the 
  literature as the \emph{action functional for BF-theory with cosmological term}
  is the integral of an invariant polynomial for $\gg$-valued differential
  forms where $\gg$ is a Lie 2-algebra.
  Namely, let $\gg_{(2)} = (\hh \stackrel{t}{\to} \gg)$ be any
   strict Lie
  2-algebra as in \ref{Lie infty algebras}. Let 
  \(
    P = \langle \cdot, \cdot \rangle
  \)
  be an invariant bilinear form on $\gg$, hence a degree 2 invariant polynomial 
  on $\gg$. 
  According to proposition 
  \ref{invariant polynomials lift to crossed module},
  $P$ therefore also is an invariant polynomial on $\gg_{(2)}$.
  
  Now for $(A,B)$ a $\gg_{(2)}$-valued differential form on $X$,
  as in the example in \ref{Lie infty-algebra valued forms},
  \(
    \xymatrix{
      \Omega^\bullet(Y)
      &&
      \mathrm{W}(\gg_{(2)})
      \ar[ll]_{((A,B),(\beta,H))}
    }
    \,,
  \)
  one finds 
  \(
    \xymatrix{
      \Omega^\bullet(Y)
      &&
      \mathrm{W}(\gg_{(2)})
      \ar[ll]_{((A,B),(\beta,H))}
      &&
      \mathrm{inv}(\gg_{(2)})
      \ar@{_{(}->}[ll]
      \ar@/^2pc/[llll]^{P \mapsto \langle \beta, \beta \rangle}      
    }
  \)
so that the corresponding characteristic form is the 4-form
  \(
    P(\beta,H) = \langle \beta \wedge \beta\rangle
    =
    \langle
      (F_A + t(B)) \wedge (F_A + t(B)) 
    \rangle
    \,.
  \)
  Collecting terms as
  \(
    P(\beta,H) 
    =
    \underbrace{
    \langle
      F_A \wedge F_A
    \rangle
    }_{\mbox{Pontryagin term}}
    +  
    2   
    \underbrace{
    \langle
      t(B) \wedge F_A 
    \rangle
    }_{\mbox{BF-term}}
    +
    \underbrace{
      \langle
        t(B)\wedge t(B)
      \rangle
    }_{\mbox{``cosmological constant''}}
  \)
  we recognize the Lagrangian for topological Yang-Mills theory and BF theory with
  cosmological term.

  For $X$ a compact 4-manifold, the corresponding action functional
  \(
    S : \Omega^\bullet(X,\gg_{(2)}) \to \mathbb{R}
  \)
  sends $\gg_{(2)}$-valued 2-forms to the intgral of this 4-form
  \(
    (A,B)
    \mapsto
    \int_X 
    \left(
    \langle
      F_A \wedge F_A
    \rangle
    +  
    2   
    \langle
      t(B) \wedge F_A  
    \rangle
    +
      \langle
        t(B)\wedge t(B)
      \rangle
    \right)
   \,.
  \)
  The first term here is usually not considered an intrinsic part of
  BF-theory, but its presence does not affect the critical points of
  $S$.

  The critical points of $S$, i.e. the $\gg_{(2)}$-valued
  differential forms on $X$ that satisfy the equations of motion defined by
  the action $S$,
  are given by the equation
  \(
    \beta := F_A + t(B) = 0
    \,.
  \)
  Notice that this implies
  \(
    d_A t(B) = 0
  \)
  but does not constrain the full 3-curvature
  \(
    H = d_A B
  \)
  to vanish. 
  In other words, the critical points of $S$ are precisely the 
  \emph{fake flat} $\gg_{(2)}$-valued forms which precisely
  integrate to strict parallel transport 2-functors 
  \cite{GirelliPfeiffer,SWII,BS}.

  While the 4-form $\langle \beta \wedge \beta \rangle$
  looks similar to the Pontrjagin 4-form 
  $\langle F_A \wedge F_A\rangle$ for an ordinary connection
  1-form $A$, one striking difference is that
  $\langle \beta \wedge \beta\rangle$ is, in general, not closed.
  Instead, according to equation 
  \ref{differential of BF-like invariant polynomial}, we have
  \(
    d \langle \beta \wedge \beta\rangle
    = 
    2
    \langle
      \beta \wedge t(H) 
    \rangle
    \,.
  \)

  \paragraph{Remark.}
  Under  the equivalence \cite{BCSS} of the skeletal
  String Lie 2-algebra to its strict version, recalled in
  proposition \ref{string Lie 2-algebra as strict Lie 2-algebra},
  the characteristic forms for strict Lie 2-algebras apply also
  to one of our central objects of interest here,
  the String 2-connections. But a little care needs to
  be exercised here, because the strict version of the
  String Lie 2-algebra is no longer finite dimensional.

 \paragraph{Remark.} Our interpretation above of BF-theory as 
  a gauge theory for Lie 2-algebras is not unrelated to, 
  but different from the one considered in 
  \cite{GirelliPfeiffer, GirelliPfeifferPopescu}. There
  only the Lie 2-algebra coming from the infinitesimal 
  crossed module
  $( |\gg| \stackrel{0}{\to} \gg \stackrel{\mathrm{ad}}{\to}
    \mathrm{der}(\gg) )$ (for $\gg$ any ordinary Lie algebra
    and $|\gg|$ its underlying vector space, regarded as an 
    abelian Lie algebra) is considered, and the action is
   restricted to the term
   $\int \langle F_A \wedge B \rangle$. We can regard 
   the above discussion as a generalization of this approach
   to arbitrary Lie 2-algebras. Standard BF-theory 
   (with ``cosmological'' term) is 
   reproduced with the above Lagrangian by using the Lie
   2-algebra $\mathrm{inn}(\gg)$ corresponding to the
   infinitesimal crossed module 
   $(\gg \stackrel{\mathrm{Id}}{\to} \gg
      \stackrel{\mathrm{ad}}{\to} \mathrm{der}(\gg))$
   discussed in 
   \ref{examples for Loo algebras}.

\section{$L_\infty$-algebra Cartan-Ehresmann connections}

\label{Loo Cartan-Ehresmann connections}

We will now combine all of the above ingredients to produce
a definition of $\gg$-valued connections. As we shall explain,
the construction we give may be thought of as a generalization
of the notion of a Cartan-Ehresmann connection, which is
given by a Lie algebra-valued 1-form on the total space
of a bundle over base space satisfying two conditions:

\begin{itemize}

  \item first Cartan-Ehresmann condition: on the fibers the
    connection form restricts to a \emph{flat} canonical form

  \item second Cartan-Ehresmann condition: under vertical
    flows the connections transforms nicely, in such a
    way that its characteristic forms descend down to base
    space.

\end{itemize}

We will essentially interpret these two conditions as 
a pullback of the universal $\gg$-bundle, in its DGC-algebraic
incarnation as given in equation 
\ref{sequence which is dgc-version of universal bundle}.

The definition we give can also be seen as the Lie algebraic
image of a similar construction involving locally
trivializable transport $n$-functors 
\cite{BS,SWIII}, but this shall not be further discussed here.

\subsection{$\gg$-Bundle descent data}
\label{gg-descent object}
 
 \begin{definition}[$\gg$-bundle descent data]

 Given a Lie $n$-algebra $\gg$, a $\gg$-bundle descent object on $X$ 
 is a pair $(Y,A_{\mathrm{vert}})$ consisting of a 
 choice of surjective submersion $\pi : Y \to X$ 
 with \emph{connected} fibers (this condition will be dropped
 when we extend to $\gg$-connection descent objects in 
  \ref{g-connection objects})
 together with a morphism of dg-algebras
 \(
   \xymatrix{
     \Omega^\bullet_{\mathrm{vert}}(Y)
     &&
     \mathrm{CE}(\gg)
     \ar[ll]_{A_{\mathrm{vert}}}
   }
   \,.
 \)
 Two such descent objects are taken to be equivalent  
 \(
   (\xymatrix{
     \Omega^\bullet_{\mathrm{vert}}(Y)
     &
     \mathrm{CE}(\gg)
     \ar[l]_{A_{\mathrm{vert}}}
   })
    \sim
   (\xymatrix{
     \Omega^\bullet_{\mathrm{vert}}(Y')
     &
     \mathrm{CE}(\gg)
     \ar[l]_{A'_{\mathrm{vert}}}
   })
 \)
 precisely if their pullbacks $\pi_1^* A_{\mathrm{vert}}$ and 
 $\pi_2^* A'_{\mathrm{vert}}$
 to the common refinement
 \(
   \xymatrix{
     Y \times_X Y'
     \ar[r]^{\pi_1}
     \ar[d]_{\pi_2}
     &
     Y
     \ar[d]^\pi
     \\
     Y'
     \ar[r]^{\pi'}
     &
     X
   }
 \)
 are concordant in the sense of definition \ref{concordance of qDGCA morphisms}.
\end{definition}

Thus two such descent objects $A_{\mathrm{vert}}$, $A'_{\mathrm{vert}}$ 
on the same $Y$ are equivalent if there is $\eta^*_{\mathrm{vert}}$ such that
\(
  \xymatrix{
    \Omega^\bullet_{\mathrm{vert}}(Y)
    &
    \Omega^\bullet_{\mathrm{vert}}(Y \times I)
    \ar@<-3pt>[l]|{s^*}
    \ar@<+3pt>[l]|{t^*}
    &&
    \mathrm{CE}(\gg)
    \ar[ll]_{\eta^*_{\mathrm{vert}}}
    \ar@/_2pc/[lll]|{A_{\mathrm{vert}}}
	\ar@/^2.2pc/[lll]|{A'_{\mathrm{vert}}}
  }
  \,.
\)

Recall from the discussion in \ref{principal n-bundles in plan}
that the surjective submersions here play the role of 
open covers of $X$.

\subsubsection{Examples}

 \label{descent objects examples}

  \paragraph{Example: ordinary $G$-bundles.}

The following example is meant to illustrate how the notion of
descent data with respect to a Lie algebra $\gg$ as defined here 
can be related to the ordinary
notion of descent data with respect to a Lie group $G$.
   Consider the case where $\gg$ is an ordinary Lie (1-)algebra. 
   A $\gg$-cocycle then is a surjective submersion $\pi : Y \to X$
   together with a $\gg$-valued flat vertical 1-form 
   $A_{\mathrm{vert}}$ on $Y$.
  Assume the fiber of $\pi : Y \to X$ to be simply connected. 
  Then for any two points 
  $(y,y') \in Y \times_X Y$ in the same fiber we obtain an element
  $g(y,y') \in G$, where $G$ is the simply connected Lie group integrating
  $\gg$, by choosing any path $\xymatrix{y \ar[r]^\gamma & y'}$ in the fiber 
  connecting $y$ with $y'$ and forming the parallel transport  
  determined by  
  $A_{\mathrm{vert}}$
  along this path
  \(
    g(y,y') := P \exp(\int_\gamma A_{\mathrm{vert}})
    \,.
  \)
  By the flatness of $A_{\mathrm{vert}}$ and the assumption that the fibers of $Y$
  are simply connected
  \begin{itemize}
    \item $g : Y \times_X Y \to G$
      is well defined (does not depend on the choice of paths), and 
    \item
      satisfies the cocycle condition for $G$-bundles
      \(
        g : 
        \raisebox{20pt}{
        \hspace{5pt}
        \xymatrix{
          & y'
          \ar[dr]
          \\
          y\ar[rr]
          \ar[ur] &&  y''
        }
        }
        \hspace{5pt}
        \mapsto
        \hspace{5pt}
        \raisebox{20pt}{
        \xymatrix{
          & \bullet
          \ar[dr]^{g(y',y'')}
          \\
          \bullet \ar[rr]_{g(y,y'')}
          \ar[ur]^{g(y,y')} &&  \bullet
        }
        }
        \,.
      \)
  \end{itemize}
  Any such cocycle $g$ defines a $G$-principal bundle. Conversely, every $G$-principal
  bundle $P \to X$ gives rise to a structure like this by choosing $Y := P$ and 
  letting $A_{\mathrm{vert}}$ be the canonical invariant $\gg$-valued vertical 1-form 
  on $Y=P$.
  Then suppose $(Y,A_{\mathrm{vert}})$ and $(Y,A'_{\mathrm{vert}})$ are two such
  cocycles defined on the same $Y$, and let $(\hat Y:= Y \times I,\hat A_{\mathrm{vert}})$ be
  a concordance between them. Then, for every path
  \(
    \xymatrix{
      y \times \{0\}
      \ar[rr]^{\gamma}
      &&
      y \times \{1\}
    }
  \)
  connecting the two copies of a point $y \in Y$ over the endpoints of the interval, we again
  obtain a group element
  \(
    h(y) := P \exp(\int_\gamma \hat A_{\mathrm{vert}})
    \,.
  \)
  By the flatness of $\hat A,$ this is
  \begin{itemize}
    \item
      well defined in that it is independent of the choice of path;

     \item
       has the property that for all $(y,y') \in Y\times_X Y$ we have
       \(
            h :
            \hspace{5pt}
            \raisebox{30pt}{
            \xymatrix{
              y\times \{0\}
              \ar[dd]
              \ar[rr]
              &&
              y \times \{1\}
              \ar[dd]
              \\
              \\
              y'\times \{0\}
              \ar[rr]
              &&
              y' \times \{1\}
            }
           }
           \hspace{7pt}
             \mapsto
           \hspace{7pt}
            \raisebox{30pt}{
            \xymatrix{
              \bullet
              \ar[dd]_{g(y,y')}
              \ar[rr]^{h(y)}
              &&
              \bullet
              \ar[dd]^{g'(y,y')}
              \\
              \\
              \bullet
              \ar[rr]^{h(y')}
              &&
              \bullet
            }
           }
           \,.
       \)
       Therefore $h$ is a gauge transformation between $g$ and $g'$, as it should be.
  \end{itemize}

Note that there is no holonomy since the fibers are assumed to be simply connected
in this example.

\paragraph{Abelian gerbes, Deligne cohomology and $(b^{n-1}\uu(1))$-descent objects}

For the case that the $L_\infty$-algebra in question is shifted $\uu(1)$, i.e.
$\gg = b^{n-1}\uu(1)$, classes of $\gg$-descent objects on $X$ should coincide with
classes of ``line $n$-bundles'', i.e. with classes of abelian $(n-1)$-gerbes
on $X$,
hence with elements in $H^n(X,\mathbb{Z})$.
In order to understand this, we relate classes of $b^{n-1}\uu(1))$-descent objects
to Deligne cohomology.
We recall Deligne cohomology for a fixed surjective submersion $\pi : Y \to X$.
For comparison with some parts of the literature, the reader should choose $Y$ to be
the disjoint union of sets of a good cover of $X$. More discussion of
this point is in \ref{surjective submersions and differential forms}.

The following definition should be thought of this way:
 a collection of $p$-forms on fiberwise intersections of a surjective
 submersion $Y \to X$ are given. The 0-form part defines an $n$-bundle
 (an $(n-1)$-gerbe) itself, while the higher forms encode a connection
 on that $n$-bundle.

\begin{definition}[Deligne cohomology]
  \label{Deligne cohomology}
  Deligne cohomology can be understood as the cohomology on differential forms on the
  simplicial space $Y^\bullet$ given by a surjective submersion $\pi : Y \to X$, where the complex of forms 
  is taken to start as
  \(
    \xymatrix{
      0 \ar[r] & C^\infty(Y^{[n]},\mathbb{R}/\mathbb{Z}) \ar[r]^{d} &
      \Omega^1(Y^{[n]},\mathbb{R}) \ar[r]^d & \Omega^2(Y^{[n]},\mathbb{R}) \ar[r]^d & \cdots
    }
    \,,
  \)
  where the first differential, often denoted $d\mathrm{log}$ in the literature, 
  is evaluated by acting with the ordinary differential on any $\mathbb{R}$-valued representative
  of a $U(1) \simeq \mathbb{R}/\mathbb{Z}$-valued function.
 
  More in detail,
  given a surjective submersion $\pi : Y \to X$, we obtain the
  augmented simplicial space
  \(
     Y^\bullet
     =
     \left(
        \xymatrix{
          \cdots
          Y^{[3]}
          \ar@<+4pt>[r]^{\pi_1}
          \ar[r]|{\pi_2}
          \ar@<-4pt>[r]_{\pi_3}
          &
          Y^{[2]}
          \ar@<+3pt>[r]^{\pi_1}
          \ar@<-3pt>[r]_{\pi_2}
          &
          Y \ar[r]^\pi 
          &  
          Y^{[0]}
        }
     \right)
  \)
  of fiberwise cartesian powers of $Y$, 
  $Y^{[n]} := \underbrace{Y \times_X Y \times_X \cdots \times_X Y}_{\mbox{$n$ factors}}$,
  with $Y^{[0]} := X$.
  The double complex of differential forms
  \(
    \Omega^\bullet(Y^\bullet) = \bigoplus\limits_{n \in \mathbb{N}} \Omega^n(Y^\bullet)
    = 
    \bigoplus\limits_{n \in \mathbb{N}} \bigoplus\limits_{r,s \in \mathbb{N} \atop r+s = n}
    \Omega^r(Y^{[s]})
  \)
  on $Y^\bullet$ has the differential $d \pm \delta$ coming from the deRham differential
  $d$ and the alternating pullback operation
  \begin{eqnarray}
    \delta : \Omega^r(Y^{[s]}) &\to& \Omega^r(Y^{[s+1]})
  \nonumber\\
    \delta : \omega &\mapsto& \pi_1^* \omega - \pi_2^*\omega + \pi_3^*\omega + \cdots
      - (-1)^{s+1}~.
  \end{eqnarray}
  Here we take 0-forms to be valued in $\mathbb{R}/\mathbb{Z}$. 
  The map $\xymatrix{ \Omega^0(Y) \ar[r]^{d} & \Omega^1(Y)}$ takes any $\mathbb{R}$-valued 
  representative $f$ of an $\mathbb{R}/\mathbb{Z}$-valued form and sends that to the 
  ordinary $d f$. This operation is often denoted $\xymatrix{ \Omega^0(Y) \ar[r]^{d \mathrm{log}} & \Omega^1(Y)}$.
  Writing
  $\Omega^\bullet_k(Y^\bullet)$ for the space of forms that vanish on 
  $Y^{[l]}$ for $l < k$ we define (everything with respect to $Y$):
  \begin{itemize}
    \item A {\bf  Deligne $n$-cocycle} is a closed
     element in $\Omega^n(Y^\bullet)$;
    \item a {\bf flat Deligne $n$-cocycle} is a closed
      element in $\Omega^n_1(Y^\bullet)$;
    \item a {\bf Deligne coboundary} is an element in 
      $(d \pm \delta)\Omega^\bullet_1(Y^\bullet)$ (i.e. no component in $Y^{[0]} = X$);
    \item a {\bf shift of connection} is an element in 
      $(d \pm \delta)\Omega^\bullet(Y^\bullet)$ (i.e. with possibly a contribution in $Y^{[0]} = X$).
  \end{itemize}
\end{definition}

The 0-form part of a Deligne cocycle is like the transition function of a 
$U(1)$-bundle. Restricting to this part yields a group homomorphism
\(
  [\cdot]
  :
  \xymatrix{
    H^n(\Omega^\bullet(Y^\bullet))
    \ar@{->>}[r]
    &
    H^n(X,\mathbb{Z})
  }
\) 
to the integral cohomology on $X$. (Notice that the degree on the right is indeed as given,
using the total degree on the double comples $\Omega^\bullet(Y^\bullet)$ as given.)

Addition
of a Deligne coboundary is a gauge transformation. Using the fact \cite{Murray}
that the ``fundamental complex''
\(
  \xymatrix{
     \Omega^r(X)
     \ar[r]^\delta
     &
     \Omega^r(Y)     
     \ar[r]^\delta
     &
     \Omega^r(Y^{[2]})     
     \cdots
  }
\)
is exact for all $r \geq 1$, one sees that Deligne cocycles with the same class
in $H^n(X,\mathbb{Z})$
 differ by  elements in 
$(d\pm \delta)\Omega^\bullet(Y^\bullet)$.
Notice that they do not, in general, differ by an element in $\Omega^\bullet_1(Y^\bullet)$:
two Deligne cochains which differ by an element in $(d\pm \delta)\Omega^\bullet_1(Y^\bullet)$
describe equivalent line $n$-bundles \emph{with equivalent connections}, while those that
differ by something in  $(d\pm \delta)\Omega^\bullet_0(Y^\bullet)$ describe equivalent
line $n$-bundles with possibly inequivalent connections on them.

Let
\(
  v : 
  \Omega^\bullet(Y^\bullet)
  \to 
  \Omega^\bullet_{\mathrm{vert}}(Y)
\)
be the map which sends each Deligne $n$-cochain $a$ with respect to $Y$ 
to the vertical part of its $(n-1)$-form on $Y^{[1]}$
\(
  \nu : a \mapsto a|_{\Omega_{\mathrm{vert}}^{n-1}(Y^{[1]})}\,.
\) 
(Recall the definition \ref{vertical deRham complex} of $\Omega^\bullet_{\mathrm{vert}}(Y)$.)
Then we have

\begin{proposition}
  If two Deligne $n$-cocycles $a$ and $b$ over $Y$ have the same class    
  in $H^n(X,\mathbb{Z}),$ then the classes of $\nu(a)$ and $\nu(b)$
  coincide. 
\end{proposition}
\proof
  As mentioned above, 
  $a$ and $b$ have the same class in $H^n(X,\mathbb{Z})$ 
  if and only if they differ by an element in 
  $(d \pm \delta)(\Omega^\bullet(Y^\bullet))$.
  This means that on $Y^{[1]}$ they differ by an element of the 
  form
  \(
    d \alpha + \delta \beta = d\alpha + \pi^* \beta
    \,.
  \)
  Since $\pi^* \beta$ is horizontal, this is exact in 
  $\Omega^\bullet_{\mathrm{vert}}(Y^{[1]})$.
\endofproof

\begin{proposition}
  If the $(n-1)$-form parts $B, B' \in \Omega^{n-1}(Y)$
  of two Deligne $n$-cocycles differ by a $d \pm \delta$-exact
  part, then the two Deligne cocycles have the same class
  in $H^{n}(X,\mathbb{Z})$.
\end{proposition}
\proof
  If the surjective submersion is not yet contractible, we 
  pull everything back to a contractible refinement, as
  described in \ref{surj subm examples}. So assume without
  restriction of generality that all $Y^{[n]}$ are 
  contractible. This implies that 
  $H^\bullet_{\mathrm{deRham}}(Y^{[n]}) = H^0(Y^{[n]})$, which
  is a vector space spanned by the connected components of 
  $Y^{[n]}$.
   Now assume
  \(
    B - B' = d \beta  + \delta \alpha
  \)
  on $Y$.
We can immediately see that this implies that the real classes
  in $H^n(X,\mathbb{R})$ coincide:
  the Deligne cocycle property says 
  \(
     d(B - B') = \delta (H - H')
  \)
  hence, by the exactness of the deRham complex we have now,
  \(
   \delta (H - H') = \delta (d \alpha)
  \)
  and by the exactness of $\delta$ we get  $[H] = [H']$.
   
  To see that also the integral classes coincide we use
  induction over $k$ in $Y^{[k]}$.
For instance on $Y^{[2]}$ we have
  \(
    \delta (B - B') = d (A - A')
  \)
  and hence
  \(
    \delta d\beta = d(A - A')
    \,.
  \)
  Now using again the exactness of the deRham differential $d$
  this implies
  \(
    A - A' = \delta \beta + d \gamma
    \,.
  \)
  This way we work our way up to $Y^{[n]}$, where it then
  follows that the 0-form cocycles are coboundant, hence 
  that they have the same class in $H^n(X,\mathbb{Z})$.
\endofproof

\begin{proposition}
  $b^{n-1}\uu(1)$-descent objects with respect to a given surjective submersion
  $Y$ are in bijection with closed vertical $n$-forms on $Y$:
  \(
    \left\{
      \xymatrix{
        \Omega^\bullet_{\mathrm{vert}}(Y)
        &&
        \mathrm{CE}(b^{n-1}\uu(1))
        \ar[ll]_{A_{\mathrm{vert}}}
      }    
    \right\}
    \hspace{5pt}
      \leftrightarrow
    \hspace{5pt}
    \left\{
      A_{\mathrm{vert}} \in \Omega^n_\mathrm{vert}(Y)
      \,,
      \hspace{5pt}
      d A_{\mathrm{vert}} = 0
    \right\}
    \,.
  \)
  Two such $b^{n-1}\uu(1)$ descent objects on $Y$ are equivalent precisely if
  these forms represent the same cohomology class
  \(
    (A_{\mathrm{vert}}
    \sim
    A'_{\mathrm{vert}})
    \hspace{4pt}
      \Leftrightarrow
    \hspace{4pt}
    [A_{\mathrm{vert}}] = [A'_{\mathrm{vert}}]
    \in 
    H^n(\Omega^\bullet_{\mathrm{vert}}(Y))
    \,.
  \)
\end{proposition}
\proof
   The first statement is a direct consequence of the definition
   of $b^{n-1}\uu(1)$ in \ref{Lie infty algebras}. The second statement
   follows from proposition \ref{concordance and chain homotopy}
   using the reasoning as in proposition
   \ref{characteristic classes are invariant under gauge transformations}.
\endofproof

Hence two Deligne cocycles with the same class in $H^n(X,\mathbb{Z})$ indeed
specify the same class of $b^{n-1}\uu(1)$-descent data.

\subsection{Connections on $\gg$-bundles: the extension problem}

  \label{g-connection objects}

  It turns out that a useful way to conceive of the curvature on a non-flat
  $\gg$ $n$-bundle is, essentially, as the \emph{$(n+1)$-bundle with connection
  obstructing}
  the existence of a flat connection on the original $\gg$-bundle.
  This superficially trivial  statement
  is crucial for our way of coming to grips with non-flat higher bundles with connection.

\label{differential g-cocycles}

\begin{definition}[descent object for $\gg$-connection]
  Given $\gg$-bundle descent object
  \(
    \xymatrix{
       \Omega^\bullet_{\mathrm{vert}}(Y)
       &&
       \mathrm{CE}(\gg)
       \ar[ll]_{A_{\mathrm{vert}}}
    }
  \)
  as above, 
  a $\gg$-connection on it  
    is a completion of this morphism to a diagram 
  \(
    \xymatrix{
       \Omega^\bullet_{\mathrm{vert}}(Y)
       &&
       \mathrm{CE}(\gg)
       \ar[ll]_{A_{\mathrm{vert}}}
       \\
       \\
       \Omega^\bullet(Y)
       \ar@{->>}[uu]_{i^*}
       &&
       W(\gg)
       \ar@{->>}[uu]
       \ar[ll]_{(A,F_A)}^<{\ }="s"
       \\
       \\
       \Omega^\bullet(X)
       \ar@{^{(}->}[uu]_{\pi^*}
       &&
       \mathrm{inv}(\gg)
       \ar@{^{(}->}[uu]
       \ar[ll]^{\{K_i\}}_>{\ }="t"
       %
    }
    \,.
  \)
As before, two $\gg$-connection descent objects are
  taken to be equivalent, if their pullbacks to a common refinement
  are concordant.    
  \label{descent object for gg-connection}
\end{definition} 

The top square can always be completed: any representative $A \in \Omega^\bullet(Y)$ 
of $A_{\mathrm{vert}} \in \Omega^\bullet_{\mathrm{vert}}(Y)$ 
will do. 
The curvature $F_A$  is then uniquely fixed by the dg-algebra homomorphism property.
The existence of the top square then says that we have a 1-form on a total space
which resticts to a canonical flat 1-form on the firbers.
The commutativity of the lower square means that for all invariant polynomials
$P$ of $\gg,$ the form $P(F_A)$ on $Y$ is a form pulled back from $X$ and is the
differential of a form $cs$ that vanishes on vertical vector fields
\(
  \label{equation}
  P(F_A) = \pi^* K
  \,.
\)
The completion of the bottom square is hence an extra condition: it demands
that $A$ has been chosen such that its curvature $F_A$ has the property that 
the form $P(F_A) \in \Omega^\bullet(Y)$ for all invariant polynomials $P$
are lifted from base space, up to that exact part.

\begin{itemize}
  \item
    The commutativity of the top square generalizes the 
   {\bf first Cartan-Ehresmann condition}: the connection form on the
   total space restricts to a nice form on the fibers.
  \item
    The commutativity of the lower square generalizes the
      {\bf second Cartan-Ehresmann condition}: the connection form
      on the total space has to behave in such a way that the invariant
      polynomials applied to its curvature descend down to the base space.
\end{itemize} 

The pullback 
\(
  f^*(Y,(A,F_A))   = (Y', (f^* A, f^* F_A))
\)
of a $\gg$-connection descent object 
$(Y, (A,F_A))$
on a surjective submersion 
$Y$ along a morphism 
\(
  \raisebox{20pt}{
  \xymatrix{
    Y'
    \ar[dr]_{\pi'}
    \ar[rr]^f
    &&
    Y
    \ar[dl]^{\pi}
    \\
    & X
  }
  }
\)
is the $\gg$-connection descent object depicted in figure 
\ref{pullback of gg-connection descent}.

\begin{figure}[h]
  $$
    \xymatrix{
       \Omega^\bullet_{\mathrm{vert}}(Y')
       &
       \Omega^\bullet_{\mathrm{vert}}(Y)
       \ar[l]_{f^*}
       &&
       \mathrm{CE}(\gg)
       \ar[ll]_{A_{\mathrm{vert}}}
       \ar@/^2.2pc/[lll]^<<<<<<<<{f^* A_{\mathrm{vert}}}
       \\
       \\
       \Omega^\bullet(Y')
       \ar@{->>}[uu]_{i'^*}
       &
       \Omega^\bullet(Y)
       \ar[l]_{f^*}
       \ar@{->>}[uu]_<<<<<<<<{i^*}|>>>>>>>{\makebox(10,10){}}
       &&
       W(\gg)
       \ar@{->>}[uu]
       \ar[ll]_{(A,F_A)}^<{\ }="s"
       \ar@/^2.2pc/[lll]^<<<<<<<<{(f^* A, F_{f^*A}}
       \\
       \\
       \Omega^\bullet(X)
       \ar@{^{(}->}[uu]_{\pi'^*}
       &
       \Omega^\bullet(X)
       \ar[l]_{\mathrm{Id}}
       \ar@{^{(}->}[uu]_<<<<<<<<{\pi^*}|>>>>>>>{\makebox(10,10){}}
       &&
       \mathrm{inv}(\gg)
       \ar@{^{(}->}[uu]
       \ar[ll]^{\{K_i\}}_>{\ }="t"
       \ar@/^2.2pc/[lll]^<<<<<<<<{\{K_i\}}
       %
    }
    \,.
  $$
  \caption{
    \label{pullback of gg-connection descent}
    {\bf Pullback of a $\gg$-connection descent object} 
    $(Y,(A,F_A))$ along a morphism $f : Y' \to Y$ of surjective
    submersions, to $f^*(Y,(A,F_A)) = (Y', (f^*A, F_{f^* A}))$.
  }  
\end{figure}
Notice that the characteristic forms remain unaffected by such a pullback.
This way, any two $\gg$-connection descent objects may be
pulled back to a common surjective submersion. 
A concordance between two $\gg$-connection descent objects on the same
surjective submersion is depicted in figure \ref{concordance for connection}.

\begin{figure}[h]
$$
  \xymatrix{
    \Omega^\bullet_{\mathrm{vert}}(Y)
    &
    \Omega^\bullet_{\mathrm{vert}}(Y)\otimes \Omega^\bullet(I)
    \ar@<-3pt>[l]|{s^*}
    \ar@<+3pt>[l]|{t^*}
    &&
    \mathrm{CE}(\gg)
    \ar[ll]_{\eta^*_{\mathrm{vert}}}
    \ar@/_2pc/[lll]|{A_{\mathrm{vert}}}
	\ar@/^2.2pc/[lll]|<<<<<<<<{A'_{\mathrm{vert}}}
    \\
    \\
    \Omega^\bullet(Y)
    \ar@{->>}[uu]
    &
    \Omega^\bullet(Y)\otimes \Omega^\bullet(I)
    \ar@{->>}[uu]|>>>>>>>{\makebox(12,12){}}
    \ar@<-3pt>[l]|{s^*}
    \ar@<+3pt>[l]|{t^*}
    &&
    \mathrm{W}(\gg)
    \ar[ll]_{\eta^*}
    \ar@{->>}[uu]
    \ar@/_2pc/[lll]_<<<<<<<<<{(A,F_A)}|>>>>>>>>>>>{\makebox(12,12){}}
    \ar@/^2.2pc/[lll]|<<<<<<<<{(A',F_{A'})}
    \\
    \\
    \Omega^\bullet(X)
    \ar@{^{(}->}[uu]
    &
    \Omega^\bullet(X)\otimes \Omega^\bullet(I)
    \ar@{^{(}->}[uu]|>>>>>>>{\makebox(12,12){}}
    \ar@<-3pt>[l]|{s^*}
    \ar@<+3pt>[l]|{t^*}
    &&
    \mathrm{inv}(\gg)
    \ar[ll]
    \ar@{^{(}->}[uu]
    \ar@/_2pc/[lll]_<<<<<<<<<{\{K_i\}}|>>>>>>>>>>>{\makebox(12,12){}}
    \ar@/^2.2pc/[lll]|{\{K'_i\}}
  }
$$
 \caption{
    \label{concordance for connection}
    {\bf Concordance between $\gg$-connection descent objects}
    $(Y,(A,F_A))$ and $(Y,(A',F_{A'}))$ defined on the same
    surjective submersion $\pi : Y \to X$. Concordance between descent objects not on the
    same surjective submersion is reduced to this case by pulling
    both back to a common refinement, as in figure 
    \ref{pullback of gg-connection descent}.
 }
\end{figure}

Suppose $(A,F_A)$ and $(A',F_{A'})$ are descent data for $\gg$-bundles
with connection over the same $Y$ (possibly after having pulled them back
to a common refinement). Then a concordance between them is a diagram
as in figure \ref{concordance for connection}.

\begin{definition}[equivalence of $\gg$-connections]
  \label{equivalence of gg-connections}
  We say that two $\gg$-connection descent objects are \emph{equivalent as 
  $\gg$-connection descent objects} if they are connected by a \emph{vertical} 
  concordance namely one for which the derivation part of $\eta^*$
  (according to proposition \ref{concordance and chain homotopy})
  vanishes on the shifted copy $\gg^*[1] \hookrightarrow \mathrm{W}(\gg)$.  
\end{definition}

We have have a closer look at concordance and equivalence of $\gg$-connection descent objects 
in \ref{section on characteristic classes}.

\subsubsection{Examples.}

\label{examples for connection descent objects}

\paragraph{Example (ordinary Cartan-Ehresmann connection).}
Let $P \to X$ be a principal $G$-bundle and consider the descent object obtained
by setting $Y = P$ and letting $A_{\mathrm{vert}}$ be the canonical
invariant vertical flat 1-form on fibers $P$. 
Then finding the morphism
\(
 \xymatrix{
  \Omega^\bullet(Y)
  &&
  \mathrm{W}(\gg)
  \ar[ll]_{(A,F_A)}
 }
\)
such that the top square commutes amounts to finding a 1-form on the total space of the
bundle which restricts to the canonical 1-form on the fibers. This is the first of the
two conditions on a Cartan-Ehresmann connection. Then requiring the lower square to commute
implies requiring that the $2n$-forms $P_i(F_A)$, formed from the curvature 2-form $F_A$
and the degree $n$-invariant polynomials $P_i$ of $\gg$, have to descend to $2n$-forms $K_i$ 
on the base $X$.
But that is precisely the case when $P_i(F_A)$ is invariant under flows along vertical vector
fields. 
Hence it is true when $A$ satisfies the second condition of a Cartan-Ehresmann connection,
the one that says that the connection form transforms
nicely under vertical flows.

\paragraph{Further examples} appear in \ref{obstruction examples}.

\subsection{Characteristic forms and characteristic classes}
  \label{section on characteristic classes}

\begin{definition}
  For any $\gg$-connection descent object $(Y,(A,F_A))$
  we say that 
  in
  \(
    \xymatrix@R=10pt{
       \Omega^\bullet_{\mathrm{vert}}(Y)
       &&
       \mathrm{CE}(\gg)
       \ar[ll]_{A_{\mathrm{vert}}}
       \\
       \\
       \Omega^\bullet(Y)
       \ar@{->>}[uu]_{i^*}
       &&
       W(\gg)
       \ar@{->>}[uu]
       \ar[ll]_{(A,F_A)}^<{\ }="s"
       \\
       \\
       \Omega^\bullet(X)
       \ar@{^{(}->}[uu]_{\pi^*}
       &&
       \mathrm{inv}(\gg)
       \ar@{^{(}->}[uu]
       \ar[ll]|{\{K_i\}}_>{\ }="t"
       \\
       H^\bullet_{\mathrm{dR}}(X)
       &&
       H^\bullet(\mathrm{inv}(\gg))
       \ar[ll]^{\{[K_i]\}}
       %
    }
  \)    
  the $\{K_i\}$ are the \emph{characteristic forms}, while their
  deRham classes $[K_i] \in H^\bullet_{\mathrm{deRham}}(X)$
  are the \emph{characteristic classes} of $(Y,(A,F_A))$.
\end{definition}

\begin{proposition}
  If two $\gg$-connection descent objects $(Y,(A,F_A))$ and 
  $(Y',(A',F_{A'}))$ are related by a concordance as in figure 
  \ref{pullback of gg-connection descent} and figure
  \ref{concordance for connection} 
  then they have the same 
  characteristic classes:
  \(
    (Y,(A,F_A)) \sim (Y',(A',F_{A'}))
    \hspace{4pt}
    \Rightarrow
    \hspace{4pt}
    \{[K_i]\} = \{[K'_i]\}
    \,.
  \)
  \label{characteristic classes of connection descent are invariant}
\end{proposition}
\proof
  We have seen that pullback does not
  change the characteristic forms. It follows from proposition 
  \ref{characteristic classes are invariant under gauge transformations}
  that the characteristic classes are invariant under concordance.
\endofproof

\begin{corollary}
  \label{same characteristic forms for equivalent gg-connections}
  If two $\gg$-connection descent objects are equivalent according to
  definition \ref{equivalence of gg-connections}, then they even have the
  same \emph{characteristic forms}.
\end{corollary}
\proof
  For concordances between equivalent $\gg$-connection descent objects the 
  derivation part of $\eta^*$ is vertical and therefore vanishes on
  $\mathrm{inv}(\gg) = \mathrm{W}(\gg)_{\mathrm{basic}}$.
\endofproof

\paragraph{Remark. (shifts of $\gg$-connections)}

 We observe that, by the very definition of $\mathrm{W}(\gg)$, any shift in the connection
 $A$,
 $$
   A \mapsto A' = A + D \in \Omega^\bullet(Y,\gg)
 $$
 can be understood as a transformation
 $$
   \xymatrix{
     \Omega^\bullet(Y)
     &&
     \mathrm{W}(\gg)
     \ar@/_2pc/[ll]_{(A,F_A)}^{\ }="s"
     \ar@/^2pc/[ll]^{(A + D,F_{A+D})}_{\ }="t"
     \ar@{=>}^\rho "s"; "t"
   }
 $$
 with the property that $\rho$ vanishes on the \emph{non}-shifted copy 
 $\gg^* \hookrightarrow \mathrm{W}(\gg)$ and
 is nontrivial only on the shifted copy $\gg^*[1] \hookrightarrow \mathrm{W}(\gg)$: 
 in that case for $a \in \gg^*$
 any element in the unshifted copy, we have
 $$
   (A + D)(a) = A(a) + [d,\rho](a) = A(a) +  \rho( d_{\mathrm{W}(\gg)} a)
   = A(a) + \rho( d_{\mathrm{CE}(\gg)}a + \sigma a ) = A(a) + \rho(\sigma a).
 $$
 and hence $D(a) = \rho(\sigma)$, which uniquely fixes $\rho$ in terms of $D$ and
 vice versa.
 
 Therefore concordances which are not purely vertical describe homotopies between
 $\gg$-connection descent objects in which the connection is allowed to 
 vary.

\paragraph{Remark (gauge transformations versus shifts of the connections).}
We therefore obtain the following picture.
\begin{itemize}
 \item 
  Vertical concordances relate gauge equivalent $\gg$-connections
  (compare definition \ref{gauge transformation of g-valued forms} 
  of gauge transformations of $\gg$-valued forms)
  \(
   \label{pure gauge concordances}
   \mbox{
     \begin{tabular}{c}
       $(A,F_A)$ and $(A',F_{A'})$ are
       \\ 
       equivalent as $\gg$-connections
     \end{tabular}
   }   
   \hspace{10pt}
    \Leftrightarrow
   \hspace{10pt}
   \exists
   \hspace{5pt}  
   \raisebox{70pt}{
  \xymatrix{
    \Omega^\bullet_{\mathrm{vert}}(Y)
    &
    \Omega^\bullet_{\mathrm{vert}}(Y)\otimes \Omega^\bullet(I)
    \ar@<-3pt>[l]|{s^*}
    \ar@<+3pt>[l]|{t^*}
    &&
    \mathrm{CE}(\gg)
    \ar[ll]_{\eta^*_{\mathrm{vert}}}
    \ar@/_2pc/[lll]|{A_{\mathrm{vert}}}
	\ar@/^2.2pc/[lll]|<<<<<<<<{A'_{\mathrm{vert}}}
    \\
    \\
    \Omega^\bullet(Y)
    \ar@{->>}[uu]
    &
    \Omega^\bullet(Y)\otimes \Omega^\bullet(I)
    \ar@{->>}[uu]|>>>>>>>{\makebox(12,12){}}
    \ar@<-3pt>[l]|{s^*}
    \ar@<+3pt>[l]|{t^*}
    &&
    \mathrm{W}(\gg)
    \ar[ll]_{\eta^*}
    \ar@{->>}[uu]
    \ar@/_2pc/[lll]_<<<<<<<<<{(A,F_A)}|>>>>>>>>>>>{\makebox(12,12){}}
    \ar@/^2.2pc/[lll]|<<<<<<<<{(A',F_{A'})}
    \\
    \\
    \Omega^\bullet(X)
    \ar@{^{(}->}[uu]
    &
    \Omega^\bullet(X)\otimes \Omega^\bullet(I)
    \ar@{^{(}->}[uu]|>>>>>>>{\makebox(12,12){}}
    \ar@<-3pt>[l]|{s^*}
    \ar@<+3pt>[l]|{t^*}
    &&
    \mathrm{inv}(\gg)
    \ar[ll]
    \ar@{^{(}->}[uu]
  }
  }
  \,;
\)
\item
  Non-vertical concordances relate $\gg$-connection descent objects whose underlying
  $\gg$-descent object -- the underlying $\gg$-$n$-bundles -- are equivalent, but which
  possibly differ in the choice of connection on these $\gg$-bundles:
\(
  \label{not-necessarily pure gauge concordances}
   \mbox{
     \begin{tabular}{c}
       $(A',F_{A'})$ is obtained from
       \\ $(A,F_{A})$ by a
       \\ 
       shift of $\gg$-connections
     \end{tabular}
   }   
   \hspace{10pt}
    \Leftrightarrow
   \hspace{10pt}
   \exists
  \raisebox{70pt}{
  \xymatrix{
    \Omega^\bullet_{\mathrm{vert}}(Y)
    &
    \Omega^\bullet_{\mathrm{vert}}(Y)\otimes \Omega^\bullet(I)
    \ar@<-3pt>[l]|{s^*}
    \ar@<+3pt>[l]|{t^*}
    &&
    \mathrm{CE}(\gg)
    \ar[ll]_{\eta^*_{\mathrm{vert}}}
    \\
    \\
    \Omega^\bullet(Y)
    \ar@{->>}[uu]
    &
    \Omega^\bullet(Y)\otimes \Omega^\bullet(I)
    \ar@{->>}[uu]
    \ar@<-3pt>[l]|{s^*}
    \ar@<+3pt>[l]|{t^*}
    &&
    \mathrm{W}(\gg)
    \ar[ll]_{\eta^*}
    \ar@{->>}[uu]
    \ar@/_2pc/[lll]_<<<<<<<<<{(A,F_A)}|>>>>>>>>>>>{\makebox(12,12){}}
    \ar@/^2.2pc/[lll]|<<<<<<<<{(A',F_{A'})}
    \\
    \\
    \Omega^\bullet(X)
    \ar@{^{(}->}[uu]
    &
    \Omega^\bullet(X)\otimes \Omega^\bullet(I)
    \ar@{^{(}->}[uu]|>>>>>>>{\makebox(12,12){}}
    \ar@<-3pt>[l]|{s^*}
    \ar@<+3pt>[l]|{t^*}
    &&
    \mathrm{inv}(\gg)
    \ar[ll]
    \ar@{^{(}->}[uu]
    \ar@/_2pc/[lll]_<<<<<<<<<{\{K_i\}}|>>>>>>>>>>>{\makebox(12,12){}}
    \ar@/^2.2pc/[lll]|{\{K'_i\}}
  }
  }
\)
\end{itemize}


\paragraph{Remark.} This in particular shows that for a given $\gg$-descent object
$$
  \xymatrix{
    \Omega^\bullet_{\mathrm{vert}}(Y) 
    &&
    \mathrm{CE}(\gg)
    \ar[ll]_{A_{\mathrm{vert}}}
  }
$$
the corresponding characteristic classes obtained by choosing a connection $(A,F_A)$
does not depend on continuous variations of that choice of connection.

In the case of ordinary Lie (1-)algebras $\gg$ it is well known that any two connections
on the same bundle may be continuously connected by a path of connections: the
space of 1-connections is an affine space modeled on $\Omega^1(X,\gg)$. If we had
an analogous statement for $\gg$-connections for higher $L_\infty$-algebras,
we could strengthen the above statement.

\subsubsection{Examples}

\paragraph{Ordinary characteristic classes of $\gg$-bundles}

Let $\gg$ be an ordinary Lie algebra and
$(Y,(A,F_A))$ be a $\gg$-descent object corresponding
to an ordinary Cartan-Ehresmann connection as in 
\ref{examples for connection descent objects}. Using the
fact \ref{ordinary Lie cohomology is reproduced} we know that 
$\mathrm{inv}(\gg)$ contains all the ordinary invariant
polynomials $P$ on $\gg$. Hence the characteristic classes
$[P(F_A)]$ are precisely the standard characteristic classes
(in deRham cohomology) of the $G$-bundle with connection.

\paragraph{Characteristic classes of $b^{n-1}\uu(1)$-bundles.-}

For $\gg = b^{n-1}\uu(1)$ we have, according to
proposition \ref{cohomology of shifted u(1)},
$\mathrm{inv}(b^{n-1}\uu(1)) = \mathrm{CE}(b^n\uu(1))$ and hence a single degree 
$n+1$ characteristic class: the curvature itself.

This case we had already discussed in the context of Deligne cohomology in 
\ref{descent objects examples}. In particular, notice that in definition
\ref{Deligne cohomology} we had already encountered the distinction between
homotopies of $L_\infty$-algebra that are or are not pure gauge transformations,
in that they
do or do not shift the connection: what is called a \emph{Deligne coboundary} in
definition \ref{Deligne cohomology} corresponds to an equivalence of
$b^{n-1}\uu$-connection descent objects as in \ref{pure gauge concordances}, 
while what is called a \emph{shift of connection}
there corresponds to a concordance that involves a shift 
as in \ref{not-necessarily pure gauge concordances}.

\subsection{Universal and generalized $\gg$-connections}

We can generalize the discussion of $\gg$-bundles with connection
on spaces $X$, by 
\begin{itemize}
  \item
   allowing all occurrences of the algebra
of differential forms to be replaced with more general 
differential graded algebras; this amounts to admitting
generalized smooth spaces as in \ref{differential forms on spaces of maps};
  \item
    by allowing all 
Chevalley-Eilenberg and Weil algebras of $L_\infty$-algebras to 
be replaced by DGCAs which may be nontrivial in degree 0. 
This amounts to allowing not just structure $\infty$-groups
but also structure $\infty$-groupoids.
\end{itemize}

\begin{definition}[generalized $\gg$-connection descent objects]
  \label{generalized g-connection}
  Given any $L_\infty$-algebra $\gg$, and given any DGCA $A$,
  we say a $\gg$-connection descent object for $A$ is
  \begin{itemize}
    \item
      a surjection $\xymatrix{F & P \ar@{->>}[l]_{i^*}}$
      such that $A \simeq P_{\mathrm{basic}}$;
    \item
      a choice of horizontal morphisms in the diagram
      \(
     \raisebox{50pt}{
    \xymatrix{
       F
       &&
       \mathrm{CE}(\gg)
       \ar[ll]_{A_{\mathrm{vert}}}
       \\
       \\
       P
       \ar@{->>}[uu]_{i^*}
       &&
       W(\gg)
       \ar@{->>}[uu]
       \ar[ll]_{(A,F_A)}^<{\ }="s"
       \\
       \\
       A
       \ar@{^{(}->}[uu]
       &&
       \mathrm{inv}(\gg)
       \ar@{^{(}->}[uu]
       \ar[ll]^{\{K_i\}}_>{\ }="t"
       %
    }
    }
    \,;
  \)
  \end{itemize}  
\end{definition}
The notion of equivalence of these descent objects is as before.

\paragraph{Horizontal forms}

Given any algebra surjection
$$
  \xymatrix{
    F
    \\
    \\
    P
    \ar[uu]_{i^*}
  }
$$
we know from definition \ref{vertical derivations} what the ``vertical directions'' on $P$ are.
After we have chosen a $\gg$-connection on $P$, we obtain also notion of \emph{horizontal} elements
in $P$:

\begin{definition}[horizontal elements]
  \label{horizontal elements}
  Given a $\gg$-connection $(A,F_A)$ on $P$, the algebra of horizontal elements
  $$
    \mathrm{hor}_A(P) \subset P
  $$
  of $P$ with respect to this connection are those elements \emph{not} in the
  ideal generated by the image of $A$.
\end{definition}
Notice that $\mathrm{hor}_A(P)$ is in general just a graded-commutative algebra, not 
a differential algebra. Accordingly the inclusion $\mathrm{hor}_A(P) \subset P$
is meant just as an inclusion of algebras.

\subsubsection{Examples.}

\paragraph{The universal $\gg$-connection.}
The tautological example is actually of interest:
for any $L_\infty$-algebra $\gg$, 
there is a \emph{canonical} $\gg$-connection descent
object on $\mathrm{inv}(\gg)$. This comes from choosing
\(
  (\xymatrix{F & P \ar@{->>}[l]_{i^*}})
  :=
  (\xymatrix{\mathrm{CE}(\gg) & \mathrm{W}(\gg) \ar@{->>}[l]_{i^*}})
\)
and then taking the horizontal morphisms to be all identities,
as shown in figure \ref{universal g-connection}:

\begin{figure}
 \(
     \raisebox{50pt}{
    \xymatrix{
       \mathrm{CE}(\gg)
       &&
       \mathrm{CE}(\gg)
       \ar[ll]_{\mathrm{Id}}
       \\
       \\
       \mathrm{W}(\gg)
       \ar@{->>}[uu]_{i^*}
       &&
       W(\gg)
       \ar@{->>}[uu]
       \ar[ll]_{\mathrm{Id}}^<{\ }="s"
       \\
       \\
       \mathrm{inv}(\gg)
       \ar@{^{(}->}[uu]
       &&
       \mathrm{inv}(\gg)
       \ar@{^{(}->}[uu]
       \ar[ll]^{\mathrm{Id}}_>{\ }="t"
       %
    }
    }
  \)
  \caption{
    \label{universal g-connection}
    {\bf The universal $\gg$-connection} descent object.
  }
\end{figure}

We can then finally give an \emph{intrinsic} interpretation of
the decomposition of the generators of the Weil algebra $\mathrm{W}(\gg)$
of any $L_\infty$-algebra into elemenets in $\gg^*$ and elements in the shifted
copy $\gg^*[1]$, which is crucial for various of our constructions 
(for instance for the vanishing condition in \ref{vanishing condition on homotopies}):

\begin{proposition}
  The horizontal elements of $\mathrm{W}(\gg)$ with respect to the 
  univeral $\gg$-connection $(A,F_A)$ on $\mathrm{W}(\gg)$ are precisely those 
  generated entirely from the shifted copy $\gg^*[1]$:
  $$
    \mathrm{hor}_{A}(\mathrm{W}(\gg)) = \wedge^\bullet(\gg^*[1]) \subset \mathrm{W}(\gg)
    \,.
  $$
\end{proposition}

\paragraph{Line $n$-bundles on classifying spaces}

\begin{proposition}
  \label{line n-bundles on classifying spaces}
  Let $\gg$ be any $L_\infty$-algebra and 
  $P \in \mathrm{inv}(\gg)$ any closed invariant polynomial
  on $\gg$ of degree $n+1$. Let $\mathrm{cs} := \tau P$ 
  be the transgression element and $\mu := i^* \mathrm{cs}$
  the cocycle that $P$ 
  transgresses to 
  according to proposition 
  \ref{every closed inv polynomial comes from transgression}. 
  Then we canonically obtain 
  a $b^{n-1}\uu(1)$-connection descent object in 
  $\mathrm{inv}(\gg)$:
 \(
     \raisebox{50pt}{
    \xymatrix{
       \mathrm{CE}(\gg)
       &&
       \mathrm{CE}(b^{n-1}\uu(1))
       \ar[ll]_{\mu}
       \\
       \\
       \mathrm{W}(\gg)
       \ar@{->>}[uu]_{i^*}
       &&
       W(b^{n-1}\uu(1))
       \ar@{->>}[uu]
       \ar[ll]_{(\mathrm{cs},P)}^<{\ }="s"
       \\
       \\
       \mathrm{inv}(\gg)
       \ar@{^{(}->}[uu]
       &&
       **[r]\mathrm{inv}(b^{n}\uu(1)) =\mathrm{CE}(b^{n-1}\uu(1))
       \ar@{^{(}->}[uu]
       \ar[ll]^{P}_>{\ }="t"
       %
    }
    }
  \)
\end{proposition}

\paragraph{Remark.} For instance for $\gg$ an ordinary 
semisimple Lie algebra and $\mu$ its canonical 3-cocylce,
we obtain a descent object for a Lie 3-bundle which plays the
role of what is known as the canonical 2-gerbe on the
classifying space $BG$ of the simply connected group $G$
integrating $\gg$ \cite{CJMSW}. 
From the above and using \ref{examples for Loo-valued forms}
we read off that its connection 3-form is the canonical
Chern-Simons 3-form. 
We will see this again in \ref{example for oo-configuration spaces}, 
where we show that
the 3-particle (the 2-brane) coupled to the above $\gg$-connection
descent object indeed reproduces Chern-Simons theory.

\section{Higher String- and Chern-Simons $n$-bundles: the lifting problem}
\label{lifting problem}

  We discuss the general concept of weak cokernels of morphisms of 
  $L_\infty$-algebras. Then we apply this to the special problem
  of lifts of differential $\gg$-cocycles through String-like
  extensions.
  
\subsection{Weak cokernels of $L_\infty$-morphisms}
\label{weak cokernels of Lie infty-algebras}

After introducing the notion of a mapping cone of qDGCAs, the main point
here is proposition \ref{weak inverse for weak cokernel}, which establishes
the existence of the weak inverse $f^{-1}$ that was mentioned in 
\ref{plan: the lifting problem}. It will turn out to be that very weak inverse which 
picks up the information about the existence or non-existence of 
the lifts discussed in \ref{lifts of differential cocycles}.
We can define the weak cokernel for \emph{normal $L_\infty$-subalgebras}:
\begin{definition}[normal $L_\infty$-subalgebra]
  \label{normal subalgebras}
  We say a Lie $\infty$-algebra $\hh$ is a normal sub $L_\infty$-algebra of 
  the $L_\infty$-algebra $\gg$ if there is a morphism
  \(
    \xymatrix{
      \mathrm{CE}(\hh)
      &&
      \mathrm{CE}(\gg)
      \ar@{->>}[ll]_{t^*}
    }
  \)
  which the property that
  \begin{itemize}
    \item
      on $\gg^*$ it restricts to a surjective linear map $\xymatrix{\hh^* & \gg^*\ar@{->>}[l]_{t_1^*}}$;
    \item
      if $a \in \mathrm{ker}(t^*)$ then $ d_{\mathrm{CE}(\gg)} a  \in \wedge^\bullet(\mathrm{ker}(t_1^*))$.
  \end{itemize}
\end{definition}

\begin{proposition}
  For $\hh$ and $\gg$ ordinary Lie algebras, the above notion of normal sub $L_\infty$-algebra
  coincides with the standard notion of normal sub Lie algebras.
\end{proposition}
\proof
  If $a \in \mathrm{ker}(t^*)$ then for any $x,y \in \gg$ the condition says that
  $(d_{\mathrm{CE}(\gg)} a)(x \vee y) = -a(D[x\vee y]) = -a([x,y])$ 
  vanishes when $x$ or $y$ are in the image of $t$. But $a([x,y])$ vanishes when
  $[x,y]$ is in the image of $t$. Hence the condition says that if at least one of $x$ and $y$ 
  is in the image of $t$, then their bracket is.
\endofproof

\begin{definition}[mapping cone of qDGCAs; crossed module of normal sub $L_\infty$-algebras]
  Let 
  $\xymatrix{
       t : \hh \hookrightarrow \gg
    }$
   be an inclusion of a normal sub $L_\infty$-algebra $\hh$ into $\gg$.
  The mapping cone of $t^*$ is the qDGCA whose underlying graded algebra is
  \(
    \wedge^\bullet( \gg^* \oplus \hh^*[1])
  \)
  and whose differential $d_{t}$ is such that it acts on 
  generators schematically as
  \(
     d_{t} =
     \left(
       \begin{array}{cc}
         d_{\gg} & 0
         \\
         t^* & d_{\hh}
       \end{array}
     \right)
     \,.
  \)
  \label{mapping cone, detailed def}
 \end{definition}
In more detail, $d_{t^*}$ is defined as follows.
We write $\sigma t^*$ for the degree +1 derivation on 
$\wedge^\bullet(  \gg^* \oplus \hh^*[1])$ which
acts on $\gg^*$ as $t^*$ followed by a shift in degree 
and which acts on $\hh^*[1]$ as 0.
Then, for any $a \in \gg^*$, we have
\(
  d_{t} a := d_{\mathrm{CE}(\gg)} a + \sigma t^*(a)
  \,.
\)
and
\(
  d_{t} \sigma t^*(a) := -\sigma t^*(d_{\mathrm{CE}(\gg)} a )
    = - d_{t} d_{\mathrm{CE}(\gg)} a
  \,.
\)
Notice that the last equation 
\begin{itemize}
  \item 
    defines $d_t$ on all of $\hh^*[1]$ since $t^*$ is surjective;
  \item
    is well defined in that it agrees for $a$ and $a'$ if $t^*(a) = t^*(a')$,
    since $t$ is normal.
\end{itemize}

\begin{proposition}
  The differential $d_{t}$ defined this way indeed satisfies
  $(d_{t})^2 = 0$.
  \label{nilpotency for mapping cone}
\end{proposition}
\proof
  For $a \in \gg^*$ we have
  \(
    d_{t} d_{t} a = d_{t}( d_{\mathrm{CE}(\gg)}a + \sigma t^*(a))
      = \sigma t^*(d_{\mathrm{CE}(\gg)}a) - \sigma t^*(d_{\mathrm{CE}(\gg)}a) = 0
    \,.
  \)
  Hence $(d_{t})^2$ vanishes on $\wedge^\bullet(\gg^*)$.
  Since
  \(
    d_{t} d_{t} \sigma t^*(a) = - d_{t} d_{t} d_{\mathrm{CE}(\gg)} a
  \)
  and since $d_{\mathrm{CE}(\gg)} a \in \wedge^\bullet(\gg^*)$ this implies
  $(d_{t})^2 = 0$.
\endofproof

We write
$\mathrm{CE}(\hh \stackrel{t}{\hookrightarrow} \gg)
   :=
   \left(
     \wedge^\bullet( \gg^* \oplus  \hh^*[1])
     ,
     d_{t}
   \right)$
for the resulting qDGCA and
  $(\hh \stackrel{t}{\hookrightarrow} \gg)$
  for the corresponding $L_\infty$-algebra.

The next proposition asserts that 
$\mathrm{CE}(\hh \stackrel{t}{\hookrightarrow} \gg)$
is indeed a (weak) kernel of $t^*$.

\begin{proposition}
  There is a canonical morphism
  $\xymatrix{
     \mathrm{CE}(\gg)
     &
     \mathrm{CE}(\hh \stackrel{t}{\hookrightarrow} \gg)
     \ar[l]
   }$
  with the property that 
  \(
   \xymatrix{
     \mathrm{CE}(\hh)
     &
     \mathrm{CE}(\gg)
     \ar[l]_{t^*}
     &
     \mathrm{CE}(\hh \stackrel{t}{\hookrightarrow} \gg)
     \ar[l]
     \ar@/^2pc/[ll]^0_{\ }="t"
     \ar@{<=}^\tau "t"+(0,4); "t"
   }
   \,.
  \)
\end{proposition}
\proof
  On components, this morphism  is the identity on $\gg^*$
  and 0 on $\hh^*[1]$. One checks that this respects the
  differentials.
  The homotopy to the 0-morphism sends
  \(
    \tau : \sigma t^*(a) \mapsto t^*(a)
    \,.
  \)
  Using definition \ref{transformation of qDGCA morphisms}
  one checks that then indeed
  \(
    [d,\tau] : a \mapsto \tau (d_{\mathrm{CE}(\gg)} a + \sigma t^* a)
      = a 
  \)
  and
  \(
    [d,\tau] : \sigma t^* a \mapsto d_{\mathrm{CE}(\gg)} a
                   + \tau(-\sigma t^* (d_{\mathrm{CE}(\gg)} a))
        = 0
        \,.
  \)
  Here the last step makes crucial use of 
  the condition 
  \ref{vanishing condition on homotopies} which demands that
  \(
    \tau(d_{\mathrm{W}(\hh \stackrel{t}{\hookrightarrow} \gg)}
      \sigma t^* a - 
      d_{\mathrm{CE}(\hh \stackrel{t}{\hookrightarrow} \gg)}
      \sigma t^* a)
      = 0
  \)
  and  the formula
  (\ref{formula for chain homotopy}) which induces precisely
  the right combinatorial factors.
\endofproof

Notice that not only is $\mathrm{CE}(\hh \stackrel{t}{\hookrightarrow} \gg)$
in the kernel of $t^*$, it is indeed the universal object with this
property, hence is \emph{the} kernel of $t^*$ (of course up to equivalence).

\begin{proposition}
Let 
$\xymatrix{
    \mathrm{CE}(\hh)
    &
    \mathrm{CE}(\gg)
    \ar@{->>}[l]_{t^*}
    &
    \mathrm{CE}(\ff)
    \ar@{_{(}->}[l]_{u^*}
  }$
be a sequence of qDGCAs with $t^*$ normal, as above, and with the property that 
$u^*$
restricts, on the underlying vector spaces of generators, to the
kernel of the linear map underlying $t^*$.
Then there is a unique morphism 
$f : \mathrm{CE}(\ff) \to \mathrm{CE}(\hh \stackrel{t}{\hookrightarrow} \gg)$ 
such that
\(
  \xymatrix{
    \mathrm{CE}(\hh)
    &
    \mathrm{CE}(\gg)
    \ar[l]_{t^*}
    &
    \ar[l]
    \mathrm{CE}(\hh \stackrel{t}{\hookrightarrow} \gg)
    \\
    &
    \mathrm{CE}(\ff)
    \ar[u]_{u^*}
    \ar@{-->}[ur]_f
  }
  \,.
\)
\label{universal property of weak kernel}
\end{proposition}
\proof
  The morphism $f$ has to be in components the same as 
  $\mathrm{CE}(\gg) \leftarrow \mathrm{CE}(\ff)$.
  By the assumption that this is in the kernel of $t^*$, the differentials
  are respected. 
\endofproof

\paragraph{Remark.} There should be a generalization of the entire
discussion where $u^*$ is not restricted to be 
the kernel of $t^*$ on generators. However, for our application here, this
simple situation is all we need.

\begin{proposition}
For a string-like extension $\gg_\mu$ from definition \ref{definition of string-like extension},
the morphism
\(
  \xymatrix{
    \mathrm{CE}(b^{n-1}\uu(1))
    &&
    \mathrm{CE}(\gg_\mu)
    \ar@{->>}[ll]_{t^*}
   }
\)
is normal in the sense of definition \ref{normal subalgebras}.
\end{proposition}

\begin{proposition}
In the case that the sequence
\(
  \xymatrix{
    \mathrm{CE}(\hh)
    &&
    \mathrm{CE}(\gg)
    \ar@{->>}[ll]_{t^*}
    &&
    \mathrm{CE}(\ff)
    \ar@{_{(}->}[ll]_{u^*}
  } 
\)
above is a String-like extension
\(
  \xymatrix{
    \mathrm{CE}(b^{n-1}\uu(1))
    &&
    \mathrm{CE}(\gg_\mu)
    \ar@{->>}[ll]_{t^*}
    &&
    \mathrm{CE}(\gg)
    \ar@{_{(}->}[ll]_{u^*}
  } 
\)
from proposition \ref{the string-like extension sequence}
or the corresponding Weil-algebra version
\(
  \xymatrix{
    \mathrm{W}(b^{n-1}\uu(1))
    \ar@{-}[d]^=
    &&
    \mathrm{W}(\gg_\mu)
    \ar@{-}[d]^=
    \ar@{->>}[ll]_{t^*}
    &&
    \mathrm{W}(\gg)
    \ar@{-}[d]^=
    \ar@{_{(}->}[ll]_{u^*}
    \\
    \mathrm{CE}(\mathrm{inn}(b^{n-1}\uu(1)))
    &&
    \mathrm{CE}(\mathrm{inn}(\gg_\mu))
    \ar@{->>}[ll]_{t^*}
    &&
    \mathrm{CE}(\mathrm{inn}(\gg))
    \ar@{_{(}->}[ll]_{u^*}
  } 
\)
  the 
  morphisms
  $f : \mathrm{CE}(\ff) \to \mathrm{CE}(\hh \stackrel{t}{\hookrightarrow} \gg)$
  and
  $\hat f : \mathrm{W}(\ff) \to \mathrm{W}(\hh \stackrel{t}{\hookrightarrow} \gg)$
  have weak inverses
  $f^{-1} : \mathrm{CE}(\hh \stackrel{t}{\hookrightarrow} \gg) \to \mathrm{CE}(\ff)$ 
  and
  $\hat f^{-1} : \mathrm{W}(\hh \stackrel{t}{\hookrightarrow} \gg) \to \mathrm{W}(\ff)$ 
  \,, respectively.
     \,.
    
  \label{weak inverse for weak cokernel}
\end{proposition}
\proof
    We first construct a morphism $f^{-1}$ and then show that it is
    weakly inverse to $f$. The statement for $\hat f$ the follows from the
    functoriality of forming the Weil algebra, proposition \ref{functoriality of W}.
    Start by choosing a splitting of the vector space $V$ underlying $\gg^*$
    as
    \(
      V = \mathrm{ker}(t^*) \oplus V_1
      \,.
    \)
    This is the non-canonical choice we need to make.
    Then take the component map of $f^{-1}$
    to be the identity on $\mathrm{ker}(t^*)$ and 0 on $V_1$.
    Moreover, for $a \in V_1$ set
    \(
      f^{-1} : \sigma t^*(a) \mapsto  -(d_{\mathrm{CE}(\gg)} a)|_{\wedge^\bullet\mathrm{ker}(t^*)}
      \,,
      \label{crucial weak inverse}
    \)
    where the restriction is again with respect to the chosen splitting of $V$.
We check that this assignment, extended as an algebra homomorphism,
   does respect the differentials.
    
        For $a \in \mathrm{ker}(t^*)$ we have
    \(
      \raisebox{30pt}{
      \xymatrix{
        a 
        \ar@{|->}[r]^{d_t}
        \ar@{|->}[d]_{f^{-1}}
        &
        d_{\mathrm{CE}(\gg)} a
        \ar@{|->}[d]_{f^{-1}}
        \\
        a
        \ar@{|->}[r]^<<<<<<{d_{\mathrm{CE}(\ff)}}
        &
        d_{\mathrm{CE}(\gg)} a
      }
      }
    \)
    using the fact that $t^*$ is normal.
    For $a \in V_1$ we have
    \(
      \raisebox{30pt}{
      \xymatrix{
        a 
        \ar@{|->}[r]^{d_t}
        \ar@{|->}[d]_{f^{-1}}
        &
        d_{\mathrm{CE}(\gg)} a + \sigma t^*(a)
        \ar@{|->}[d]_{f^{-1}}
        \\
        0
        \ar@{|->}[r]^<<<<<<{d_{\mathrm{CE}(\ff)}}
        &
        (d_{\mathrm{CE}(\gg)}a)|_{\wedge^\bullet\mathrm{ker}(t^*)}
        -
        (d_{\mathrm{CE}(\gg)}a)|_{\wedge^\bullet\mathrm{ker}(t^*)}        
      }
      }
      \,.
    \)
    and
    \(
      \raisebox{30pt}{
      \xymatrix{
        \sigma t^*(a) 
        \ar@{|->}[r]^{d_t}
        \ar@{|->}[d]_{f^{-1}}
        &
        -\sigma t^*(d_{\mathrm{CE}(\gg)}a)
        \ar@{|->}[d]_{f^{-1}}
        \\
        -(d_{\mathrm{CE}(\gg)}a)|_{\wedge^\bullet \mathrm{ker}(t^*)}
        \ar@{|->}[r]^{d_{\mathrm{CE}(\ff)}}
        & 
        - d_{\mathrm{CE}(\ff)}((d_{\mathrm{CE}(\gg)}a)|_{\mathrm{ker}(t^*)})
      }
      }
      \,.
      \label{a subtle condition on the weak inverse}
    \) 
    This last condition happens to be satisfied for the examples stated in the 
    proposition. The details for that are discussed in 
    \ref{weak cokernel examples} below. 
    By the above, $f^{-1}$ is indeed a morphism of qDGCAs.
    
    Next we check that $f^{-1}$ is a weak inverse of $f$. Clearly
    \(
      \xymatrix{
        \mathrm{CE}(\ff) 
        &
        \ar[l]
        \mathrm{CE}(\hh \stackrel{t}{\hookrightarrow} \gg)
        &
        \ar[l]
        \mathrm{CE}(\ff) 
      }
    \)
    is the identity on $\mathrm{CE}(\ff)$.
    What remains is to construct a homotopy
    \(
      \xymatrix{
        \mathrm{CE}(\hh \stackrel{t}{\hookrightarrow} \gg)
        &
        \mathrm{CE}(\ff)
        \ar[l]
        &
        \mathrm{CE}(\hh \stackrel{t}{\hookrightarrow} \gg)        
        \ar[l]
        \ar@/^2.3pc/[ll]^{\mathrm{Id}}_{\ }="t"
        \ar@{=>}^\tau "t"+(0,6); "t"
      }
      \,.
    \)
    One checks that this is accomplished by taking $\tau$ to
    act on $\sigma V_1$ as $\tau : \sigma V_1 \stackrel{\simeq}{\to} V_1$
    and extended suitably.
\endofproof

\subsubsection{Examples}

\label{weak cokernel examples}

\paragraph{Weak cokernel for the String-like extension.}

Let our sequence
\(
  \xymatrix{
    \mathrm{CE}(\hh)
    &&
    \mathrm{CE}(\gg)
    \ar@{->>}[ll]_{t^*}
    &&
    \mathrm{CE}(\ff)
    \ar@{_{(}->}[ll]_{u^*}
  } 
\)
be a String-like extension
\(
  \xymatrix{
    \mathrm{CE}(b^{n-1}\uu(1))
    &&
    \mathrm{CE}(\gg_\mu)
    \ar@{->>}[ll]_{t^*}
    &&
    \mathrm{CE}(\gg)
    \ar@{_{(}->}[ll]_{u^*}
  } 
\)
from proposition \ref{the string-like extension sequence}. Then the mapping cone
Chevalley-Eilenberg algebra
\(
  \mathrm{CE}(b^{n-1}\uu(1) \hookrightarrow \gg_\mu)
\)
is
\(
  \wedge^\bullet(
    \gg^* \oplus \mathbb{R}[n] \oplus \mathbb{R}[n+1]
  )
\)
with differential given by
\begin{eqnarray}
  d_{t}|_{\gg^*} &=& d_{\mathrm{CE}(\gg)}
\\
  d_t|_{\mathbb{R}[n]} &=& -\mu + \sigma 
\\
  d_t|_{\mathbb{R}[n+1]} &=& 0 \,.
 \end{eqnarray}
(As always, $\sigma$ is the canonical degree shifting isomorphism on generators
extended as a derivation.)
The morphism 
\(
  \xymatrix{
    \mathrm{CE}(\gg)
    &&
    \mathrm{CE}(b^{n-1}\uu(1) \hookrightarrow \gg_\mu)
    \ar[ll]_{f^{-1}}^\simeq
  }
\)
acts as
\begin{eqnarray}
  f^{-1}|_{\gg^*} &=& \mathrm{Id}
\\
  f^{-1}|_{\mathbb{R}[n]} &=& 0
\\
  f^{-1}|_{\mathbb{R}[n+1]} &=& \mu \,.
  \end{eqnarray}
To check the condition in equation \ref{a subtle condition on the weak inverse}
explicitly in this case, let $b \in \mathbb{R}[n]$ and write $b := t^* b$
for simplicity (since $t^*$ is the identity on $\mathbb{R}[n]$). Then 
\(
  \xymatrix{
    \sigma b
    \ar@{|->}[r]^{d_t}
    \ar[d]^{f^{-1}}
    &
    0
    \ar[d]^{f^{-1}}
    \\
    \mu
    \ar@{|->}[r]^{d_{\mathrm{CE}(\gg)}}
    &
    0    
  }
\)
does commute.

\paragraph{Weak cokernel for the String-like extension in terms of the Weil algebra.}

We will also need the analogous discussion not for the Chevalley-Eilenberg algebras,
but for the corresponding Weil algebras.
To that end consider now the sequence
\(
  \xymatrix{
    \mathrm{W}(b^{n-1}\uu(1))
    &&
    \mathrm{W}(\gg_\mu)
    \ar@{->>}[ll]_{t^*}
    &&
    \mathrm{W}(\gg)
    \ar@{_{(}->}[ll]_{u^*}
  } 
  \,.
\)
This is handled most conveniently by inserting the isomorphism
\(
  \mathrm{W}(\gg_\mu) \simeq \mathrm{CE}(\mathrm{cs}_P(\gg))
\)
from proposition \ref{Chern and Chern-Simons}
as well as the identitfcation
\(
  \mathrm{W}(\gg) = \mathrm{CE}(\mathrm{inn}(\gg))
\)
such that we get
\(
  \xymatrix{
    \mathrm{CE}(\mathrm{inn}(b^{n-1}\uu(1)))
    &&
    \mathrm{CE}(\mathrm{cs}_P(\gg))
    \ar@{->>}[ll]_{t^*}
    &&
    \mathrm{CE}(\mathrm{inn}(\gg))
    \ar@{_{(}->}[ll]_{u^*}
  } 
  \,.
\)
Then we find that the mapping cone algebra $\mathrm{CE}(b^{n-1}\uu(1) \hookrightarrow
   \mathrm{cs}_P(\gg))$
is
\(
  \wedge^\bullet( \gg^* \oplus \gg^*[1] \oplus 
     (\mathbb{R}[n] \oplus \mathbb{R}[n+1]) 
     \oplus  
     (\mathbb{R}[n+1] \oplus \mathbb{R}[n+2]))
  \,.
  \)
Write $b$ and $c$ for the canonical basis elements of 
$\mathbb{R}[n] 
\oplus \mathbb{R}[n+1]$, then the differential is characterized by

\begin{eqnarray}
  d_t|_{\gg^* \oplus \gg^*} &=& d_{\mathrm{W}(\gg)}
\\
  d_t &:& b \mapsto c - \mathrm{cs} + \sigma b
\\
  d_t &:& c \mapsto P + \sigma c
\\
  d_t &:& \sigma b \mapsto - \sigma c
  \\
  d_t &:& \sigma c \mapsto 0
  \,.
\end{eqnarray}
Notice above the relative sign between $\sigma b$ and $\sigma c$. 
This implies that the canonical injection
\(
  \xymatrix{ 
    \mathrm{CE}(b^{n-1}\uu(1) \hookrightarrow
    \mathrm{cs}_P(\gg))
    &&
    \mathrm{W}(b^n \uu(1))
    \ar[ll]_i
  }
  \label{injection of W(bnu) into cs mapping cone}
\)
also carries a sign: if we denote the degree $n+1$ and $n+2$ generators
of $\mathrm{W}(b^n \uu(1))$ by $h$ and $dh$, then
\begin{eqnarray}
  i &:& h \mapsto \sigma b
\\
  i &:& dh \mapsto - \sigma c
  \,.
\end{eqnarray}
This sign has no profound structural role, but we need to carefully keep track
of it, for instance in order for our examples in \ref{obstruction examples} 
to come out right.
The morphism
\(
  \xymatrix{
     \mathrm{CE}(b^{n-1}\uu(1) \hookrightarrow
    \mathrm{cs}_P(\gg))
    &&
    \mathrm{W}(\gg)
    \ar[ll]_>>>>>>>>>>>>>{f^{-1}}^{\simeq}
  }
\)
acts as
\begin{eqnarray}
  f^{-1}|_{\gg^* \oplus \gg^*[1]} &=& \mathrm{Id}
\\
  f^{-1} : \sigma b &\mapsto& \mathrm{cs}
\\
  f^{-1} : \sigma c &\mapsto& - P
  \,.
\end{eqnarray}
Again, notice the signs, as they follow from the general prescription
in proposition \ref{weak inverse for weak cokernel}.
We again check explicitly equation (\ref{a subtle condition on the weak inverse}):
\(
  \raisebox{20pt}{
  \xymatrix{
    \sigma b
    \ar@{|->}[r]^{d_t}
    \ar@{|->}[d]^{f^{-1}}
    &
    - \sigma c
    \ar@{|->}[d]^{f^{-1}}
    \\
    \mathrm{cs}
    \ar@{|->}[r]^{d_{\mathrm{W}(\gg)}}
    &
    P 
  }
  }
  \,.
\)

\subsection{Lifts of $\gg$-descent objects through String-like extensions}

We need the above general theory for the special case where
we have the mapping cone
$\mathrm{CE}(b^{n-1}\uu(1) \hookrightarrow \gg_\mu)$
as the weak kernel of the left morphism in a String-like extension
\(
  \xymatrix{
    \mathrm{CE}(b^{n-1}\uu(1))
    &&
    \mathrm{CE}(\gg_\mu)
    \ar@{->>}[ll]
    &&
    \mathrm{CE}(\gg)
    \ar@{_{(}->}[ll]
  }
\)
coming from an $(n+1)$ cocycle $\mu$ on an ordinary Lie algebra $\gg$.
In this case $\mathrm{CE}(b^{n-1}\uu(1) \hookrightarrow \gg_\mu)$ looks like
\(
  \mathrm{CE}(b^{n-1}\uu(1) \hookrightarrow \gg_\mu)
  =
  (\wedge^\bullet(\gg^* \oplus \mathbb{R}[n] \oplus \mathbb{R}[n+1]), d_{t}  )
  \,.
\)
By chasing this through the above definitions, we find
\begin{proposition}
The morphism
\(
  f^{-1} : \mathrm{CE}(b^{n-1}\uu(1) \hookrightarrow \gg_\mu)
   \to \mathrm{CE}(\gg)
\)
acts 
as the identity on $\gg^*$
\(
  f^{-1}|_{\gg^*} = \mathrm{Id}\,,
\) 
vanishes on $\mathbb{R}[n]$
\(
  f^{-1}|_{\mathbb{R}[n]} : b \mapsto  0,
\) 
and satisfies
\(
  f^{-1}|_{\mathbb{R}[n+1]} : \sigma t^* b \mapsto \mu
  \,.
\)

\label{weak inverse for string extension}
\end{proposition}

Therefore we find the $(n+1)$-cocycle 
\(
  \xymatrix{
    \Omega_{\mathrm{vert}}^\bullet(Y)
    &&
    \mathrm{CE}(b^{n}\uu(1))
    \ar[ll]_{\hat A_{\mathrm{vert}}}
  }
\)
obstructing the lift of a 
$\gg$-cocycle
\(
  \xymatrix{
    \Omega_{\mathrm{vert}}^\bullet(Y)
    &&
    \mathrm{CE}(\gg)
    \ar[ll]_{A_{\mathrm{vert}}}
  },
\)
according to \ref{plan: the lifting problem} given by 
  \(
   \xymatrix@C=8pt{
   &&&&&&
   \mathrm{CE}(b^{n-1}\uu(1) \hookrightarrow \gg_\mu)
   \ar[dllll]
   \ar[dll]|{f^{-1}}
   &&
   \mathrm{CE}(b^n \uu(1))
   \ar@{_{(}->}[ll]_<<<j
   \ar[ddlllll]|{\hat A_{\mathrm{vert}}}
   \\
     \mathrm{CE}(b^{n-1}\uu(1))
     &&
     \mathrm{CE}(\gg_\mu)
     \ar@{->>}[ll]_<<<<{i^*}
     \ar@{..>}[dr]
     &&
     \mathrm{CE}(\gg)     
     \ar@{_{(}->}[ll]
      \ar[dl]|{A_{\mathrm{vert}}}
     \\
     &&&\Omega^\bullet_{\mathrm{vert}}(Y)
   },
 \)
to be the $(n+1)$-form
\(
  \mu(A_{\mathrm{vert}})  \in \Omega_{\mathrm{vert}}^{n+1}(Y)
  \,.
\)

\begin{proposition} 
Let $A_{\mathrm{vert}} \in \Omega_{\mathrm{vert}}^1(Y,\gg)$ 
be the cocycle 
of a $G$-bundle $P \to X$ for $\gg$ semisimple and let 
$\mu = \langle \cdot, [\cdot,\cdot]\rangle$ be the canonical
3-cocycle. Then $\gg_\mu$ is the standard String Lie 3-algebra
and the obstruction to lifting $P$ to a String 2-bundle,
i.e. lifitng to a $\gg_\mu$-cocycle, is the Chern-Simons 3-bundle
with cocycle given by the vertical 3-form
\(
  \langle A_{\mathrm{vert}} \wedge [A_{\mathrm{vert}} \wedge A_{\mathrm{vert}}] \rangle
  \in 
  \Omega^3_{\mathrm{vert}}(Y)
  \,.
\)
\label{bundle cocycle obstruction for lift through string extension}
\end{proposition}

In the following we will express these obstruction in a more familiar
way in terms of their characteristic classes. In order to do that,
we first need to generalize the discussion to differential $\gg$-cocycle.
But that is now straightforward.

\subsubsection{Examples}

 \label{examples for lifting discussion}

The continuation of the discussion of \ref{examples for Lie 00-algebra cohomology} 
to coset spaces
gives a classical illustration of the lifting construction considered here.

\paragraph{Cohomology of coset spaces.}

The above relation between the cohomology of groups and that of their
Chevalley-Eilenberg qDGCAs generalizes to coset spaces.
This  also illustrates the constructions which are discussed later in 
\ref{lifting problem}.
Consider the case of an ordinary extension of (compact connected) Lie groups:
\(
 1\to H \to G\to G/H\to 1
 \)
or even the same sequence in which $G/H$ is only a homogeneous space and not itself a group.
For a closed connected subgroup $t:H\hookrightarrow G$, there is the  induced map $Bt:BH\to BG$ 
and a commutative diagram 
\(
  \xymatrix{
    \mathrm{W}(\gg)
    \ar[rr]^{dt^*}
    &&
    \mathrm{W}(\hh)
    \\
    \\
    \wedge^\bullet P_G
    \ar[rr]^{dt^*}
    \ar@{^{(}->}[uu]
    &&
    \wedge^\bullet P_H
    \ar@{^{(}->}[uu]
  }.
\)
By analyzing the fibration sequence
\( 
G/H\to EG/H\simeq BH\to BG, 
\)
Halperin and Thomas \cite{halperin-thomas} show
there is a morphism 
\(
 \wedge^\bullet (P_G\oplus Q_H)\to \Omega^\bullet(G/K)
\)
 inducing an isomorphism in cohomology.
 It is not hard to see that their morphism factors through
 \(
 \wedge^\bullet(\gg^*\oplus \hh^*[1]).
 \)
In general, the homogeneous space $G/H$ itself is not a group, but in
case of an extension
$H\to G\to K$, we also have $BK$ and the sequences $K\to BH\to BG$ and
$BH\to BG\to BK$.  Up to homotopy equivalence, the fiber of the bundle  $BH\to BG$ is $K$ and that of $BG\to BK$ is $BH.$
In particular, consider
an extension of $\gg$ by a String-like Lie $\infty$-algebra
 \(
   \xymatrix{
     \mathrm{CE}(b^{n-1}\uu(1))
     &&
     \mathrm{CE}(\gg_\mu)
     \ar@{->>}[ll]_i
     &&
     \mathrm{CE}(\gg)     
     \ar@{_{(}->}[ll]
   }\,.
 \)
Regard $\gg$ now as the quotient $\gg_\mu/b^{n-1}\uu(1)$ and recognize that corresponding to $BH$
we have $b^{n}\uu(1)$. Thus we have a quasi-isomorphism 
\( 
\mathrm{CE}(b^{n-1}\uu(1)\hookrightarrow\gg_\mu)\simeq \mathrm{CE}(\gg)
\)
 and hence a morphism 
 \(
 \mathrm{CE}(b^{n}\uu(1))\to\mathrm{CE}(\gg). 
 \)
Given a $\gg$-bundle cocycle
 \(
   \xymatrix{
     &
     \mathrm{CE}(\gg)
     \ar[dl]^{A_{\mathrm{vert}}}
     \\
     \Omega^\bullet_{\mathrm{vert}}(Y)
   }
 \)
 and given an extension of $\gg$ by a String-like Lie $\infty$-algebra
 \(
   \xymatrix{
     \mathrm{CE}(b^{n-1}\uu(1))
     &&
     \mathrm{CE}(\gg_\mu)
     \ar@{->>}[ll]_i
     &&
     \mathrm{CE}(\gg)     
     \ar@{_{(}->}[ll]
   }
 \)
 we ask if it is possible to \emph{lift the cocycle} through this 
 extension, i.e. to find a dotted arrow in
 \(
   \label{the lift to be found}
   \xymatrix@C=8pt{
     \mathrm{CE}(b^{n-1}\uu(1))
     &&
     \mathrm{CE}(\gg_\mu)
     \ar@{->>}[ll]
     \ar@{..>}[dr]
     &&
     \mathrm{CE}(\gg)     
     \ar@{_{(}->}[ll]
      \ar[dl]^{A_{\mathrm{vert}}}
     \\
     &&&\Omega^\bullet_{\mathrm{vert}}(Y)
  }
  \,.
 \)
 In general this is not possible. Indeed, consider the map $A'_{\mathrm{vert}}$ given by $ \mathrm{CE}(b^n \uu(1))\to CE(\gg)$
 composed with $A_{\mathrm{vert}}.$
 The nontriviality of the $b^n\uu(1)$-cocycle $A'_{\mathrm{vert}}$
 is the obstruction to constructing the desired lift.

\subsection{Lifts of $\gg$-connections through String-like extensions}
\label{lifts of differential cocycles}

In order to find the obstructing characteristic classes, we would like to extend the
above lift \ref{the lift to be found}
of $\gg$-descent objects to a lift of $\gg$-connection descent objects
extending them,
according to \ref{differential g-cocycles}.
Hence we would like first to \emph{extend} $A_{\mathrm{vert}}$ to $(A,F_A)$
\(
  \xymatrix{
    \mathrm{CE}(b^{n-1}\uu(1))
    \ar@{<-}[r]
    \ar@{<-}[dd]
    &
    \mathrm{CE}(\gg_\mu)
    \ar@{<-}[rr]
    \ar@{<-}[dd]
    &&
    \mathrm{CE}(\gg)
    \ar@{<-}[dd]
    \\
     && \Omega^\bullet_{\mathrm{vert}}(Y)
       \ar@{<-}[ur]|{A_{\mathrm{vert}}}
       \ar@{<-}[dd]
    \\
    \mathrm{W}(b^{n-1}\uu(1))
    \ar@{<-}[r]
    \ar@{<-}[dd]
    &
    \mathrm{W}(\gg_\mu)
    \ar@{<-}[rr]|{\hspace{12pt}}
    \ar@{<-}[dd]
    &&
    \mathrm{W}(\gg)
    \ar@{<-}[dd]
    \\
     && \Omega^\bullet(Y) 
       \ar@{<-}[ur]|{(A,F_A)}
       \ar@{<-}[dd]
    \\
    \mathrm{inv}(b^{n-1}\uu(1))
    \ar@{<-}[r]
    &
    \mathrm{inv}(\gg_\mu)
    \ar@{<-}[rr]|{\hspace{12pt}}
    &&
    \mathrm{inv}(\gg)
    \\
    &&
      \Omega^\bullet(X)
        \ar@{<-}[ur]|{\{K_i\}}
  }
\)
and then \emph{lift} the resulting $\gg$-connection descent object $(A,F_A)$ to a
$\gg_\mu$-connection object $(\hat A, F_{\hat A})$
\(
  \xymatrix{
    \mathrm{CE}(b^{n-1}\uu(1))
    \ar@{<-}[r]
    \ar@{<-}[dd]
    &
    \mathrm{CE}(\gg_\mu)
    \ar@{<-}[rr]
    \ar@{<-}[dd]
    &&
    \mathrm{CE}(\gg)
    \ar@{<-}[dd]
    \\
     && \Omega^\bullet_{\mathrm{vert}}(Y)
       \ar@{<-}[ur]|{A_{\mathrm{vert}}}
       \ar@{<-}[dd]
       \ar@{<--}[ul]|{\hat A_{\mathrm{vert}}}
    \\
    \mathrm{W}(b^{n-1}\uu(1))
    \ar@{<-}[r]
    \ar@{<-}[dd]
    &
    \mathrm{W}(\gg_\mu)
    \ar@{<-}[rr]|{\hspace{12pt}}
    \ar@{<-}[dd]
    &&
    \mathrm{W}(\gg)
    \ar@{<-}[dd]
    \\
     && \Omega^\bullet(Y) 
       \ar@{<-}[ur]|{(A,F_A)}
       \ar@{<-}[dd]
       \ar@{<--}[ul]|{(\hat A_{\mathrm{vert}},F_{\hat A_{\mathrm{vert}}})}
    \\
    \mathrm{inv}(b^{n-1}\uu(1))
    \ar@{<-}[r]
    &
    \mathrm{inv}(\gg_\mu)
    \ar@{<-}[rr]|{\hspace{12pt}}
    &&
    \mathrm{inv}(\gg)
    \\
    &&
      \Omega^\bullet(X)
        \ar@{<-}[ur]|{\{K_i\}}
       \ar@{<--}[ul]|{\{\hat K_i\}}
  }
  \,.
\)
The situation is essentially an obstruction problem as before, only that
instead of single morphisms, we are now lifting an entire sequence of morphisms.
As before, we measure the obstruction to the existence of the lift 
by precomposing everything with the a map from a weak cokernel:
$$
 
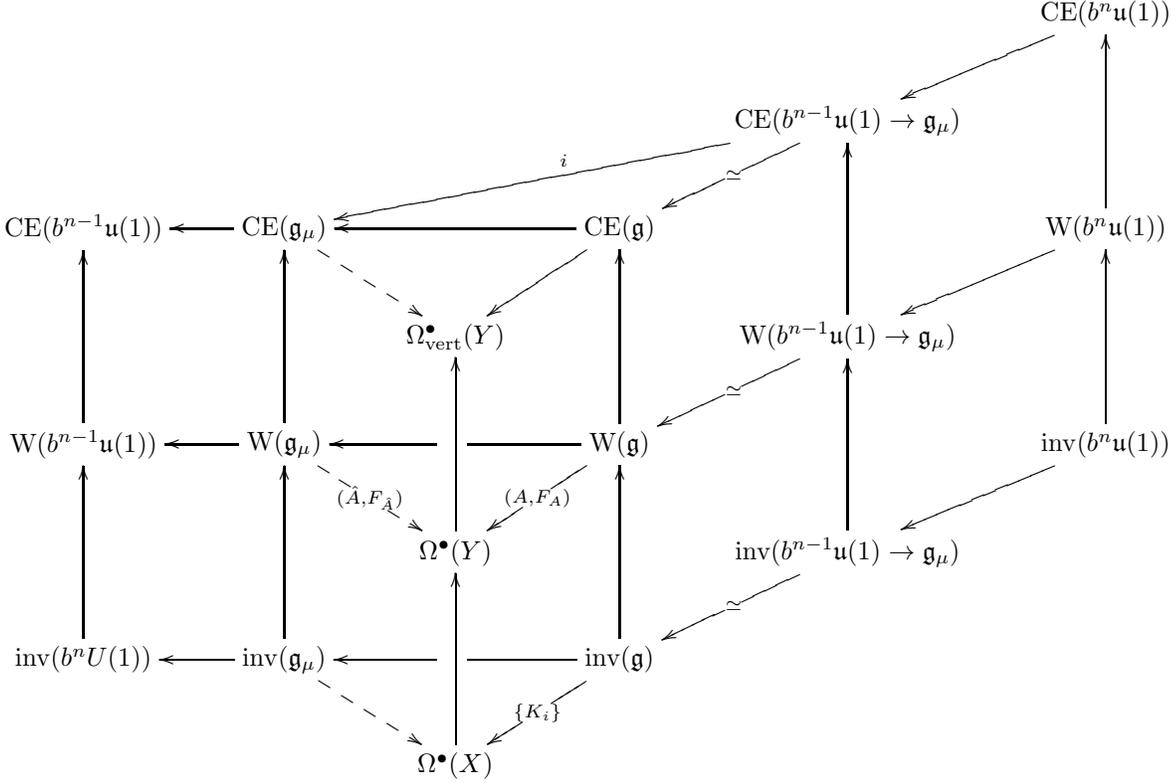
\hspace{-2cm}
  \xymatrix{
    &&&&&
    \mathrm{CE}(b^{n}\uu(1))
    \ar@{<-}[dd]
    \\
    &&&&
    \mathrm{CE}(b^{n-1}\uu(1)\to \gg_\mu)
    \ar@{<-}[ur]
    \ar@{<-}[dd]
    \\
    \mathrm{CE}(b^{n-1}\uu(1))
    \ar@{<-}[r]
    \ar@{<-}[dd]
    &
    \mathrm{CE}(\gg_\mu)
    \ar@{<-}[rr]
    \ar@{<-}[dd]
    \ar@{<-}[urrr]^i
    &&
    \mathrm{CE}(\gg)
    \ar@{<-}[dd]
    \ar@{<-}[ur]|{\simeq}
    &&
    \mathrm{W}(b^{n}\uu(1))
    \ar@{<-}[dd]
    \\
     && \Omega^\bullet_{\mathrm{vert}}(Y)
       \ar@{<-}[ur]
        \ar@{<--}[ul]
       \ar@{<-}[dd]
     &&
      \mathrm{W}(b^{n-1}\uu(1)\to \gg_\mu)
      \ar@{<-}[ur]
      \ar@{<-}[dd]
    \\
    \mathrm{W}(b^{n-1}\uu(1))
    \ar@{<-}[r]
    \ar@{<-}[dd]
    &
    \mathrm{W}(\gg_\mu)
    \ar@{<-}[rr]|{\hspace{12pt}}
    \ar@{<-}[dd]
    &&
    \mathrm{W}(\gg)
    \ar@{<-}[ur]|{\simeq}
    \ar@{<-}[dd]
    && 
    \mathrm{inv}(b^n \uu(1))
    \\
     && \Omega^\bullet(Y) 
       \ar@{<-}[ur]|{(A,F_A)}
       \ar@{<--}[ul]|{(\hat A,F_{\hat A})}
       \ar@{<-}[dd]
      &&
      \mathrm{inv}(b^{n-1}\uu(1)\to \gg_\mu)
      \ar@{<-}[ur]
    \\
    \mathrm{inv}(b^n U(1))
    \ar@{<-}[r]
    &
    \mathrm{inv}(\gg_\mu)
    \ar@{<-}[rr]|{\hspace{12pt}}
    &&
    \mathrm{inv}(\gg)
    \ar@{<-}[ur]|{\simeq }
    \\
    &&
      \Omega^\bullet(X)
        \ar@{<-}[ur]|{\{K_i\}}
        \ar@{<--}[ul]
        }
  $$

The result is a $b^n \uu(1)$-connection object. We will call 
(the class of) this the
{\bf generalized Chern-Simons $(n+1)$-bundle} obstructing the lift.

\begin{figure}[h]
$$
 \hspace{-2cm}
  \xymatrix{
    &&&&&
    \mathrm{CE}(b^{n}\uu(1))
    \ar@{<-}[dd]
    \ar@/^2.3pc/[dddlll]|<<<<<<<<<<<<<<<{\mathrm{CS}(A)_{\mathrm{vert}}}
    \\
    &&&&
    \mathrm{CE}(b^{n-1}\uu(1)\to \gg_\mu)
    \ar@{<-}[ur]
    \ar@{<-}[dd]|<<<<<<<<<<<{\makebox(8,8){}}
    \\
    \mathrm{CE}(b^{n-1}\uu(1))
    \ar@{<-}[r]
    \ar@{<-}[dd]
    &
    \mathrm{CE}(\gg_\mu)
    \ar@{<-}[rr]
    \ar@{<-}[dd]
    \ar@{<-}[urrr]^i
    &&
    \mathrm{CE}(\gg)
    \ar@{<-}[dd]|<<<<<<<<<{\makebox(8,8){}}
    \ar@{<-}[ur]|{\simeq}
    &&
    \mathrm{W}(b^{n}\uu(1))
    \ar@{<-}[dd]
    \ar@/^2.3pc/[dddlll]|<<<<<<<<<<<<<<<{(\mathrm{CS}(A),P(F_A))}
    \\
     && \Omega^\bullet_{\mathrm{vert}}(Y)
       \ar@{<-}[ur]
        \ar@{<--}[ul]
       \ar@{<-}[dd]
     &&
      \mathrm{W}(b^{n-1}\uu(1)\to \gg_\mu)
      \ar@{<-}[ur]
    \ar@{<-}[dd]|<<<<<<<<<<<{\makebox(8,8){}}
    \\
    \mathrm{W}(b^{n-1}\uu(1))
    \ar@{<-}[r]
    \ar@{<-}[dd]
    &
    \mathrm{W}(\gg_\mu)
    \ar@{<-}[rr]|{\hspace{12pt}}
    \ar@{<-}[dd]
    &&
    \mathrm{W}(\gg)
    \ar@{<-}[ur]|{\simeq}
    \ar@{<-}[dd]|<<<<<<<<<{\makebox(8,8){}}
    && 
    \mathrm{inv}(b^n \uu(1))
    \ar@/^2.3pc/[dddlll]|<<<<<<<<<<<<<<<{\{ P(F_A)\}}
    \\
     && \Omega^\bullet(Y) 
       \ar@{<-}[ur]|{(A,F_A)}
       \ar@{<--}[ul]|{(\hat A,F_{\hat A})}
       \ar@{<-}[dd]
      &&
      \mathrm{inv}(b^{n-1}\uu(1)\to \gg_\mu)
      \ar@{<-}[ur]
    \\
    \mathrm{inv}(b^n U(1))
    \ar@{<-}[r]
    &
    \mathrm{inv}(\gg_\mu)
    \ar@{<-}[rr]|{\hspace{12pt}}
    &&
    \mathrm{inv}(\gg)
    \ar@{<-}[ur]|{\simeq }
    \\
    &&
      \Omega^\bullet(X)
        \ar@{<-}[ur]|{\{K_i\}}
        \ar@{<--}[ul]
  }
$$
\caption{
  The {\bf generalized Chern-Simons $b^n\uu(1)$-bundle
  that obstructs the lift of a given $\gg$-bundle
  to a $\gg_\mu$-bundle}, or rather the descent object representing it.
}
\end{figure}

In order to construct the lift it is convenient, for similar reasons as in the
proof of proposition \ref{exactness of P in inv(g_mu)}, to work with
$\mathrm{CE}(\mathrm{cs}_P(\gg))$ instead of the isomorphic $\mathrm{W}(\gg_\mu)$,
using the isomorphism from proposition \ref{Chern and Chern-Simons}.
Furthermore, using the identity
\(
  \mathrm{W}(\gg) = \mathrm{CE}(\mathrm{inn}(\gg))
  \)
mentioned in \ref{Lie infty algebras}, we can hence consider instead of
\(
  \xymatrix{ 
     \mathrm{W}(b^{n-1})
     &&
     \mathrm{W}(\gg_\mu)
     \ar@{->>}[ll]
     &&
     \mathrm{W}(\gg)
     \ar@{_{(}->}[ll]
  }
\)
the sequence
\(
  \xymatrix{ 
     \mathrm{CE}(\mathrm{inn}(b^{n-1}))
     &&
     \mathrm{CE}(\mathrm{cs}_P(\gg))
     \ar@{->>}[ll]
     &&
     \mathrm{CE}(\mathrm{inn}(\gg))
     \ar@{_{(}->}[ll]
  }
  \,.
\)
Fortunately, this still satisfies the assumptions of proposition 
\ref{universal property of weak kernel}. So in complete analogy, 
we find the extension of proposition 
\ref{weak inverse for string extension}
from $\gg$-bundle cocyces to differential $\gg$-cocycles:

\begin{proposition}
  The morphism
  \(
    f^{-1} : \mathrm{CE}(
        \mathrm{inn}(b^{n-1}\uu(1)) \hookrightarrow 
    \mathrm{CE}(\mathrm{cs}_P(\gg))
    \to 
    \mathrm{CE}(\mathrm{inn}(\gg))
  \)
  constructed as in proposition \ref{weak inverse for string extension}
  acts as the identity on $\gg^* \oplus \gg^*[1]$
  \(
    f^{-1}|_{\gg^* \oplus \gg^*[1]} = \mathrm{Id}
  \)  
  and satisfies
  \(
    f^{-1}|_{\mathbb{R}[n+2]}
    : 
    c
    \mapsto 
    P
    \,.
  \)
\end{proposition}

This means that, as an extension of proposition
\ref{bundle cocycle obstruction for lift through string extension}, we find
the differential $b^n\uu(1)$ $(n+1)$-cocycle 
\(
  \xymatrix{
    \Omega^\bullet(Y)
    &&
    \mathrm{W}(b^{n}\uu(1))
    \ar[ll]_{\hat A}
  }
\)
obstructing the lift of a 
differential $\gg$-cocycle
\(
  \xymatrix{
    \Omega^\bullet(Y)
    &&
    \mathrm{W}(\gg)
    \ar[ll]_{(A,F_A)}
  },
\)
according to the above discussion
  \(
   \xymatrix@C=8pt{
   &&&&&&
   \mathrm{CE}(\mathrm{inn}(b^{n-1}\uu(1)) \hookrightarrow \mathrm{inn}(\gg_\mu))
   \ar[dllll]
   \ar[dll]|{f^{-1}}
   &&
   \mathrm{W}(b^n \uu(1))
   \ar@{_{(}->}[ll]_<<<j
   \ar[ddlllll]|{(\hat A, F_{\hat A})}
   \\
     \mathrm{W}(b^{n-1}\uu(1))
     &&
     \mathrm{W}(\gg_\mu)
     \ar@{->>}[ll]_{i^*}
     \ar@{..>}[dr]
     &&
     \mathrm{W}(\gg)     
     \ar@{_{(}->}[ll]
      \ar[dl]|{(A,F_A)}
     \\
     &&&\Omega^\bullet(Y)
   },
 \)
to be the connection $(n+1)$-form
\(
  \hat A = \mathrm{CS}(A) \in \Omega^{n+1}(Y)
\)
with the corresponding curvature $(n+2)$-form
\(
  F_{\hat A} = P(F_A) \in \Omega^{n+2}(Y)
  \,.
\)

Then we finally find, in particular,
\begin{proposition}
  For $\mu$ a cocycle on the ordinary Lie algebra $\gg$ in transgression
  with the invariant polynomial $P$, the obstruciton to lifting a
  $\gg$-bundle cocycle through the String-like extension determined
  by $\mu$ is the characteristic class given by $P$.
  \label{coc}
\end{proposition}

\paragraph{Remark.} Notice that, so far, all our statements about characteristic
classes are in deRham cohomology. Possibly our construction actually obtains for
integral cohomology classes, but if so, we have not extracted that yet. A more detailed
consideration of this will be the subject of \cite{SSS2}.

\subsubsection{Examples}

\label{obstruction examples}

\paragraph{Chern-Simons 3-bundles obstructing lifts of $G$-bundles to 
$\mathrm{String}(G)$-bundles.}

Consider, on a base space $X$ for some 
semisimple Lie group $G,$ with Lie algebra $\gg$ 
a principal 
$G$-bundle $\pi : P \to X$. Identify our surjective submersion with
the total space of this bundle 
\(
Y := P
 \,.
 \) 
Let $P$ be equipped with a connection, $(P,\nabla)$, 
realized in terms of an Ehresmann connection 1-form
\(
  A \in \Omega^1(Y,\gg)
\)
with curvature
\(
  F_A \in \Omega^2(Y,\gg)
\)
i.e. a dg-algebra morphism
\(
  \xymatrix{
    \Omega^\bullet(Y)
    &&
    \mathrm{W}(\gg)
    \ar[ll]_{(A,F_A)}
  }
\)
satisfying the two Ehresmann conditions. By the discussion in 
\ref{examples for connection descent objects} this yields a 
$\gg$-connection descent object $(Y,(A,F_A))$ in our sense. 

We would like to compute the obstruction to lifting this $G$-bundle
to a String 2-bundle, i.e. to lift the $\gg$-connection descent 
object to a $\gg_\mu$-connection descent object, for
\(
  0 \to b \uu(1) \to \gg_\mu \to \gg \to 0
\)
the ordinary String extension from definition \ref{string Lie 2-algebra}.   
By the above discussion in \ref{lifts of differential cocycles}, the 
obstruction is the (class of the) $b^2\uu(1)$-connection descent object
$(Y,(H_{(3)}, G_{(4)}))$ whose connection and curvature are given by
the composite
\(
 \hspace{-1cm}
  \raisebox{30pt}{
  \xymatrix{
    &
    &&
    \mathrm{W}(b^{2}\uu(1))
    \ar@/^2.6pc/[dddlll]|{(H_{(3)},G_{(4)})}
    \\
     &&
      (\mathrm{W}(b\uu(1)) \to \mathrm{CE}(\mathrm{cs}_P(\gg)))
      \ar@{<-}[ur]
    \\
    &
    \mathrm{W}(\gg)
    \ar@{<-}[ur]|{\simeq}
    && 
    \\
    \Omega^\bullet(Y) 
       \ar@{<-}[ur]|{(A,F_A)}
  }
  }
  \,,
\)
where, as discussed above, we are making use of the isomorphism $\mathrm{W}(\gg_\mu) \simeq 
\mathrm{CE}(\mathrm{cs}_P(\gg))$ from proposition \ref{Chern and Chern-Simons}.
The crucial aspect of this composite is the isomorphism
\(
  \xymatrix{
     \mathrm{W}(\gg)
     &&
     (\mathrm{W}(b \uu(1))
     \to
     \mathrm{CE}_P(\gg))
     \ar[ll]_<<<<<<<<{f^{-1}}^<<<<<<<<\simeq
  }
\)
from proposition \ref{weak inverse for weak cokernel}. This is where the
obstruction data is picked up. The important formula governing this is
equation \ref{crucial weak inverse}, which describes how the shifted elements
coming from $\mathrm{W}(b \uu(1))$ in the mapping cone
$     (\mathrm{W}(b \uu(1))
     \to
     \mathrm{CE}_P(\gg))$ are mapped to $\mathrm{W}(\gg)$.
     
Recall that $\mathrm{W}(b^2 \uu(1)) = \mathrm{F}(\mathbb{R}[3])$
is generated from elements $(h,dh)$ of degree 3 and 4, respectively,
that $\mathrm{W}(b \uu(1)) = \mathrm{F}(\mathbb{R}[2])$ is
generated from elements $(c, dc)$ of degree 2 and 3, respectively, and that
$\mathrm{CE}(\mathrm{cs}_P(\gg))$ is generated from 
$\gg^* \oplus \gg^*[1]$ together with elements $b$ and $c$ of degree 2
and 3, respectively, with
\(
  d_{\mathrm{CE}(\mathrm{cs}_P(\gg))} b = c - \mathrm{cs}
\)
and
\(
  d_{\mathrm{CE}(\mathrm{cs}_P(\gg))} c = P
  \,,
\)
where $\mathrm{cs} \in \wedge^3 (\gg^* \oplus \gg^*[1])$ is the transgression
element interpolating between the cocycle 
$\mu = \langle \cdot,  [\cdot,\cdot]\rangle \in \wedge^3 (\gg^*)$ and the
invariant polynomial $P = \langle \cdot, \cdot \rangle \in \wedge^2 (\gg^*[1])$.
Hence the map $f^{-1}$ acts as
\(
  f^{-1} : \sigma b \mapsto 
     -(d_{\mathrm{CE}(\mathrm{cs}_P(\gg))} b)|_{\wedge^\bullet (\gg^* \oplus \gg^*[1])} 
     = + \mathrm{cs}
\)
and
\(
  f^{-1} : \sigma c \mapsto 
    - (d_{\mathrm{CE}(\mathrm{cs}_P(\gg))} c)|_{\wedge^\bullet (\gg^* \oplus \gg^*[1])}
    = -P
  \,.
\)
Therefore the above composite $(H_{(3)}, G_{(4)})$ maps the generators $(h,dh)$
of $\mathrm{W}(b^2 \uu(1))$ as
\(
 \hspace{-1cm}
  \raisebox{30pt}{
  \xymatrix{
    &
    &&
    h
    \ar@{|->}@/^2.6pc/[dddlll]|{(H_{(3)},G_{(4)})}
    \\
     &&
      \makebox(130,12){$\sigma b$}
      \ar@{<-|}[ur]
    \\
    &
    \mathrm{cs}
    \ar@{<-|}[ur]|{\simeq}
    && 
    \\
    \mathrm{CS}_P(F_A) 
       \ar@{<-|}[ur]|{(A,F_A)}
  }
  }
\)
and
\(
 \hspace{-1cm}
  \raisebox{30pt}{
  \xymatrix{
    &
    &&
    dh
    \ar@{|->}@/^2.6pc/[dddlll]|{(H_{(3)},G_{(4)})}
    \\
     &&
      \makebox(130,12){-$\sigma c$}
      \ar@{<-|}[ur]
    \\
    &
    P
    \ar@{<-|}[ur]|{\simeq}
    && 
    \\
    P(F_A) 
       \ar@{<-|}[ur]|{(A,F_A)}
  }
  }
  \,.
\)
Notice  the signs here, as discussed around equation 
\ref{injection of W(bnu) into cs mapping cone}.
We then have that the connection 3-form of the Chern-Simons 3-bundle given by 
our obstructing $b^2\uu(1)$-connection descent object is the Chern-Simons
form
\(
  \label{3-fom connection of CS bundle}
  H_{(3)} = -\mathrm{CS}(A,F_A) = 
  -
  \langle
    A \wedge d A
  \rangle
  -
  \frac{1}{3}
  \langle
    A \wedge [A \wedge A]
  \rangle
  \in \Omega^3(Y)
\)
of the original Ehresmann connection 1-form $A$, and its 4-form curvature is therefore
the corresponding 4-form
\(
  G_{(4)} = - P(F_A) = \langle F_A \wedge F_A\rangle
  \in \Omega^4(Y)
  \,.
\)
This descends down to $X$, where it constitutes the characteristic form
which classifies the obstruction. Indeed, noticing that 
$\mathrm{inv}(b^2\uu(1)) = \wedge^\bullet (\mathbb{R}[4])$,
we see that (this works the same for all line $n$-bundles, i.e., for all
$b^{n-1}\uu(1)$-connection descent objects) the characteristic forms of the
obstructing Chern-Simons 3-bundle
\(
 \hspace{-1cm}
  \raisebox{30pt}{
  \xymatrix{
    &
    &&
    \mathrm{inv}(b^{2}\uu(1))
    \ar@/^2.6pc/[dddlll]|{\{G_{(4)}\}}
    \\
     &&
      \mathrm{inv}(b\uu(1) \to \gg_\mu)
      \ar@{<-}[ur]
    \\
    &
    \mathrm{inv}(\gg)
    \ar@{<-}[ur]|{\simeq}
    && 
    \\
    \Omega^\bullet(X) 
       \ar@{<-}[ur]|{\{K_i\}}
  }
  }
\)
consist only and precisely of this curvature 4-form: the second
Chern-form of the original $G$-bundle $P$.

\section{$L_\infty$-algebra parallel transport}

 \label{partra and sigma model}

One of the main points about a connection is that it allows to do parallel transport. 
Connections on ordinary bundles give rise to a notion of parallel transport along
curves, known as holonomy if these curves are closed.

Higher connections on $n$-bundles should yield a way to obtain a notion of parallel
transport over $n$-dimensional spaces. In physics, this assignment plays the
role of the gauge coupling term in the non-kinetic part of the action functional: 
the action functional of the charged particle is essentially its parallel
transport with respect to an ordinary (1-)connection, while the action functional
of the string contains the parallel transport of a 2-connection (the Kalb-Ramond field).
Similarly the action functional of the membrane contains the parallel transport of a
3-connection (the supergravity ``$C$-field'').

There should therefore be a way to assign to any one of our $\gg$-connection descent objects
for $\gg$ any Lie $n$-algebra
\begin{itemize}
  \item a prescription for parallel transport over $n$-dimensional spaces;
  \item a configuration space for the $n$-particle coupled to that transport;
  \item a way to transgress the transport to an action functional on that configuration space;
  \item a way to obtain the corresponding quantum theory.
\end{itemize}

Each point separately deserves a separate discussion, but in the remainder we shall quickly
give an impression for how each of these points is addressed in our context.

\subsection{$L_\infty$-parallel transport}

 \label{parallel transport}

In this section we indicate briefly how our notion of $\gg$-connections
give rise to a notion of parallel transport over $n$-dimensional
spaces. 
The abelian case (meaning here that $\gg$ is an $L_\infty$ algebra such
that $\mathrm{CE}(\gg)$ has trivial differential) 
is comparatively easy to discuss. It is in fact
the only case considered in most of the literature. Nonabelian parallel
$n$-transport in the integrated picture for $n$ up to 2 is discussed in
\cite{BS,SW,SWII,SWIII}. There is a close relation between all
differential concepts we develop here and the 
corresponding integrated concepts,
but here we will not attempt to give a comprehensive discussion
of the translation.

  Given an $(n-1)$-brane 
  (``$n$-particle'') whose $n$-dimensional 
  worldvolume is modeled on the 
  smooth parameter space $\Sigma$ (for instance 
  $\Sigma = T^2$ for the closed string) and which propagates
  on a target space $X$ in that its configurations are given by
  maps
  \(
    \phi : \Sigma \to X
  \)
  hence by dg-algebra morphisms
  \(
    \xymatrix{
      \Omega^\bullet(\Sigma)    
      &&
      \Omega^\bullet(X)
      \ar[ll]_{\phi^*}
    }
  \)
  we can couple it to
  a $\gg$-descent connection object $(Y, (A,F_A))$ over $X$
  pulled back to $\Sigma$
  if $Y$ is such that
  for every map 
  \(
    \phi : \Sigma \to X
  \)
  the pulled back surjective submersion has a global section
  \(
    \raisebox{20pt}{
    \xymatrix{
      & \phi^*Y
      \ar[d]^\pi
      \\
      \Sigma
      \ar[r]^{\mathrm{Id}}
      \ar[ur]^{\hat \phi}
      &
      \Sigma
    }
    }
    \,.
  \) 

\begin{definition}[parallel transport]
  Given a $\gg$-descent object $(Y,(A,F_A))$ 
  on a target space $X$ and a parameter space $\Sigma$
  such that for all maps $\phi : \Sigma \to X$ the pullback 
  $\phi^* Y$ has a global section,
  we obtain a map
  \(
    \mathrm{tra}_{(A)}
    :
    \mathrm{Hom}_{\mathrm{DGCA}}(
      \Omega^\bullet(X),
      \Omega^\bullet(\Sigma)
    )
    \to
    \mathrm{Hom}_{\mathrm{DGCA}}(\mathrm{W}(\gg), \Omega^\bullet(\Sigma))
  \)
  by precomposition with 
  \(
    \xymatrix{
      \Omega^\bullet(Y)
      &&
      \mathrm{W}(\gg)
      \ar[ll]_{(A,F_A)}
    }
    \,.
  \)
\end{definition}
  This is essentially the parallel transport of the $\gg$-connection 
  object $(Y,(A,F_A))$. 
  A full discussion is beyond the scope of 
  this article, but for the 
  special case that our $L_\infty$-algebra is $(n-1)$-fold 
  shifted $\uu(1)$, $\gg = b^{n-1}\uu(1)$, 
  the elements in 
  \(
   \mathrm{Hom}_{\mathrm{dgca}}(\mathrm{W}(\gg), \Omega^\bullet(\Sigma))
   =
   \Omega^\bullet(\Sigma, b^{n-1}\uu(1))
   \simeq
   \Omega^n(\Sigma)
  \)
  are in bijection with $n$-forms on $\Sigma$. Therefore
  they can be integrated over $\Sigma$. Then the functional 
  \(
    \int_\Sigma \mathrm{tra}_A : 
    \mathrm{Hom}_{\mathrm{dgca}}(\Omega^\bullet(Y), \Omega^\bullet(\Sigma))
    \to 
    \mathbb{R}
  \)
  is the full parallel transport of $A$.

  \begin{proposition}
    The map $\mathrm{tra}_{(A)}$ is indeed well defined, in that
    it depends at most on the homotopy class of the choice of global section
    $\hat \phi$ of $\phi$.
  \end{proposition}
  \proof
    Let $\hat \phi_1$ and $\hat \phi_2$ be two global sections of 
    $\phi^* Y$.
    Let $\hat \phi : \Sigma \times I \to \phi^*Y$ be a homotopy 
    between them,
    i.e. such that
    $\hat \phi|_0 = \hat \phi_1$ and 
    $\hat \phi|_1 = \hat \phi_2$. Then the difference in 
    the parallel transport using $\hat \phi_1$ and $\hat \phi_2$
    is the integral of the pullback of the curvature 
    form of the $\gg$-descent object over $\Sigma \times I$.
    But that vanishes, due to the commutativity of
    \(
      \xymatrix{  
          & \Omega^\bullet(\phi^*Y)
          \ar[ddl]|{\hat \phi^*}
          &&
          \mathrm{W}(\gg)
          \ar[ll]_{(A,F_A)}
          \\
          \\
          \Omega^\bullet(\Sigma \times I)
          &
          \Omega^\bullet(\Sigma)
          \ar@{_{(}->}[l]
          \ar@{^{(}->}[uu]
          \ar[l]_{\phi^*}
          &&
          **[r]\mathrm{inv}(b^{n-1}\uu(1)) = b^n \uu(1)
          \ar@/^2pc/[lll]^{0}
          \ar@{^{(}->}[uu]
          \ar[ll]_{K}
      }
    \)
    The composite of the morphisms on the top boundary of 
    this diagram send the single degree $(n+1)$-generator
    of $\mathrm{inv}(b^{n-1}\uu(1)) = \mathrm{CE}(b^n \uu(1))$
    to the curvature form of the $\gg$-connection descent
    object pulled back to $\Sigma$. 
    It is equal to the composite of the horizontal morphisms 
    along the bottom boundary by the definition of $\gg$-descent objects. 
    These vanish,
    as there is no nontrivial $(n+1)$-form on the 
    $n$-dimensional $\Sigma$.
  \endofproof

\subsubsection{Examples.}

\label{examples for parallel transport}

\paragraph{Chern-Simons and higher Chern-Simons action functionals}

\begin{proposition}
  \label{CS functional from transport}
  For $G$ simply connected, 
  the parallel transport coming from the Chern-Simons
  3-bundle discussed in \ref{obstruction examples} for
  $\gg = \mathrm{Lie}(G)$
  reproduces the familiar Chern-Simons action functional \cite{Freed}
  \(
    \int_\Sigma 
    \left(
       \langle A \wedge dA\rangle
       +
       \frac{1}{3}
       \langle A \wedge [A \wedge A] \rangle
    \right)
  \)
  over 3-dimensional $\Sigma$.
\end{proposition}
\proof
  Recall from \ref{obstruction examples} that we can build the
  connection descent object for the Chern-Simons connection on 
  the surjective submersion $Y$ coming from the total space $P$ of the
  underlying $G$-bundle $P \to X$. Then $\phi^* Y = \phi^* P$ 
  is simply the pullback of
  that $G$-bundle to $\Sigma$. For $G$ simply connected, $B G$ is 3-connected
  and hence any $G$-bundle on $\Sigma$ is trivializable. Therefore the
  required lift $\hat \phi$ exists and we can construct the above diagram.
  By equation \ref{3-fom connection of CS bundle} one sees that the 
  integral which gives the parallel transport is indeed precisely the
  Chern-Simons action functional.
\endofproof

Higher Chern-Simons $n$-bundles, coming from obstructions
to fivebrane lifts or still higher lifts, similarly induce
higher dimensional generalizations of the Chern-Simons
action functional.

\paragraph{BF-theoretic functionals}

From proposition \ref{transgressive elements from gg to hhtogg}
it follows that we can similarly obtain the action functional
of BF theory, discussed in \ref{characteristic forms}, as 
the parallel transport of the 4-connection descent object
which arises as the obstruction to lifting a 2-connection
descent object for a strict Lie 2-algebra 
$(\hh \stackrel{t}{\to} \gg)$ through the string-like extension
\(
  b^2 \uu(1) \to 
   (\hh \stackrel{t}{\to} \gg)_{d_{\mathrm{CE}(\hh \stackrel{t}{\to} \gg)}\mu}
   \to 
   (\hh \stackrel{t}{\to} \gg)
\)
for $\mu$ the 3-cocycle on $\mu$ which transgresses to the
invariant polynomial $P$ on $\gg$ which appears in the
BF-action functional.

\subsection{Transgression of $L_\infty$-transport}

An important operation on parallel transport is
its \emph{transgression} to mapping spaces.
This is familiar from simple examples, where for instance
$n$-forms on some space transgress to $(n-1)$-forms on
the corresponding loop space. We should think of the
$n$-form here as a $b^{n-1}\uu(1)$-connection 
which transgresses to an $b^{n-2}\uu(1)$ connection on 
loop space. 

This modification of the structure $L_\infty$-algebra under
transgression is crucial. In \cite{SWII} it is
shown that for parallel transport $n$-functors ($n=2$ there),
the operation of transgression is a very natural one,
corresponding to acting on the transport functor with 
an inner hom operation. As shown there, this operation automatically
induces the familiar pull-back followed by a fiber integration
on the corresponding differential form data, and also 
automatically takes care of the modification of the structure
Lie $n$-group.

The analogous construction in the differential world of $L_\infty$
algebras we state now, without here going into details about
its close relation to \cite{SWII}.

\begin{definition}[transgression of $\gg$-connections]
  For any $\gg$-connection descent object
  \(
    \xymatrix{
       F
       &&
       \mathrm{CE}(\gg)
       \ar[ll]_{A_{\mathrm{vert}}}
       \\
       \\
       P
       \ar@{->>}[uu]_{i^*}
       &&
       W(\gg)
       \ar@{->>}[uu]
       \ar[ll]_{(A,F_A)}^<{\ }="s"
       \\
       \\
       P_{\mathrm{basic}}
       \ar@{^{(}->}[uu]_{\pi^*}
       &&
       \mathrm{inv}(\gg)
       \ar@{^{(}->}[uu]
       \ar[ll]^{\{K_i\}}_>{\ }="t"
       %
    }
  \)
  and any smooth space $\mathrm{par}$, we can form the
  image of the above diagram under the functor
  \(
    \mathrm{maps}(-, \Omega^\bullet(\mathrm{par}))
    :
    \mathrm{DGCAs} \to \mathrm{DGCAs}
  \)
  from definition \ref{forms on maps functor} to 
  obtain the generalized $\gg$-connection descent object
  (according to definition \ref{generalized g-connection})
  \(
    \raisebox{50pt}{
    \xymatrix{
       \mathrm{maps}(F,\Omega^\bullet(\mathrm{par})
       &&
       \mathrm{maps}(\mathrm{CE}(\gg),\Omega^\bullet(\mathrm{par})
       \ar[ll]_{\mathrm{tg}_{\mathrm{par}}(A_{\mathrm{vert})}}
       \\
       \\
       \mathrm{maps}(P, \Omega^\bullet(\mathrm{par})
       \ar@{->>}[uu]_{\mathrm{tg}_{\mathrm{par}} i^*}
       &&
       \mathrm{maps}(W(\gg), \Omega^\bullet(\mathrm{par})       
       \ar@{->>}[uu]
       \ar[ll]_{\mathrm{tg}_{\mathrm{par}}(A,F_A)}^<{\ }="s"
       \\
       \\
       \mathrm{maps}(P_{\mathrm{basic}}, \Omega^\bullet(\mathrm{par})              
       \ar@{^{(}->}[uu]_{\mathrm{tg}_{\mathrm{par}}(\pi^*)}
       &&
       \mathrm{maps}(\mathrm{inv}(\gg), \Omega^\bullet(\mathrm{par})       
       \ar@{^{(}->}[uu]
       \ar[ll]_{\mathrm{tg}_{\mathrm{par}}(\{K_i\})}_>{\ }="t"
       %
    }
    }
    \,.
  \)
  This new $\mathrm{maps}(\mathrm{CE}(\gg),\Omega^\bullet(\mathrm{par}))$-connection
  descent object 
  we call the transgression of the original one to $\mathrm{par}$.
\end{definition}

The operation of transgression is closely related to that 
of integration. 

\subsubsection{Examples}

 \paragraph{Transgression of $b^{n-1}\uu(1)$-connections.}

  Let $\gg$ be an $L_\infty$-algebra of the form shifted $\uu(1)$, $\gg = b^{n-1}\uu(1)$.
  By proposition \ref{cohomology of shifted u(1)} the Weil algebra $\mathrm{W}(b^{n-1}\uu(1))$
  is the free DGCA on a single degree $n$-generator $b$ with differential $c := db$.
  Recall from \ref{examples for Loo-valued forms} that a DGCA morphism 
  $\mathrm{W}(b^{n-1}\uu(1)) \to \Omega^\bullet(Y)$ is just an $n$-form on $Y$.
  For every point $y \in \mathrm{par}$ and for every multivector
  $v \in \wedge^n T_y \mathrm{par}$
  we get a 0-form on the smooth space
  \(
    \mathrm{maps}(\mathrm{W}(b^{n-1}\uu(1)),\Omega^\bullet(\mathrm{par}))
  \)
  of all $n$-forms on $\mathrm{par}$, which we denote
  \(
    A(v) \in \Omega(\mathrm{maps}(\mathrm{W}(b^{n-1}\uu(1)),\Omega^\bullet(\mathrm{par})))
    \,.
  \)
  This is the 0-form on this space of maps obtained from the element 
  $b \in \mathrm{W}(b^{n-1}\uu(1))$ and the current $\delta_y$ (the ordinary
  delta-distribution on 0-forms) according to proposition \ref{forms from currents}.
  Its value on any any $n$-form $\omega$ is the value of that form evaluated on $v$.
  
  Since this, and its generalizations which we discuss in 
  \ref{example for oo-configuration spaces},
  is crucial for making contact with standard constructions in physics, it may be
  worthwhile to repeat that statement more explicitly in terms of components:
   Assume that $\mathrm{par} = \mathbb{R}^k$ and for any point $y$ let $v$ be the 
  unit in $\wedge T^n \mathbb{R}^n \simeq \mathbb{R}$. Then $A(v)$ is the 0-form on
  the space of forms which sends any form 
  $\omega = \omega_{\mu_1\mu_2 \dots \mu_n} dx^{\mu_1} \wedge \cdots \wedge dx^{\mu_n}$
  to its component
  \(
    A(v) : \omega \mapsto \omega(y)_{12, \cdots n}
    \,.
  \)
  This implies that when a $b^{n-1}\uu(1)$-connection is transgressed to the space
  of maps from an $n$-dimensional parameter space $\mathrm{par}$, it becomes a map
  that pulls back functions on the space of $n$-forms on $\mathrm{par}$ to the 
  space of functions on maps from parameter space into target space. But such 
  pullbacks correspond to \emph{functions} (0-forms) 
  on the space of maps $\mathrm{par} \to \mathrm{tar}$
  with values in the space of $n$-forms on $\mathrm{tra}$.

\subsection{Configuration spaces of $L_\infty$-transport}

 \label{configuration spaces}

  With the notion of $\gg$-connections and their parallel 
  transport and transgression in hand, we can say what it means
  to 
  \emph{couple an $n$-particle/$(n-1)$-brane to a $\gg$-connection}.

 \begin{definition}[the charged $n$-particle/$(n-1)$-brane]
   \label{charged n-particle}
    We say a charged $n$-particle/$(n-1)$-brane is 
    a tuple $(\mathrm{par},(A,F_A))$ consisting of
    \begin{itemize}
      \item {\bf parameter space} $\mathrm{par}$: a smooth space
      \item {\bf a background field} $(A,F_A)$: a $\gg$-connection 
               descent object involving
        \begin{itemize}
           \item {\bf target space} $\mathrm{tar}$: the smooth space
             that the $\gg$-connection $(A,F_A)$ lives over;
           \item {\bf space of phases} $\mathrm{phas}$: the smooth 
              space such that $\Omega^\bullet(\mathrm{phas}) 
                \simeq \mathrm{CE}(\gg)$
        \end{itemize}
    \end{itemize}   
    From such a tuple we form
    \begin{itemize}
      \item {\bf configuration space} 
       $\mathrm{conf} = \mathrm{hom}_{S^\infty}(\mathrm{par},\mathrm{tar})$;
      \item
        the {\bf action functional} $\exp(S) := 
           \mathrm{tg}_{\mathrm{par}}$: the transgression of the 
           background field to configuration space.
    \end{itemize}
 \end{definition}
   The configuration space thus  
    defined automatically comes equipped with a notion of
   vertical derivations as described in 
  \ref{vertical flows and basic forms}.
\(
  \xymatrix{
      \mathrm{maps}(F,\Omega^\bullet(\mathrm{par}))
      &&
      \mathrm{maps}(F,\Omega^\bullet(\mathrm{par}))
      \ar@/^1.5pc/[ll]^{[d,\rho']}_{\ }="t1"
      \ar@/_1.5pc/[ll]_{0}^{\ }="s1"
      \\
      \\
      \mathrm{maps}(P,\Omega^\bullet(\mathrm{par}))
      \ar@{->>}[uu]
      &&
      \mathrm{maps}(P,\Omega^\bullet(\mathrm{par}))
      \ar@{->>}[uu]
      \ar@/^1.5pc/[ll]^{[d,\rho]}_{\ }="t2"
      \ar@/_1.5pc/[ll]_{0}^{\ }="s2"
      \ar@{=>}^{\rho'} "s1"; "t1"
      \ar@{=>}^{\rho} "s2"; "t2"
  }
  \,.
\)
These form 
\begin{itemize}
  \item {\bf the gauge symmetries} $\gg_{\mathrm{gauge}}$:
    an $L_\infty$-algebra.
\end{itemize}
These act on the horizontal elements of configuration space, which form
\begin{itemize}
  \item the {\bf anti-fields and anti-ghosts}
\end{itemize}
in the language of BRST-BV-quantization \cite{TeitelboimHenneaux}.

We will not go into further details of this here, except for spelling out,
as the archetypical example, some details
of the computation of the configuration space of ordinary gauge theory.

\subsubsection{Examples}

\label{example for oo-configuration spaces}

\paragraph{Configuration space of ordinary gauge theory.}

We compute here the the configuration space of ordinary 
gauge theory on a manifold $\mathrm{par}$ with respect
to an ordinary Lie algebra $\gg$.
A configuration of such a theory is a $\gg$-valued differential
form on $\mathrm{par}$, hence, 
according to \ref{Lie infty-algebra valued forms},
an element in 
$\mathrm{Hom}_{\mathrm{DGCAs}}(\mathrm{W}(\gg),
\Omega^\bullet(\mathrm{par}))$. So we are interested in 
understanding the smooth space
\(
  \mathrm{maps}(\mathrm{W}(\gg),\Omega^\bullet(\mathrm{par})) =: \Omega^\bullet(\mathrm{par},\gg)
\)
according to definition \ref{space of maps between two DGCAs},
and the differential graded-commutative algebra
\(
  \mathrm{maps}(\mathrm{W}(\gg),\Omega^\bullet(\mathrm{par}))
  =: \Omega^\bullet(\Omega^\bullet(\mathrm{par},\gg))
\)
of differential forms on it.

To make contact with the physics literature, we 
describe everything in components.
So let $\mathrm{par} = \mathbb{R}^n$ 
and let $\{x^\mu\}$ be the canonical set of coordinate functions 
on $\mathrm{par}$.
Choose a basis $\{t_a\}$ of $\gg$ and let $\{t^a\}$ be the 
corresponding dual basis of 
$\gg^*$. Denote by
\(
  \delta_y \iota_{\frac{\partial}{\partial x^\mu}}
\)
the delta-current on $\Omega^\bullet(\mathrm{par})$, according to definition \ref{currents},
 which sends 
a 1-form $\omega $ to
\(
  \omega_\mu(y) := \omega(\frac{\partial}{\partial x^\mu})(y)
  \,.
\)

\subparagraph{Summary of the structure of forms on configuration space of ordinary gauge theory.}

Recall that the Weil algebra $\mathrm{W}(\gg)$ is generated 
from the $\{t^a\}$ in degree
1 and the $\sigma t^a$ in degree 2, with the differential defined by
\begin{eqnarray}
  d t^a &=& -\frac{1}{2}C^a{}_{bc} t^b \wedge t^c + \sigma t^a
\\
  d (\sigma t^a) &=& - C^a{}_{bc}t^b \wedge (\sigma t^c)
  \,.
\end{eqnarray}

We will find that 
$\mathrm{maps}(\mathrm{W}(\gg),\Omega^\bullet(\mathrm{par}))$ 
does look
pretty much entirely like this, only that all generators are now 
forms on $\mathrm{par}$. See table \ref{list of forms config space of ordinary gauge theory}.

\begin{table}[h]
  \begin{tabular}{l|l}
     {\bf fields} & 
             $\left\{ 
             A_\mu^a(y), (F_A)_{\mu \nu}(y)  \in \Omega^0(\Omega(\mathrm{par},\gg))
                   \hspace{4pt}|\hspace{4pt} 
                    y \in \mathrm{par}, \mu,\nu \in \{1,\cdots, 
                    \mathrm{dim}(\mathrm{par}), a \in \{1,\ldots, \mathrm{dim}(\gg)\}\}
              \right\}  
              $
     \\
     \\
     {\bf ghosts} 
        & 
           $\left\{ 
             c^a(y) \in \Omega^1(\Omega(\mathrm{par},\gg))
                \hspace{4pt}|\hspace{4pt} 
                y \in \mathrm{par}, a \in \{1,\ldots, \mathrm{dim}(\gg)\}\}
              \right\}  
           $
     \\
     \\
     {\bf antifields} 
        & 
           $\left\{ 
             \iota_{(\delta A^a_\mu(y))} 
               \in \mathrm{Hom}(\Omega^1(\Omega(\mathrm{par},\gg)),\mathbb{R})
                   \hspace{4pt}|\hspace{4pt} 
                    y \in \mathrm{par}, \mu \in \{1,\cdots, 
                    \mathrm{dim}(\mathrm{par}), a \in \{1,\ldots, \mathrm{dim}(\gg)\}\}
              \right\}  
           $
     \\
     \\
     {\bf anti-ghosts} 
        & 
           $\left\{ 
             \iota_{(\beta^a(y))} 
               \in \mathrm{Hom}(\Omega^2(\Omega(\mathrm{par},\gg)),\mathbb{R})
                   \hspace{4pt}|\hspace{4pt} 
                    y \in \mathrm{par},  
                    \mathrm{dim}(\mathrm{par}), a \in \{1,\ldots, \mathrm{dim}(\gg)\}\}
              \right\}  
           $
  \end{tabular}
     \caption{
        \label{list of forms config space of ordinary gauge theory}
       {\bf The BRST-BV field content of gauge theory} obtained from our almost internal
         hom of dg-algebras, definition \ref{forms on maps functor}.
         The dgc-algebra $\mathrm{maps}(\mathrm{W}(\gg),\Omega^\bullet(\mathrm{par}))$ is the
         algebra of differential forms on a smooth space of maps from $\mathrm{par}$ to the 
         smooth space underlying $\mathrm{W}(\gg)$.     In the above table $\beta$ is a certain
        2-form that one finds in this algebra of forms on the space of $\gg$-valued forms.
     }     
\end{table}

\subparagraph{Remark.}

Before looking at the details of the computation,
recall from from \ref{differential forms on spaces of maps} 
that an $n$-form $\omega$ in 
$\mathrm{maps}(\mathrm{W}(\gg),\Omega^\bullet(\mathrm{par}))$
is an assignment
\(
  \xymatrix{
     U
     \ar[dd]_{\phi}
     &&
     \mathrm{Hom}_{\mathrm{DGCA}s}(\mathrm{W}(\gg),\Omega^\bullet(\mathrm{par} \times U))
     \ar[rr]^<<<<<<<<{\omega_U}
     &&
     \Omega^\bullet(U)
     \\
     \\
     V
     &&
     \mathrm{Hom}_{\mathrm{DGCA}s}(\mathrm{W}(\gg),\Omega^\bullet(\mathrm{par} \times V))
     \ar[rr]^<<<<<<<<{\omega_V}
     \ar[uu]_{\phi^*}
     &&
     \Omega^\bullet(V)
     \ar[uu]_{\phi^*}
  }
\)
of forms on $U$ to $\gg$-valued forms on $\mathrm{par} \times U$ 
for all plot domains $U$ (subsets of $\mathbb{R} \cup \mathbb{R}^2 \cup \cdots$ for us), 
natural in $U$.
We	 concentrate on those $n$-forms $\omega$ which arise in the way 
of proposition \ref{forms from currents}.

\subparagraph{0-Forms.}

The 0-forms on the space of $\gg$-value forms are constructed as 
in proposition
\ref{forms from currents} from an element $t^a \in \gg^*$ and a current 
$\delta_y \iota_{\frac{\partial}{\partial x^\mu}}$ using
\(
  t^a \delta_y \iota_{\frac{\partial}{\partial x^\mu}}
\)
and from an element $\sigma t^a \in \gg^*[1]$ and a current
\(
  \delta_y \iota_{\frac{\partial}{\partial x^\mu}}\iota_{\frac{\partial}{\partial x^\nu}}
  \,.
\)
This way we obtain the families of functions (0-forms) on the space of $\gg$-valued forms:
\(
  A^a_{\mu}(y) : (\Omega^\bullet(\mathrm{par} \times U) 
     \leftarrow \mathrm{W}(\gg) : A) \mapsto 
     (u \mapsto \iota_{\frac{\partial}{\partial x^\mu}}A(t^a)(y,u))
\)
and
\(
  F^a_{\mu\nu}(y) : (\Omega^\bullet(\mathrm{par} \times U) \leftarrow \mathrm{W}(\gg) : F_A) 
    \mapsto 
   (u \mapsto \iota_{\frac{\partial}{\partial x^\mu}}\iota_{\frac{\partial}{\partial x^\nu}}F_A(\sigma t^a)(y,u))
\)
which pick out the corresponding components of the $\gg$-valued 1-form and of its curvature 2-form,
respectively.
These are the \emph{fields} of ordinary gauge theory.

\subparagraph{1-Forms.}

A 1-form on the space of $\gg$-valued forms is obtained 
from either starting with a degree 1 element and contracting with a degree 0 delta-current
\(
  t^a \delta_y
\)
or starting with a degree 2 element and contracting with a degree 1 delta current:
\(
  (\sigma t^a)\delta_y \frac{\partial}{\partial x^\mu}
  \,.
\)
To get started,  consider first the case where $U = I$ is the interval. Then a DGCA morphism
\(
  (A,F_A) : \mathrm{W}(\gg) \to \Omega^\bullet(\mathrm{par}) \otimes \Omega^\bullet(I)
\)
can be split into its components proportional to $dt\in \Omega^\bullet(I)$ and those
not containing $dt$. 
We can hence write the general $\gg$-valued 1-form on $\mathrm{par} \times I$ as
\(
  (A,F_A) : t^a \mapsto A^a(y,t)  + g^a(y,t) \wedge dt 
\)
and the corresponding curvature 2-form as
$$
  (A,F_A)
  :
  \sigma t^a \mapsto (d_{\mathrm{par}} + d_t)(A^a(y,t)  + g^a(y,t) \wedge dt )
    + \frac{1}{2}C^a{}_{bc}(A^a(y,t)  + g^a(y,t) \wedge dt)\wedge (A^b(y,t)  + g^b(y,t) \wedge dt)
$$
\(
  = F_A^a(y,t) + 
    ( \partial_t A^a(y,t) + d_{\mathrm{par}} g^a(y,t) + [g,A]^a ) \wedge dt
    \,.
\)
By contracting this again with the current $\delta_y \frac{\partial}{\partial x^\mu}$
we obtain the 1-forms
\(
  t \mapsto g^a(y,t) d t
\)
and
\(
  t \mapsto ( \partial_t A_\mu^a(y,t) + \partial_\mu g^a(y,t) + [g,A_\mu]^a ) d t
\)
on the interval. 
We will identify the first one with the component of the 1-forms on the space of $\gg$-valued forms
on $\mathrm{par}$ called the \emph{ghosts}  and the second one with the 1-forms which are killed by the
objects called the \emph{anti-fields}.

To see more of this structure, consider now $U = I^2$, the unit square. 
Then a DGCA morphism
\(
  (A,F_A) : \mathrm{W}(\gg) \to \Omega^\bullet(\mathrm{par}) 
     \otimes \Omega^\bullet(I^2)
\)
can be split into its components proportional to $dt^1, dt^2\in \Omega^\bullet(I^2)$.
We hence can write the general $\gg$-valued 1-form on $Y \times I$ as
\(
  (A,F_A) : t^a \mapsto A^a(y,t)  + g^a_i(y,t) \wedge dt^i 
  \,, 
\)
and the corresponding curvature 2-form as
$$
  (A,F_A)
  :
  \sigma t^a \mapsto (d_Y + d_{I^2})(A^a(y,t)  + g^a_i(y,t) \wedge dt^i )
$$
$$
    + \frac{1}{2}C^a{}_{bc}(A^a(y,t)  + g^a_i(y,t) \wedge dt^i )
     \wedge (A^b(y,t)  + g^b_i(y,t) \wedge dt^i )
$$
$$
  = F_A^a(y,t) + 
    ( \partial_{t^i} A^a(y,t) + d_Y g_i^a(y,t) + [g_i,A]^a ) \wedge dt^i
$$
\(
  + (\partial_i g^a_j + [g_i ,g_j]^a )dt^i \wedge dt^j 
  \,.
\)
By contracting this again with the current $\delta_y \frac{\partial}{\partial x^\mu}$
we obtain the 1-forms
\(
  t \mapsto g^a_i(y,t) d t^i
\)
and
\(
  t \mapsto ( \partial_t A_\mu^a(y,t) + \partial_\mu g_i^a(y,t) + [g_i,A_\mu]^a ) d t^i
\)
on the unit square.
These are again the local values of our
\(
  c^a(y) \in \Omega^1(\Omega^\bullet(\mathrm{par},\gg))
\)
and
\(
  \delta A^a_\mu(Y) \in \Omega^1(\Omega^\bullet(\mathrm{par},\gg))
  \,.
\)
The second 1-form vanishes in directions 
in which the variation of the $\gg$-valued 1-form $A$ is
a pure gauge transformation induced by the 
function $g^a$ which is measured by the first 1-form.
Notice that it is the sum of the exterior derivative of the 
0-form $A^a_\mu(y)$ with another
term. 
\(
  \delta A^a_\mu(y) = d (A^a_\mu(y)) + \delta_g A^a_\mu(y)
  \,.
\)
The first term on the right measures the change of the connection, 
the second subtracts the 
contribution to this change due to gauge transformations. So the 1-form 
$\delta A^a_\mu(y)$ on the space of $\gg$-valued 
forms vanishes along all 
directions along which the form $A$ is modfied 
purely by a gauge transformation.
The $\delta A^a_\mu(y)$ are the 1-forms the pairings dual to which 
will be the \emph{antifields}.

\subparagraph{2-Forms.}

We have already seen the 2-form appear on the standard square. We call this 2-form
\(
  \beta^a \in \Omega^2(\Omega^\bullet(\mathrm{par},\gg))
  \,,
\)
corresponding on the unit square to the assignment
\(
  \beta^a : (\Omega^\bullet(\mathrm{par} \times I^2) 
    \leftarrow \mathrm{W}(\gg) : A) 
   \mapsto 
   (\partial_i g_j^a + [g_i, g_j]^a )dt^i \wedge dt^j 
   \,.
\)
There is also a 2-form coming from $(\sigma t^a) \delta_y$.
Then one immediately sees that our forms on the space 
of $\gg$-valued forms satisfy
the relations
\begin{eqnarray}
  d c^a(y) &=& -\frac{1}{2}C^a{}_{bc}c^b(y) \wedge c^c(y) + \beta^a(y)
  \\
  d \beta^a(y) &=& - C^a{}_{bc}c^b(y)\wedge \beta^c(y)
  \,.
\end{eqnarray}
The 2-form $\beta$ on the space of $\gg$-valued forms is what is being 
contracted by the horizontal pairings called the \emph{antighosts}.
We see, in total, that $\Omega^\bullet(\Omega^\bullet(\mathrm{par},\gg))$ 
is the Weil algebra
of a DGCA, which is obtained from 
the above formulas by setting $\beta = 0$ and
$\delta A = 0$. This DGCA is the algebra of the 
gauge groupoid, that where the
only morphisms present are gauge transformations.

The computation we have just performed are over $U = I^2$. However, it should be
clear how this extends to the general case.

\paragraph{Chern-Simons theory.}

One can distinguish two ways to set up Chern-Simons theory. In one approach
one regards principal $G$-bundles on abstract 3-manifolds, in the other approach
one fixes a given principal $G$-bundle $P \to X$ on some base space $X$, and
pulls it back to  3-manifolds equipped with a map into $X$. 
Physically, the former case is thought of as Chern-Simons theory proper,
while the latter case arises as the gauge coupling part of the membrane
propagating on $X$.
One tends to
want to regard the first case as a special case of the second, obtained by letting
$X = BG$ be the classifying space for $G$-bundles and $P$ the 
universal $G$-bundle on that.

In our context this is realized by proposition \ref{line n-bundles on classifying spaces},
which gives the canonical Chern-Simons 3-bundle on $BG$ in terms of a $b^{2}\uu(1)$-connection
descent object on $\mathrm{W}(\gg)$.
Picking some 3-dimensional parameter space manifold $\mathrm{par}$, we can transgress this 
$b^{2}\uu(1)$-connection to the configuration space 
$\mathrm{maps}(\mathrm{W}(\gg),\Omega^\bullet(\mathrm{par}))$, which we learned is the 
configuration space of ordinary gauge theory. 
$$
     \raisebox{75pt}{
    \xymatrix{
       \mathrm{CE}(\gg)
       &&
       \mathrm{CE}(b^{n-1}\uu(1))
       \ar[ll]_{\mu}
       \\
       \\
       \mathrm{W}(\gg)
       \ar@{->>}[uu]_{i^*}
       &&
       W(b^{n-1}\uu(1))
       \ar@{->>}[uu]
       \ar[ll]_{(\mathrm{cs},P)}^<{\ }="s"
       \\
       \\
       \mathrm{inv}(\gg)
       \ar@{^{(}->}[uu]
       &&
       \mathrm{inv}(b^{n}\uu(1))
       \ar@{^{(}->}[uu]
       \ar[ll]^{P}_>{\ }="t"
       %
    }
    }
    \hspace{10pt}
      \mapsto
    \hspace{10pt}
     \raisebox{75pt}{
    \xymatrix@C=9pt{
       \mathrm{maps}(\mathrm{CE}(\gg),\Omega^\bullet(\mathrm{par})
       &&
       \mathrm{maps}(\mathrm{CE}(b^{n-1}\uu(1)),\Omega^\bullet(\mathrm{par}))
       \ar[ll]_{\mathrm{tg}_{\mathrm{par}}\mu}
       \\
       \\
       \mathrm{maps}(\mathrm{W}(\gg),\Omega^\bullet(\mathrm{par}))
       \ar@{->>}[uu]_{i^*}
       &&
       \mathrm{maps}(W(b^{n-1}\uu(1)),\Omega^\bullet(\mathrm{par}))
       \ar@{->>}[uu]
       \ar[ll]_{\mathrm{tg}_{\mathrm{par}}(\mathrm{cs},P)}^<{\ }="s"
       \\
       \\
       \mathrm{maps}(\mathrm{inv}(\gg),\Omega^\bullet(\mathrm{par}))
       \ar@{^{(}->}[uu]
       &&
       \mathrm{maps}(\mathrm{CE}(b^{n-1}\uu(1)),\Omega^\bullet(\mathrm{par}))
       \ar@{^{(}->}[uu]
       \ar[ll]_{\mathrm{tg}_{\mathrm{par}}P}_>{\ }="t"
       %
    }
    }
    \,.
$$
Proposition \ref{CS functional from transport} says that the 
transgressed connection is the Chern-Simons
action functional.

Further details of this should be discussed elsewhere.

\paragraph{Transgression of $p$-brane structures to loop space}

It is well known that obstructions to String structures on a space $X$ --
for us: Chern-Simons 3-bundles as in \ref{lifting problem} -- can be conceived 

\begin{itemize}

\item 
either in terms of 
a 3-bundle on $X$ classified by a four class on $X$ obstructing the
lift of a 1-bundle on $X$ to a 2-bundle;

\item
 or in terms of a 2-bundle on $L X$ classified by a 3-class on $L X$
 obstructing the lift of a 1-bundle on $L X$ to another 1-bundle, principal
 for a Kac-Moody central extension of the loop group.
\end{itemize}

In the second case, one is dealing with the transgression
of the first case to loop space.

The relation between the two points of views is carefully described in
\cite{Kuribayashi}. Essentially, the result is that \emph{rationally}
both obstructions are equivalent.

\paragraph{Remark.} Unfortunately, there is no universal agreement on
the convention of the direction of the operation called transgression.
Both possible conventions are used in the iterature relevant for
our purpose here. For instance \cite{BrylinskiMcLaughlin} say transgression
for what \cite{Asada} calls the inverse of transgression (which, in turn,
should be called suspension).

We will demonstrate in the context of $L_\infty$-algebra connections
how Lie algebra $(n+1)$-cocycles related to $p$-brane structures on $X$
transgress to loop Lie algebra $n$-cocycles on loop space. One can understand
this also as an alternative proof of the strictification theorem of the
String Lie 2-algebra (proposition \ref{string Lie 2-algebra as strict Lie 2-algebra}),
but this will not be further discussed here.

So let $\gg$ be an ordinary Lie algebra, $\mu$ an $(n+1)$-cocycle on it
in transgression with an invariant polynomial $P$, where the transgression
is mediated by the transgression element $\mathrm{cs}$ as described 
in \ref{Lie infty-algebra cohomology}.

According to proposition \ref{line n-bundles on classifying spaces} 
the corresponding universal obstruction structure
is the $b^n \uu(1)$-connection
\(
  \raisebox{80pt}{
  \xymatrix{
    \mathrm{CE}(\gg)
    &&
    \mathrm{CE}(b^n\uu(1))
    \ar[ll]_{\mu}
    \\
    \\
    \mathrm{W}(\gg)
    \ar@{->>}[uu]
    &&
    \mathrm{W}(b^n\uu(1))
    \ar[ll]_{(\mathrm{cs},P)}
    \ar@{->>}[uu]
    \\
    \\
    \mathrm{inv}(\gg)
    \ar@{_{(}->}[uu]
    &&
    **[r]\mathrm{inv}(b^{n}\uu(1)) = \mathrm{CE}(b^{n+1}\uu(1))
    \ar[ll]_{\{P\}}
    \ar@{_{(}->}[uu]
  }
  }
\)
to be thought of as the universal higher Chern-Simons $(n+1)$-bundle with connection
on the classifying space of the simply connected Lie group integrating $\gg$.

We transgress this to loops by applying the functor
$\mathrm{maps}(-, \Omega^\bullet(S^1))$ from definition \ref{forms on maps functor}
to it, which can be thought of as computing for
all DGC algebras the DGC algebra of differential forms on the space of maps from the
circle into the space that the original DGCA was the algebra of differential forms of:
\(
  \raisebox{80pt}{
  \xymatrix{
    \mathrm{maps}(\mathrm{CE}(\gg),\Omega^\bullet(S^1))
    &&
    \mathrm{maps}(\mathrm{CE}(b^n\uu(1)),\Omega^\bullet(S^1))
    \ar[ll]_{\mathrm{tg}_{S^1} \mu}
    \\
    \\
    \mathrm{maps}(\mathrm{W}(\gg),\Omega^\bullet(S^1))
    \ar@{->>}[uu]
    &&
    \mathrm{maps}(\mathrm{W}(b^n\uu(1)),\Omega^\bullet(S^1))
    \ar[ll]_{\mathrm{tg}_{S^1}(\mathrm{cs},P)}
    \ar@{->>}[uu]
    \\
    \\
    \mathrm{maps}(\mathrm{inv}(\gg),\Omega^\bullet(S^1))
    \ar@{_{(}->}[uu]
    &&
    \mathrm{maps}(\mathrm{CE}(b^{n+1}\uu(1)),\Omega^\bullet(S^1))
    \ar[ll]_{\{\mathrm{tg}_{S^1} P\}}
    \ar@{_{(}->}[uu]
  }
  }
  \,.
\)
We want to think of the result as a $b^{n-1}\uu(1)$-bundle. This we can achieve
by pulling back along the inclusion
\(
  \mathrm{CE}(b^{n-1}\uu(1)) \hookrightarrow 
    \mathrm{maps}(\mathrm{CE}(b^n\uu(1)),\Omega^\bullet(S^1))
\)
which comes from the integration current $\int_{S^1}$ on $\Omega^\bullet(S^1)$
according to proposition \ref{forms from currents}. 

(This restriction to the integration current can be understood from looking at
the basic forms of the loop bundle descent object, which induces 
\emph{integration without integration} essentially  in the sense of 
\cite{Kauffmann}. But this we shall not further go into here.)

We now show that the transgressed cocycles $\mathrm{tg}_{S^1}\mu$
are the familiar cocycles on loop algebras, as appearing for instance
in Lemma 1 of \cite{Asada}. For simplicity of exposition, we shall consider
explicitly just the case where $\mu = \langle \cdot, [\cdot, \cdot ]\rangle$
is the canonical 3-cocycle on a Lie algebra with bilinear invariant form
$\langle \cdot, \cdot \rangle$.

\begin{proposition}
  \label{transgression to get Kac-Moody cocycle}
  The transgressed cocycle in this case is the 2-cocycle of the
  Kac-Moody central extension of the loop Lie algebra $\Omega \gg$
  \(
    \mathrm{tg}_{S^1} \mu : (f,g)
    \mapsto
    \int_{S^1} \langle f(\sigma), g'(\sigma)\rangle d\sigma
    + 
    \mbox{(a coboundary)}
  \)
  for all $f,g \in \Omega \gg$\,.
\end{proposition}
\proof
We compute $\mathrm{maps}(\mathrm{CE}(\gg),\Omega^\bullet(S^1))$
as before from proposition \ref{forms from currents} along the same lines
as in the above examples:
for $\{t_a\}$ a basis of $\gg$ and $U$ any test domain, a DGCA homomorphism
\(
  \phi : \mathrm{CE}(\gg) \to \Omega^\bullet(S^1)\otimes \Omega^\bullet(U)
\)
sends
\(
  \xymatrix{
    t^a 
    \ar@{|->}[rr]^{\phi}
    \ar@{|->}[dd]^{d_{\mathrm{CE}(\gg)}}
    && c^a + A^a \theta
    \ar@{|->}[dd]^{d_{S^1} + d_U}
    \\
    \\
    -\frac{1}{2}
    C^a{}_{bc}t^b \wedge t^c
    \ar@{|->}[rr]^{\phi}
    &&
    {
      {
        \theta \wedge (\frac{\partial}{\partial \sigma}c^a) + d_U c^a + 
             d_U A^a \wedge \theta
      }
      \atop
      {= -\frac{1}{2}C^a{}_{bc} c^b \wedge c^c}
      - C^a{}_{bc} c^b \wedge A^b \wedge \theta
    }
  }
  \,.
\)
Here $\theta \in \Omega^1(S^1)$ is the canonical 1-form on $S^1$ and
$\frac{\partial}{\partial \sigma}$ the canonical  vector field;
moreover $c^a\in \Omega^0(S^1)\otimes\Omega^1(U)$
and $A^a \theta \in \Omega^1(S^1) \otimes \Omega^0(U)$. 

By contracting with $\delta$-currents on $S^1$
we get 1-forms $c^a(\sigma)$, $\frac{\partial}{\partial \sigma} c^a(\sigma)$ 
and 0-forms $A^a(\sigma)$ 
for all $\sigma \in S^1$ on $\mathrm{maps}(\mathrm{CE}(\gg),\Omega^\bullet(S^1))$
satisfying
\(
  \label{ghost algebra on the circle}
  d_{\mathrm{maps}(\cdots)} c^a(\sigma) + \frac{1}{2} C^a{}_{bc} c^b(\sigma) \wedge c^c(\sigma) = 0
\)
and
\(
  \label{differential for field on the circle}
  d_{\mathrm{maps}(\cdots)} A^a(\sigma) - C^a{}_{bc}A^b(\sigma) \wedge c^c(\sigma) 
     = \frac{\partial}{\partial \sigma} c^a(\sigma)
  \,.
\)
Notice the last term appearing here, which is the crucial one responsible for
the appearance of derivatives in the loop cocycles, as we will see now.

So $A^a(\sigma)$ (a ``field'') is the function on 
(necessarily flat) $\gg$-valued 1-forms on $S^1$ which sends each such 
1-form for its $t^a$-component along $\theta$ at $\sigma$, while $c^a(\sigma)$ (a ``ghost'') is the 1-form 
which sends
each tangent vector field to the space of flat $\gg$-valued forms to the gauge transformation
in $t^a$ direction which it induces on the given 1-form at $\sigma \in S^1$.

Notice that the transgression of our 3-cocycle
\(
  \mu = \mu_{abc}t^a \wedge t^b \wedge t^c = 
   C_{abc}t^a \wedge t^b \wedge t^c \in H^3(\mathrm{CE}(\gg))
\)
is
\(
  \mathrm{tg}_{S^1} \mu 
  = 
  \int_{S^1} C_{abc} A^a(\sigma)  c^b(\sigma)\wedge c^c(\sigma) \; d\sigma
  \hspace{5pt}
  \in
  \hspace{5pt}
  \Omega^2(\Omega^1_{\mathrm{flat}}(S^1, \gg)
  \,.
\)
We can rewrite this using the identity
\(
  \label{differential for Ac term}
  d_{\mathrm{maps}(\cdots)} \left(
     \int_{S^1} P_{ab} A^a(\sigma) c^b(\sigma) d \sigma
  \right)
  =
  \int_{S^1} P_{ab} \left(\partial_\sigma c^a(\sigma)\right) \wedge c^b(\sigma) 
  + \frac{1}{2}\int_{S^1} C_{abc}A^a(\sigma) c^b(\sigma) \wedge c^c(\sigma)
  \,,
\)
which follows from \ref{ghost algebra on the circle} and \ref{differential for field on the circle},
as
\(
  \mathrm{tg}_{S^1}\mu = \int_{S^1} P_{ab} \left(\partial_\sigma c^a(\sigma)\right) \wedge c^b(\sigma)
   + d_{\mathrm{maps}(\cdots)}(\cdots)
   \,.
\)

Then notice that 
\begin{itemize}
  \item
    equation \ref{ghost algebra on the circle}
    is the Chevalley-Eilenberg algebra of the loop algebra $\Omega \gg$;
  \item
   the term $\int_{S^1} P_{ab} (\partial_\sigma c^a(\sigma))\wedge c^b(\sigma)$ 
   is the familiar 2-cocycle on the loop algebra obtained
   from transgression of the 3-cocycle $\mu = \mu_{abc}t^a \wedge t^b \wedge t^c = 
   C_{abc}t^a \wedge t^b \wedge t^c$
   \,.
\end{itemize}
\endofproof

\appendix

\section{Appendix: Explicit formulas for 2-morphisms of $L_\infty$-algebras}

\label{explicit formulas for homotopies}

To the best of our knowledge, the only place in the literature where 
2-morphisms between 1-morphisms of $L_\infty$-algebras have been spelled
out in detail is \cite{BaC}, which gives a definition of 
2-morphisms for Lie 2-algebras, i.e. for $L_\infty$-algebras concentrated
in the lowest two degrees.
Our definition \ref{transformation of qDGCA morphisms} provides
an algorithm for 
computing 2-morphisms between morphisms of arbitrary (finite dimensional)
$L_\infty$-algebras. We had already demonstrated in \ref{homotopies and concordances} 
one application of that algorithm, showing explicitly how it allows to compute
transgression elements (Chern-Simons forms).

For completeness, we demonstrate that the formulas given in 
\cite{BaC} for the special case of Lie 2-algebras also follow
as a special case from our general definition
\ref{transformation of qDGCA morphisms}.
This is of relevance to our discussion of the String Lie 2-algebra,
since the equivalence of its strict version with its weak skeletal
version, mentioned in our proposition \ref{string Lie 2-algebra as strict Lie 2-algebra},
has been established in \cite{BCSS} using these very formulas. 
First we quickly recall the relevant definitions from 
\cite{BaC,BCSS}:
A ``2-term'' $L_\infty$-algebra is an $L_\infty$-algebra concentrated in the
lowest two degrees.
A morphism
\(
\varphi : \gg \to \hh
\)
of 2-term $L_\infty$-algebras $\gg$ and $\hh$
is a pair of maps
\begin{eqnarray}
  \phi_0 &:& \gg_1 \to \hh_1
\\
  \phi_1 &:&  \gg_2 \to \hh_2
\end{eqnarray}
together with a skew-symmetric map
\(
  \phi_2 : \gg_1 \otimes \gg_1 \to \hh_2
\)
satisfying
\(
  \phi_0(d(h)) = d(\phi_1(h))
\)
as well as
\begin{eqnarray}
  d(\phi_2(x,y)) &=& \phi_0(l_2(x,y)) - l_2(\phi_0(x),\phi_0(y))  
\\
  \phi_2(x,dh) &=& \phi_1(l_2(x,h)) - l_2(\phi_0(x),\phi_1(h))
\end{eqnarray}
and finally
\begin{eqnarray}
  l_3(\phi_0(x),\phi_0(y),\phi_0(z)) - 
  \phi_1(l_3(x,y,z)) &=&
  \phi_2(x,l_2(y,z))
  +
  \phi_2(y,l_2(z,x))
  +
  \phi_2(z,l_2(x,y)) +
  \nonumber\\
&&  l_2(\phi_0(x),\phi_2(y,z))
  +l_2(\phi_0(y),\phi_2(z,x))
  +
  l_2(\phi_0(z),\phi_2(x,y))
  \nonumber
  \,.
\end{eqnarray}
for all $x,y,z \in \gg_1$ and $h \in \gg_2$.
This follows directly from the requirement that morphisms of $L_\infty$-algebras
be homomorphisms of the corresponding codifferential coalgebras, according to
definition \ref{Loo algebra}. The not quite so obvious aspect are the analogous
formulas for 2-morphisms:

\begin{definition}[Baez-Crans]
  \label{2-morphisms of Lie 2-algebras}
  A 2-morphism 
  \(
    \xymatrix{
      \gg
      \ar@/^2pc/[rr]^{\phi}_{\ }="s"
      \ar@/_2pc/[rr]_{\psi}^{\ }="t"
      &&
      \hh
      \ar@{=>}^\tau "s"; "t"
    } 
  \)
  of 1-morphisms of 2-term $L_\infty$-algebras is a linear map
  \(
    \tau : \gg_1 \to \hh_2
  \)
  such that
  \(
    \label{2-term Loo 2-morphism, first equation}
    \psi_0 - \phi_0 = t_W \circ \tau
  \)
  \(
    \label{2-term Loo 2-morphism, second equation}
    \psi_1 - \phi_1 = \tau \circ t_v
  \)
  and
  \(
    \label{2-term Loo 2-morphism, third equation}
    \phi_2(x,y)
    -
    \psi_2(x,y)
    =
    l_2(\phi_0(x),\tau(y))
    +
    l_2(\tau(x),\psi_0(y))
    -
    \tau(l_2(x,y))
  \)
\end{definition}

Notice that $[d,\tau] := d_\hh \circ \tau + \tau \circ d_\gg$
and that it restricts to $d_\hh \circ \tau$ on $\gg_1$ and
to $\tau \circ d_\gg$ on $\gg_2$.

\begin{proposition}
  For finite dimensional $L_\infty$-algebras,
  definition \ref{2-morphisms of Lie 2-algebras} is equivalent to 
  the restriction of our definition \ref{transformation of qDGCA morphisms}
  to 2-term $L_\infty$-algebras.
\end{proposition}
\proof
Let $\gg = \gg_1 \oplus \gg_2$ and $\hh = \hh_1 \oplus \hh_2$ be any two 2-term
$L_\infty$-algebras.
Then take
\(
   \psi, \phi : \gg \to \hh
\)
to be any two $L_\infty$ morphisms 
with 
\(
  \xymatrix{
    \mathrm{CE}(\gg)
    &&
    \mathrm{CE}(\hh)
    \ar[ll]_{\psi^* , \phi^*}
  }
\)
the corresponding DGCA morphisms.
We would like to describe the collection of all 2-morphisms 
\(
  \xymatrix{
    \mathrm{CE}(\gg)
    &&
    \mathrm{CE}(\hh)
    \ar@/_2pc/[ll]_{\phi^*}^{\ }="s"
    \ar@/^2pc/[ll]^{\psi^*}_{\ }="t"
    \ar@{=>}^\tau "s"; "t"
  }
\)
according to definition \ref{transformation of qDGCA morphisms}.
We do this in terms of a basis.
With $\{t^a\}$ a basis for $\hh_1$ and $\{b^i\}$ a basis for $\hh_2$,
and accordingly $\{t'^a\}$ and $\{b'^i\}$ a basis of $\gg_1$ and $\gg_2$,
respectively,
this comes from a map
\(
  \tau^* : \hh^*_1 \oplus \hh^*_2 \oplus \hh^*_1[1] \oplus \hh^*_2[1] 
  \to \wedge^\bullet \gg^* 
\)
of degree -1 which acts on these basis elements as
\(
  \tau^* : b^i \mapsto \tau^i{}_a t'^a
\)
and
\(
  \tau^* : a^a \mapsto 0
\)
for some coefficients $\{\tau^i{}_a\}$.
Now the crucial requirement \ref{vanishing condition on homotopies} 
of definition \ref{transformation of qDGCA morphisms}
is that \ref{2-morphism by pull-back to free qDGCA}
 \emph{vanishes} 
when restricted 
\(
\hspace{-1cm}
         \xymatrix{
           & \mathrm{CE}(\hh) 
           \ar@/_1pc/[dl]_{\phi^*}
           \\
			\mathrm{CE}(\gg) 
            && 
            \mathrm{W}(\hh)
				\ar@{->>}[ul]^{\ }="s"
            \ar@{->>}[dl]
            &
            \hh^*[1]
            \ar@{_{(}->}[l]
		   \\
		   & \mathrm{CE}(\hh) 
           \ar@/^1pc/[ul]^{\psi^*}_{\ }="t"
		   \ar@{=>} "s"; "t"
         }
       \)
to generators in the shifted copy of the Weil algebra.
This implies the following.
For $\tau^*$ to vanish on all $\sigma t^a$ we find that its value on
$d_{\mathrm{W}(\hh)} t^a = -\frac{1}{2}C^a{}_{bc}t^a \wedge t^b - t^a{}_i b^i + \sigma t^a$ 
is fixed to be
\(
  \tau^* : d_{\mathrm{W}(\gg)}  t^a \mapsto -t^a{}_i \tau^i{}_b t'^b
\)
and on $d_{\mathrm{W}(\hh)} b^i = -\alpha^i{}_{aj} t^a \wedge b^j + c^i$ to be
\(
  \tau^*(db^i) = \tau^*(-\alpha^i{}_{aj} t^a \wedge b^j)
  \,.
\)
The last expression needs to be carefully evaluated using formula 
\ref{formula for chain homotopy}.
Doing so we get
\(
  [d,\tau^*] : t^a \mapsto  - t^a{}_i \tau^i{}_b t'^b 
\)
and 
\(
  [d,\tau^*] : b^i \mapsto 
     -\frac{1}{2}\tau^i{}_{a} C'^a{}_{bc} t'^b t'^c 
     - \tau^i{}_a t'^a{}_j b'^j
      + \alpha^i{}_{aj} \frac{1}{2}(\phi + \psi)^a{}_{b}\tau^j{}_c t'^b t'^c
  \,.
\)
Then the expression
\(
  \phi^* - \psi^* = [d,\tau^*]
\)
is equivalent to the following ones
\begin{eqnarray}
  (\psi^a{}_b - \phi^a{}_b) t'^b &=& t^a{}_i \tau^i{}_b t'^b
\\
  (\psi^i{}_j - \phi^i{}_j) b'^j &=& \tau^i{}_a t'^a{}_j b'^j
\\
  \frac{1}{2}(\phi^i{}_{ab} - \psi^i{}_{ab}) t'^a t'^b 
  &=&
     -\frac{1}{2}\tau^i{}_{a} C'^a{}_{bc} t'^b t'^c 
      + \alpha^i{}_{aj} \frac{1}{2}(\phi + \psi)^a{}_b\tau^j{}_c t'^b t'^c  
  \,.
\end{eqnarray}
The first two equations express the fact that $\tau$ is a chain homotopy
with respect to $t$ and $t'$.
The last equation is equivalent to
\begin{eqnarray}
  \phi_2(x,y) - \psi_2(x,y)
  &=&
  -
  \tau([x,y])
  +
  [q(x) + \frac{1}{2}t(\tau(x)),\tau(y)]
  -
  [q'(y) - \frac{1}{2}t(\tau(y)),\tau(x)]
  \nonumber\\
  &=&
  -\tau([x,y])
  +
  [q(x),\tau(y)]
  +
  [\tau(x),q'(y)]  
\end{eqnarray}
This is indeed the Baez-Crans condition on a 2-morphism.
\endofproof

\vspace{1cm}

\noindent {\bf \large Acknowledgements}\\

We thank Danny Stevenson for many helpful comments.
Among other things, the discussion of the cohomology of the string Lie 2-algebra
has been greatly influenced by conversation with him. 
We acknowledge stimulating remarks by John Baez 
on characteristic classes of string bundles. H.S. thanks Matthew 
Ando and U.S. thanks Stephan Stolz and Peter 
Teichner for discussions at an early stage of this work.

U.S. thankfully acknowledges invitations to the University of Oxford
by Nigel Hitchin;
to the conference ``Lie algebroids and Lie groupoids in 
differential geometry''
in Sheffield by Kirill Mackenzie and Ieke Moerdijk;
to the Erwin Schr{\"o}dinger institute in the context of the program
``Poisson $\sigma$-models, Lie algebroids, deformations and higher analogues'';
to ``Categories in Geometry and Physics'' by Zoran {\v S}koda and Igor
Bakovi{\v c};
and to Yale University which led to collaboration with H.S. 
and to useful discussions with Mikhail Kapranov. U.S. thanks Todd Trimble for 
discussion of DGCAs,
and their relation to smooth spaces, thanks Mathieu Dupont 
and Larry Breen for discussion about weak cokernels and 
thanks Simon Willerton and Bruce Bartlett for helpful disucssion
of the notion of transgression. H.S. thanks Akira Asada and 
Katsuhiko Kuribayashi for useful correspondence on their work
on (higher) string structures.

We thank the contributors to the weblog
\emph{The $n$-Category Caf{\'e}} for much interesting, 
often very helpul and sometimes outright invaluable discussion concerning
the ideas presented here and plenty of related issues. 

We have gratefully received useful comments on earlier versions of 
this document from Gregory Ginot and David Roberts. Evan Jenkins and
Todd Trimble have provided useful help with references to literature on smooth
function algebras. 

While this article was being written, the preprints \cite{Yekutieli}, \cite{KS}
and \cite{GinotStienon} (based on \cite{GinotXu}) 
appeared, which are closely related to our discussion here
in that they address the issues of characteristic classes and extensions of
$n$-bundles in one way or another. We are thankful to 
Amnon Yekutieli and to Gregory Ginot for 
pointing out their work to us and for further discussion.

Finally, we thank Bertfried Fauser and the
organizers of the conference ``Recent developments in QFT'' 
in Leipzig.

U. S. thankfully acknowledges financial support from the
DAAD German-Croatian bilateral project ``Nonabelian cohomology and applications''.


\vspace{1cm}

{\small

\begin{tabular}{ccc}

\begin{tabular}{l}
Hisham Sati\\
Department of Mathematics\\Yale University\\10 Hillhouse Avenue\\
New Haven, CT 06511
\end{tabular}
&
\begin{tabular}{l}
Urs Schreiber\\
Fachbereich Mathematik\\Schwerpunkt Algebra und Zahlentheorie\\Universit\"at
Hamburg\\Bundesstra\ss e 55\\D--20146 Hamburg 
\end{tabular}
&
\begin{tabular}{l}
Jim Stasheff\\
Department of Mathematics\\University of Pennsylvania\\
David Rittenhouse Lab. \\
209 South 33rd Street \\
Philadelphia, PA 19104-6395
\end{tabular}

\end{tabular}
}

\end{document}